\documentclass[12pt,bezier]{article}

\usepackage{booktabs}
\usepackage[UTF8, fontset=fandol]{ctex} 
\usepackage{enumitem}
\usepackage{pifont}
\usepackage{floatrow}
\usepackage{caption}
\usepackage{mathrsfs}
\usepackage[fleqn]{amsmath}
\usepackage{amsfonts,amsthm,amssymb,mathrsfs,bbding}
\usepackage{txfonts}
\usepackage{tikz}
\usepackage{graphics,multicol}
\usepackage{graphicx}
\usepackage{geometry}
\usepackage{color}

\usepackage[colorlinks=true,anchorcolor=blue,filecolor=blue,linkcolor=blue,urlcolor=blue,citecolor=black]{hyperref}
\usepackage{caption}
\usepackage{tikz}
\usepackage{cite}
\usepackage{latexsym,bm}
\usepackage{indentfirst}
\usepackage{mathtools}
\usepackage{amsmath}
\evensidemargin=\oddsidemargin\topmargin=-1.5cm

\newtheorem{thm}{Theorem}[section]

\newtheorem{lem}{Lemma}[section]
\newtheorem{cor}{Corollary}[section]

\newtheoremstyle{factstyle}
  {}{}{\itshape}{}{\bfseries}{)}{ }
  {(F\thmnumber{#2}}

\theoremstyle{factstyle}
\newtheorem{fact}{F}[section]

\theoremstyle{definition}

\addtocounter{section}{0}

\usetikzlibrary{decorations.pathreplacing}
\begin{document}
\title{Planar graphs whose second largest eigenvalue does not exceed $1$\footnote{This work is supported
by the NSFC (No. 12171154).}}
\author{{\Large Zihao Geng$^{a}$,\,\, Zhiwen Wang$^{a}$,\,\, Wen-Jun Li$^{b}$,\,\, Xue-Ying Li$^{c}$,}\\
{\Large Qianqian Chen$^{d}$\footnote{
This work was supported by the NSFC (No. 12601690), the Yancheng Applied Basic
Research Program (No. YCBK2026136).},\,\, Ji-Ming Guo$^{a}$\footnote{
Corresponding author.}
\setcounter{footnote}{-1}
\footnote{
\emph{E-mail addresses:} gzh635516986@163.com (Z.H. Geng), walkerwzw@163.com (Z.W. Wang), leewj1375@163.com (W.J. Li), lixueying19990704@163.com (X.Y. Li), qianqian\_chen@yeah.net (Q.Q. Chen), jimingguo@hotmail.com (J.-M. Guo).} }\\[2mm]
\parbox[t]{0.96\textwidth}{%
\centering\fontsize{10}{12}\selectfont\itshape
\setlength{\parindent}{0pt}\setlength{\parskip}{4pt}
\textsuperscript{a}\,Department of Mathematics, East China University of Science and Technology,\\
Shanghai 200237, China\endgraf
\textsuperscript{b}\,Changzhou Institute of Technology, Changzhou 213032, China\endgraf
\textsuperscript{c}\,Suzhou Chien-Shiung Institute of Technology, Taicang 215411, China\endgraf
\textsuperscript{d}\,School of Mathematics and Statistics, Yancheng Teachers University,\\
Yancheng 224002, China}}

\date{}
\maketitle

{\flushleft\large\bf Abstract}\\
The problem posed by Hoffman, namely, to characterize graphs whose second largest adjacency eigenvalue is at most \(1\), was first studied by D.~M.~Cvetkovi\'{c} in 1982 [\textit{Publ. Inst. Math. (Beograd)} \textbf{31} (1982), 15--20].
In this paper, we completely characterize all connected planar graphs whose second largest eigenvalue does not exceed $1$.\\
\vspace{0.1cm}

\noindent {\bf Keywords:} Planar graph; adjacency matrix; the second largest eigenvalue\\[0.5\baselineskip]

\noindent \textbf{AMS Classification:} 05C50

\section{Introduction}

In this paper, all graphs are simple. 
For a graph $G$, the eigenvalues of $G$ are defined as those of its adjacency matrix $A(G)$. 
Let $\lambda_2(G)$ denote the second largest eigenvalue of  $G$.
As usual, $P_n$, $C_n$, and $K_n$ denote the path, cycle, and complete graph on $n$ vertices, respectively.

The characterization of graphs with second largest adjacency eigenvalue at most $1$ was originally proposed as a problem by Hoffman.
The problem has been extensively studied in spectral graph theory.
The first result in this direction appeared in a paper of D.~M.~Cvetkovi\'{c} in 1982~\cite{Cv}.
Many graph classes with $\lambda_2 \leq 1$ have been characterized, including trees~\cite{Neu}, unicyclic graphs~\cite{Xu1}, bicyclic graphs~\cite{Guo}, tricyclic graphs~\cite{Li}, bipartite graphs~\cite{PetM}, regular graphs and coronas~\cite{Zor}, biregular graphs with three distinct eigenvalues~\cite{CGGK16}, connected graphs with three distinct eigenvalues~\cite{Ch}, line graphs~\cite{Pet}, generalized line graphs~\cite{Petr}, generalized $\theta$-graphs~\cite{Gao}, $(K_{1,3},K_5-e)$-free graphs~\cite{Liu}, $(K_4-e)$-free graphs~\cite{Liu1}, connected claw-free graphs~\cite{Chen}, regular graphs with $n\geq2r$~\cite{Wang}, and $K_4$-minor-free graphs~\cite{CK}.
Xu and Shao~\cite{Xu} characterized non-bipartite graphs with girth at least $4$ and $\lambda_2<1$, while Stani\'{c}~\cite{Stan} characterized nested split graphs with $\lambda_2<1$.
Zhang et al.~\cite{Zhang} classified $2$-partially distance-regular graphs for which every local graph has second largest eigenvalue at most $1$.

A graph $G$ is said to be planar if it does not contain a subdivision of $K_{3,3}$ or $K_5$ as a subgraph.
Recently, the problem has also been studied for planar graphs. In particular, Li et al. \cite{Lis} characterized outerplanar graphs with the second largest eigenvalue at most $1$.
In our previous work \cite{Geng},  we determined all connected planar graphs that are $T_{1,2}$-free and whose second largest eigenvalue does not exceed $1$, where $T_{1,2}$ denotes the tree of diameter $3$ with exactly three pendant vertices. 

However, these results apply only to specific subclasses of planar graphs, and a complete characterization of connected planar graphs with $\lambda_2\leq 1$ is still lacking.

In this paper, we completely characterize all connected planar graphs whose second largest eigenvalue does not exceed $1$. 
Since trees, unicyclic graphs, bicyclic graphs, tricyclic graphs, and outerplanar graphs are all subclasses of planar graphs, our result extends and unifies several previous characterizations of graphs with $\lambda_2 \leq 1$.

The proof is based on a detailed structural analysis combined with eigenvalue techniques and forbidden subgraph arguments.

To state our main result, we first introduce several families of graphs 
that will serve as extremal configurations. These graphs are denoted by 
$A_1,\dots,A_{91}$, $H_1,\dots,H_{153}$, together with two additional graphs
$D_1$ and $D_2$ shown in Figs.~\ref{fig-3}, \ref{fig-4}, and \ref{F1}, respectively. 
Let $d(G)$ denote the diameter of $G$.
The clique number $\omega(G)$ of $G$ is the maximum number of vertices in a clique of $G$.
Let $\mathcal{G}$ be the set of connected planar graphs whose second largest eigenvalue is at most $1$, and let
$\mathcal{G}_i=\{G\in\mathcal{G}:\omega(G)=i\text{ and }d(G)\leq 3\}$.


\begin{thm}\label{thm-6}
$G \in \mathcal{G}$ if and only if $G$ is a connected induced subgraph of a graph in
$\{A_{1},\ldots,A_{91},H_{1},\ldots,H_{153},D_1,D_2\}$.
\end{thm}

The rest of the paper is organized as follows. 
In Section $2$, we present some preliminary results. 
In Sections $3$ and $4$, we characterize the graphs according to their clique numbers under the additional assumption that $G$ has diameter less than $4$. 
In Section $5$, we handle the remaining case where $G$ has diameter at least $4$.
Finally, in Section $6$, we complete the proof of Theorem $1.1$.

\section{Preliminaries}\label{s-2}

Throughout this section, suppose that $G$ has $n$ vertices and $m$ edges, with vertex set $V(G)$ and edge set $E(G)$.
For two vertices $x_{1},x_{2}\in V(G)$, write $x_{1}\sim x_{2}$ if they are adjacent and $x_{1}\nsim x_{2}$ otherwise.
For a vertex set $X\subseteq V(G)$, write
$N(X):=\{u\in V(G)\setminus X:\ \exists x\in X \text{ such that } u\sim x\}$,
and call $N(X)$ the set of neighbors of $X$.
In particular, for a vertex $v\in V(G)$, we write $N(v)=N(\{v\})$.
A matching is a set of edges in which no two edges share a common vertex; that is, the edges are pairwise vertex-disjoint.
The degree of a vertex $x$ is denoted by $d(x)$.
The girth $g(G)$ of a graph $G$ is the length of a shortest cycle in $G$, with $g(G)=\infty$ if $G$ is acyclic.
For $U\subseteq V(G)$, let $G[U]$ denote the induced subgraph of $G$ induced by $U$, and for simplicity, write $E(G[U])=E(U)$.
For graphs $H$ and $G$, we write $H\subset G$ if $H$ is an induced subgraph of $G$ and write $H\not\subset G$ otherwise.
For two disjoint subsets $V_1,V_2\subseteq V(G)$, let $E(V_1,V_2)=\{uv\in E(G):u\in V_1,\ v\in V_2\}$.
When $V_1=\{u\}$, we write $E(u,V_2)$ instead of $E(\{u\},V_2)$.
Moreover, $\left|E(V_1,V_2)\right|$ denotes the number of edges between $V_{1}$ and $V_{2}$.
Let $G[E(V_1,V_2)]$ denote the subgraph of $G$ with vertex set $V_{1}\cup V_{2}$ and edge set $E(V_1,V_2)$.
Whenever we write $\left|E(A,B)\right|\geq 1$, it is understood that both $A$ and $B$ are nonempty.

The union of two vertex-disjoint graphs $G_1$ and $G_2$, denoted by $G_1\cup G_2$, is the graph with $V(G_1\cup G_2)=V(G_1)\cup V(G_2)$ and $E(G_1\cup G_2)=E(G_1)\cup E(G_2)$.
For a positive integer $a$, we denote by $aG$ the union of $a$ copies of a graph $G$.
The join of two disjoint graphs $G_1$ and $G_2$ is denoted by $G_1 \nabla G_2$ and is the graph obtained from their disjoint union by adding all edges between $V(G_1)$ and $V(G_2)$.

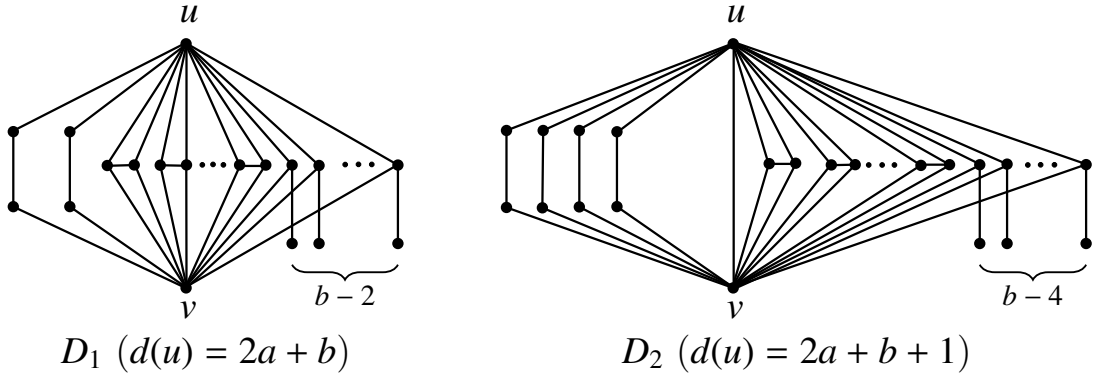
\begin{figure}[H]
\centering
\begin{tikzpicture}[x=0.70mm, y=0.80mm, inner xsep=0pt, inner ysep=0pt, outer xsep=0pt, outer ysep=0pt]
\path[line width=0mm] (74.29,204.25) rectangle +(208.72,65.32);
\definecolor{L}{rgb}{0,0,0}
\definecolor{F}{rgb}{0,0,0}
\path[line width=0.30mm, draw=L, fill=F] (109.93,260.01) circle (0.60mm);
\path[line width=0.30mm, draw=L, fill=F] (109.93,219.66) circle (0.60mm);
\path[line width=0.45mm, draw=L, dash pattern=on 2.00mm off 1.00mm] (105.32,251.02);
\draw(113.62,207.33) node[anchor=base]{\fontsize{14.23}{17.07}\selectfont $D_1$ ($d(u)=2a+b$)};
\draw(225.21,207.33) node[anchor=base]{\fontsize{14.23}{17.07}\selectfont $D_2$ ($d(u)=2a+b+1$)};
\path[line width=0.15mm, draw=L, fill=F] (116.98,240.00) circle (0.30mm);
\path[line width=0.15mm, draw=L, fill=F] (115.02,240.07) circle (0.30mm);
\path[line width=0.15mm, draw=L, fill=F] (113.00,240.00) circle (0.30mm);
\draw(108.75,263.64) node[anchor=base west]{\fontsize{14.23}{17.07}\selectfont $u$};
\draw(108.75,214.76) node[anchor=base west]{\fontsize{14.23}{17.07}\selectfont $v$};
\path[line width=0.30mm, draw=L, fill=F] (100.14,239.98) circle (0.60mm);
\path[line width=0.30mm, draw=L] (100.14,240.12) -- (109.95,260.15);
\path[line width=0.30mm, draw=L] (100.28,239.98) -- (109.82,219.95);
\path[line width=0.30mm, draw=L, fill=F] (120.04,239.98) circle (0.60mm);
\path[line width=0.30mm, draw=L] (120.04,239.98) -- (109.82,260.01);
\path[line width=0.30mm, draw=L] (120.17,239.98) -- (109.82,219.81);
\path[line width=0.30mm, draw=L, fill=F] (124.99,239.96) circle (0.60mm);
\path[line width=0.30mm, draw=L] (124.96,239.94) -- (109.84,260.18);
\path[line width=0.30mm, draw=L] (125.09,239.82) -- (119.97,239.82);
\path[line width=0.30mm, draw=L] (124.84,239.82) -- (109.84,219.82);
\path[line width=0.30mm, draw=L, fill=F] (95.04,239.89) circle (0.60mm);
\path[line width=0.30mm, draw=L] (94.96,239.98) -- (109.79,260.00);
\path[line width=0.30mm, draw=L] (95.04,239.98) -- (109.71,219.97);
\path[line width=0.30mm, draw=L, fill=F] (129.96,239.97) circle (0.60mm);
\path[line width=0.30mm, draw=L] (129.96,239.97) -- (110.09,259.84);
\path[line width=0.30mm, draw=L] (129.96,239.97) -- (109.96,219.98);
\path[line width=0.30mm, draw=L, fill=F] (149.95,239.97) circle (0.60mm);
\path[line width=0.30mm, draw=L] (149.95,239.97) -- (109.96,260.09);
\path[line width=0.30mm, draw=L] (150.07,240.09) -- (109.96,219.85);
\path[line width=0.15mm, draw=L, fill=F] (139.99,240.05) circle (0.30mm);
\path[line width=0.15mm, draw=L, fill=F] (145.01,240.07) circle (0.30mm);
\path[line width=0.15mm, draw=L, fill=F] (142.50,240.08) circle (0.30mm);
\path[line width=0.30mm, draw=L] (129.96,227.00) -- (129.96,239.97);
\path[line width=0.30mm, draw=L] (149.95,227.00) -- (149.95,239.97);
\path[line width=0.30mm, draw=L] (94.97,239.84) -- (100.22,239.97);
\path[line width=0.30mm, draw=L, fill=F] (87.97,245.47) circle (0.60mm);
\path[line width=0.30mm, draw=L] (87.97,245.47) -- (87.97,233.10);
\path[line width=0.30mm, draw=L, fill=F] (87.97,233.10) circle (0.60mm);
\path[line width=0.30mm, draw=L] (87.97,233.10) -- (109.96,219.73);
\path[line width=0.30mm, draw=L] (87.97,245.34) -- (109.84,259.84);
\path[line width=0.30mm, draw=L, fill=F] (213.29,260.01) circle (0.60mm);
\path[line width=0.30mm, draw=L, fill=F] (213.29,219.66) circle (0.60mm);
\path[line width=0.45mm, draw=L, dash pattern=on 2.00mm off 1.00mm] (208.68,251.02);
\draw(212.11,263.64) node[anchor=base west]{\fontsize{14.23}{17.07}\selectfont $u$};
\draw(212.11,214.76) node[anchor=base west]{\fontsize{14.23}{17.07}\selectfont $v$};
\path[line width=0.30mm, draw=L] (220.02,240.02) -- (213.18,260.01);
\path[line width=0.30mm, draw=L] (219.99,240.06) -- (213.18,219.81);
\path[line width=0.30mm, draw=L] (225.09,239.92) -- (213.20,260.18);
\path[line width=0.30mm, draw=L] (225.28,240.16) -- (220.15,240.16);
\path[line width=0.30mm, draw=L] (225.02,240.02) -- (213.20,219.82);
\path[line width=0.15mm, draw=L, fill=F] (238.89,239.98) circle (0.30mm);
\path[line width=0.15mm, draw=L, fill=F] (243.75,239.98) circle (0.30mm);
\path[line width=0.15mm, draw=L, fill=F] (241.43,239.98) circle (0.30mm);
\path[line width=0.30mm, draw=L, fill=F] (191.33,245.47) circle (0.60mm);
\path[line width=0.30mm, draw=L] (191.33,245.47) -- (191.33,233.10);
\path[line width=0.30mm, draw=L, fill=F] (191.33,233.10) circle (0.60mm);
\path[line width=0.30mm, draw=L] (191.33,233.10) -- (213.32,219.73);
\path[line width=0.30mm, draw=L] (191.33,245.34) -- (213.20,259.84);
\path[line width=0.30mm, draw=L, fill=F] (170.46,245.68) circle (0.60mm);
\path[line width=0.30mm, draw=L] (170.46,245.68) -- (170.46,232.96);
\path[line width=0.30mm, draw=L, fill=F] (170.46,232.96) circle (0.60mm);
\path[line width=0.30mm, draw=L] (170.46,232.96) -- (213.24,219.67);
\path[line width=0.30mm, draw=L] (170.46,245.68) -- (213.44,260.14);
\path[line width=0.30mm, draw=L, fill=F] (177.39,245.68) circle (0.60mm);
\path[line width=0.30mm, draw=L] (177.39,245.68) -- (177.24,233.05);
\path[line width=0.30mm, draw=L, fill=F] (177.24,232.90) circle (0.60mm);
\path[line width=0.30mm, draw=L] (177.24,232.90) -- (213.20,219.83);
\path[line width=0.30mm, draw=L] (177.39,245.68) -- (213.35,259.95);
\path[line width=0.30mm, draw=L, fill=F] (184.22,245.68) circle (0.60mm);
\path[line width=0.30mm, draw=L] (184.22,245.68) -- (184.22,233.05);
\path[line width=0.30mm, draw=L, fill=F] (184.22,233.05) circle (0.60mm);
\path[line width=0.30mm, draw=L] (184.22,233.05) -- (213.05,219.98);
\path[line width=0.30mm, draw=L] (184.22,245.83) -- (213.35,259.80);
\path[line width=0.30mm, draw=L] (213.38,259.93) -- (213.38,219.71);
\path[line width=0.30mm, draw=L, fill=F] (248.76,239.91) circle (0.60mm);
\path[line width=0.30mm, draw=L] (248.76,239.91) -- (213.39,259.90);
\path[line width=0.30mm, draw=L] (248.90,240.03) -- (213.39,219.96);
\path[line width=0.30mm, draw=L, fill=F] (254.16,240.03) circle (0.60mm);
\path[line width=0.30mm, draw=L] (254.16,240.03) -- (213.39,260.21);
\path[line width=0.30mm, draw=L] (254.37,240.03) -- (213.39,219.86);
\path[line width=0.30mm, draw=L] (248.69,240.03) -- (254.26,240.03);
\path[line width=0.30mm, draw=L, fill=F] (259.87,240.02) circle (0.60mm);
\path[line width=0.30mm, draw=L] (259.87,240.11) -- (213.40,260.25);
\path[line width=0.30mm, draw=L] (260.05,240.02) -- (213.40,219.78);
\path[line width=0.15mm, draw=L, fill=F] (268.95,239.98) circle (0.30mm);
\path[line width=0.15mm, draw=L, fill=F] (273.98,240.01) circle (0.30mm);
\path[line width=0.15mm, draw=L, fill=F] (271.37,240.01) circle (0.30mm);
\path[line width=0.30mm, draw=L, fill=F] (279.96,240.04) circle (0.60mm);
\path[line width=0.30mm, draw=L] (280.04,240.09) -- (213.46,260.26);
\path[line width=0.30mm, draw=L] (280.04,240.09) -- (213.35,219.81);
\path[line width=0.30mm, draw=L] (259.87,227.00) -- (259.87,240.02);
\path[line width=0.30mm, draw=L] (279.96,227.00) -- (279.96,240.04);
\path[line width=0.30mm, draw=L, fill=F] (77.29,245.49) circle (0.60mm);
\path[line width=0.30mm, draw=L] (77.29,245.49) -- (77.29,233.09);
\path[line width=0.30mm, draw=L, fill=F] (77.29,233.09) circle (0.60mm);
\path[line width=0.30mm, draw=L] (77.29,233.09) -- (109.99,219.60);
\path[line width=0.30mm, draw=L] (77.29,245.63) -- (109.99,260.21);
\path[line width=0.30mm, draw=L, fill=F] (149.95,227.00) circle (0.60mm);
\path[line width=0.30mm, draw=L, fill=F] (279.96,227.00) circle (0.60mm);
\path[line width=0.30mm, draw=L, fill=F] (259.87,227.00) circle (0.60mm);
\path[line width=0.30mm, draw=L, fill=F] (110.01,239.99) circle (0.60mm);
\path[line width=0.30mm, draw=L, fill=F] (105.05,239.94) circle (0.60mm);
\path[line width=0.30mm, draw=L] (105.05,239.94) -- (110.01,239.85);
\path[line width=0.30mm, draw=L] (110.01,239.85) -- (110.01,259.93);
\path[line width=0.30mm, draw=L] (110.01,259.93) -- (105.09,239.81);
\path[line width=0.30mm, draw=L] (105.09,239.81) -- (109.97,220.14);
\path[line width=0.30mm, draw=L] (110.01,219.91) -- (109.97,239.85);
\path[line width=0.30mm, draw=L, fill=F] (135.05,239.90) circle (0.60mm);
\path[line width=0.30mm, draw=L] (135.05,239.90) -- (109.97,260.01);
\path[line width=0.30mm, draw=L] (135.14,239.72) -- (110.01,219.69);
\path[line width=0.30mm, draw=L] (135.05,239.90) -- (135.05,227.00);
\path[line width=0.30mm, draw=L, fill=F] (135.05,227.00) circle (0.60mm);
\path[line width=0.30mm, draw=L, fill=F] (129.96,227.00) circle (0.60mm);
\path[line width=0.30mm, draw=L, fill=F] (231.83,239.99) circle (0.60mm);
\path[line width=0.30mm, draw=L, fill=F] (236.34,240.03) circle (0.60mm);
\path[line width=0.30mm, draw=L] (236.34,240.03) -- (231.88,239.95);
\path[line width=0.30mm, draw=L] (231.88,239.95) -- (213.34,259.79);
\path[line width=0.30mm, draw=L] (213.34,259.79) -- (236.43,240.07);
\path[line width=0.30mm, draw=L] (236.43,240.07) -- (213.30,219.82);
\path[line width=0.30mm, draw=L] (213.30,219.90) -- (231.71,239.99);
\path[line width=0.30mm, draw=L, fill=F] (220.09,240.20) circle (0.60mm);
\path[line width=0.30mm, draw=L, fill=F] (225.09,240.27) circle (0.60mm);
\path[line width=0.30mm, draw=L, fill=F] (265.03,240.12) circle (0.60mm);
\path[line width=0.30mm, draw=L] (265.03,240.12) -- (213.53,259.91);
\path[line width=0.30mm, draw=L] (265.10,239.92) -- (213.32,219.44);
\path[line width=0.30mm, draw=L] (265.03,240.12) -- (265.03,227.00);
\path[line width=0.30mm, draw=L, fill=F] (265.03,227.00) circle (0.60mm);
\tikzset{range brace/.style={decorate,decoration={brace,mirror,amplitude=2mm},line width=.25mm,preaction={draw=white,line width=1.1mm}},range label/.style={fill=white,inner sep=1pt,font=\small}}
\draw[range brace] (129.96,223.5) -- (149.95,223.5) node[midway,below=2.4mm,range label] {$b-2$};
\draw[range brace] (259.87,223.5) -- (279.96,223.5) node[midway,below=2.4mm,range label] {$b-4$};
\end{tikzpicture}%
\caption{The graphs $D_{1}$ and $D_{2}$.}\label{F1}
\end{figure}

Let $a$ and $b$ be two sufficiently large integers.
Denote by $D_1$ (or $D_2$, respectively) the graph consisting of $\{u,v\}\nabla aK_2$ and $bP_2$,
where $u$ and $v$ are non-adjacent (or adjacent, respectively),
for $(b-2)P_2$ (or $(b-4)P_2$, respectively), one of the vertices in each $P_2$ is adjacent to both $u$ and $v$,
and for the remaining $2P_2$ (or $4P_2$, respectively), one vertex in each $P_2$ is adjacent to $u$,
and the other is adjacent to $v$ (see Fig.~\ref{F1}).

\begin{lem}[\cite{Rap}]\label{lem-2}
Let $A$ be a real symmetric $n \times n$ matrix and let $B$ be a principal submatrix of $A$ of order $k$. If $\lambda_1 \geq \cdots \geq \lambda_n$ and $\mu_1 \geq \cdots \geq \mu_{k}$ are the eigenvalues of $A$ and $B$, respectively, then
$$\lambda_i \geq \mu_i \geq \lambda_{n-k+i}, \qquad i=1,2,\ldots,k.$$
\end{lem}

The following corollary follows immediately from Lemma \ref{lem-2}.

\begin{cor}\label{cor-subgraph}
If $H$ is an induced subgraph of $G$ with $\lambda_2(H)>1$, then $\lambda_2(G)>1$.
\end{cor}

Suppose that $A$ is a symmetric real matrix whose rows and columns are indexed by
$X = \{1,\ldots,n\}$. Let ${X_1,\ldots,X_m}$ be a partition of $X$.
Let $A=(A_{i,j})$ be partitioned according to ${X_1,\ldots,X_m}$,
where $A_{i,j}$ denotes the submatrix (block) of $A$ formed by rows in $X_i$ and the columns
in $X_j$.
Let $b_{i,j}$ denote the average row sum of $A_{i,j}$.
Then the matrix $B=(b_{i,j})$ is called the quotient matrix of $A$.
If the row sum of each block $A_{i,j}$ is constant then the partition is called equitable.
For an equitable partition, the following standard result is useful.

\begin{lem}[\cite{Brouwer}]\label{equitable}
Let $A$ be a symmetric matrix with an equitable partition $\pi$. Then the spectrum of the quotient matrix $B_{\pi}$ is contained in the spectrum of $A$.
\end{lem}

The graphs $R_{i}$ $(i=1,2,\ldots,98)$ are shown in Fig.~\ref{fig-2}.
By direct computation, we obtain the following lemma.

\begin{lem}\label{lem-1}
For each graph $R_{i}$, $i=1,2,\ldots,98$, we have $\lambda_{2}(R_{i})> 1$.
\end{lem}

\begin{lem}\label{lem-3}
If $G \in \left\{A_{1},\ldots, A_{91}, H_{1},\ldots, H_{153}, {D}_1, {D}_2\right\}$, then $\lambda_{2}(G)\leq 1$.
\end{lem}
\begin{proof}
Direct computation gives $\lambda_2(A_i)\leq 1$ for $i=1,\ldots,91$ and $\lambda_2(H_i)\leq 1$ for $i=1,\ldots,153$.
Choose a sufficiently large integer $c\geq\max\{a,b\}$. For $i=1,2$, the graph $D_i$ with parameters $a,b$ is an induced subgraph of the corresponding graph with parameters $c,c$. Thus, by Lemma \ref{lem-2}, it suffices to consider the case $a=b$, with $a$ sufficiently large.
After several elementary row transformations of $A(D_1)-I$, the resulting matrix contains at least $2a$ zero rows.
Similarly, we get at least $2a-2$ (or 1, respectively) zero rows after the elementary row transformation of $A(D_1)+I$ (or $A(D_1)+2I$, respectively).
Thus, $1$, $-1$, and $-2$ are eigenvalues of $A(D_1)$ with multiplicities at least $2a$, $2a-2$, and $1$, respectively.
Consider the equitable partition $\pi$ of $V(D_1)$ given by $U_1=\{u,v\}$, $U_2=\{w:d(w)=1\}$, $U_3=\{w:d(w)=3\text{ and }w\text{ is adjacent to some vertex in }U_2\}$, $U_4=\{w:d(w)=3\text{ and }w\text{ is adjacent to no vertex in }U_2\}$, and $U_5=V(D_1)\setminus(U_1\cup U_2\cup U_3\cup U_4)$.
Then the corresponding quotient matrix of $A(D_1)$ is
$$
B_{\pi}^1=
\left(
\begin{array}{ccccc}
  0 & 0 & a-2 & 2a & 2 \\
  0 & 0 & 1 & 0 & 0 \\
  2 & 1 & 0 & 0 & 0 \\
  2 & 0 & 0 & 1 & 0 \\
  1 & 0 & 0 & 0 & 1
\end{array}
\right).
$$

The characteristic polynomial of $B_{\pi}^1$ is
$$\Phi(B_{\pi}^1,\lambda)=(\lambda-1)^2(\lambda^3-(6a-1)\lambda-(4a+2))\triangleq (\lambda-1)^2g(\lambda).$$
Combining these eigenvalue multiplicities with Lemma \ref{equitable} gives
$$\Phi(A(D_1),\lambda)=(\lambda-1)^{2a}(\lambda+1)^{2a-2}(\lambda+2)g(\lambda).$$
For sufficiently large $a$, we have $g(-a)<0<g(-2)$, $g(-1)>0>g(1)$, and $g(1)<0<g(a)$, where $g(-2)=8a-12$, $g(-1)=2a-4$, and $g(1)=-10a$. Hence the three roots of $g$ lie in $(-a,-2)$, $(-1,1)$, and $(1,a)$, respectively, and therefore $\lambda_2(D_1)=1$.

For $D_2$, the same argument shows that $1$, $-1$, and $-3$ are eigenvalues of $A(D_2)$ with multiplicities at least $2a$, $2a-2$, and $1$, respectively.
Consider the analogous equitable partition $\pi$ of $V(D_2)$.
Then the corresponding quotient matrix of $A(D_2)$ is
$$
B_{\pi}^2=
\left(
\begin{array}{ccccc}
  1 & 0 & a-4 & 2a & 4 \\
  0 & 0 & 1 & 0 & 0 \\
  2 & 1 & 0 & 0 & 0 \\
  2 & 0 & 0 & 1 & 0 \\
  1 & 0 & 0 & 0 & 1
\end{array}
\right).
$$
The characteristic polynomial of $B_{\pi}^2$ is
$$\Phi(B_{\pi}^2,\lambda)=(\lambda-1)^2(\lambda^3-\lambda^2-(6a-3)\lambda-(4a+3))\triangleq (\lambda-1)^2h(\lambda).$$
Combining these eigenvalue multiplicities with Lemma \ref{equitable} gives
$$\Phi(A(D_2),\lambda)=(\lambda-1)^{2a}(\lambda+1)^{2a-2}(\lambda+3)h(\lambda).$$
For sufficiently large $a$, we have $h(-a)<0<h(-3)$, $h(-1)>0>h(1)$, and $h(1)<0<h(a)$, where $h(-3)=14a-48$, $h(-1)=2a-8$, and $h(1)=-10a$. Hence the three roots of $h$ lie in $(-a,-3)$, $(-1,1)$, and $(1,a)$, respectively, and therefore $\lambda_2(D_2)=1$.
\end{proof}

Let $K_4-e$ denote the graph obtained from $K_4$ by deleting an edge.
By  Theorem 2.3 of~\cite{Liu1} and Lemma \ref{lem-3}, we obtain the following lemma.

\begin{lem}\label{lem-5}
If $G \in \mathcal{G}$ and $G$ is $K_{4}-e$-free, then $G$ is a connected induced subgraph of a graph in $\left\{A_{48},\ldots,A_{91},D_1,D_2\right\}$.
\end{lem}

\section{Planar graphs with $\omega(G)=4$ and $\lambda_{2}(G)\leq 1$}\label{s-3}

Let $K_{4}$ be a maximal clique of $G \in \mathcal{G}_4$, and $V(K_{4})=\{x_{0}, x_{1}, x_{2}, x_{3}\}$.
Denote by $T$ the set $V(G)\backslash V(K_{4})$.
Note that a vertex of $T$ is adjacent to at most $3$ vertices of $K_{4}$.
Thus, $T = T_{0} \cup T_{1} \cup T_{2} \cup T_{3}$, where $T_{i}$ $(i=0,1,2,3)$ is the set of vertices of $T$
which are adjacent to exactly $i$ vertices of $K_{4}$.
Also, we may suppose that $T_{k}=\bigcup\limits_{0\leq i_{1}< \cdots < i_{k}\leq 3} T_{x_{i_{1}}\ldots x_{i_{k}}}$ $(0\leq k \leq 3)$,
where $T_{x_{i_{1}}\ldots x_{i_{k}}}$ is the set of vertices which are adjacent precisely to the vertices $x_{i_{1}}, \ldots, x_{i_{k}}$.

For the notation used below, let $i,j,k,z$ be distinct integers in $\{0,1,2,3\}$, and label the vertices of each class as follows:
$T_0=\{t_r\mid 1\le r\le |T_0|\}$,
$T_{x_i}=\{p_{i,r}\mid 1\le r\le |T_{x_i}|\}$,
$T_{x_j}=\{p_{j,r}\mid 1\le r\le |T_{x_j}|\}$,
$T_{x_k}=\{p_{k,r}\mid 1\le r\le |T_{x_k}|\}$,
$T_{x_z}=\{p_{z,r}\mid 1\le r\le |T_{x_z}|\}$,
$T_{x_ix_j}=\{q_{ij,r}\mid 1\le r\le |T_{x_ix_j}|\}$,
$T_{x_ix_k}=\{q_{ik,r}\mid 1\le r\le |T_{x_ix_k}|\}$,
$T_{x_ix_z}=\{q_{iz,r}\mid 1\le r\le |T_{x_ix_z}|\}$,
$T_{x_jx_k}=\{q_{jk,r}\mid 1\le r\le |T_{x_jx_k}|\}$,
$T_{x_jx_z}=\{q_{jz,r}\mid 1\le r\le |T_{x_jx_z}|\}$,
$T_{x_kx_z}=\{q_{kz,r}\mid 1\le r\le |T_{x_kx_z}|\}$,
$T_{x_ix_jx_z}=\{w_{ijz,r}\mid 1\le r\le |T_{x_ix_jx_z}|\}$,
$T_{x_ix_kx_z}=\{w_{ikz,r}\mid 1\le r\le |T_{x_ix_kx_z}|\}$,
$T_{x_ix_jx_k}=\{w_{ijk,r}\mid 1\le r\le |T_{x_ix_jx_k}|\}$, and
$T_{x_jx_kx_z}=\{w_{jkz,r}\mid 1\le r\le |T_{x_jx_kx_z}|\}$.
If the upper bound is zero, the corresponding vertex class is understood to be empty.

Before proving Theorem \ref{thm-31}, the main theorem of this section, we establish several structural properties of graphs in $\mathcal{G}_4$. 
These properties will play a crucial role in the proof.

\begin{lem}\label{lem-9}
Let $K_{4}$ be a clique of $G \in \mathcal{G}_4$, $V(K_{4})=\{x_{0},x_{1},x_{2},x_{3}\}$, and  $i,j,k,z$ be distinct integers in $\{0,1,2,3\}$.
Then the following statements hold:
\begin{enumerate}[label=(\arabic*),leftmargin=2em,itemsep=2pt]

\item\label{lem-9:item-1} $\left|E(T_{x_{i}})\right|=0$, $\left|T_{x_{i}}\right| \leq 4$,  $\left|E(T_{3})\right|=0$ and $\left|T_{x_{i}x_{j}x_{k}}\right| \leq 1$;

\item\label{lem-9:item-2} If $w\in T$, then $\left|N(w)\cap T_{x_{i}x_{j}}\right| \leq 1$; 
    if $u\in T_{x_{i}}\cup T_{x_{j}}$ and $q_{ij,1}\in T_{x_{i}x_{j}}$, then $u\nsim q_{ij,1}$;

\item\label{lem-9:item-3} $G[T_{x_{i}x_{j}}]=aK_{2}\cup (\left|T_{x_{i}x_{j}}\right|-2a)K_{1}$ ($a\geq 0$);

\item\label{lem-9:item-5} If $\left|T_{3}\right|\geq 1$ ($\left|T_{x_{i}x_{j}x_{k}}\right|=1$), then $\left|E(T_{x_{i}},T_{x_{i}x_{j}x_{k}})\right| \leq 1$, $\left|E(T_{x_{z}},T_{x_{i}x_{j}x_{k}})\right|= 0$, $\left|T_{x_{z}}\right| \leq  1$, $\left|E(T_{1})\right| =0$ and $\left|E(T_{2},T_{0})\right| =0$. Furthermore, if $\left|T_{3}\right| \geq 2$, then  $\left|T_{0}\right|=0$  and $\left|E(T_{1},T_{3})\right|= 0$;

\item\label{lem-9:item-6} $\left|E(T_{0},T_{1})\right|= 0$ and $G[E(T_{0},T_{2}\cup T_{3})]\subset aK_{2}$ ($a\geq 0$). In particular, if $q_{ij,1},q_{ij,2}\in T_{x_{i}x_{j}}$ and $q_{ij,1} \sim q_{ij,2}$, then $\left|E(\{q_{ij,1},q_{ij,2}\},T_{0})\right| =0$;

\item\label{lem-9:item-11} $\left|E(T_{x_{i}x_{j}},T_{x_{j}x_{k}})\right| \leq 1$ and $\left|E(T_{x_{i}x_{j}},T_{x_{k}x_{z}})\right| =0$. Moreover, if $\left|E(T_{x_{i}x_{j}},T_{x_{j}x_{k}})\right| =1$, then $\left|T_{x_{k}x_{z}}\cup T_{x_{i}x_{z}}\right|=0$;

\item\label{lem-9:item-12} If $\left|T_{x_{i}x_{j}}\right| \geq 2$, then $\left|T_{x_{k}x_{z}}\cup T_{x_{i}x_{k}x_{z}}\cup T_{x_{j}x_{k}x_{z}}\right|=0$;

\item\label{lem-9:item-13} If  $u\in T_{x_{i}x_{j}x_{k}}\cup T_{x_{k}}\cup (T_{2}-T_{x_{i}x_{j}}- T_{x_{k}x_{z}})\triangleq U$, $q_{ij,1},q_{ij,2}\in T_{x_{i}x_{j}}$ and $q_{ij,1}\sim q_{ij,2}$, then 
$u$ is  exactly adjacent to  the  one of $q_{ij,1}$ and $q_{ij,2}$, say $q_{ij,1}$. Furthermore, if $v(\neq u)\in U$, then $v\sim q_{ij,1}$ and $v\nsim q_{ij,2}$;

\item\label{lem-9:item-15} $\left|E(T_{x_{i}x_{j}},T_{x_{i}x_{k}x_{z}}\cup T_{x_{j}x_{k}x_{z}})\right|=0$ and $\left|E(T_{x_{i}x_{j}},T_{x_{i}x_{j}x_{k}}\cup T_{x_{i}x_{j}x_{z}})\right| \leq 1$;

\item\label{lem-9:item-16} If $\left|T_{x_{i}x_{j}}\right|\geq 1$ and $\left|T_{x_{k}}\cup T_{x_{z}}\right|\geq 1$, then $\left|E(T_{1})\right| =0$, $\left|T_{x_{i}}\right|\leq 2$ and $\left|E(T_{2},T_{0})\right| =0$;

\item\label{lem-9:item-17} If $\left|E(T_{0},T_{x_{i}x_{j}x_{k}})\right| = 1$, then $\left|T_{x_{i}x_{z}}\cup T_{x_{j}x_{z}}\cup T_{x_{k}x_{z}}\right|=0$; 
if $\left|T_{x_{i}x_{j}}\right| \geq\left|E(T_{0},T_{x_{i}x_{j}x_{k}})\right| = 1$, then $\left|E(T_{x_{i}}\cup T_{x_{j}},T_{x_{i}x_{j}x_{k}})\right|=0$.
\end{enumerate}
\end{lem}
\begin{proof}

Statements \ref{lem-9:item-1}--\ref{lem-9:item-3} follow from (1), (2), (4), (6), and (7) of Lemma 5.1 in \cite{Geng}.
We prove the remaining statements by contradiction.


If $\left|E(T_{x_{i}},T_{x_{i}x_{j}x_{k}})\right| \geq 2$, then  $R_{36}\subset G$ by Lemma \ref{lem-9}\ref{lem-9:item-1}.
By Lemma \ref{lem-1} and Corollary \ref{cor-subgraph}, $\lambda_2(G)> 1$, a contradiction.
If $\left|E(T_{x_{z}},T_{x_{i}x_{j}x_{k}})\right|\geq 1$, then $G$ contains a subdivision of $K_{5}$, a contradiction.
If $\left|T_{x_{z}}\right| \geq  2$, then $R_{32}\subset G$ by Lemma \ref{lem-9}\ref{lem-9:item-1},  a contradiction.
To prove that $\left|E(T_{1})\right| =0$, suppose $\left|E(T_{1})\right| \geq 1$. 
Then there are two vertices $p_{i,1}\in T_{x_{i}}$, $p_{j,1}\in T_{x_{j}}$ and $p_{i,1}\sim p_{j,1}$ by  Lemma \ref{lem-9}\ref{lem-9:item-1}.
If $w_{ijk,1}\in T_{x_{i}x_{j}x_{k}}$, then $R_{9} \subset G$ (if $E(w_{ijk,1},\{p_{i,1},p_{j,1}\})=\varnothing$, then deleting $x_{j}$) or $R_{12} \subset G$ (if $E(w_{ijk,1},\{p_{i,1},p_{j,1}\})=\{w_{ijk,1}p_{j,1}\}$, then deleting $x_{j}$) or $R_{18} \subset G$ (if $\left|E(w_{ijk,1},\{p_{i,1},p_{j,1}\})\right|=2$, then deleting $x_{j}$), a contradiction.
If $w_{ikz,1}\in T_{x_{i}x_{k}x_{z}}$, then either $R_{6} \subset G$ (if $\left|E(w_{ikz,1},\{p_{i,1},p_{j,1}\})\right|=0$, then deleting $x_{j}$) or  $G$ contains a subdivision of $K_{5}$ (if $\left|E(w_{ikz,1},\{p_{i,1},p_{j,1}\})\right|\geq 1$), a contradiction.
Next, to prove that $\left|E(T_{2},T_{0})\right| =0$,  assume that  $q_{ij,1}\in T_{x_{i}x_{j}}$ and $\left|E(q_{ij,1},T_{0})\right| \geq 1$.
If $w_{ijk,1}\in T_{x_{i}x_{j}x_{k}}$, then either $R_{9} \subset G$ (if $\left|E(w_{ijk,1},q_{ij,1})\right|=0$, then deleting $x_{j}$) or $R_{14} \subset G$ (if $\left|E(w_{ijk,1},q_{ij,1})\right|=1$, then deleting $x_{j}$) by Lemma \ref{lem-9}\ref{lem-9:item-6}, a contradiction.
If $w_{ikz,1}\in T_{x_{i}x_{k}x_{z}}$, then either $R_{6} \subset G$ (if $\left|E(w_{ikz,1},q_{ij,1})\right|=0$, then deleting $x_{j}$) or $G$ contains a subdivision of $K_{5}$ (if $\left|E(w_{ikz,1},q_{ij,1})\right|=1$) by Lemma \ref{lem-9}\ref{lem-9:item-6}, a contradiction.
We now suppose that $\left|T_{3}\right| \geq 2$.
Let $\left|T_{x_{i}x_{j}x_{k}}\right| =\left|T_{x_{i}x_{j}x_{z}}\right| =1$ by Lemma \ref{lem-9}\ref{lem-9:item-1}.
If $t_{1}\in T_{0}$ and $\left|E(t_{1},T_{x_{i}x_{j}x_{k}})\right|= 1$, then $R_{14}\subset G$ (deleting $x_{i}$) by Lemma \ref{lem-9}\ref{lem-9:item-6}, a contradiction. 
Thus,  $\left|T_{0}\right|=0$ by symmetry and Lemma \ref{lem-9}\ref{lem-9:item-6}.
Suppose that $y\in T_{1}$ and $\left|E(y,T_{3})\right|\geq 1$.
Since $G$ contains no subdivision of $K_{5}$, $\left|E(y,T_{3})\right|= 1$.
If $y\in T_{x_{i}}$ and $\left|E(y,T_{x_{i}x_{j}x_{z}})\right|= 1$, then $R_{14}\subset G$ (deleting $x_{i}$), a contradiction.
If $y\in T_{x_{z}}$ and $\left|E(y,T_{x_{i}x_{j}x_{z}})\right|= 1$, then $R_{18}\subset G$ (deleting $x_{i}$), a contradiction.
Thus, by symmetry and Lemma \ref{lem-9}\ref{lem-9:item-5}, $\left|E(T_{1},T_{3})\right|= 0$.
Thus, \ref{lem-9:item-5} holds.

If $\left|E(T_{0},T_{1})\right|\neq 0$, then $R_{15}\subset G$, a contradiction.
Let $t_{1}\in T_{0}$ and $v\in T_{2}\cup T_{3}$.
To prove that $G[E(T_{0},T_{2}\cup T_{3})]\subset aK_{2}$ ($a\geq 0$), it suffices to show that  $\left|N(v)\cap T_{0}\right| \leq 1$ and $\left|N(t_{1})\cap (T_{2}\cup T_{3})\right| \leq 1$.
We first suppose that $\left|N(v)\cap T_{0}\right|  \geq2$ and $t_{2},t_{3}\in N(v)\cap T_{0}$.
If  $v\in T_{x_{i}x_{j}}$, then either $R_{4} \subset G$ (if $t_{2}\nsim t_{3}$, then deleting $x_{i}$) or $R_{5} \subset G$ (if $t_{2}\sim t_{3}$, then deleting $x_{i}$), a contradiction.
If  $v\in T_{x_{i}x_{j}x_{k}}$, then either $R_{13} \subset G$ (if $t_{2}\sim t_{3}$, then deleting $x_{i}$) or $R_{32} \subset G$ (if $t_{2}\nsim t_{3}$), a contradiction.
Thus, $\left|N(v)\cap T_{0}\right| \leq 1$.
It remains to prove that  $\left|N(t_{1})\cap (T_{2}\cup T_{3})\right|  \leq1$.
Suppose that  $\left|N(t_{1})\cap (T_{2}\cup T_{3})\right|  \geq2$ and let $y,z\in N(t_{1})\cap (T_{2}\cup T_{3})$.
We first consider the case when $y,z\in T_{2}$.
If $y,z\in T_{x_{i}x_{j}}$, then either $R_{8} \subset G$ (if $y\nsim z$, then deleting $x_{i}$) or $R_{13} \subset G$ (if $y\sim z$, then deleting $x_{i}$), a contradiction;
if $y\in T_{x_{i}x_{j}}$ and $z\in T_{x_{j}x_{k}}$, then either $R_{12} \subset G$ (if $y\nsim z$, then deleting $x_{i}$) or $R_{18} \subset G$ (if $y\sim z$, then deleting $x_{i}$), a contradiction;
if $y\in T_{x_{i}x_{j}}$ and $z\in T_{x_{k}x_{z}}$, then $G$ contains a subdivision of $K_{5}$, a contradiction.
We next consider the case when $y,z\in T_{3}$.
Then $G$ contains a subdivision of $K_{5}$ by  Lemma \ref{lem-9}\ref{lem-9:item-1}, a contradiction.
Finally, let $y\in T_{2}$ and $z\in T_{3}$.
If $y\in T_{x_{i}x_{j}}$ and $z\in  T_{x_{i}x_{j}x_{k}}$, then either $R_{12} \subset G$ (if $y\nsim z$, then deleting $x_{i}$) or $R_{18} \subset G$ (if $y\sim z$, then deleting $x_{i}$), a contradiction;
if $y\in T_{x_{i}x_{j}}$ and $z\in T_{x_{i}x_{k}x_{z}}$, then $G$ contains a subdivision of $K_{5}$, a contradiction.
Thus, $\left|N(t_{1})\cap (T_{2}\cup T_{3})\right| \leq 1$.
If  $\left|E(\{q_{ij,1},q_{ij,2}\},T_{0})\right| \geq 1$, then  $R_{10} \subset G$ (deleting $x_{j}$), a contradiction.
Thus, \ref{lem-9:item-6} holds.

If $\left|E(T_{x_{i}x_{j}},T_{x_{j}x_{k}})\right| \geq 2$ or $\left|E(T_{x_{i}x_{j}},T_{x_{k}x_{z}})\right| \geq 1$, then  $G$ contains a subdivision of $K_{5}$ or $K_{3,3}$ by  Lemma \ref{lem-9}\ref{lem-9:item-2}, a contradiction.
If $\left|E(T_{x_{i}x_{j}},T_{x_{j}x_{k}})\right| =1$ and $\left|T_{x_{k}x_{z}}\right|\geq 1$, then either $R_{18} \subset G$ (if $\left|E(T_{x_{j}x_{k}},T_{x_{k}x_{z}})\right|=0$, then deleting $x_{i}$) or $G$ contains a subdivision of $K_{5}$ (if $\left|E(T_{x_{j}x_{k}},T_{x_{k}x_{z}})\right| \geq 1$), a contradiction.
Thus, by symmetry, \ref{lem-9:item-11} holds.

If $\left|T_{x_{k}x_{z}}\right|\geq 1$, then either $R_{13} \subset G$ (if $\left|E(T_{x_{i}x_{j}})\right| \geq 1$, then deleting $x_{i}$) or $R_{34} \subset G$ (if $\left|E(T_{x_{i}x_{j}})\right| =0$) by Lemma \ref{lem-9}\ref{lem-9:item-11}, a contradiction.
If $\left|T_{x_{i}x_{k}x_{z}}\cup T_{x_{j}x_{k}x_{z}}\right|\geq 1$, then $R_{19}\subset G$ (if $\left|E(T_{x_{i}x_{j}})\right| \geq1$) or  $R_{36}\subset G$ (if $\left|E(T_{x_{i}x_{j}})\right|=0$) or $G$ contains a subdivision of $K_{5}$ (if $\left|E(T_{x_{i}x_{j}}, T_{x_{i}x_{k}x_{z}}\cup T_{x_{j}x_{k}x_{z}})\right| \geq1$), a contradiction.
Thus, \ref{lem-9:item-12} holds.

If $\left|E(u,\{q_{ij,1},q_{ij,2}\})\right|=0$, then either $R_{10}\subset G$ (if $u\in T_{x_{k}}\cup T_{2}-T_{x_{i}x_{j}}$, then deleting $x_{j}$) or $R_{16}\subset G$ (if $u\in T_{x_{i}x_{j}x_{k}}$, then deleting $x_{j}$), a contradiction.
Thus, $\left|E(u,\{q_{ij,1},q_{ij,2}\})\right| =1$ by Lemma \ref{lem-9}\ref{lem-9:item-2}.
Suppose that $v (\neq u)\in U$, $u\sim q_{ij,1}$ and $v\sim q_{ij,2}$.
If $u,v\in T_{x_{j}x_{k}}$, then $G$ contains a subdivision of $K_{3,3}$ (on vertex set $\{x_i,x_j,x_k,x_z,q_{ij,1},q_{ij,2},u,v\}$), a contradiction.
If $u\in T_{x_{j}x_{k}}$ and $v\in T_{2}-T_{x_{j}x_{k}}$, then  $G$ contains a subdivision of $K_{3,3}$ by Lemma \ref{lem-9}\ref{lem-9:item-12}, a contradiction.
If $u\in T_{x_{j}x_{k}}$ and $v\in T_{x_{i}x_{j}x_{k}}\cup T_{x_{k}}$, then  $G$ contains a subdivision of $K_{3,3}$, a contradiction.
If $u\in T_{x_{i}x_{j}x_{k}}$ and $v\in T_{x_{k}}$, then  $G$ contains a subdivision of $K_{5}$, a contradiction.
Thus,  by Lemma \ref{lem-9}\ref{lem-9:item-12}, $u,v$ are both adjacent to  the same  one of $\{q_{ij,1},q_{ij,2}\}$.
Thus, \ref{lem-9:item-13} holds.

If $\left|E(T_{x_{i}x_{j}},T_{x_{i}x_{k}x_{z}}\cup T_{x_{j}x_{k}x_{z}})\right|\geq 1$, then $G$ contains a subdivision of $K_{5}$, a contradiction.
Suppose that $\left|E(T_{x_{i}x_{j}},T_{x_{i}x_{j}x_{k}}\cup T_{x_{i}x_{j}x_{z}})\right| > 1$.
If $q_{ij,1},q_{ij,2}\in T_{x_{i}x_{j}}$, $w_{ijk,1}\in T_{x_{i}x_{j}x_{k}}$ and $w_{ijk,1}\sim q_{ij,1},q_{ij,2}$, then $G$ contains a subdivision of $K_{3,3}$ (deleting $x_{z}$), a contradiction.
If $q_{ij,1},q_{ij,2}\in T_{x_{i}x_{j}}$, $w_{ijk,1}\in T_{x_{i}x_{j}x_{k}}$, $w_{ijz,1}\in T_{x_{i}x_{j}x_{z}}$, $q_{ij,1}\sim w_{ijk,1}$ and $q_{ij,2}\sim w_{ijz,1}$, then either $R_{16} \subset G$ (if $q_{ij,1}\nsim q_{ij,2}$, then deleting $x_{j}$ and $x_{k}$) or $G$ contains a subdivision of $K_{5}$ (if $q_{ij,1}\sim q_{ij,2}$), a contradiction.
If $q_{ij,1}\in T_{x_{i}x_{j}}$, $w_{ijk,1}\in T_{x_{i}x_{j}x_{k}}$, $w_{ijz,1}\in T_{x_{i}x_{j}x_{z}}$ and $q_{ij,1}\sim w_{ijk,1},w_{ijz,1}$, then $G$ contains a subdivision of $K_{5}$, a contradiction.
Thus, by Lemma \ref{lem-9}\ref{lem-9:item-1}, \ref{lem-9:item-15} holds.

Suppose that $q_{ij,1}\in T_{x_{i}x_{j}}$ and $p_{k,1}\in T_{x_{k}}$.
If $p_{i,1}\in T_{x_{i}}$ and $p_{k,1}\sim p_{i,1}$, then $R_{12} \subset G$ (if $q_{ij,1}\sim p_{k,1}$, then deleting $x_{k}$  or if $q_{ij,1}\nsim p_{k,1}$, then deleting $x_{z}$) by Lemma \ref{lem-9}\ref{lem-9:item-2}, a contradiction.
If $p_{z,1}\in T_{x_{z}}$ and $p_{k,1}\sim p_{z,1}$, then either $R_{6} \subset G$ (if $\left|E(q_{ij,1},\{p_{k,1},p_{z,1}\})\right|=0$, then deleting $x_{k}$)  or $G$ contains a subdivision of $K_{5}$ (if $\left|E(q_{ij,1},\{p_{k,1},p_{z,1}\})\right|\geq 1$), a contradiction.
Thus, by symmetry and Lemma \ref{lem-9}\ref{lem-9:item-1}, $\left|E(T_{1})\right| =0$.
If $\left|T_{x_{i}}\right|\geq 3$, then either $R_{38} \subset G$ (if $q_{ij,1}\sim p_{k,1}$, then deleting $x_{z}$) or $R_{51} \subset G$ (if $q_{ij,1}\nsim p_{k,1}$) by Lemma \ref{lem-9}\ref{lem-9:item-1} and \ref{lem-9:item-2}, a contradiction.
If $\left|E(T_{0},q_{ij,1})\right|\geq 1$, then either $R_{3} \subset G$ (if $q_{ij,1}\nsim p_{k,1}$, then deleting $x_{j}$) or $R_{7} \subset G$ (if $q_{ij,1}\sim p_{k,1}$, then deleting $x_{j}$) by Lemma \ref{lem-9}\ref{lem-9:item-6}, a contradiction.
Thus, \ref{lem-9:item-16} holds.

If  $\left|T_{x_{k}x_{z}}\right|\geq 1$, then either $R_{14} \subset G$ (if $\left|E(T_{x_{k}x_{z}},T_{x_{i}x_{j}x_{k}})\right|=0$, then deleting $x_{i}$) by Lemma \ref{lem-9}\ref{lem-9:item-6} or $G$ contains a subdivision of $K_{5}$ (if $\left|E(T_{x_{k}x_{z}},T_{x_{i}x_{j}x_{k}})\right|\geq 1$), a contradiction.
If   $\left|E(T_{x_{i}},T_{x_{i}x_{j}x_{k}})\right|\geq 1$, then either $R_{25} \subset G$ (if $\left|E(T_{x_{i}x_{j}x_{k}},T_{x_{i}x_{j}})\right|\geq 1$, then deleting $x_{i}$) or $R_{28} \subset G$ (if $\left|E(T_{x_{i}x_{j}x_{k}},T_{x_{i}x_{j}})\right|=0$, then deleting $x_{k}$) by Lemma \ref{lem-9}\ref{lem-9:item-1}\ref{lem-9:item-2} and \ref{lem-9:item-6}, a contradiction.
Thus,  by symmetry, \ref{lem-9:item-17} holds.
\end{proof}

\input{generated_figures/R_figures.tex}

\begin{lem}\label{lem-10}
Suppose that $G \in \mathcal{G}_4$, $K_{4}$ is a  clique of $G \in \mathcal{G}_4$ and $V(K_{4})=\{x_{0},x_{1},x_{2},x_{3}\}$.
If $\left|T_{2}\right|=0$, then $G$ is an induced subgraph of $\{H_{1},H_{3},H_{5},H_{6},H_{34},H_{64},
A_{14}, A_{91}, {D}_2\}$.
\end{lem}
\begin{proof}

Let $i, j, k, z$ ($0\leq i, j, k, z\leq 3$) be different integers.
If $\left|T\right|= 0$, then $G=K_{4}\subset D_2$, the result holds.
We distinguish the following three cases.

\noindent$\textbf{Case 1.}$ $\left|T_{x_{i}}\right| \geq 1$, $\left|T_{x_{j}}\right| \geq 1$ and $\left|T_{x_{k}}\right| \geq 1$.

Let $p_{i,1}\in T_{x_{i}}$, $p_{j,1}\in T_{x_{j}}$ and $p_{k,1}\in T_{x_{k}}$.
We first prove that $\left|E(T_{1})\right| =0$.
If $\left|E(\{p_{i,1},p_{j,1},p_{k,1}\})\right| =1$ ($p_{i,1}\sim p_{j,1}$ and $p_{k,1}\nsim p_{i,1},p_{j,1}$), then  $R_{3}\subset G$ (deleting $x_{i}$), a contradiction.
If  $\left|E(\{p_{i,1},p_{j,1},p_{k,1}\})\right| =2$ ($p_{i,1}\sim p_{j,1},p_{k,1}$ and  $p_{j,1}\nsim p_{k,1}$), then  $R_{7}\subset G$ (deleting $x_{j}$), a contradiction.
If  $\left|E(\{p_{i,1},p_{j,1},p_{k,1}\})\right| =3$, then  $R_{11}\subset G$ (deleting $x_{j}$), a contradiction.
Thus, by Lemma \ref{lem-9}\ref{lem-9:item-1}, $\left|E(T_{1})\right| =0$.
Then $\left|T_{x_{g}}\right|\leq 2$ ($0\leq g \leq 3$) for $R_{48} \not\subset G$.
Thus, if $\left|T_{3}\right| =0$, then $G\subset A_{91}$ by  Lemma \ref{lem-9}\ref{lem-9:item-6}.
If $2\leq \left|T_{3}\right| \leq 4$,  then $G\subset H_{34}-v\subset H_{34}$ (if $\left|T_{3}\right| = 2$ and  $v$ is a vertex of degree $2$ in $H_{34}$, as shown in Fig.~\ref{fig-4}) or $G\subset H_{5}$ (if $\left|T_{3}\right| = 3$) or $G\subset H_{6}$ (if $\left|T_{3}\right| = 4$)  by Lemma \ref{lem-9}\ref{lem-9:item-1} and \ref{lem-9:item-5}.

Suppose that $\left|T_{3}\right| =1$ by Lemma \ref{lem-9}\ref{lem-9:item-1}.
If   $\left|E(T_{1},T_{3})\right| =0$, then $G\subset A_{14}$ by  Lemma \ref{lem-9}\ref{lem-9:item-5} and \ref{lem-9:item-6}. 
We now suppose that $\left|E(T_{1},T_{3})\right| \geq 1$.
By Lemma \ref{lem-9}\ref{lem-9:item-5} and $R_{14}\not\subset G$ (if $\left|T_{x_{i}x_{j}x_{z}}\right| =1$ and $\left|E(T_{x_{i}x_{j}x_{z}},p_{i,1})\right| \geq 1$, then deleting $x_{j}$), $\left|T_{3}-T_{x_{i}x_{j}x_{k}}\right| =0$.
Thus, $\left|T_{x_{i}x_{j}x_{k}}\right| =1$.
Let $w_{ijk,1}\in T_{x_{i}x_{j}x_{k}}$.
By Lemma \ref{lem-9}\ref{lem-9:item-5} and without loss of generality, let $p_{j,1}\sim w_{ijk,1}$.
Then  $\left|T_{x_{z}}\right| =0$; otherwise, $R_{14}\subset G$ (deleting $x_{i}$),  a contradiction.
Then $\left|E(T_{0},T_{x_{i}x_{j}x_{k}})\right| =0$ by Lemma \ref{lem-9}\ref{lem-9:item-6}; otherwise, either $R_{22}\subset G$ (if $p_{i,1}\nsim w_{ijk,1}$, then deleting $x_{j}$) or $R_{25}\subset G$ (if $p_{i,1}\sim w_{ijk,1}$, then deleting $x_{j}$),  a contradiction.
Thus, $\left|T_{0}\right| =0$ by Lemma \ref{lem-9}\ref{lem-9:item-6}.
If $\left|T_{x_{i}}\right| =2$, then $\left|E(T_{x_{i}},T_{x_{i}x_{j}x_{k}})\right| = 1$; otherwise, $R_{25}\subset G$ (deleting $x_{k}$) by Lemma \ref{lem-9}\ref{lem-9:item-5},  a contradiction.
Thus,  $G\subset H_{3}-\{u,v\}\subset H_{3}$ (where  $u,v$ are  vertices of degree $2$ in $H_{3}$, as shown in Fig.~\ref{fig-4}).

\noindent$\textbf{Case 2.}$ $\left|T_{x_{i}}\right| \geq 1$, $\left|T_{x_{j}}\right| \geq 1$, and $\left|T_{x_{k}}\cup T_{x_{z}}\right| =0$.

We first suppose that $\left|E(T_{1})\right| =0$. 
If $\left|T_{x_{i}}\right| \leq 2$ and $\left|T_{x_{j}}\right| \leq 2$, an argument similar to that in Case $1$ gives the same conclusion.
Without loss of generality, we now suppose $\left|T_{x_{i}}\right|\geq 3$.
By Lemma \ref{lem-9}\ref{lem-9:item-1} and $R_{49} \not\subset G$, $\left|T_{x_{i}}\right|= 3$.
Since $R_{21} \not\subset G$,  $\left|T_{x_{j}}\right|= 1$.
By Lemma \ref{lem-9}\ref{lem-9:item-5}, we get $\left|T_{x_{j}x_{k}x_{z}}\right| = 0$.
If $\left|T_{x_{i}x_{j}x_{k}}\right| =1$ and $\left|E(T_{x_{i}},T_{x_{i}x_{j}x_{k}})\right| =0$, then either $R_{25}\subset G$ (if $\left|E(T_{x_{j}},T_{x_{i}x_{j}x_{k}})\right| =1$) or $R_{53}\subset G$ (if $\left|E(T_{x_{j}},T_{x_{i}x_{j}x_{k}})\right| =0$), a contradiction.
If $\left|T_{x_{i}x_{j}x_{k}}\right| =1$ and $\left|E(T_{x_{i}},T_{x_{i}x_{j}x_{k}})\right| \geq 1$, then  either $R_{25}\subset G$ (if $\left|E(T_{x_{j}},T_{x_{i}x_{j}x_{k}})\right| =1$) or $R_{57}\subset G$ (if $\left|E(T_{x_{j}},T_{x_{i}x_{j}x_{k}})\right| =0$) by Lemma \ref{lem-9}\ref{lem-9:item-5}, a contradiction.
Thus, $\left|T_{x_{i}x_{j}x_{k}}\right| =0$ by Lemma \ref{lem-9}\ref{lem-9:item-1}.
By symmetry, $\left|T_{x_{i}x_{j}x_{z}}\right| =0$.
If $\left|T_{x_{i}x_{k}x_{z}}\right| =0$, then $G\subset D_{2}$.
By Lemma \ref{lem-9}\ref{lem-9:item-1}, let $\left|T_{x_{i}x_{k}x_{z}}\right| =1$.
By  Lemma \ref{lem-9}\ref{lem-9:item-5} and $R_{14} \not\subset G$ (deleting $x_{k}$), $\left|E(T_{1}, T_{x_{i}x_{k}x_{z}})\right| =0$.
By Lemma \ref{lem-9}\ref{lem-9:item-6} and $d(G)< 4$,  $\left|T_{0}\right| =0$.
Thus, $G\subset H_{1}$. 

We now suppose that $\left|E(T_{1})\right| \neq 0$. 
Let $p_{i,1}\in T_{x_{i}}$, $p_{j,1}\in T_{x_{j}}$ and $p_{i,1}\sim p_{j,1}$ by Lemma \ref{lem-9}\ref{lem-9:item-1}.
Then $\left|T_{3}\right| =0$ by Lemma \ref{lem-9}\ref{lem-9:item-5}.
Then  $\left|E(T_{x_{i}}-\{p_{i,1}\},p_{j,1})\right| = 0$; otherwise, $R_{4} \subset G$ (deleting $x_{i}$), a contradiction.
Since $R_{49} \not\subset G$, at most three vertices in $T_{x_{j}}$ are not adjacent to a vertex of $T_{x_{i}}$.
Since $R_{21} \not\subset G$, at most two vertices in $T_{x_{j}}$ are not adjacent to two vertices of $T_{x_{i}}$.
Thus, by symmetry and Lemma \ref{lem-9}\ref{lem-9:item-1}, $G\subset D_{2}$.

\noindent$\textbf{Case 3.}$ $\left|T_{x_{j}}\cup T_{x_{k}}\cup T_{x_{z}}\right| =0$. 

If $\left|T_{x_{i}}\right| = 0$, then $G\subset H_{6}$ by Lemma \ref{lem-9}\ref{lem-9:item-1}.
By an argument analogous to that in Case $1$, the same conclusion holds for  $1 \leq \left|T_{x_{i}}\right| \leq 2$.
From now on, we will only consider $\left|T_{x_{i}}\right| \geq 3$. 
If $\left|T_{3}\right|=0$, then $G\subset H_{1}$ by Lemma \ref{lem-9}\ref{lem-9:item-1}.
Let $\left|T_{3}-T_{x_{j}x_{k}x_{z}}\right|\geq 1$ by Lemma \ref{lem-9}\ref{lem-9:item-5}.
Since $R_{54} \not\subset G$, at most three vertices in $T_{x_{i}}$ are not adjacent to a vertex of $T_{3}$.
If $\left|T_{3}\right|= 1$, then $G\subset H_{1}$ by Lemma \ref{lem-9}\ref{lem-9:item-1}\ref{lem-9:item-5} and \ref{lem-9:item-6}. 
If $\left|T_{3}\right|\geq 2$, then  $G\subset H_{64}-\{u,v\}\subset H_{64}$ (where $u,v$ are  vertices of degree $2$ in $H_{64}$, as shown in Fig.~\ref{fig-4}) by Lemma \ref{lem-9}\ref{lem-9:item-1} and \ref{lem-9:item-5}. 
\end{proof}

\input{generated_figures/A_figures.tex}

\input{generated_figures/H_figures.tex}

\begin{lem}\label{lem-11}
Suppose that $G \in \mathcal{G}_4$, $K_{4}$ is a  clique of $G \in \mathcal{G}_4$ and $V(K_{4})=\{x_{0},x_{1},x_{2},x_{3}\}$.
If  $\left|T\right|=\left|T_{2}\right|\geq 1$, then $G$ is an induced subgraph of $\{H_{7},\ldots,H_{13},  H_{45}, H_{46}, 
A_{2}, {D}_2\}$.
\end{lem}

\begin{proof}
Let  $i, j, k, z$ ($0\leq i, j, k, z\leq 3$) be different integers.
If $T_{2}=T_{x_{i}x_{j}}$, then $G\subset D_{2}$ by Lemma \ref{lem-9}\ref{lem-9:item-3}.
We now suppose that at least two of the sets $T_{x_gx_h}$ ($0\leq g<h\leq3$) are nonempty.
Choose distinct $i,j$ such that $\left|T_{x_ix_j}\right|\geq1$.
Then
\begin{equation}\label{eq1}
\left|T_{x_ix_j}\right|\geq1
\quad\text{and}\quad
\left|T_2-T_{x_ix_j}\right|\geq1.
\tag{1}
\end{equation}
By  Lemma \ref{lem-9}\ref{lem-9:item-3}\ref{lem-9:item-12} and \ref{lem-9:item-13},  if $\left|E(T_{x_{i}x_{j}})\right| \geq 2$, then $G$ contains a subdivision of $K_{3,3}$, a contradiction.
Thus, for any $i,j$ ($0\leq i< j\leq 3$), $\left|E(T_{x_{i}x_{j}})\right| \leq 1$.
Thus, without loss of generality, we  distinguish the following two cases.

\noindent$\textbf{Case 1.}$  For some $i,j$ ($0\leq i< j\leq 3$), $\left|E(T_{x_{i}x_{j}})\right| =1$.

By Lemma \ref{lem-9}\ref{lem-9:item-12}\ref{lem-9:item-13} and  $G$ contains no subdivision of $K_{3,3}$, $\left|T_{x_{k}x_{z}}\right|=0$ and  $\left|T_{x_{j}x_{k}}\right|\leq 1$.
Without loss of generality,  we assume that  $\left|T_{x_{j}x_{k}}\right|=1$ by \eqref{eq1}.
Then $\left|T_{x_{i}x_{z}}\right|=\left|T_{x_{j}x_{z}}\right|=0$ by Lemma \ref{lem-9}\ref{lem-9:item-13} and  $G$ contains no subdivision of $K_{3,3}$.
Then $\left|T_{x_{i}x_{j}}\right| \leq 5$; otherwise, $R_{37}\subset G$ (deleting $x_{j}$ and $x_{z}$), a contradiction.
Thus, if $\left|T_{x_{i}x_{k}}\right|=0$, then $G\subset H_{45}$ by  Lemma \ref{lem-9}\ref{lem-9:item-11} and \ref{lem-9:item-13} (In fact, $G\subset H_{45}^{'}$, where  $H_{45}^{'}$ is obtained from $H_{45}$ by deleting four pendant vertices).
Let $\left|T_{x_{i}x_{k}}\right|=1$ and $q_{ik,1}\in T_{x_{i}x_{k}}, q_{jk,1}\in T_{x_{j}x_{k}}$.
If $q_{ik,1} \nsim q_{jk,1}$, then $G\subset H_{9}$ by   Lemma \ref{lem-9}\ref{lem-9:item-11} and \ref{lem-9:item-13}; otherwise, $R_{47} \subset G$ (if $\left|T_{x_{i}x_{j}}\right| \geq 5$, then deleting $x_{i}$ and $x_{k}$).
If $q_{ik,1} \sim q_{jk,1}$, then $G\subset H_{7}$ by   Lemma \ref{lem-9}\ref{lem-9:item-11}  and \ref{lem-9:item-13}; otherwise, $R_{25} \subset G$ (if $\left|T_{x_{i}x_{j}}\right| \geq 3$, then deleting $x_{j}$ and a vertex of $T_{x_{i}x_{j}}$ with degree $\geq 3$).

\noindent$\textbf{Case 2.}$  For any $i,j$ ($0\leq i< j\leq 3$), $\left|E(T_{x_{i}x_{j}})\right| =0$.

We proceed with the following fact. 

\begin{fact}\label{f-2.1}
If $\left|E(T_{2})\right| \geq 1$ and $\left|T_{x_{i}x_{j}}\right| \geq 2$, then $\left|E(T_{x_{i}x_{j}},T_{2}-T_{x_{i}x_{j}})\right| \geq 1$.
\end{fact}

\begin{proof}
If $\left|E(T_{x_{i}x_{j}},T_{2}-T_{x_{i}x_{j}})\right|=0$, then  either $R_{1}\subset G$ (if $\left|E(T_{x_{j}x_{k}}, T_{x_{j}x_{z}})\right| \geq 1$, then deleting $x_{j}$ and $x_{k}$) or $R_{33}\subset G$ (if $\left|E(T_{x_{j}x_{k}}, T_{x_{i}x_{k}})\right| \geq 1$, then deleting $x_{z}$) by Lemma \ref{lem-9}\ref{lem-9:item-11} and \ref{lem-9:item-12}, a contradiction.
\end{proof}

For any $g,h$ ($0\leq g< h\leq 3$), $\left|T_{x_{g}x_{h}}\right| \leq 4$ by \ref{lem-9:item-11} and \ref{lem-9:item-12} of Lemma \ref{lem-9}; otherwise, either $R_{37} \subset G$ (if $\left|E(T_{x_{g}x_{h}}, T_{2}-T_{x_{g}x_{h}})\right| \geq 1$, then deleting $x_{j}$ and $x_{z}$) or $R_{46} \subset G$ (if $\left|E(T_{x_{g}x_{h}}, T_{2}-T_{x_{g}x_{h}})\right| =0$, then deleting $x_{j}$), a contradiction.
Let $|T_{x_i x_j}| = \max\limits_{0\leq g< h\leq 3} |T_{x_{g}x_{h}}|$.
Without loss of generality, we  distinguish the following four subcases.

\noindent$\textbf{Subcase 2.1.}$  $\left|T_{x_{i}x_{j}}\right|=4$.

Then $\left|T_{x_{j}x_{k}}\right|\leq 1$ by Lemma \ref{lem-9}\ref{lem-9:item-11}; otherwise, $R_{21} \subset G$ (deleting $x_{j}$ and $x_{z}$), a contradiction.
By  \eqref{eq1} and Lemma \ref{lem-9}\ref{lem-9:item-12}, let $\left|T_{x_{j}x_{k}}\right| = 1$.

We first suppose that $\left|E(T_{2})\right| =0$.
Then  $\left|T_{x_{i}x_{k}}\cup T_{x_{j}x_{z}}\right| = 0$; otherwise, either  $R_{45} \subset G$ (if  $\left|T_{x_{j}x_{z}}\right| \neq 0$, then deleting $x_{j}$) or $R_{61} \subset G$ (if  $\left|T_{x_{i}x_{k}}\right| \neq 0$, then deleting $x_{z}$), a contradiction.
Thus, $G\subset H_{10}-q_{ik,1} \subset H_{10}$ (where $q_{ik,1}$ is a vertex of $T_{x_{i}x_{k}}$ in $H_{10}$, as shown in Fig.~\ref{fig-4}).

Suppose that $\left|E(T_{2})\right| \geq 1$.
Without loss of generality, we assume that  $\left|E(T_{x_{i}x_{j}},T_{x_{j}x_{k}})\right| = 1$ by Lemma \ref{lem-9}\ref{lem-9:item-11} and (F\ref{f-2.1}).
Let $q_{ij,1} \in T_{x_{i}x_{j}}, q_{jk,1}\in T_{x_{j}x_{k}}$ and $q_{ij,1}\sim q_{jk,1}$.
Since  $G$ contains no subdivision of $K_{3,3}$, $\left|E(T_{x_{i}x_{k}},T_{x_{i}x_{j}}-\{q_{ij,1}\})\right| = 0$.
Then  $\left|E(T_{x_{i}x_{k}},\{q_{ij,1},q_{jk,1}\})\right| = 0$; otherwise, $R_{23} \subset G$ (if $\left|E(T_{x_{i}x_{k}},q_{ij,1})\right| = 1$ and $\left|E(T_{x_{i}x_{k}},q_{jk,1})\right| = 1$, then deleting $x_{j}$ and $x_{k}$) or $R_{37} \subset G$ (if $\left|E(T_{x_{i}x_{k}},q_{ij,1})\right| = 0$ and $\left|E(T_{x_{i}x_{k}},q_{jk,1})\right| = 1$, then deleting $x_{j}$ and $x_{k}$) or $R_{38} \subset G$ (if $\left|E(T_{x_{i}x_{k}},q_{ij,1})\right| = 1$ and $\left|E(T_{x_{i}x_{k}},q_{jk,1})\right| = 0$, then deleting $x_{j}$ and $x_{z}$), a contradiction.
By an argument similar to the one above, $\left|E(T_{x_{j}x_{z}},T_{x_{i}x_{j}}\cup T_{x_{j}x_{k}})\right|=0$.
Thus,  $G\subset H_{10}$ by Lemma \ref{lem-9}\ref{lem-9:item-11}.

\noindent$\textbf{Subcase 2.2.}$  $\left|T_{x_{i}x_{j}}\right|=3$.

Suppose that $\left|E(T_{2})\right| =0$.
If  for some $c,v$ ($c=i$ or $j$ and  $v=k$ or $z$), $\left|T_{x_{c}x_{v}}\right|\geq 2$,  then $R_{21} \subset G$ (if  $\left|T_{x_{j}x_{k}}\right| \geq 2$, then deleting $x_{j}$ and $x_{z}$) by Lemma \ref{lem-9}\ref{lem-9:item-12}, a contradiction.
If  for any $c,v$ ($c=i$ or $j$ and  $v=k$ or $z$), $\left|T_{x_{c}x_{v}}\right|\leq 1$,  then $G\subset H_{13}$ by Lemma \ref{lem-9}\ref{lem-9:item-12}.

Suppose that $\left|E(T_{2})\right|  \geq1$.
Without loss of generality, we assume that  $\left|E(T_{x_{i}x_{j}},T_{x_{j}x_{k}})\right| = 1$ by Lemma \ref{lem-9}\ref{lem-9:item-11} and (F\ref{f-2.1}).
Let $q_{ij,1}\in T_{x_{i}x_{j}}, q_{jk,1}\in T_{x_{j}x_{k}}$ and $q_{ij,1}\sim q_{jk,1}$.
Then $\left|T_{x_{j}x_{k}}\right| \leq 2$ by Lemma \ref{lem-9}\ref{lem-9:item-11}; otherwise, $R_{21} \subset G$ (deleting $x_{j}$ and $x_{z}$), a contradiction.
Then $\left|E(T_{x_{i}x_{k}}, (T_ {x_{j}x_{z}}\cup T_{x_{i}x_{j}}\cup T_{x_{j}x_{k}})-\{q_{ij,1}, q_{jk,1}\})\right| = 0$ by Lemma \ref{lem-9}\ref{lem-9:item-11} and  $G$ contains no subdivision of $K_{3,3}$.
Then $\left|E(T_{x_{i}x_{k}}, q_{jk,1})\right| = 0$; otherwise, $R_{33}\subset G$ (deleting $q_{ij,1}$ and $x_{z}$), a contradiction.
Then $\left|E(T_{x_{j}x_{z}}, \{q_{jk,1},q_{ij,1}\})\right| = 0$; otherwise, $G$ contains a subdivision of $K_{5}$, a contradiction.
Then $\left|E(T_{x_{j}x_{z}}, (T_{x_{i}x_{j}}-\{q_{ij,1}\})\cup (T_{x_{j}x_{k}}-\{q_{jk,1}\}))\right| = 0$ by Lemma \ref{lem-9}\ref{lem-9:item-11}; otherwise, $R_{2} \subset G$ (if $\left|E(T_{x_{j}x_{z}}, T_{x_{i}x_{j}}-\{q_{ij,1}\})\right| = 1$, then deleting $x_{j}$ and $x_{z}$ or  if $\left|E(T_{x_{j}x_{z}}, T_{x_{j}x_{k}}-\{q_{jk,1}\})\right| = 1$, then deleting $x_{i}$ and $x_{j}$).
Thus,  $\left|T_{x_{j}x_{z}}\right|\leq 1$; otherwise, $R_{21}\subset G$ (deleting $x_{j}$ and $x_{k}$), a contradiction.
If $\left|E(T_{x_{i}x_{k}}, q_{ij,1})\right| = 0$, then $\left|T_{x_{i}x_{k}}\right|\leq 1$ for $R_{21}\not\subset G$.
Thus, either $G\subset H_{11}$ (if  $\left|E(T_{x_{i}x_{k}}, q_{ij,1})\right| = 1$) or  $G\subset H_{12}$ (if  $\left|E(T_{x_{i}x_{k}}, q_{ij,1})\right| = 0$) by Lemma \ref{lem-9}\ref{lem-9:item-11}.

\noindent$\textbf{Subcase 2.3.}$  $\left|T_{x_{i}x_{j}}\right|=2$.

By an argument similar to that in Subcase 2.2, if $\left|E(T_{2})\right|\geq1$, then either $G\subset H_{11}$ or $G\subset H_{12}$.
Then we suppose that $\left|E(T_{2})\right| =0$.
Without loss of generality, let  $1\leq \left|T_{x_{j}x_{k}}\right|\leq 2$ and $|T_{x_j x_k}| = \max\limits_{T_{x_{g}x_{h}} \subseteq T_2 \setminus T_{x_i x_j}} |T_{x_{g}x_{h}}|$ by Lemma \ref{lem-9}\ref{lem-9:item-12}.
If $\left|T_{x_{j}x_{k}}\right| = 1$, then $G\subset H_{13}$  by Lemma \ref{lem-9}\ref{lem-9:item-12}.
If $\left|T_{x_{j}x_{k}}\right| = 2$, then $G\subset A_{2}-\{u,v\}\subset A_{2}$ (if $\left|T_{x_{i}x_{k}}\right| = 2$ and $u,v$ are  vertices of degree $1$ in $A_{2}$, as shown in Fig.~\ref{fig-3}) or $G\subset H_{46}-\{u,v,u',v'\}\subset H_{46}$ (if $\left|T_{x_{j}x_{z}}\right| = 2$ and  $u,v,u',v'$ are  vertices of degree $1$ in $H_{46}$, as shown in Fig.~\ref{fig-4})
or $G\subset H_{12}-\{q_{ij,1}\}\subset H_{12}$ (if $\left|T_{x_{i}x_{k}}\right| \leq 1$ and  $\left|T_{x_{j}x_{z}}\right| \leq 1$, and $q_{ij,1}$ is a  vertex of $T_{x_{i}x_{j}}$ with degree $3$ in $H_{12}$, as shown in Fig.~\ref{fig-4})  by Lemma \ref{lem-9}\ref{lem-9:item-12}.

\noindent$\textbf{Subcase 2.4.}$  $\left|T_{x_{i}x_{j}}\right|=1$. 

If  $\left|E(T_{2})\right| =0$,  then  $G\subset H_{8}$.
Without loss of generality, let $\left|E(T_{2})\right| \geq 1$ and $\left|E(T_{x_{i}x_{j}}, T_{x_{j}x_{k}})\right| =1$ by Lemma \ref{lem-9}\ref{lem-9:item-11}.
Then $\left|E(T_{x_{j}x_{z}}, T_{x_{i}x_{j}}\cup T_{x_{j}x_{k}})\right| =0$ for G contains no subdivision of $K_{3,3}$.
Then $\left|E(T_{x_{i}x_{k}}, T_{x_{i}x_{j}}\cup T_{x_{j}x_{k}})\right| \leq 2$ by  Lemma \ref{lem-9}\ref{lem-9:item-11}.
Thus, $G\subset H_{7}$ (if  $\left|E(T_{x_{i}x_{k}}, T_{x_{i}x_{j}}\cup T_{x_{j}x_{k}})\right| =2$)
or  $G\subset H_{11}$ (if  $\left|E(T_{x_{i}x_{k}}, T_{x_{i}x_{j}}\cup T_{x_{j}x_{k}})\right| =1$)
or  $G\subset H_{12}$ (if  $\left|E(T_{x_{i}x_{k}}, T_{x_{i}x_{j}}\cup T_{x_{j}x_{k}})\right| =0$) by Lemma \ref{lem-9}\ref{lem-9:item-11}.
\end{proof}

\begin{lem}\label{lem-12}
Suppose that $G \in \mathcal{G}_4$, $K_{4}$ be a clique of $G \in \mathcal{G}_4$ and $V(K_{4})=\{x_{0},x_{1},x_{2},x_{3}\}$.
If  $\left|T\right|>\left|T_{2}\right|\geq 1$, then $G$ is an induced subgraph of $\{H_{2},\allowbreak  H_{3},\allowbreak H_{6},\allowbreak  H_{14},\allowbreak \ldots,\allowbreak H_{70},\allowbreak 
A_{3},\allowbreak \ldots,\allowbreak A_{6},\allowbreak A_{9},\allowbreak   A_{11},\allowbreak {D}_2\}$.
\end{lem}

\begin{proof}
Let  $i, j, k, z$ ($0\leq i, j, k, z\leq 3$) be different integers.
For the induced subgraphs below, use $t_r\in T_0$, $p_{i,r}\in T_{x_i}$, $p_{j,r}\in T_{x_j}$, $p_{k,r}\in T_{x_k}$, $p_{z,r}\in T_{x_z}$ and $q_{ij,r}\in T_{x_ix_j}$.
When the corresponding classes are nonempty, write $T_{x_ix_jx_k}=\{w_{ijk,1}\}$, $T_{x_ix_jx_z}=\{w_{ijz,1}\}$, $T_{x_ix_kx_z}=\{w_{ikz,1}\}$ and $T_{x_jx_kx_z}=\{w_{jkz,1}\}$ by Lemma \ref{lem-9}\ref{lem-9:item-1}.
The labels are chosen locally; different subscripts in the same class denote distinct vertices, and only existing vertices are used.
The labels in each parenthesis are chosen to satisfy the stated conditions, with symmetric choices understood.
By Lemma \ref{lem-9}\ref{lem-9:item-3}\ref{lem-9:item-12} and \ref{lem-9:item-13}, together with the fact that no subdivision of
$K_{3,3}$ is an  subgraph of G, we obtain the following fact.

\begin{fact}\label{f-l34.1}
If $\left|T-(T_{x_{i}}\cup T_{x_{j}}\cup T_{x_{i}x_{j}}\cup T_{0})\right|\geq 1$, then $\left|E(T_{x_{i}x_{j}})\right| \leq 1$.
\end{fact}

Without loss of generality, we  distinguish the following two cases.

\noindent$\textbf{Case 1.}$  For some $i,j$ ($0\leq i< j\leq 3$), $\left|T_{2}\right|=\left|T_{x_{i}x_{j}}\right| \geq1$.

\begin{fact}\label{f-l34.2}
If $\left|T_{x_{i}x_{j}}\right| \geq 5$, then $\left|T_{x_{k}}\cup T_{x_{z}}\right|= \left|E(T_{1},T_{2})\right|= 0$  and $G[T_{1}]\subset 4K_{2}$.
\end{fact}
\begin{proof}
Suppose that $\left|T_{x_{i}x_{j}}\right| \geq 5$.
If $\left|T_{x_{k}}\cup T_{x_{z}}\right|\geq 1$, then $R_{62}\subset G$ by Lemma \ref{lem-9}\ref{lem-9:item-2} and \ref{lem-9:item-13}, a contradiction.
Then $\left|E(T_{1},T_{2})\right|= 0$  by Lemma \ref{lem-9}\ref{lem-9:item-2}.
By reasoning similar to that in Case $2$ of Lemma \ref{lem-10}, if  $\left|T_{x_{i}x_{j}}\right| \geq 5$, then  $G[T_{1}]\subset 4K_{2}$.
\end{proof}

If $\left|T_{x_{i}x_{j}}\right| \geq 9$, then $\left|T_{x_{i}x_{j}x_{k}}\cup T_{x_{i}x_{j}x_{z}}\right|= 0$ by Lemma \ref{lem-9}\ref{lem-9:item-2} and \ref{lem-9:item-13}; otherwise, $R_{75}\subset G$, a contradiction.
Thus, if  $\left|T_{x_{i}x_{j}}\right| \geq 9$, then $G\subset D_{2}$  by (F\ref{f-l34.2}) and Lemma \ref{lem-9}\ref{lem-9:item-2}\ref{lem-9:item-3}\ref{lem-9:item-6} and \ref{lem-9:item-12}.
Next, we consider the case $\left|T_{x_{i}x_{j}}\right|\leq 8$, which is divided into the following five subcases.

\noindent$\textbf{Subcase 1.1.}$ $5\leq \left|T_{x_{i}x_{j}}\right|\leq 8$.

By  Lemma \ref{lem-9}\ref{lem-9:item-12} and (F\ref{f-l34.2}), $\left|T_{x_{k}}\cup T_{x_{z}}\cup T_{x_{i}x_{k}x_{z}}\cup T_{x_{j}x_{k}x_{z}}\cup T_{x_{k}x_{z}}\right|=\left|E(T_{1},T_{2})\right| = 0$ and $G[T_{1}]\subset 4K_{2}$.
Thus, if $\left|T_{3}\right|=0$, then $G \subset D_2$.
Let $1 \leq \left|T_{3}\right| \leq 2$ by   Lemma \ref{lem-9}\ref{lem-9:item-1}.
Then $\left|E(T_{0},T_{3})\right| =0$ by   Lemma \ref{lem-9}\ref{lem-9:item-2}\ref{lem-9:item-6} and \ref{lem-9:item-13}; otherwise, $R_{62} \subset G$ (if $\left|E(T_{0}, T_{x_{i}x_{j}x_{k}})\right|  \geq1$, then deleting $x_{z}$ and a vertex of $T_{x_{i}x_{j}}$ adjacent to $T_{x_{i}x_{j}x_{k}}$), a contradiction.
Thus, $\left|T_{0}\right| =0$ by Lemma \ref{lem-9}\ref{lem-9:item-5} and \ref{lem-9:item-6}.
By  Lemma \ref{lem-9}\ref{lem-9:item-1} and \ref{lem-9:item-12}, $\left|T_{3}\right|\leq 2$.

We first suppose that $\left|T_{3}\right|=2$.
Then $\left|E(T_{x_{i}x_{j}})\right| = 0$ by  Lemma \ref{lem-9}\ref{lem-9:item-13} and  $G$ contains no subdivision of $K_{5}$.
Then $\left|T_{x_{i}}\right| \leq 2$ and $\left|T_{x_{j}}\right| \leq 2$ by  Lemma \ref{lem-9}\ref{lem-9:item-5}; otherwise, $R_{68} \subset G$, a contradiction.
Then $\left|T_{x_{i}x_{j}}-(N(T_{x_{i}x_{j}x_{k}}\cup T_{x_{i}x_{j}x_{z}})\cap T_{x_{i}x_{j}})\right| \leq 6$; otherwise, $R_{76} \subset G$, a contradiction.
If $\left|T_{x_{i}}\cup T_{x_{j}}\right| \geq 1$, then at most five vertices in $T_{x_{i}x_{j}}$  are adjacent to neither vertex of $T_{x_{i}x_{j}x_{k}}\cup T_{x_{i}x_{j}x_{z}}$ by  Lemma \ref{lem-9}\ref{lem-9:item-5}; otherwise, $R_{74} \subset G$, a contradiction.
Thus,  either $G\subset H_{14}$ (if $\left|T_{x_{i}}\cup T_{x_{j}}\right| =0$) or $G\subset H_{28}$ (if $\left|T_{x_{i}}\cup T_{x_{j}}\right| \geq 1$) by  Lemma \ref{lem-9}\ref{lem-9:item-5} and \ref{lem-9:item-15}.

We next suppose that $\left|T_{3}\right|=1$.
Without loss of generality, let $\left|T_{x_{i}x_{j}x_{k}}\right|=1$.
Then $\left|E(T_{1}, T_{3})\right| \leq 1$ by Lemma \ref{lem-9}\ref{lem-9:item-2} and \ref{lem-9:item-5}; otherwise, $R_{61} \subset G$ (deleting $x_{k}$), a contradiction.
By an argument similar to the proof of Lemma \ref{lem-10}, either $T_{1}=T_{x_{i}}$ and $G[T_1\cup T_3]\subset 3K_{1}\cup K_{2}$ (see $H_{1}$ in Fig.~\ref{fig-4}), or $\left|T_{x_{g}}\right|\leq2$ for $g\in\{i,j,k\}$, $\left|T_{x_{z}}\right|\leq1$, and $\left|E(T_{1},T_{3})\right|=0$ (see $A_{14}$ in Fig.~\ref{fig-3}).
By reasoning similar to the case $\left|T_{3}\right|=2$, if $G[T_1\cup T_3]\subset A_{14}[T_1\cup T_3]$, then either $G\subset H_{14}$ or $G\subset H_{28}$.
Let $G[T_1\cup T_3]\subset H_{1}[T_1\cup T_3]$ and $G[T_1\cup T_3]\not\subset A_{14}[T_1\cup T_3]$.
Then either $\left|T_{x_{i}}\right|\geq 3$ or $\left|E(T_{x_{i}},T_{3})\right| =1$.
Then $\left|E(T_{2})\right| = 1$ by Lemma \ref{lem-9}\ref{lem-9:item-2}\ref{lem-9:item-15} and (F\ref{f-l34.1}); otherwise, $R_{46} \subset G$ (if $\left|E(T_{x_{i}},T_{3})\right| =1$ and $\left|E(T_{2},T_{3})\right| =0$, then deleting $x_{i}$ and $x_{z}$), $R_{47} \subset G$ (if $\left|E(T_{x_{i}},T_{3})\right| =1$ and $\left|E(T_{2},T_{3})\right| =1$, then deleting $x_{i}$ and $x_{z}$), $R_{68} \subset G$ (if $\left|E(T_{x_{i}},T_{3})\right| =0$, $\left|T_{x_{i}}\right| \geq 3$ and $\left|E(T_{2},T_{3})\right| =0$), or $R_{69} \subset G$ (if $\left|E(T_{x_{i}},T_{3})\right| =0$, $\left|T_{x_{i}}\right| \geq 3$ and $\left|E(T_{2},T_{3})\right| =1$), a contradiction.
Then $\left|T_{x_{i}x_{j}}\right| =5$ by Lemma \ref{lem-9}\ref{lem-9:item-2}\ref{lem-9:item-2}\ref{lem-9:item-13}; otherwise, either $R_{47} \subset G$ (if $\left|E(T_{x_{i}},T_{3})\right| =1$) or $R_{69} \subset G$ (if $\left|T_{x_{i}}\right| \geq 3$), a contradiction.
Thus,  $G\subset H_{27}$ by Lemma \ref{lem-9}\ref{lem-9:item-2}\ref{lem-9:item-2}\ref{lem-9:item-13}.

\noindent$\textbf{Subcase 1.2.}$ $\left|T_{x_{i}x_{j}}\right|=4$.

Lemma \ref{lem-9}\ref{lem-9:item-2} and  \ref{lem-9:item-12},  $\left|T_{x_{i}x_{k}x_{z}}\cup T_{x_{j}x_{k}x_{z}}\cup T_{x_{k}x_{z}}\right|=\left|E(T_{x_{i}}\cup T_{x_{j}},T_{2})\right| = 0$.
Then, by reasoning similar to that in Subcase 1.1, if $\left|T_{x_{k}}\cup T_{x_{z}}\right|=\left|E(T_{0},T_{3})\right| = 0$, then  $G$ is an induced subgraph of $\{H_{14},H_{27},H_{28}, D_{2}\}$.

We first suppose that $\left|E(T_{0},T_{x_{i}x_{j}x_{k}})\right| \geq 1$.
Then $\left|E(T_{2})\right| = 1$ 
by  Lemma \ref{lem-9} \ref{lem-9:item-6}\ref{lem-9:item-15} and (F\ref{f-l34.1}); otherwise, either $R_{62} \subset G$ (if $\left|E(T_{x_{i}x_{j}},T_{x_{i}x_{j}x_{k}})\right| =0$, then deleting $x_{z}$) or $R_{64} \subset G$ (if $\left|E(T_{x_{i}x_{j}},T_{x_{i}x_{j}x_{k}})\right| =1$, then deleting $x_{z}$), a contradiction.
Then $\left|T_{x_{i}x_{j}x_{z}}\right| = 0$ by Lemma \ref{lem-9}\ref{lem-9:item-13} and  $G$ contains no subdivision of $K_{5}$.
Then $\left|T_{x_{i}}\cup T_{x_{j}}\right| =0$ by Lemma \ref{lem-9}\ref{lem-9:item-2}\ref{lem-9:item-6} and \ref{lem-9:item-13}; otherwise, either $R_{22} \subset G$ (if $\left|E(T_{x_{i}},T_{x_{i}x_{j}x_{k}})\right| =1$, then deleting $x_{i}$, $x_{k}$ and $x_{z}$) or $R_{58} \subset G$ (if $\left|E(T_{x_{i}},T_{x_{i}x_{j}x_{k}})\right| =0$, then deleting $x_{z}$ and a vertex of $T_{x_{i}x_{j}}$ with degree is $4$), a contradiction.
Similarly,   $\left|T_{x_{k}}\cup T_{x_{z}}\right| =0$.
Thus,  $G=H_{18}$ by Lemma \ref{lem-9}\ref{lem-9:item-2}\ref{lem-9:item-5}\ref{lem-9:item-6} and \ref{lem-9:item-13}.

We now suppose that $\left|E(T_{0},T_{3})\right| =0$ and $\left|T_{x_{k}}\cup T_{x_{z}}\right|\geq 1$.
Without loss of generality, suppose $\left|T_{x_{k}}\right| \geq 1$ and $p_{k,1}\in T_{x_{k}}$.
Then $\left|E(T_{2})\right| =1$  by  Lemma \ref{lem-9}\ref{lem-9:item-2} and (F\ref{f-l34.1}); otherwise, either $R_{62} \subset G$ (if $\left|E(p_{k,1},T_{x_{i}x_{j}})\right| =0$) or  $R_{41} \subset G$ (if $\left|E(p_{k,1},T_{x_{i}x_{j}})\right| =1$, then deleting $x_{z}$), a contradiction.
Then $\left|T_{x_{k}}\right| =1$ by  Lemma \ref{lem-9}\ref{lem-9:item-13}; otherwise, $R_{21} \subset G$ (deleting $x_{j}$ and $x_{k}$), a contradiction.
Then $\left|T_{x_{z}}\right| =0$  by  Lemma \ref{lem-9}\ref{lem-9:item-13} and  $G$ contains no subdivision of $K_{5}$.
Then $\left|T_{x_{i}}\cup T_{x_{j}}\right| =0$ by  Lemma \ref{lem-9}\ref{lem-9:item-2}\ref{lem-9:item-16}; otherwise, $R_{58} \subset G$, a contradiction. 
Then $\left|T_{3}\right| =0$ by Lemma \ref{lem-9}\ref{lem-9:item-13} and  $G$ contains no subdivision of $K_{5}$.
Then  $\left|T_{0}\right| = 0$  by Lemma \ref{lem-9}\ref{lem-9:item-6}; otherwise, $R_{3} \subset G$ (deleting $x_{j}$), a contradiction.
Thus,  $G=H_{15}$.

\noindent$\textbf{Subcase 1.3.}$ $\left|T_{x_{i}x_{j}}\right|=3$.

By Lemma \ref{lem-9}\ref{lem-9:item-2} and  \ref{lem-9:item-12},  $\left|T_{x_{i}x_{k}x_{z}}\cup T_{x_{j}x_{k}x_{z}}\cup T_{x_{k}x_{z}}\right|=\left|E(T_{x_{i}}\cup T_{x_{j}},T_{2})\right| = 0$.
We first suppose that $\left|T_{x_k}\cup T_{x_z}\right|=\left|E(T_0,T_3)\right|=0$.
If $\left|T_3\right|=2$, then $G$ is an induced subgraph of $\{H_{14},H_{27},H_{28}\}$ by  Lemma \ref{lem-9}\ref{lem-9:item-5}\ref{lem-9:item-13} and \ref{lem-9:item-15} and $R_{54},R_{53}\not\subset G$.
By an argument similar to the proof of Lemma \ref{lem-10},  if $\left|E(T_{1})\right|\geq 1$, then $G\subset D_2$ and if $\left|E(T_1)\right|=0$ and $\left|T_3\right|\leq1$, then $G[T_1\cup T_3\cup K_4]\subset H_{1}[T_1\cup T_3\cup K_4]$ (choosing the copy of $K_4$ with $T_{2}=\varnothing$ in $H_{1}$) or $G[T_1\cup T_{3}\cup K_4]\subset A_{14}[T_1\cup T_3\cup K_4]$ or $G[T_1\cup T_{3}\cup K_4]\subset H_{3}[T_1\cup T_3\cup K_4]$ (choosing the copy of $K_4$ with $T_{2}=\varnothing$ in $H_{3}$).
By reasoning similar to that in Subcase 1.1, if   $G[T_1\cup T_3\cup K_4]\subset H_{1}[T_1\cup T_3\cup K_4]$ (or $A_{14}[T_1\cup T_3\cup K_4]$), then  $G$ is an induced  subgraph of $\{H_{14},H_{27},H_{28}, D_2\}$. 
If  $G[T_1\cup T_3\cup K_4]\subset H_3[T_1\cup T_3\cup K_4]$, $G[T_1\cup T_3\cup K_4]\not\subset H_1[T_1\cup T_3\cup K_4]$,  $G[T_1\cup T_3\cup K_4]\not\subset A_{14}[T_1\cup T_3\cup K_4]$ and $\left|T_{3}\right| \geq 1$, then  $\left|E(T_{2})\right| =1$  by  Lemma \ref{lem-9} \ref{lem-9:item-13}\ref{lem-9:item-15} and (F\ref{f-l34.1}); otherwise, $G[\{x_i,x_j,x_z,q_{ij,1},q_{ij,2},q_{ij,3},p_{i,1},p_{j,1},w_{ijk,1}\}]\cong R_{61}$ (if $p_{i,1},p_{j,1}\sim w_{ijk,1}$ and $\left|E(T_2,w_{ijk,1})\right|=0$) or $G[\{x_i,x_j,x_k,x_z,q_{ij,1},q_{ij,2},q_{ij,3},p_{i,1},p_{j,1},w_{ijk,1}\}]\cong R_{66}$ (if $p_{i,1},p_{j,1}\sim w_{ijk,1}$ and $\left|E(T_2,w_{ijk,1})\right|=1$) or $G[\{x_i,x_j,x_k,q_{ij,1},q_{ij,2},q_{ij,3},p_{i,1},p_{j,1},w_{ijk,1}\}]\cong R_{63}$ (if $p_{i,1}\sim w_{ijk,1}$, $p_{j,1}\nsim w_{ijk,1}$ and $\left|E(T_2,w_{ijk,1})\right|=0$) or $G[\{x_i,x_j,x_z,q_{ij,1},q_{ij,2},q_{ij,3},p_{i,1},p_{j,1},w_{ijk,1}\}]\cong R_{63}$ (if $p_{i,1}\sim w_{ijk,1}$, $p_{j,1}\nsim w_{ijk,1}$ and $\left|E(T_2,w_{ijk,1})\right|=1$), a contradiction.
If $G[T_1\cup T_3\cup K_4]\subset H_3[T_1\cup T_3\cup K_4]$, $G[T_1\cup T_3\cup K_4]\not\subset H_1[T_1\cup T_3\cup K_4]$, and $G[T_1\cup T_3\cup K_4]\not\subset A_{14}[T_1\cup T_3\cup K_4]$, then $G\subset H_{26}$.
Next, we consider $\left|E(T_{0},T_{3})\right| \geq 1$ or $\left|T_{x_{k}}\cup T_{x_{z}}\right|\geq 1$.

We first suppose that $\left|E(T_{0},T_{3})\right| \geq 1$.
By reasoning similar to that in Subcase 1.2, if \\ $\left|E(T_{0},T_{x_{i}x_{j}x_{k}})\right| \geq 1$ and  $\left|E(T_{x_{i}x_{j}})\right| = 1$, then  $\left|T_{x_{i}x_{j}x_{z}}\right| =\left|T_{x_{k}}\right| =\left|T_{x_{z}}\right| = 0$.
If $\left|E(T_{0},T_{x_{i}x_{j}x_{k}})\right| \geq 1$,  $\left|E(T_{x_{i}x_{j}})\right| = 1$ and  $\left|T_{x_{i}}\cup T_{x_{j}}\right| =0$,  then $G\subset H_{18}$  by Lemma \ref{lem-9}\ref{lem-9:item-1}\ref{lem-9:item-5}\ref{lem-9:item-6} and \ref{lem-9:item-13}.
If $\left|E(T_{0},T_{x_{i}x_{j}x_{k}})\right| \geq 1$,  $\left|E(T_{2})\right| = 1$ and  $\left|T_{x_{i}}\cup T_{x_{j}}\right| \geq 1$, then $\left|E(T_{1},T_{x_{i}x_{j}x_{k}})\right| =0$,  $\left|T_{x_{i}}\right|\leq 2$ and $\left|T_{x_{j}}\right| \leq 2$;  otherwise, $G[\{x_i,x_j,x_z,q_{ij,3},p_{i,1},w_{ijk,1},t_1\}]\cong R_{28}$ ($p_{i,1}\sim w_{ijk,1},\ q_{ij,1}\sim q_{ij,2}$) or $G[\{x_i,x_j,x_k,q_{ij,3},w_{ijk,1},t_1,p_{i,1},p_{i,2},p_{i,3}\}]\cong R_{51}$ ($|T_{x_i}|\geq3,\ E(T_{x_i},w_{ijk,1})=\varnothing,\ q_{ij,1}\sim q_{ij,2}$), a contradiction.
Thus, if  $\left|E(T_{0},T_{x_{i}x_{j}x_{k}})\right| \geq 1$,  $\left|E(T_{2})\right| = 1$ and  $\left|T_{x_{i}}\cup T_{x_{j}}\right| \geq 1$,  then $G\subset H_{31}$.
We now supppse that $\left|E(T_{0},T_{x_{i}x_{j}x_{k}})\right| \geq 1$ and $\left|E(T_{2})\right| = 0$.
Then, $\left|T_{x_{i}x_{j}x_{z}}\right| = 0$; otherwise, $G[\{x_i,x_k,x_z,w_{ijk,1},w_{ijz,1},t_1\}]\cong R_{14}$ ($t_1\sim w_{ijk,1}$), a contradiction.
Then $\left|T_{x_{i}}\right| = \left|T_{x_{j}}\right| = 0$ by Lemma \ref{lem-9}\ref{lem-9:item-6}; otherwise, $G[\{x_i,x_j,x_z,q_{ij,3},p_{i,1},w_{ijk,1},t_1\}]\cong R_{28}$ ($p_{i,1}\sim w_{ijk,1},\ q_{ij,3}\nsim w_{ijk,1}$) or $G[\{x_i,x_j,x_k(x_z),q_{ij,1},q_{ij,2},q_{ij,3},p_{i,1},w_{ijk,1},t_1\}]\cong R_{58}$ ($p_{i,1}\nsim w_{ijk,1},\ |E(T_2,w_{ijk,1})|=0 (1)$), a contradiction.
Since $d(G) < 4$,  $\left|T_{x_{z}}\right|  = 0$.
If $\left|E(T_{2},T_{x_{i}x_{j}x_{k}})\right| =0$, then  $\left|T_{x_{k}}\right| = 0$ by Lemma \ref{lem-9}\ref{lem-9:item-6}; otherwise, $G[K_4 \cup \{q_{ij,1},p_{k,1},w_{ijk,1}\}]\cong R_{35}$ ($p_{k,1}\sim w_{ijk,1}$) or $G[K_4 \cup \{q_{ij,1},q_{ij,2},q_{ij,3},p_{k,1},w_{ijk,1}\}]\cong R_{64}$ ($p_{k,1}\nsim w_{ijk,1}$), a contradiction.
Thus, if  $\left|E(T_{2},T_{x_{i}x_{j}x_{k}})\right| =0$, then $G\subset H_{18}$.
We now suppose that $\left|E(T_{2},T_{x_{i}x_{j}x_{k}})\right| =1$.
Since $G$ contains no subdivision of $K_{3,3}$ and $d(G) < 4$, $\left|E(T _{x_{k}},T_{2})\right| =0$.
Then  $\left|E(T_{x_{k}},T_{x_{i}x_{j}x_{k}})\right| =0$ and $\left|T_{x_{k}}\right| \leq 1$  by Lemma \ref{lem-9}\ref{lem-9:item-6}; otherwise, $G[\{x_i,x_k,x_z,w_{ijk,1},p_{k,1},q_{ij,2},q_{ij,3}\}]\cong R_{25}$ ($q_{ij,1}\sim w_{ijk,1},\ p_{k,1}\sim w_{ijk,1}$) or $G[\{x_i,x_k,x_z,w_{ijk,1},q_{ij,1},p_{k,1},p_{k,2}\}]\cong R_{25}$ ($q_{ij,1}\sim w_{ijk,1},\ p_{k,1},p_{k,2}\nsim w_{ijk,1}$), a contradiction.
Thus, $G\subset H_{20}$.

Suppose that $\left|E(T_{0},T_{3})\right| =0$ and $\left|T_{x_{k}}\cup T_{x_{z}}\right|\geq 1$.
Then  $\left|T_{0}\right| =0$ by Lemma \ref{lem-9}\ref{lem-9:item-6}; otherwise, either $R_{3} \subset G$ (if  $\left|E(T_{2},T_{x_{k}})\right|=0$, then deleting $x_{j}$) or $R_{7} \subset G$ (if $\left|E(T_{2},T_{x_{k}})\right|=1$, then deleting $x_{j}$), a contradiction.
By arguments similar to those used in the preceding paragraph,  if $\left|T_{3}\right| \geq 1$, then either $G\subset H_{18}$ or  $G\subset H_{20}$.
Thus, let $\left|T_{3}\right| =0$.
Note that $\left|E(T_{2})\right| \leq 1$.
By reasoning similar to that in Subcase 1.2, if  $\left|E(T_{2})\right|=1$, then $G[K_{4}\cup T_{x_{k}}\cup T_{x_{z}}\cup T_{x_{i}x_{j}}]\subset H_{15}$.
Thus, if $\left|E(T_{2})\right|=1$, then $G\subset H_{25}$ by Lemma \ref{lem-9}\ref{lem-9:item-16}.
We now supppse that $\left|E(T_{2})\right|=0$.
Then  $\left|T_{x_{i}}\cup T_{x_{j}}\right|=0$; otherwise, $G[\{x_i,x_j,x_k,q_{ij,1},q_{ij,2},q_{ij,3},p_{i,1},p_{k,1}\}]\cong R_{40}$ ($p_{k,1}\sim q_{ij,1}$) or $G[\{x_i,x_j,x_k,x_z,q_{ij,1},q_{ij,2},q_{ij,3},p_{i,1},p_{k,1}\}]\cong R_{58}$ ($E(p_{k,1},T_2)=\varnothing$), a contradiction.
Then $\left|T_{x_{k}}\right|\leq 1$ by  Lemma \ref{lem-9}\ref{lem-9:item-1} and \ref{lem-9:item-2}; otherwise, $G[\{x_i,x_j,x_k,q_{ij,2},q_{ij,3},p_{k,1},p_{k,2}\}]\cong R_{24}$ ($E(\{q_{ij,2},q_{ij,3}\},\{p_{k,1},p_{k,2}\})=\varnothing$) or $G[\{x_i,x_j,x_k,q_{ij,1},q_{ij,2},p_{k,1},p_{k,2}\}]\cong R_{29}$ ($p_{k,1},p_{k,2}\sim q_{ij,1}$), a contradiction.
Thus, $G\subset H_{17}$ for $d(G)<4$; otherwise, either $G[K_4 \cup \{q_{ij,1},q_{ij,2},q_{ij,3},p_{k,1},p_{z,1}\}]\cong R_{56}$ ($E(\{p_{k,1},p_{z,1}\},T_2)=\varnothing$) or $G$ contains a subdivision of $K_{5}$, a contradiction.

\noindent$\textbf{Subcase 1.4.}$ $\left|T_{x_{i}x_{j}}\right|=2$.

By Lemma \ref{lem-9}\ref{lem-9:item-2} and  \ref{lem-9:item-12},  $\left|T_{x_{i}x_{k}x_{z}}\cup T_{x_{j}x_{k}x_{z}}\cup T_{x_{k}x_{z}}\right|=\left|E(T_{x_{i}}\cup T_{x_{j}},T_{2})\right| = 0$.
By reasoning similar to that in Subcase 1.3, if $\left|T_{x_{k}}\cup T_{x_{z}}\right|=\left|E(T_{0},T_{3})\right| = 0$, then  $G$ is an induced subgraph of $\{H_{14},H_{26},H_{27},H_{28}, D_2\}$. 

Suppose that $\left|E(T_{0},T_{x_{i}x_{j}x_{k}})\right| \geq 1$.
By   Lemma \ref{lem-9} \ref{lem-9:item-2}\ref{lem-9:item-6} and \ref{lem-9:item-5}, $\left|T_{3}\right|  = \left|T_{x_{i}x_{j}x_{k}}\right|=1$, $\left|E(T_{0},T_{x_{i}x_{j}x_{k}})\right| = 1$, $\left|E(T_{1})\right| = 0$ and $\left|E(T_{2},T_{0})\right| = 0$.
By reasoning similar to that in Subcase 1.3,  if $\left|E(T_{2})\right| = 1$, then $G$ is an induced subgraph of $\{H_{18},H_{31}\}$. 
Then let $\left|E(T_{2})\right| = 0$ by (F\ref{f-l34.1}).
By reasoning similar to that in Subcase 1.3,  we get $\left|T_{x_{z}}\right|  =\left|E(T_{x_{k}},T_{2}\cup T_{3})\right|= 0$ and $\left|T_{x_{k}}\right| \leq 1$. 
By reasoning similar to that in Subcase 1.3,  if   $\left|E(T_{2},T_{x_{i}x_{j}x_{k}})\right| =0$, then $G\subset H_{18}$.
If $\left|E(T_{2},T_{x_{i}x_{j}x_{k}})\right| =1$, then  $\left|T_{x_{i}}\right|,\left|T_{x_{j}}\right|  \leq 2$ by   Lemma \ref{lem-9} \ref{lem-9:item-17}; otherwise, $R_{53} \subset G$, a contradiction. 
If $\left|E(T_{2},T_{x_{i}x_{j}x_{k}})\right| =1$, then  $\left|E(T_{x_{k}},T_{x_{i}x_{j}})\right| =0$; otherwise, $G$ contains a subdivision of $K_{5}$ or $K_{3,3}$, a contradiction. 
If $\left|E(T_{2},T_{x_{i}x_{j}x_{k}})\right| =1$, then  $\left|E(T_{x_{k}},T_{x_{i}x_{j}x_{k}})\right| =0$; otherwise, $R_{24} \subset G$ (deleting $x_{k}$ and a vertex of $T_{x_{i}x_{j}}$ with degree $3$), a contradiction. 
Thus, if  $\left|E(T_{2},T_{x_{i}x_{j}x_{k}})\right| =1$, then $G\subset H_{34}$ by   Lemma \ref{lem-9} \ref{lem-9:item-6} and \ref{lem-9:item-17}.

Suppose that $\left|E(T_{0},T_{3})\right| =0$ and $\left|T_{x_{k}}\cup T_{x_{z}}\right|\geq 1$.
Then  $\left|T_{0}\right| =0$ by Lemma \ref{lem-9}\ref{lem-9:item-6}; otherwise, either $G[\{x_i,x_k,x_z,q_{ij,1},p_{k,1},t_1\}]\cong R_{3}$ ($t_1\sim q_{ij,1},\ p_{k,1}\nsim q_{ij,1}$) or $G[\{x_i,x_k,x_z,q_{ij,1},p_{k,1},t_1\}]\cong R_{7}$ ($t_1\sim q_{ij,1},\ p_{k,1}\sim q_{ij,1}$), a contradiction.
Then $\left|E(T_{x_{k}}\cup T_{x_{z}},T_{3})\right| =0$ by  Lemma \ref{lem-9}\ref{lem-9:item-2} and \ref{lem-9:item-5}; otherwise, either $G[K_4 \cup \{w_{ijk,1},p_{k,1},q_{ij,2}\}]\cong R_{35}$ ($p_{k,1}\sim w_{ijk,1},\ q_{ij,2}\nsim w_{ijk,1}$) (if $\left|E(T_{2},T_{x_k})\right|=0$) or $G$ contains a subdivision of $K_{3,3}$ (if $\left|E(T_{2},T_{x_k})\right|\geq 1$), a contradiction. 
If $\left|T_{3}\right| \geq 1$, then $\left|E(T_{x_{k}}\cup T_{x_{z}},T_{2})\right| =0$ for $G$ contains no subdivision of $K_{3,3}$.
By symmetry and arguments similar to those used in the preceding paragraph, if $\left|T_{x_{i}x_{j}x_{k}}\right| =1$ and $\left|E(T_{2},T_{3})\right| =1$ (or $\left|E(T_{2})\right| = 1$ ), then  $G$ is an induced subgraph of $\{H_{18},H_{31},H_{34}\}$. 
We now suppose that $\left|T_{x_{i}x_{j}x_{k}}\right| =1$ and $\left|E(T_{2})\right| = \left|E(T_{2},T_{3})\right| =0$.
Then $\left|T_{x_{k}}\right| \leq 1$ for $R_{24}\not\subset G$.
By Lemma \ref{lem-9}\ref{lem-9:item-5}, $\left|T_{x_{z}}\right| \leq 1$ and $\left|E(T_{1})\right| =0$.
Thus, if  $\left|T_{x_{i}}\cup T_{x_{j}}\right| =  0$, then $G\subset H_{20}$.
If  $\left|T_{x_{z}}\right|  = \left|T_{x_{k}}\right| =1$, then  $\left|T_{x_{i}}\cup T_{x_{j}}\right| =  0$; otherwise, $G[K_4 \cup \{q_{ij,1},q_{ij,2},p_{i,1},p_{k,1},p_{z,1}\}]\cong R_{52}$, a contradiction.
Thus, if $\left|T_{x_{z}}\right|  = \left|T_{x_{k}}\right| =1$, then $G\subset H_{20}$.
If  $\left|T_{x_{z}}\right|   = 1$ and $\left|T_{x_{i}}\cup T_{x_{j}}\right| \geq 1$, then $\left|E(T_{x_{i}}\cup T_{x_{j}},T_{3})\right|=0$; otherwise, $G[\{x_i,x_j,x_z,p_{i,1},p_{z,1},w_{ijk,1}\}]\cong R_{14}$ ($p_{i,1}\sim w_{ijk,1}$), a contradiction.
If  $\left|T_{x_{z}}\right|   = 1$ and $\left|T_{x_{i}}\cup T_{x_{j}}\right| \geq 1$, then $\left|T_{x_{i}x_{j}x_{z}}\right| =0$; otherwise, $G[\{x_i,x_j,x_z,q_{ij,1},q_{ij,2},p_{i,1},p_{z,1},w_{ijk,1},w_{ijz,1}\}]\cong R_{58}$, a contradiction.
Thus, if  $\left|T_{x_{z}}\right|   = 1$, $\left|T_{x_{k}}\right| = 0$ and $\left|T_{x_{i}}\cup T_{x_{j}}\right| \geq 1$, then $G\subset H_{31}$ by Lemma \ref{lem-9}\ref{lem-9:item-16}.
By symmetry, if  $\left|T_{x_{k}}\right|   = 1$ and $\left|T_{x_{i}}\cup T_{x_{j}}\right| \geq 1$, then $\left|T_{x_{i}x_{j}x_{z}}\right| =0$.
If $\left|T_{x_{k}}\right|   = 1$ and $\left|T_{x_{i}}\right|   \geq 2$ (or  $\left|T_{x_{j}}\right|   \geq 2$), then $\left|E(T_{x_{i}},T_{3})\right| = 1$ (or $\left|E(T_{x_{j}},T_{3})\right| = 1$) by Lemma \ref{lem-9}\ref{lem-9:item-5}; otherwise, $G[\{x_i,x_k,x_z,q_{ij,1},q_{ij,2},p_{i,1},p_{i,2},p_{k,1},w_{ijk,1}\}]\cong R_{47}$ ($p_{i,1},p_{i,2}\nsim w_{ijk,1}$), a contradiction.
Thus, if $\left|T_{x_{z}}\right|   = 0$ and $\left|T_{x_{k}}\right| = 1$  by Lemma \ref{lem-9}\ref{lem-9:item-16}, then $G\subset H_{30}$.
By the above discussion and symmetry, if $\left|T_{3}\right| \geq 1$, then $G$ is an induced subgraph of $\{H_{18},H_{20},H_{30},H_{31},H_{34}\}$. 
Then let $\left|T_{3}\right|=0$.
If $\left|E(T_{2})\right|=1$ and  $\left|T_{x_{k}}\right|\geq 1$, then   $\left|T_{x_{z}}\right|=0$  by Lemma \ref{lem-9}\ref{lem-9:item-13}; otherwise, $G$ contains a subdivision of $K_{5}$, a contradiction.
Thus, if $\left|E(T_{2})\right|=1$ and $\left|T_{x_{k}}\right|=1$, then $G\subset H_{25}$ by Lemma \ref{lem-9}\ref{lem-9:item-2}\ref{lem-9:item-13} and \ref{lem-9:item-16}.
If $\left|E(T_{2})\right|=1$ and $\left|T_{x_{k}}\right|\geq 2$, then $\left|T_{x_{i}}\cup T_{x_{j}}\right|= 0$ by Lemma \ref{lem-9}\ref{lem-9:item-1}\ref{lem-9:item-2}\ref{lem-9:item-13} and \ref{lem-9:item-16}; otherwise, $G[\{x_i,x_j,x_z,q_{ij,1},p_{i,1},p_{k,1},p_{k,2}\}]\cong R_{22}$ ($q_{ij,1}\sim q_{ij,2},\ q_{ij,1}\sim p_{k,1},p_{k,2}$), a contradiction.
If $\left|E(T_{2})\right|=1$, then $\left|T_{x_{k}}\right|\leq  2$ by Lemma \ref{lem-9}\ref{lem-9:item-1} and \ref{lem-9:item-13}; otherwise, $G[\{x_i,x_j,x_z,q_{ij,1},p_{k,1},p_{k,2},p_{k,3}\}]\cong R_{23}$ ($q_{ij,1}\sim p_{k,1},p_{k,2},p_{k,3}$), a contradiction.
Thus, if $\left|E(T_{2})\right|=1$ and $\left|T_{x_{k}}\right|= 2$, then $G\subset H_{16}$ by Lemma \ref{lem-9}\ref{lem-9:item-1} and \ref{lem-9:item-13}.
We now suppose that $\left|E(T_{2})\right|=0$.
Without loss of generality, let $\left|T_{x_{k}}\right|\geq 1$.
If $\left|E(T_{x_{k}},T_{2})\right| \geq 2$, then every vertex in $T_{x_{k}}$ is adjacent to exactly one vertex in $T_{x_{i}x_{j}}$ by Lemma \ref{lem-9}\ref{lem-9:item-2}; otherwise, $G$ contains a subdivision of $K_{3,3}$, a contradiction.
Then $\left|T_{x_{k}}\right|\leq 2$; otherwise, $R_{23} \subset G$, a contradiction.
If $\left|T_{x_{k}}\right|= 2$, then $\left|E(T_{x_{k}},T_{2})\right| =1$; otherwise, either $G[\{x_i,x_j,x_k,q_{ij,1},q_{ij,2},p_{k,1},p_{k,2}\}]\cong R_{24}$ ($E(T_{x_k},T_2)=\varnothing$) or $G[\{x_i,x_j,x_z,q_{ij,1},q_{ij,2},p_{k,1},p_{k,2}\}]\cong R_{24}$ ($p_{k,1},p_{k,2}\sim q_{ij,1}$), a contradiction.
Thus, if $\left|T_{x_{k}}\right|= 2$, then $\left|T_{x_{i}}\cup T_{x_{j}}\right|=0$ by Lemma \ref{lem-9}\ref{lem-9:item-2} and \ref{lem-9:item-16}; otherwise, $G[\{x_i,x_j,x_k,q_{ij,2},p_{i,1},p_{k,1},p_{k,2}\}]\cong R_{22}$ ($p_{k,1}\sim q_{ij,1}$), a contradiction.
If $\left|T_{x_{k}}\right|= 2$, then $\left|E(T_{x_{z}},T_{x_{i}x_{j}})\right|=0$ by Lemma \ref{lem-9}\ref{lem-9:item-16}; otherwise, either $P_{6} \subset G$ or $G$ contains a subdivision of $K_{5}$, a contradiction.
If $\left|T_{x_{k}}\right|= 2$, then $\left|T_{x_{z}}\right|\leq 1$; otherwise, $G[\{x_i,x_z,q_{ij,1},p_{k,1},p_{z,1},p_{z,2}\}]\cong R_{1}$ ($p_{k,1}\sim q_{ij,1}$), a contradiction.
Thus, if $\left|T_{x_{k}}\right|= 2$, then $G\subset H_{19}$.
If $\left|T_{x_{k}}\right|\leq 1$ and $\left|T_{x_{z}}\right|\leq 1$,  then $\left|E(T_{2},T_{x_{k}}\cup T_{x_{z}})\right|\leq 1$ by Lemma \ref{lem-9}\ref{lem-9:item-16}; otherwise, either $P_{6} \subset G$ or $G$ contains a subdivision of $K_{5}$, a contradiction.
If $\left|T_{x_{k}}\right|= 1$, $\left|T_{x_{z}}\right|= 1$ and $\left|E(T_{2},T_{x_{k}}\cup T_{x_{z}})\right|= 0$, then  $\left|T_{x_{i}}\cup T_{x_{j}}\right|=0$  by Lemma \ref{lem-9}\ref{lem-9:item-2} and \ref{lem-9:item-16}; otherwise, $G[K_4 \cup \{q_{ij,1},q_{ij,2},p_{i,1},p_{k,1},p_{z,1}\}]\cong R_{52}$, a contradiction.
Thus, if $\left|T_{x_{k}}\right|\leq 1$ and $\left|T_{x_{z}}\right|\leq 1$,  then $G\subset H_{29}$ by Lemma \ref{lem-9}\ref{lem-9:item-2} and \ref{lem-9:item-16}.

\noindent$\textbf{Subcase 1.5.}$ $\left|T_{x_{i}x_{j}}\right|=1$.

We distinguish the following two subsubcases.

\noindent$\textbf{Subsubcase 1.5.1.}$ $\left|T_{x_{i}x_{k}x_{z}}\cup T_{x_{j}x_{k}x_{z}}\right|= 0$.

By reasoning similar to that in Subcase 1.4, if $\left|T_{x_{k}}\cup T_{x_{z}}\right|=\left|E(T_{0},T_{3})\right| = 0$, then  $G$ is an induced subgraph of $\{H_{14},H_{26},H_{27},H_{28}, D_2\}$. 

Suppose that $\left|E(T_{0},T_{x_{i}x_{j}x_{k}})\right| \geq 1$.
Then  $\left|T_{x_{z}}\cup T_{x_{i}x_{j}x_{z}}\right|  =0$ by Lemma \ref{lem-9}\ref{lem-9:item-5}; otherwise, $P_{6} \subset G$, a contradiction.
Then  $\left|T_{x_{i}}\right|  \leq 2$ and $\left|T_{x_{j}}\right|  \leq 2$ by Lemma \ref{lem-9}\ref{lem-9:item-1}\ref{lem-9:item-2}\ref{lem-9:item-5}\ref{lem-9:item-6} and \ref{lem-9:item-17}; otherwise, either $G[\{x_i,x_j,x_k,q_{ij,1},p_{i,1},p_{i,2},p_{i,3},w_{ijk,1},t_1\}]\cong R_{51}$ ($q_{ij,1}\nsim w_{ijk,1}, t_1\sim w_{ijk,1}$) or $G[\{x_i,x_j,x_z,q_{ij,1},p_{i,1},p_{i,2},p_{i,3},w_{ijk,1},t_1\}]\cong R_{51}$ ($q_{ij,1}\sim w_{ijk,1}, t_1\sim w_{ijk,1}$), a contradiction.
Then $\left|E(T_{x_{k}}, T_{x_{i}x_{j}})\right|=0$; otherwise, $G$ contains a subdivision of $K_{3,3}$, a contradiction.
If  $\left|E(T_{x_{k}},T_{x_{i}x_{j}x_{k}})\right|\geq  1$, then $\left|E(T_{x_{i}x_{j}x_{k}},T_{x_{i}x_{j}})\right|=1$ by Lemma \ref{lem-9}\ref{lem-9:item-6}; otherwise, $G[\{x_i,x_j,x_k,x_z,q_{ij,1},p_{k,1},w_{ijk,1}\}]\cong R_{35}$ ($p_{k,1}\sim w_{ijk,1},\ q_{ij,1}\nsim w_{ijk,1}$), a contradiction.
If $\left|E(T_{x_{k}},T_{x_{i}x_{j}x_{k}})\right|\geq 1$, then $\left|E(T_{x_{k}},T_{x_{i}x_{j}x_{k}})\right|=1$ by Lemma \ref{lem-9}\ref{lem-9:item-5} and \ref{lem-9:item-6}; otherwise, $G[\{x_i,x_j,x_k,x_z,p_{k,1},p_{k,2},w_{ijk,1}\}]\cong R_{36}$ ($p_{k,1},p_{k,2}\sim w_{ijk,1}$), a contradiction.
If $\left|E(T_{x_{k}},T_{x_{i}x_{j}x_{k}})\right|= 1$, then $\left|T_{x_{k}}\right|\leq 2$ by Lemma \ref{lem-9}\ref{lem-9:item-5} and \ref{lem-9:item-6}; otherwise, $G[\{x_i,x_j,x_k,q_{ij,1},p_{k,2},p_{k,3},w_{ijk,1}\}]\cong R_{32}$ ($p_{k,1}\sim w_{ijk,1}$), a contradiction.
If $\left|E(T_{x_{k}},T_{x_{i}x_{j}x_{k}})\right|= 1$, then $\left|T_{x_{i}}\cup T_{x_{j}}\right|=0$ by Lemma \ref{lem-9}\ref{lem-9:item-2}\ref{lem-9:item-5}\ref{lem-9:item-6} and \ref{lem-9:item-17}; otherwise, $G[\{x_i,x_j,x_z,p_{i,1},p_{k,1},w_{ijk,1},t_1\}]\cong R_{22}$ ($p_{i,1}\nsim w_{ijk,1},\ p_{k,1}\sim w_{ijk,1},\ t_1\sim w_{ijk,1}$), a contradiction.
Thus, if  $\left|E(T_{x_{k}},T_{x_{i}x_{j}x_{k}})\right|=1$, then $G\subset H_{24}$ by Lemma \ref{lem-9}\ref{lem-9:item-5} and \ref{lem-9:item-6}.
Let $\left|E(T_{x_{k}},T_{x_{i}x_{j}x_{k}})\right|=0$.
Then $\left|T_{x_{k}}\right|\leq 2$  by Lemma \ref{lem-9}\ref{lem-9:item-5} and \ref{lem-9:item-6}; otherwise, $G[\{x_i,x_j,x_k,q_{ij,1},p_{k,1},p_{k,2},p_{k,3}\}]\cong R_{23}$, a contradiction.
If $\left|T_{x_{k}}\right|\leq 1$, then  $G\subset H_{34}$ (whether $\left|E(T_{x_{i}x_{j}x_{k}} ,T_{x_{i}x_{j}})\right|=0$ or $\left|E(T_{x_{i}x_{j}x_{k}} ,T_{x_{i}x_{j}})\right|=1$) by Lemma \ref{lem-9}\ref{lem-9:item-2}\ref{lem-9:item-5}\ref{lem-9:item-6} and \ref{lem-9:item-17}.
If $\left|T_{x_{k}}\right|=2$, then  $\left|E(T_{x_{i}x_{j}x_{k}} ,T_{x_{i}x_{j}})\right|=0$ by Lemma \ref{lem-9}\ref{lem-9:item-5} and \ref{lem-9:item-6}; otherwise, $G[\{x_i,x_j,x_k,q_{ij,1},p_{k,1},p_{k,2},w_{ijk,1}\}]\cong R_{32}$ ($q_{ij,1}\sim w_{ijk,1}$), a contradiction.
If $\left|T_{x_{k}}\right|=2$, then  $\left|T_{x_{i}}\cup T_{x_{j}}\right| =0$ by Lemma \ref{lem-9}\ref{lem-9:item-2}\ref{lem-9:item-5}\ref{lem-9:item-6} and \ref{lem-9:item-17}; otherwise, $G[\{x_i,x_j,x_k,q_{ij,1},p_{i,1},p_{k,1},p_{k,2}\}]\cong R_{22}$, a contradiction.
Thus, if $\left|T_{x_{k}}\right|=2$, then  $G= A_{11}$.

Suppose that $\left|E(T_{0},T_{3})\right| =0$ and $\left|T_{x_{k}}\cup T_{x_{z}}\right|\geq 1$.
Without loss of generality, let $\left|T_{x_{k}}\right|\geq 1$.
By reasoning similar to that in Subcase 1.3,  $\left|T_{0}\right| =0$.
We distinguish the following two subsubsubcases.

\noindent$\textbf{Subsubsubcase 1.5.1.1}$ $\left|T_{3}\right| =0$.

Then $\left|E(T_{x_{k}},T_{x_{i}x_{j}})\right|\leq 2$; otherwise, $G[\{x_i,x_j,x_z,q_{ij,1},p_{k,1},p_{k,2},p_{k,3}\}]\cong R_{23}$ ($q_{ij,1}\sim p_{k,1},p_{k,2},p_{k,3}$), a contradiction.
If $\left|E(T_{x_{k}},T_{x_{i}x_{j}})\right|=2$, then $\left|T_{x_{i}}\cup T_{x_{j}}\right|=0$ by Lemma \ref{lem-9}\ref{lem-9:item-2} and \ref{lem-9:item-16}; otherwise, $G[\{x_i,x_j,x_z,q_{ij,1},p_{i,1},p_{k,1},p_{k,2}\}]\cong R_{22}$ ($q_{ij,1}\sim p_{k,1},p_{k,2}$), a contradiction.
If $\left|E(T_{x_{k}},T_{x_{i}x_{j}})\right|=2$, then $\left|T_{x_{z}}\right|=0$ by Lemma \ref{lem-9}\ref{lem-9:item-16}; otherwise, either $P_{6} \subset G$ or $G$ contains a subdivision of $K_{5}$, a contradiction.
Thus, if $\left|E(T_{x_{k}},T_{x_{i}x_{j}})\right|=2$, then $G\subset H_{21}$ by Lemma \ref{lem-9}\ref{lem-9:item-1}.
If $\left|E(T_{x_{k}},T_{x_{i}x_{j}})\right|=1$, then  $\left|E(T_{x_{z}},T_{x_{i}x_{j}})\right|=0$ and $\left|T_{x_{z}}\right| \leq 1$; otherwise, either $G[\{x_k,x_z,q_{ij,1},p_{k,1},p_{z,1},p_{z,2}\}]\cong R_{1}$ ($q_{ij,1}\sim p_{k,1},\ |T_{x_z}|\geq2$) or $G$ contains a subdivision of $K_{5}$, a contradiction.
If $\left|E(T_{x_{k}},T_{x_{i}x_{j}})\right|=1$, then $\left|T_{x_{k}}\right|\leq 3$ by Lemma \ref{lem-9}\ref{lem-9:item-1}; otherwise, $G[\{x_i,x_j,x_k,q_{ij,1},p_{k,2},p_{k,3},p_{k,4}\}]\cong R_{23}$ ($N(q_{ij,1})\cap T_{x_k}=\{p_{k,1}\}$), a contradiction.
If $\left|E(T_{x_{k}},T_{x_{i}x_{j}})\right|=1$ and $\left|T_{x_{k}}\right|\leq 2$, then $G\subset H_{33}$   by Lemma \ref{lem-9}\ref{lem-9:item-2} and \ref{lem-9:item-16}.
If $\left|E(T_{x_{k}},T_{x_{i}x_{j}})\right|=1$ and $\left|T_{x_{k}}\right|=3$, then $\left|T_{x_{i}}\cup T_{x_{j}}\right|=0$    by Lemma \ref{lem-9}\ref{lem-9:item-2} and \ref{lem-9:item-16}; otherwise, $G[\{x_i,x_j,x_k,q_{ij,1},p_{i,1},p_{k,2},p_{k,3}\}]\cong R_{22}$ ($N(q_{ij,1})\cap T_{x_k}=\{p_{k,1}\}$), a contradiction.
Thus, if $\left|E(T_{x_{k}},T_{x_{i}x_{j}})\right|=1$ and $\left|T_{x_{k}}\right|=3$, then $G\subset H_{22}$ by Lemma \ref{lem-9}\ref{lem-9:item-16}.
If $\left|E(T_{x_{k}}\cup T_{x_{z}},T_{x_{i}x_{j}})\right|=0$, then $\left|T_{x_{k}}\right|\leq 2$, $\left|T_{x_{z}}\right|\leq 2$ by Lemma \ref{lem-9}\ref{lem-9:item-1}; otherwise, $G[\{x_i,x_j,x_k,q_{ij,1},p_{k,1},p_{k,2},p_{k,3}\}]\cong R_{23}$, a contradiction.
If $\left|E(T_{x_{k}}\cup T_{x_{z}},T_{x_{i}x_{j}})\right|=0$, $\left|T_{x_{k}}\right|\leq 1$ and $\left|T_{x_{z}}\right|\leq 1$, then $G\subset H_{33}$ by Lemma \ref{lem-9}\ref{lem-9:item-2} and \ref{lem-9:item-16}.
If $\left|E(T_{x_{k}}\cup T_{x_{z}},T_{x_{i}x_{j}})\right|=0$ and $\left|T_{x_{k}}\right|=2$ (or $\left|T_{x_{z}}\right|=2$), then $\left|T_{x_{i}}\cup T_{x_{j}}\right|=0$ by Lemma \ref{lem-9}\ref{lem-9:item-2} and \ref{lem-9:item-16}; otherwise, $G[\{x_i,x_j,x_k,q_{ij,1},p_{i,1},p_{k,1},p_{k,2}\}]\cong R_{22}$, a contradiction.
Thus, if $\left|E(T_{x_{k}}\cup T_{x_{z}},T_{x_{i}x_{j}})\right|=0$ and $\left|T_{x_{k}}\right|=2$, then $G\subset H_{23}$ by Lemma \ref{lem-9}\ref{lem-9:item-16}.

\noindent$\textbf{Subsubsubcase 1.5.1.2.}$ $\left|T_{3}\right| \geq 1$.

Without loss of generality, let $\allowbreak \left|\allowbreak T_{x_{i}x_{j}x_{k}}\allowbreak\right|\allowbreak= \allowbreak1$ by Lemma \ref{lem-9}\ref{lem-9:item-1}.
Since  $G$ contains no subdivision of $K_{3,3}$, $\left|E(T_{x_{k}}\cup T_{x_{z}},T_{2})\right| =0$.
If $\left|E(T_{x_{k}},T_{x_{i}x_{j}x_{k}})\right|\geq 1$, then $\left|E(T_{x_{i}x_{j}x_{k}},T_{x_{i}x_{j}})\right|=1$ by Lemma \ref{lem-9}\ref{lem-9:item-5}; otherwise, $G[\{x_i,x_j,x_k,x_z,q_{ij,1},p_{k,1},w_{ijk,1}\}]\cong R_{35}$ ($p_{k,1}\sim w_{ijk,1},\ q_{ij,1}\nsim w_{ijk,1}$), a contradiction.
By Lemma \ref{lem-9}\ref{lem-9:item-5}, we  discuss the following two configurations.

Suppose that $\left|E(T_{x_{k}},T_{x_{i}x_{j}x_{k}})\right|=1$.
Then $\left|T_{x_{z}}\right|=0$ by Lemma \ref{lem-9}\ref{lem-9:item-5}; otherwise, $G[\{x_i,x_k,x_z,q_{ij,1},p_{z,1},w_{ijk,1}\}]\cong R_{14}$ ($q_{ij,1}\sim w_{ijk,1}$), a contradiction.
Then $\left|T_{x_{i}x_{j}x_{z}}\right|=0$  by Lemma \ref{lem-9}\ref{lem-9:item-1} and \ref{lem-9:item-5}; otherwise, either $G[\{x_i,x_k,x_z,p_{k,1},w_{ijk,1},w_{ijz,1}\}]\cong R_{18}$ ($p_{k,1}\sim w_{ijk,1}$) (if  $\left|E(T_{2},T_{x_{i}x_{j}x_{z}})\right|=0$) or $G$ contains a subdivision of $K_{5}$ (if  $\left|E(T_{2},T_{x_{i}x_{j}x_{z}})\right|=1$), a contradiction.
Then $\left|T_{x_{k}}\right|\leq 2$; otherwise, $G[\{x_i,x_j,x_k,q_{ij,1},p_{k,2},p_{k,3},w_{ijk,1}\}]\cong R_{32}$ ($p_{k,1}\sim w_{ijk,1}$), a contradiction.
If $\left|T_{x_{k}}\right|=2$, then $\left|T_{x_{i}}\cup T_{x_{j}}\right|=0$  by Lemma \ref{lem-9}\ref{lem-9:item-2} and \ref{lem-9:item-5}; otherwise, $G[\{x_i,x_j,x_k,q_{ij,1},p_{i,1},p_{k,1},p_{k,2}\}]\cong R_{22}$, a contradiction.
Thus, if  $\left|T_{x_{k}}\right|=2$, then $G\subset H_{24}$  by Lemma \ref{lem-9}\ref{lem-9:item-1}\ref{lem-9:item-5} and \ref{lem-9:item-16}.
If  $\left|T_{x_{k}}\right|=1$ and $\left|T_{x_{i}}\right|=2$, then $\left|E(T_{x_{i}},T_{x_{i}x_{j}x_{k}})\right|=1$; otherwise, either $G[\{x_i,x_k,x_z,p_{i,1},p_{i,2},p_{k,1},w_{ijk,1}\}]\cong R_{25}$ ($E(T_{x_i},w_{ijk,1})=\varnothing$) or $G[\{x_i,x_j,x_k,x_z,p_{i,1},p_{i,2},w_{ijk,1}\}]\cong R_{36}$ ($p_{i,1},p_{i,2}\sim w_{ijk,1}$), a contradiction.
Thus, if $\left|T_{x_{k}}\right|=1$, then $G\subset H_{32}$  by Lemma \ref{lem-9}\ref{lem-9:item-1}\ref{lem-9:item-2}\ref{lem-9:item-5} and \ref{lem-9:item-16}.

We now suppose that $\left|E(T_{x_{k}},T_{x_{i}x_{j}x_{k}})\right|=0$ (and $\left|E(T_{x_{z}},T_{x_{i}x_{j}x_{z}})\right|=0$).
Then $\left|T_{x_{k}}\right|\leq 2$ by Lemma \ref{lem-9}\ref{lem-9:item-1}; otherwise, $G[\{x_i,x_j,x_k,q_{ij,1},p_{k,1},p_{k,2},p_{k,3}\}]\cong R_{23}$, a contradiction.
If $\left|T_{x_{k}}\right|=2$, then $\left|T_{x_{i}}\cup T_{x_{j}}\right|=0$  by Lemma \ref{lem-9}\ref{lem-9:item-2} and \ref{lem-9:item-16}; otherwise, $G[\{x_i,x_j,x_k,q_{ij,1},p_{i,1},p_{k,1},p_{k,2}\}]\cong R_{22}$, a contradiction.
If $\left|T_{x_{k}}\right|=2$, then we have $\left|T_{x_{i}x_{j}x_{z}}\right|=0$ by Lemma \ref{lem-9}\ref{lem-9:item-5}; otherwise, $G[\{x_i,x_j,x_k,x_z,p_{k,1},p_{k,2},w_{ijz,1}\}]\cong R_{32}$, a contradiction.
If $\left|T_{x_{k}}\right|=2$, then $\left|E(T_{x_{i}x_{j}},T_{x_{i}x_{j}x_{k}})\right|=0$; otherwise, $G[\{x_i,x_j,x_k,q_{ij,1},p_{k,1},p_{k,2},w_{ijk,1}\}]\cong R_{32}$ ($q_{ij,1}\sim w_{ijk,1}$), a contradiction.
Thus, if $\left|T_{x_{k}}\right|=2$, then $G\subset A_{11}$ by Lemma \ref{lem-9}\ref{lem-9:item-5} and \ref{lem-9:item-16}.
Let $\left|T_{x_{k}}\right|=1$.
If $\left|E(T_{x_{i}x_{j}},T_{x_{i}x_{j}x_{k}})\right|=1$, then $\left|T_{x_{z}}\right|=0$ by Lemma \ref{lem-9}\ref{lem-9:item-5} and \ref{lem-9:item-16}; otherwise, $G[\{x_i,x_k,x_z,q_{ij,1},p_{z,1},w_{ijk,1}\}]\cong R_{14}$ ($q_{ij,1}\sim w_{ijk,1}$), a contradiction.
If $\left|E(T_{x_{i}x_{j}},T_{x_{i}x_{j}x_{k}})\right|=1$, then $\left|E(T_{x_{i}}\cup T_{x_{j}},T_{x_{i}x_{j}x_{k}})\right|=0$ by Lemma \ref{lem-9}\ref{lem-9:item-2}\ref{lem-9:item-5} and \ref{lem-9:item-16}; otherwise, $G[\{x_i,x_k,x_z,q_{ij,1},p_{i,1},p_{k,1},w_{ijk,1}\}]\cong R_{28}$ ($q_{ij,1}\sim w_{ijk,1},\ p_{i,1}\sim w_{ijk,1}$), a contradiction.
If $\left|E(T_{x_{i}x_{j}},T_{x_{i}x_{j}x_{k}})\right|=1$, then $\left|E(T_{x_{i}}\cup T_{x_{j}},T_{x_{i}x_{j}x_{z}})\right|=0$ by Lemma \ref{lem-9}\ref{lem-9:item-2}\ref{lem-9:item-5} and \ref{lem-9:item-16}; otherwise, $G[\{x_i,x_j,x_k,p_{i,1},w_{ijk,1},w_{ijz,1}\}]\cong R_{20}$ ($p_{i,1}\sim w_{ijz,1},\ p_{i,1}\nsim w_{ijk,1}$), a contradiction.
Thus, if $\left|E(T_{x_{i}x_{j}},T_{x_{i}x_{j}x_{k}})\right|=1$, then $G\subset H_{31}$ by Lemma \ref{lem-9}\ref{lem-9:item-2}\ref{lem-9:item-5}\ref{lem-9:item-6} and \ref{lem-9:item-16}.
Let $\left|E(T_{x_{i}x_{j}} ,T_{x_{i}x_{j}x_{k}}\cup T_{x_{i}x_{j}x_{z}})\right|=0$.
If $\left|T_{x_{z}}\right|=1$, then $\left|E(T_{x_{i}}\cup T_{x_{j}},T_{x_{i}x_{j}x_{k}}\cup T_{x_{i}x_{j}x_{z}})\right|=0$ by Lemma \ref{lem-9}\ref{lem-9:item-2} and \ref{lem-9:item-16}; otherwise, $G[\{x_i,x_j,x_z,p_{i,1},p_{z,1},w_{ijk,1}\}]\cong R_{14}$ ($p_{i,1}\sim w_{ijk,1}$), a contradiction.
Thus, if $\left|T_{x_{z}}\right|=1$, then either $G\subset A_{12}$ (if $\left|T_{3}\right|=1$) or $G\subset H_{34}$ (if $\left|T_{3}\right|=2$) by Lemma \ref{lem-9}\ref{lem-9:item-1}\ref{lem-9:item-2}\ref{lem-9:item-5}\ref{lem-9:item-6} and \ref{lem-9:item-16}.
Let $\left|T_{x_{z}}\right|=0$.
If $\left|T_{x_{i}x_{j}x_{z}}\right|=1$, then $\left|E(T_{x_{i}}\cup T_{x_{j}},T_{x_{i}x_{j}x_{k}}\cup T_{x_{i}x_{j}x_{z}})\right|=0$; otherwise, either $G[\{x_i,x_j,x_k,p_{i,1},w_{ijk,1},w_{ijz,1}\}]\cong R_{20}$ ($p_{i,1}\sim w_{ijz,1},\ p_{i,1}\nsim w_{ijk,1}$) or  $G$ contains a subdivision of $K_{5}$, a contradiction.
Thus, if $\left|T_{x_{i}x_{j}x_{z}}\right|=1$, then $G\subset H_{34}$ by Lemma \ref{lem-9}\ref{lem-9:item-1}\ref{lem-9:item-2}\ref{lem-9:item-5} and \ref{lem-9:item-16}.
If $\left|T_{x_{i}x_{j}x_{z}}\right|=0$, $\left|T_{x_{j}}\right|=2$ and $\left|E(T_{x_{i}},T_{x_{i}x_{j}x_{k}})\right|\geq 1$, then $\left|E(T_{x_{j}},T_{x_{i}x_{j}x_{k}})\right|=1$ by Lemma \ref{lem-9}\ref{lem-9:item-2}\ref{lem-9:item-5} and \ref{lem-9:item-16}; otherwise, $G[\{x_i,x_j,x_z,p_{i,1},p_{j,1},p_{j,2},w_{ijk,1}\}]\cong R_{25}$ ($p_{i,1}\sim w_{ijk,1},\ p_{j,1},p_{j,2}\nsim w_{ijk,1}$), a contradiction.
Thus, if $\left|T_{x_{i}x_{j}x_{z}}\right|=0$, then either $G\subset H_{30}$ (if $\left|E(T_{x_{i}}\cup T_{x_{j}},T_{x_{i}x_{j}x_{k}})\right|\geq1$) or $G\subset A_{12}$ (if $\left|E(T_{x_{i}}\cup T_{x_{j}},T_{x_{i}x_{j}x_{k}})\right|=0$) by Lemma \ref{lem-9}\ref{lem-9:item-2}\ref{lem-9:item-5}\ref{lem-9:item-6}\ref{lem-9:item-16} and \ref{lem-9:item-17}.

\noindent$\textbf{Subsubcase 1.5.2.}$ $\left|T_{x_{i}x_{k}x_{z}}\cup T_{x_{j}x_{k}x_{z}}\right|\geq 1$.

Without loss of generality, let $\left|T_{x_{j}x_{k}x_{z}}\right|\geq 1$.
Then $\left|E(T_{x_{i}x_{j}},T_{3})\right|=0$; otherwise, either $G[\{x_j,x_k,x_z,q_{ij,1},w_{ijk,1},w_{jkz,1}\}]\cong R_{20}$ ($q_{ij,1}\sim w_{ijk,1}$) (if $\left|E(T_{x_{i}x_{j}},T_{x_{i}x_{j}x_{k}})\right|=1$) or  $G$ contains a subdivision of $K_{5}$ (if $\left|E(T_{x_{i}x_{j}},T_{x_{j}x_{k}x_{z}})\right|=1$), a contradiction.
Then $\left|T_{0}\right|=0$ by Lemma \ref{lem-9}\ref{lem-9:item-5} and \ref{lem-9:item-6}; otherwise, $R_{14}\subset G$, a contradiction.
Then $\left|E(T_{x_{k}}\cup T_{x_{z}}, T_{x_{j}x_{k}x_{z}}\cup T_{x_{i}x_{k}x_{z}})\right|=0$; otherwise, either $G[\{x_i,x_k,x_z,q_{ij,1},p_{k,1},w_{jkz,1}\}]\cong R_{14}$ ($p_{k,1}\sim w_{jkz,1}$) (if $\left|E(T_{x_k},T_{x_ix_j})\right|=0$) or $G$ contains a subdivision of $K_{5}$ (if $\left|E(T_{x_k},T_{x_ix_j})\right|=1$), a contradiction.
Then $\left|E(T_{x_{k}}\cup T_{x_{z}}, T_{x_{i}x_{j}})\right|=0$; otherwise, $G[\{x_i,x_j,x_z,q_{ij,1},p_{k,1},w_{jkz,1}\}]\cong R_{14}$ ($p_{k,1}\sim q_{ij,1}$), a contradiction.
Then $\left|T_{x_{k}}\right|\leq 1$ and $\left|T_{x_{z}}\right|\leq 1$; otherwise, $G[\{x_i,x_j,x_k,q_{ij,1},p_{k,1},p_{k,2},w_{jkz,1}\}]\cong R_{25}$, a contradiction.
By Lemma \ref{lem-9}\ref{lem-9:item-5},  we discuss the following three configurations.

We first suppose that $\left|T_{x_{i}}\right|=1$.
Then $\left|T_{x_{k}}\cup T_{x_{z}}\right|=0$ by Lemma \ref{lem-9}\ref{lem-9:item-2} and \ref{lem-9:item-16}; otherwise, either $G[\{x_i,x_j,x_z,q_{ij,1},p_{k,1},w_{jkz,1}\}]\cong R_{14}$ ($p_{k,1}\sim q_{ij,1}$) or $G[\{x_i,x_k,x_z,q_{ij,1},p_{i,1},p_{k,1},w_{jkz,1}\}]\cong R_{22}$ ($p_{k,1}\nsim q_{ij,1}$), a contradiction.
Then  $\left|T_{x_{i}x_{j}x_{k}}\cup T_{x_{i}x_{j}x_{z}}\right|=0$ by Lemma \ref{lem-9}\ref{lem-9:item-1}\ref{lem-9:item-2} and \ref{lem-9:item-5}; otherwise, $G[\{x_i,x_k,x_z,q_{ij,1},p_{i,1},w_{ijk,1},w_{jkz,1}\}]\cong R_{25}$, a contradiction.
If $\left|T_{x_{i}x_{k}x_{z}}\right|=1$, then we have $\left|E(T_{x_{i}}\cup T_{x_{j}}, T_{x_{j}x_{k}x_{z}}\cup T_{x_{i}x_{k}x_{z}})\right|=0$ by Lemma \ref{lem-9}\ref{lem-9:item-2} and \ref{lem-9:item-5}; otherwise, either $G[\{x_i,x_j,x_k,p_{i,1},w_{ikz,1},w_{jkz,1}\}]\cong R_{18}$ ($p_{i,1}\sim w_{ikz,1}$) or $G[\{x_i,x_j,x_k,p_{j,1},w_{ikz,1},w_{jkz,1}\}]\cong R_{18}$ ($p_{j,1}\sim w_{jkz,1}$), a contradiction.
Thus, if  $\left|T_{x_{i}x_{k}x_{z}}\right|=1$, then  $G\subset H_{24}$ by Lemma \ref{lem-9}\ref{lem-9:item-1}\ref{lem-9:item-2}\ref{lem-9:item-5} and \ref{lem-9:item-5}.
If $\left|T_{x_{i}x_{k}x_{z}}\right|=0$,  then  $\left|E(T_{x_{j}x_{k}x_{z}}, T_{x_{i}}\cup T_{x_{j}})\right|=0$ by Lemma \ref{lem-9}\ref{lem-9:item-2} and \ref{lem-9:item-5}; otherwise, $G[\{x_i,x_j,x_k,p_{i,1},p_{j,1},w_{jkz,1}\}]\cong R_{14}$ ($p_{j,1}\sim w_{jkz,1}$), a contradiction.
If $\left|T_{x_{i}x_{k}x_{z}}\right|=0$,  then  $\left|T_{x_{j}}\right|\leq 3$  by Lemma \ref{lem-9}\ref{lem-9:item-2} and \ref{lem-9:item-5}; otherwise, $G[K_4 \cup \{p_{j,1},\ldots,p_{j,4},w_{jkz,1}\}]\cong R_{54}$, a contradiction.
Thus, if $\left|T_{x_{i}x_{k}x_{z}}\right|=0$,  then $G\subset H_{1}$  by Lemma \ref{lem-9}\ref{lem-9:item-2} and \ref{lem-9:item-5}.

Suppose that $\left|T_{x_{i}}\right|=0$ and $\left|T_{x_{k}}\right|=  1$.
Note that $\left|E(T_{x_{i}x_{j}},T_{3})\right|=0$.
Then $\left|T_{x_{i}x_{j}x_{z}}\right|=0$ by Lemma \ref{lem-9}\ref{lem-9:item-5}; otherwise, $G[\{x_i,x_j,x_k,q_{ij,1},p_{k,1},w_{ijz,1},w_{jkz,1}\}]\cong R_{28}$, a contradiction.
Then  $\left|E(T_{x_{k}}, T_{x_{i}x_{j}x_{k}})\right|=0$ by Lemma \ref{lem-9}\ref{lem-9:item-1}; otherwise, $G[\{x_i,x_j,x_z,p_{k,1},w_{ijk,1},w_{jkz,1}\}]\cong R_{14}$ ($p_{k,1}\sim w_{ijk,1}$), a contradiction.
Then  $\left|T_{x_{z}}\cup T_{x_{i}x_{j}x_{k}}\right|\leq 1$ by Lemma \ref{lem-9}\ref{lem-9:item-1} and \ref{lem-9:item-5}; otherwise, either $G[\{x_i,x_j,x_z,q_{ij,1},p_{z,1},p_{z,2},w_{jkz,1}\}]\cong R_{25}$ ($|T_{x_z}|\geq2$) or $G[\{x_i,x_j,x_z,q_{ij,1},p_{z,1},w_{ijk,1},w_{jkz,1}\}]\cong R_{28}$ ($|T_{x_z}|=|T_{x_ix_jx_k}|=1$), a contradiction.
Then  $\left|T_{x_{j}}\cup T_{x_{i}x_{k}x_{z}}\right|\leq 1$ by Lemma \ref{lem-9}\ref{lem-9:item-1} and \ref{lem-9:item-5}; otherwise, either  $G[\{x_j,x_k,x_z,q_{ij,1},p_{j,1},p_{k,1},w_{ikz,1}\}]\cong R_{22}$ ($|T_{x_j}|=\ |T_{x_ix_kx_z}|=1$) or $G[K_4 \cup \{q_{ij,1},p_{j,1},p_{j,2},p_{k,1},w_{jkz,1}\}]\cong R_{57}$ ($|T_{x_j}|\geq2,\ E(\{p_{j,1},p_{j,2}\},w_{jkz,1})=\varnothing$), a contradiction.
Thus, 
$G\subset H_{3}$ (if $\left|T_{x_{i}x_{j}x_{k}}\right|=\left|T_{x_{i}x_{k}x_{z}}\right|=0$) or 
$G\subset H_{26}$ (if $\left|T_{x_{i}x_{k}x_{z}}\right|=0$ and $\left|T_{x_{i}x_{j}x_{k}}\right|=1$) or 
$G\subset H_{32}$ (if $\left|T_{x_{i}x_{k}x_{z}}\right|=1$ and $\left|T_{x_{i}x_{j}x_{k}}\right|=0$) or 
$G\subset H_{37}$ (if $\left|T_{x_{i}x_{j}x_{k}}\right|=\left|T_{x_{i}x_{k}x_{z}}\right|=1$) 
by Lemma \ref{lem-9}\ref{lem-9:item-1}\ref{lem-9:item-2}\ref{lem-9:item-5} and \ref{lem-9:item-16}.

Suppose that $\left|T_{x_{i}}\right|=\left|T_{x_{k}}\right|=\left|T_{x_{z}}\right|=0$.
Obviously, if  $\left|T_{x_{j}}\right|=0$, then $G\subset H_{35}$ by Lemma \ref{lem-9}\ref{lem-9:item-5}.
Let $\left|T_{x_{j}}\right|\geq 1$.
By an argument similar to that for the situation $\left|T_{x_{i}}\right|=1$, if $\left|T_{x_{i}x_{k}x_{z}}\right|=1$, then $G\subset H_{24}$.
Let $\left|T_{x_{i}x_{k}x_{z}}\right|=0$.
By an argument similar to that used in Case $3$ of Lemma \ref{lem-10}, we obtain that $G[T_{x_{j}}\cup T_{x_{j}x_{k}x_{z}}]\subset 3K_{1}\cup K_{2}$ (see $H_{1}$ in Fig.~\ref{fig-4}).
If $\left|E(T_{x_{j}}, T_{x_{j}x_{k}x_{z}})\right|=1$, then $\left|T_{x_{i}x_{j}x_{k}}\cup T_{x_{i}x_{j}x_{z}}\right|=0$ by Lemma \ref{lem-9}\ref{lem-9:item-1} and \ref{lem-9:item-2}; otherwise, either $G[\{x_i,x_j,x_k,p_{j,1},w_{ijk,1},w_{jkz,1}\}]\cong R_{20}$ ($p_{j,1}\sim w_{jkz,1},\ p_{j,1}\nsim w_{ijk,1}$) (if $\left|E(T_{x_{j}},T_{x_{i}x_{j}x_{k}})\right|=0$) or  $G$ contains a subdivision of $K_{5}$ (if $\left|E(T_{x_{j}},T_{x_{i}x_{j}x_{k}})\right|=1$), a contradiction.
Thus, if $\left|E(T_{x_{j}}, T_{x_{j}x_{k}x_{z}})\right|=1$, then $G\subset H_{66}$  by Lemma \ref{lem-9}\ref{lem-9:item-1} and \ref{lem-9:item-2}.
If $\left|E(T_{x_{j}}, T_{x_{j}x_{k}x_{z}})\right|=0$, then   $\left|E(T_{x_{j}}, T_{x_{i}x_{j}x_{k}}\cup T_{x_{i}x_{j}x_{z}})\right|=0$ by Lemma \ref{lem-9}\ref{lem-9:item-2}; otherwise, $G[\{x_i,x_k,x_z,p_{j,1},w_{ijk,1},w_{jkz,1}\}]\cong R_{14}$ ($p_{j,1}\sim w_{ijk,1},\ p_{j,1}\nsim w_{jkz,1}$), a contradiction.
Thus,  if $\left|E(T_{x_{j}}, T_{x_{j}x_{k}x_{z}})\right|=0$,  then $G\subset H_{64}$ by Lemma \ref{lem-9}\ref{lem-9:item-1} and \ref{lem-9:item-2}.

\noindent$\textbf{Case 2.}$ For any $i,j$ ($0\leq i< j\leq 3$), $\left|T_{2}\right|>\left|T_{x_{i}x_{j}}\right|$.

\begin{fact}\label{f-13}
If $\left|T_{x_{i}x_{j}}\right| \neq 0$, then $\left|T_{x_{k}x_{z}}\right| =0$.
\end{fact}

\begin{proof}
Suppose that $\left|T_{x_{i}x_{j}}\right| \geq1$ and $\left|T_{x_{k}x_{z}}\right|  \geq1$.
If $\left|T_{x_{i}}\right|\geq 1$, then  either $R_{12}\subset G$ (if $\left|E(T_{x_{i}}, T_{x_{k}x_{z}})\right| \geq 1$, then deleting $x_{z}$) or $R_{31}\subset G$ (if $\left|E(T_{x_{i}}, T_{x_{k}x_{z}})\right| =0$) by Lemma \ref{lem-9}\ref{lem-9:item-2} and \ref{lem-9:item-11}, a contradiction.
If $\left|T_{x_{i}x_{j}x_{k}}\right|\geq 1$, then   $R_{18}\subset G$ (if $\left|E(T_{x_{i}x_{j}}, T_{x_{i}x_{j}x_{k}})\right| \geq 1$ and $\left|E(T_{x_{k}x_{z}}, T_{x_{i}x_{j}x_{k}})\right| = 0$, then deleting $x_{j}$) or $R_{35}\subset G$ (if $\left|E(T_{x_{i}x_{j}}\cup T_{x_{k}x_{z}}, T_{x_{i}x_{j}x_{k}})\right| = 0$) or $G$ contains a subdivision of $K_{5}$ (if $\left|E(T_{x_{k}x_{z}}, T_{x_{i}x_{j}x_{k}})\right| \geq 1$) by Lemma \ref{lem-9}\ref{lem-9:item-11}, a contradiction.
If $\left|E(T_{x_{i}x_{j}},T_{0})\right|\geq 1$, then   $R_{6}\subset G$ (deleting $x_{j}$) by Lemma \ref{lem-9}\ref{lem-9:item-6} and \ref{lem-9:item-11}, a contradiction.
This completes the proof of (F\ref{f-13}) by symmetry.
\end{proof}

\begin{fact}\label{f-14}
If $\left|E(T_{x_{i}x_{j}},T_{x_{j}x_{k}})\right| =1$, then $\left|T_{x_{k}x_{z}}\cup T_{x_{i}x_{z}} \cup T_{3}\cup T_{x_{z}}\right|=0$.
\end{fact}

\begin{proof}
If $\left|T_{x_{i}x_{j}x_{z}}\right|\geq 1$, then either $R_{12}\subset G$ (if $\left|E(T_{x_{i}x_{j}x_{z}}, T_{x_{i}x_{j}}\cup T_{x_{j}x_{k}} )\right| =0$, then deleting $x_{j}$) or $G$ contains a subdivision of $K_{5}$ (if $\left|E(T_{x_{i}x_{j}x_{z}}, T_{x_{i}x_{j}}\cup T_{x_{j}x_{k}} )\right| \geq 1$), a contradiction.
If $\left|T_{x_{i}x_{j}x_{k}}\right|\geq 1$, then $G$ contains a subdivision of $K_{3,3}$, a contradiction.
If $\left|T_{x_{z}}\right|\geq 1$, then either $R_{14}\subset G$ (if $\left|E(T_{x_{z}}, T_{x_{i}x_{j}}\cup T_{x_{j}x_{k}} )\right| =0$, then deleting $x_{j}$) or $G$ contains a subdivision of $K_{5}$ (if $\left|E(T_{x_{z}}, T_{x_{i}x_{j}}\cup T_{x_{j}x_{k}} )\right| \geq 1$), a contradiction.
This completes the proof of (F\ref{f-14}) by Lemma \ref{lem-9}\ref{lem-9:item-11}.
\end{proof}

By Lemma \ref{lem-9}\ref{lem-9:item-1} and \ref{lem-9:item-16}, we immediately arrive at the following one.

\begin{fact}\label{f-15}
$\left|E(T_{1})\right| =0$.
\end{fact}

\begin{fact}\label{f-16}
$\left|E(T_{2},T_{0})\right| =0$.
\end{fact}
\begin{proof}
Otherwise, either $R_{7}\subset G$ (if $\left|E(T_{x_{i}x_{j}},T_{x_{j}x_{k}})\right| =1$, then deleting $x_{j}$) or $R_{14} \subset G$ (if $\left|E(T_{x_{i}x_{j}},T_{x_{j}x_{k}})\right| =0$, then deleting $x_{z}$)   by Lemma \ref{lem-9}\ref{lem-9:item-6}, a contradiction.
\end{proof}

Let $|T_{x_i x_j}| = \max\limits_{T_{x_{g}x_{h}} \subseteq T_2} |T_{x_{g}x_{h}}|$.
Without loss of generality, let  $\left|T_{x_{j}x_{k}}\right|\geq 1$ and $|T_{x_j x_k}| = \max\limits_{T_{x_{g}x_{h}} \subseteq T_2 \setminus T_{x_i x_j}} |T_{x_{g}x_{h}}|$ by (F\ref{f-13}).
We  distinguish the following two subcases by (F\ref{f-l34.1}).

\noindent$\textbf{Subcase 2.1.}$  For some $g,h$ ($0\leq g< h\leq 3$), $\left|E(T_{x_{g}x_{h}})\right| =1$.

If $\left|E(T_{x_{i}x_{j}})\right| =0$, then  $G$ contains a subdivision of $K_{3,3}$ by Lemma \ref{lem-9}\ref{lem-9:item-13}, a contradiction.
Thus, $\left|E(T_{x_{i}x_{j}})\right| =1$ by (F\ref{f-l34.1}).
Then  $\left|E(T_{x_{i}x_{j}},T_{x_{j}x_{k}})\right| =1$ and $\left|T_{x_{j}x_{k}}\right| =1$ by Lemma \ref{lem-9}\ref{lem-9:item-13}; otherwise, $G$ contains a subdivision of $K_{3,3}$, a contradiction.
Then $\left|T_{x_{j}x_{z}}\right| =0$ by Lemma \ref{lem-9}\ref{lem-9:item-13}; otherwise, $G$ contains a subdivision of $K_{5}$, a contradiction.
Then $\left|T_{0}\right| =0$ by (F\ref{f-14}), (F\ref{f-16}) and Lemma \ref{lem-9}\ref{lem-9:item-6}.
Let $q_{ij,1}\in T_{x_{i}x_{j}}$ and $\left|E(q_{ij,1},T_{x_{j}x_{k}})\right| =1$.
Then $\left|T_{x_{k}}\right|\leq 1$ by Lemma \ref{lem-9}\ref{lem-9:item-1}\ref{lem-9:item-2}\ref{lem-9:item-2} and \ref{lem-9:item-13}; otherwise, $G[\{x_i,\allowbreak x_j,\allowbreak x_k,\allowbreak q_{ij,1},\allowbreak q_{jk,1},\allowbreak p_{k,1},\allowbreak p_{k,2}\}]\cong R_{25}$, a contradiction.

If $\left|T_{x_{k}}\right|= 1$, then $\left|T_{x_{i}x_{j}}\right| =2$ by  Lemma \ref{lem-9}\ref{lem-9:item-2}\ref{lem-9:item-3} and  \ref{lem-9:item-13}; otherwise, $G[\{x_i,\allowbreak x_j,\allowbreak x_k,\allowbreak q_{ij,1},\allowbreak q_{ij,2},\allowbreak q_{jk,1},\allowbreak p_{k,1}\}]\cong R_{28}$, a contradiction.
If $\left|T_{x_{k}}\right|= 1$, then $\left|E(T_{x_{i}}\cup T_{x_{j}},T_{x_{i}x_{j}}\cup T_{x_{j}x_{k}}\cup T_{x_{i}x_{k}})\right| =0$ by Lemma \ref{lem-9}\ref{lem-9:item-2} and (F\ref{f-15}); otherwise, $G$ contains a subdivision of $K_{3,3}$, a contradiction.
If $\left|T_{x_{k}}\right|= 1$, then $\left|T_{x_{i}}\right| \leq 1$; otherwise, $G[\{x_i,\allowbreak x_j,\allowbreak x_z,\allowbreak q_{ij,1},\allowbreak q_{jk,1},\allowbreak p_{i,1},\allowbreak p_{i,2}\}]\cong R_{25}$ ($E(\{q_{ij,1},\allowbreak q_{jk,1},\allowbreak p_{i,1},\allowbreak p_{i,2}\})=\{q_{ij,1}q_{jk,1}\}$), a contradiction.
By an argument similar to the one above, if $\left|T_{x_{k}}\right|=1$ and $\left|T_{x_{i}x_{k}}\right|=1$, then $\left|T_{x_{j}}\right|\leq1$.
If $\left|T_{x_{k}}\right|= 1$, then $\left|E(T_{x_{i}x_{k}},T_{x_{j}x_{k}})\right| =0$ by Lemma \ref{lem-9}\ref{lem-9:item-13}; otherwise, $G$ contains a subdivision of $K_{3,3}$, a contradiction.
Thus, if $\left|T_{x_{k}}\right|= 1$, then either $G\subset H_{40}$ (if $\left|T_{x_{i}x_{k}}\right|=1$)  or $G\subset H_{50}$ (if $\left|T_{x_{i}x_{k}}\right|=0$) by (F\ref{f-14}), (F\ref{f-15}) and Lemma \ref{lem-9}\ref{lem-9:item-16}.
Let $\left|T_{x_{k}}\right|= 0$.

If $\left|E(T_{x_{i}x_{k}},T_{x_{j}x_{k}})\right| =1$, then  $\left|T_{x_{i}}\cup T_{x_{j}}\right| =0$ by Lemma \ref{lem-9}\ref{lem-9:item-2} and \ref{lem-9:item-13}; otherwise, either $R_{26} \subset G$ (if $\left|E(T_{x_{i}},T_{x_{j}x_{k}})\right| =0$) or $R_{30} \subset G$ (if $\left|E(T_{x_{i}},T_{x_{j}x_{k}})\right| =1$), a contradiction.
Thus, $\left|E(T_{x_{i}x_{k}},T_{x_{j}x_{k}})\right| =0$ by (F\ref{f-14}) and $\left|T-T_{2}\right|\geq 1$.
If $\left|T_{x_{i}x_{k}}\right|=1$, then $\left|E(T_{x_{i}}\cup T_{x_{j}},T_{x_{i}x_{j}}\cup T_{x_{j}x_{k}}\cup T_{x_{i}x_{k}})\right| =0$ by Lemma \ref{lem-9}\ref{lem-9:item-2} and  $G$ contains no subdivision of $K_{3,3}$.
Thus, if $\left|T_{x_{i}x_{k}}\right|=1$, then $\left|T_{x_{i}}\right| \leq 1$ and $\left|T_{x_{j}}\right| \leq 1$ by Lemma \ref{lem-9}\ref{lem-9:item-2}; otherwise, either $G[\{x_i,\allowbreak x_j,\allowbreak x_z,\allowbreak q_{ij,1},\allowbreak q_{jk,1},\allowbreak p_{i,1},\allowbreak p_{i,2}\}]\cong R_{25}$ ($E(\{q_{ij,1},\allowbreak q_{jk,1},\allowbreak p_{i,1},\allowbreak p_{i,2}\})=\{q_{ij,1}q_{jk,1}\}$) or $G[\{x_i,\allowbreak x_j,\allowbreak x_z,\allowbreak q_{ij,1},\allowbreak q_{ik,1},\allowbreak p_{j,1},\allowbreak p_{j,2}\}]\cong R_{25}$ ($E(\{q_{ij,1},\allowbreak q_{ik,1},\allowbreak p_{j,1},\allowbreak p_{j,2}\})=\{q_{ij,1}q_{ik,1}\}$), a contradiction.
If $\left|T_{x_{i}x_{k}}\right|=1$, then $\left|T_{x_{i}x_{j}}\right|\leq 4$ by Lemma \ref{lem-9}\ref{lem-9:item-3} and \ref{lem-9:item-2}; otherwise, $G[\{x_i,\allowbreak x_k,\allowbreak q_{ij,1},\ldots,q_{ij,4},\allowbreak q_{jk,1},\allowbreak q_{ik,1}\}]\cong R_{38}$ ($E(\{q_{ij,1},\allowbreak q_{ij,2},\allowbreak q_{ij,3},\allowbreak q_{ij,4},\allowbreak q_{jk,1},\allowbreak q_{ik,1}\})=\{q_{ij,1}q_{jk,1},\allowbreak q_{ij,1}q_{ik,1}\}$), a contradiction.
If $\left|T_{x_{i}x_{k}}\right|=1$ and $\left|T_{x_{i}x_{j}}\right|=4$, then  $\left|T_{x_{i}}\cup T_{x_{j}}\right|=0$  by Lemma \ref{lem-9}\ref{lem-9:item-2}; otherwise, $G[\{x_i,\allowbreak x_j,\allowbreak x_k,\allowbreak q_{ij,1},q_{ij,2},q_{ij,3},\allowbreak q_{jk,1},\allowbreak p_{i,1}\}]\cong R_{43}$ ($E(\{q_{ij,1},q_{ij,2},q_{ij,3},\allowbreak q_{jk,1},\allowbreak p_{i,1}\})=\{q_{ij,1}q_{jk,1}\}$), a contradiction.
Thus, if $\left|T_{x_{i}x_{k}}\right|=1$, then $\left|T_{x_{i}x_{j}}\right| \leq 3$ for $\left|T-T_{2}\right|\geq 1$.
Thus, if $\left|T_{x_{i}x_{k}}\right|=1$, then $G\subset H_{36}$ by (F\ref{f-14}) and (F\ref{f-15}).
Let $\left|T_{x_{i}x_{k}}\right|=0$.
If $\left|T_{x_{j}}\right| \geq 3$, then $\left|T_{x_{i}}\right|=0$ by  Lemma \ref{lem-9}\ref{lem-9:item-2} and (F\ref{f-15}); otherwise, either $R_{38} \subset G$ (if $\left|E(T_{x_{i}},T_{x_{j}x_{k}})\right|=1$) or $R_{51} \subset G$ (if $\left|E(T_{x_{i}},T_{x_{j}x_{k}})\right|=0$), a contradiction.
If $\left|T_{x_{j}}\right| \geq 1$, then $\left|T_{x_{i}x_{j}}\right|\leq 5$; otherwise, $G[\{x_i,\allowbreak x_k,\allowbreak q_{ij,1},\ldots,q_{ij,5},\allowbreak q_{jk,1}\}]\cong R_{37}$ ($E(\{q_{ij,1},\ldots,\allowbreak q_{ij,5},\allowbreak q_{jk,1}\})=\{q_{ij,1}q_{jk,1}\}$), a contradiction.
Thus, if  $\left|T_{x_{j}}\right| \geq 3$ or $\left|T_{x_{i}}\right|=0$, then $G\subset H_{45}$ by (F\ref{f-14}), (F\ref{f-15}) and  Lemma \ref{lem-9}\ref{lem-9:item-1}\ref{lem-9:item-2} and \ref{lem-9:item-13}.
Let $\left|T_{x_{j}}\right|  \leq 2$ and $\left|T_{x_{i}}\right|\geq 1$.
Then $\left|T_{x_{i}x_{j}}\right|\leq 3$ by (F\ref{f-15}) and  Lemma \ref{lem-9}\ref{lem-9:item-2} and \ref{lem-9:item-13}; otherwise, either $R_{42} \subset G$ (if $\left|E(T_{x_{i}},T_{x_{j}x_{k}})\right|=1$) or $R_{43} \subset G$ (if $\left|E(T_{x_{i}},T_{x_{j}x_{k}})\right|=0$), a contradiction.
If $\left|T_{x_{i}}\right|\geq 2$, then $\left|E(T_{x_{i}},T_{x_{j}x_{k}})\right| =1$ by (F\ref{f-15}) and Lemma \ref{lem-9}\ref{lem-9:item-2} and \ref{lem-9:item-13}; otherwise, $R_{25} \subset G$ (whether $\left|E(T_{x_{i}},T_{x_{j}x_{k}})\right| =0$ or $\left|E(T_{x_{i}},T_{x_{j}x_{k}})\right|  \geq2$), a contradiction.
Thus, $\left|T_{x_{i}}\right|\leq 2$; otherwise, $R_{25} \subset G$, a contradiction.
Thus, $G\subset H_{55}$ by (F\ref{f-14}), (F\ref{f-15}) and  Lemma \ref{lem-9}\ref{lem-9:item-2} and \ref{lem-9:item-13}.

\noindent$\textbf{Subcase 2.2.}$  For any $g,h$ ($0\leq g< h\leq 3$), $\left|E(T_{x_{g}x_{h}})\right| =0$.

If $\left|T_{x_{i}x_{j}}\right| \geq 5$ and $\left|T_{x_{j}x_{k}}\right| \geq 1$, then we have either $G[\{x_i,\allowbreak x_k,\allowbreak x_z,\allowbreak q_{ij,1},\ldots,q_{ij,5},\allowbreak q_{jk,1}\}]\cong R_{46}$ (if $\left|E(T_{x_{i}x_{j}},T_{x_{j}x_{k}})\right|=0$) or $G[\{x_i,\allowbreak x_k,\allowbreak q_{ij,1},\ldots,q_{ij,5},\allowbreak q_{jk,1}\}]\cong R_{37}$ (if $\left|E(T_{x_{i}x_{j}},T_{x_{j}x_{k}})\right|=1$) by  Lemma \ref{lem-9}\ref{lem-9:item-11}, a contradiction.
Thus, $\left|T_{x_{i}x_{j}}\right| \leq 4$ by (F\ref{f-13}).
Without loss of generality, we  distinguish the following four subsubcases.

\noindent$\textbf{Subsubcase 2.2.1.}$ $\left|T_{x_{i}x_{j}}\right|=4$.

By an argument similar to that used in subcase 2.1 of Lemma \ref{lem-11} and (F\ref{f-13}), $\left|T_{x_{k}x_{z}}\right|=\left|T_{x_{i}x_{z}}\right|=0$, and $\left|T_{x_{j}x_{k}}\right|=1$.
Then $\left|T_{x_{i}}\right|=0$ by  Lemma \ref{lem-9}\ref{lem-9:item-2} and \ref{lem-9:item-11}; otherwise, either $R_{37} \subset G$ (if $\left|E(T_{x_{i}},T_{x_{j}x_{k}})\right| =\left|E(T_{x_{i}x_{j}},T_{x_{j}x_{k}})\right| =0$) or  $R_{63} \subset G$ (if $\left|E(T_{x_{i}},T_{x_{j}x_{k}})\right| +\left|E(T_{x_{i}x_{j}},T_{x_{j}x_{k}})\right| \geq 1$), a contradiction.
Then $\left|T_{x_{z}}\right|=0$  by  (F\ref{f-14}) and Lemma \ref{lem-9} \ref{lem-9:item-11}; otherwise, $R_{14} \subset G$ (if $\left|E(T_{x_{z}},T_{x_{i}x_{j}}\cup T_{x_{j}x_{k}})\right|  =0$) or $R_{59} \subset G$ (if $\left|E(T_{x_{z}},T_{x_{i}x_{j}}\cup T_{x_{j}x_{k}})\right|  =1$) or $G$ contains a subdivision of $K_{5}$ (if $\left|E(T_{x_{z}},T_{x_{i}x_{j}}\cup T_{x_{j}x_{k}})\right|  =2$), a contradiction.
Then $\left|T_{x_{k}}\right|=0$ by  Lemma \ref{lem-9}\ref{lem-9:item-2} and \ref{lem-9:item-11}; otherwise, $R_{28} \subset G$ (deleting $x_{z}$), a contradiction.

If $\left|T_{x_{i}x_{k}}\right|= 1$, then  $\left|E(T_{x_{i}x_{k}},T_{x_{j}x_{k}})\right| =0$ by  Lemma \ref{lem-9}\ref{lem-9:item-11}; otherwise, $G[\{x_i,\allowbreak x_j,\allowbreak x_k,\allowbreak q_{ij,1},\allowbreak q_{ij,2},\allowbreak q_{jk,1},\allowbreak q_{ik,1}\}]\cong R_{33}$, a contradiction.
If $\left|T_{x_{i}x_{k}}\right|=1$, then $\left|E(T_{x_{i}x_{j}},T_{x_{i}x_{k}}\cup T_{x_{j}x_{k}})\right|=1$ by Lemma \ref{lem-9}\ref{lem-9:item-11}; otherwise, $R_{38}\subset G$ (if $q_{ij,1}\in T_{x_{i}x_{j}}$ and $\left|E(q_{ij,1},T_{x_{i}x_{k}}\cup T_{x_{j}x_{k}})\right|=2$), $R_{61}\subset G$ (if $\left|E(T_{x_{i}x_{j}},T_{x_{i}x_{k}}\cup T_{x_{j}x_{k}})\right|=0$), or $G$ contains a subdivision of $K_{3,3}$ (if $\left|E(q_{ij,1},T_{x_{i}x_{k}})\right|=\left|E(q_{ij,2},T_{x_{j}x_{k}})\right|=1$ and $q_{ij,1},q_{ij,2}\in T_{x_{i}x_{j}}$), a contradiction.
Thus, if $\left|T_{x_{i}x_{k}}\right|= 1$, then  $\left|T_{x_{j}}\right|=0$ by Lemma \ref{lem-9}\ref{lem-9:item-2}; otherwise, $R_{12} \subset G$ (if $\left|E(T_{x_{j}},T_{x_{i}x_{k}})\right|  =1$ and $\left|E(T_{x_{i}x_{j}},T_{x_{j}x_{k}})\right|  =1$) or $R_{44} \subset G$ (if $\left|E(T_{x_{j}},T_{x_{i}x_{k}})\right|  =0$ and $\left|E(T_{x_{i}x_{j}},T_{x_{j}x_{k}})\right|  =1$), a contradiction.
Thus, if $\left|T_{x_{i}x_{k}}\right|= 1$, then  $\left|T_{x_{i}x_{j}x_{k}}\cup T_{x_{i}x_{j}x_{z}}\right|=0$; otherwise, either $G[\{x_i,\allowbreak x_j,\allowbreak x_z,\allowbreak q_{ij,1},\allowbreak q_{jk,1},\allowbreak w_{ijz,1}\}]\cong R_{20}$ ($E(\{q_{ij,1},\allowbreak q_{jk,1},\allowbreak w_{ijz,1}\})=\{q_{ij,1}q_{jk,1}\}$)  or $G$ contains a subdivision of $K_{3,3}$ or $K_{5}$, a contradiction.
Then $\left|T-T_{2}\right|=0$ by (F\ref{f-16}) and Lemma \ref{lem-9}\ref{lem-9:item-12}, a contradiction.
Thus, $\left|T_{x_{i}x_{k}}\right|= 0$.

If $\left|E(T_{x_{i}x_{j}},T_{x_{j}x_{k}}(T_{x_{j}x_{z}}))\right| =1$, then $\left|E(T_{x_{j}x_{z}}(T_{x_{j}x_{k}}),T_{x_{i}x_{j}}\cup T_{x_{j}x_{k}})\right| =0$ for $G$ contains no subdivision of $K_{5}$.
Thus, if $\left|E(T_{x_{i}x_{j}},T_{x_{j}x_{k}}\cup T_{x_{j}x_{z}})\right| =1$, then $G\subset H_{47}$ by (F\ref{f-14}), (F\ref{f-16}) and Lemma \ref{lem-9}\ref{lem-9:item-1}\ref{lem-9:item-2}\ref{lem-9:item-6} and \ref{lem-9:item-11}.
Let  $\left|E(T_{x_{i}x_{j}},T_{x_{j}x_{k}}\cup T_{x_{j}x_{z}})\right| =0$.
Then  $\left|T_{x_{j}x_{z}}\right| =0$; otherwise, either $R_{16} \subset G$ (if $\left|E(T_{x_{j}x_{z}},T_{x_{j}x_{k}})\right|  =1$) or $R_{45} \subset G$ (if $\left|E(T_{x_{j}x_{z}},T_{x_{j}x_{k}})\right|  =0$), a contradiction.
By an argument similar to that used in case $3$ of Lemma \ref{lem-10}, $G[T_{x_{j}}\cup T_{x_{i}x_{j}x_{k}}(T_{x_{i}x_{j}x_{z}})]\subset 3K_{1}\cup K_{2}$ (see $H_{1}$ in Fig.~\ref{fig-4}).
Since $R_{20} \not\subset G$ and $G$ contains no subdivision of $K_{5}$, $\left|E(T_{x_{i}x_{j}x_{z}},T_{x_{i}x_{j}}\cup T_{x_{j}x_{k}})\right| =0$.
Since $R_{36} \not\subset G$, $\left|E(T_{x_{i}x_{j}x_{k}},T_{x_{j}x_{k}})\right| =0$.
If $\left|T_{x_{i}x_{j}x_{k}}\right| =1$, then $\left|E(T_{x_{i}x_{j}x_{k}},T_{x_{i}x_{j}})\right| =1$  by Lemma \ref{lem-9}\ref{lem-9:item-2}; otherwise, $G[\{x_i,\allowbreak x_k,\allowbreak x_z,\allowbreak q_{ij,1},\allowbreak q_{ij,2},\allowbreak q_{ij,3},\allowbreak q_{ij,4},\allowbreak q_{jk,1},\allowbreak w_{ijk,1}\}]\cong R_{47}$, a contradiction.
Then $\left|T_{0}\right| =0$ by  (F\ref{f-16}) and Lemma \ref{lem-9}\ref{lem-9:item-6}; otherwise, either $R_{14} \subset G$ (if $\left|E(T_{x_{i}x_{j}x_{z}},T_{0})\right| =1$) or  $R_{59} \subset G$ (if $\left|E(T_{x_{i}x_{j}x_{k}},T_{0})\right| =1$), a contradiction.
Thus, if $\left|T_{3}\right|\leq 1$, then either $G\subset H_{65}$ (if $\left|T_{x_{i}x_{j}x_{z}}\right|=0$) or $G\subset H_{66}$ (if $\left|T_{x_{i}x_{j}x_{k}}\right|=0$) by Lemma \ref{lem-9}\ref{lem-9:item-1}\ref{lem-9:item-2} and \ref{lem-9:item-12}.
By  Lemma \ref{lem-9}\ref{lem-9:item-1} and  \ref{lem-9:item-5},  if $\left|T_{3}\right|=2$, then  $\left| T_{x_{j}}\right| \leq 3$ and $\left|E(T_{x_{j}},T_{3})\right| =0$.
Thus, if $\left|T_{3}\right|=2$, then $G\subset H_{27}$ by Lemma \ref{lem-9}\ref{lem-9:item-1}\ref{lem-9:item-5} and \ref{lem-9:item-12}.

\noindent$\textbf{Subsubcase 2.2.2.}$ $\left|T_{x_{i}x_{j}}\right|=3$.

By reasoning similar to that in Subsubcase 2.2.1, we obtain $\left|T_{x_i}\cup T_{x_k}\cup T_{x_z}\right|=0$ and $\left| T_{x_{i}x_{k}x_{z}}\cup T_{x_{j}x_{k}x_{z}}\cup T_{x_{k}x_{z}}\cup T_{x_{i}x_{z}}\right| =0.$
By Lemma \ref{lem-9}\ref{lem-9:item-11},  we discuss the following two configurations.

Suppose that $\left|E(T_{x_{i}x_{j}},T_{x_{j}x_{k}})\right| =1$.
By an argument similar to that used in subcase $2.2$ of Lemma \ref{lem-11}, $G[K_4\cup T_{x_{i}x_{j}}\cup T_{x_{j}x_{k}}\cup T_{x_{i}x_{k}}]\subset H_{11}$.
If $\left|T_{x_{i}x_{k}}\right| \geq 1$, then $\left|T_{x_{j}}\right| =0$ by  Lemma \ref{lem-9}\ref{lem-9:item-2}; otherwise, $R_{12} \subset G$  (if $\left|E(T_{x_{i}x_{j}},T_{x_{i}x_{k}})\right| =0$ and $\left|E(T_{x_{j}},T_{x_{i}x_{k}})\right| =1$)  or $R_{42} \subset G$ (if $\left|E(T_{x_{i}x_{j}},T_{x_{i}x_{k}})\right| =1$ and $\left|E(T_{x_{j}},T_{x_{i}x_{k}})\right| =1$)  or $R_{43} \subset G$ (if $\left|E(T_{x_{i}x_{j}},T_{x_{i}x_{k}})\right| =1$ and $\left|E(T_{x_{j}},T_{x_{i}x_{k}})\right| =0$)  or $R_{44} \subset G$ (if $\left|E(T_{x_{i}x_{j}},T_{x_{i}x_{k}})\right| =0$ and $\left|E(T_{x_{j}},T_{x_{i}x_{k}})\right| =0$), a contradiction.
Thus,  $\left|T_{x_{i}x_{k}}\right| =0$ by (F\ref{f-14}) and (F\ref{f-16}); otherwise, $\left|T-T_{2}\right|=0$, a contradiction.
Then $\left|E(T_{x_{i}x_{j}}\cup T_{x_{j}x_{k}},T_{x_{j}x_{z}})\right| =0$; otherwise, $G$ contains a subdivision of $K_{5}$, a contradiction.
Thus,  $\left|T_{x_{j}x_{z}}\right| \leq 1$; otherwise, $G[\{x_i,\allowbreak x_z,\allowbreak q_{ij,1},q_{ij,2},q_{ij,3},\allowbreak q_{jz,1},\allowbreak q_{jz,2}\}]\cong R_{21}$, a contradiction.
Thus,   $G\subset H_{48}$ by (F\ref{f-14}), (F\ref{f-16}) and Lemma \ref{lem-9}\ref{lem-9:item-1}\ref{lem-9:item-2}\ref{lem-9:item-6} and \ref{lem-9:item-11}.

Suppose that $\left|E(T_{2})\right| =0$ (if $\left|E(T_{2})\right|  \geq1$, then  $G\subset H_{48}$ by symmetry and an argument similar to that used in (F\ref{f-2.1})).
Then  $\left|T_{x_{j}x_{k}}\right| \leq 1$, $\left|T_{x_{i}x_{k}}\right| \leq 1$ and  $\left|T_{x_{j}x_{z}}\right| \leq 1$; otherwise, $G[\{x_i,\allowbreak x_k,\allowbreak q_{ij,1},q_{ij,2},q_{ij,3},\allowbreak q_{jk,1},\allowbreak q_{jk,2}\}]\cong R_{21}$, a contradiction.
Then   $\left|E(T_{x_{j}x_{k}}\cup T_{x_{i}x_{k}}\cup T_{x_{j}x_{z}},T_{x_{i}x_{j}x_{k}})\right| =0$ by Lemma \ref{lem-9}\ref{lem-9:item-2}; otherwise, either $G[\{x_i,\allowbreak x_j,\allowbreak x_k,\allowbreak q_{ij,1},\allowbreak q_{ij,2},\allowbreak q_{jk,1},\allowbreak w_{ijk,1}\}]\cong R_{36}$ ($q_{jk,1}\sim w_{ijk,1}$) or $G$ contains a subdivision of $K_{5}$, a contradiction.
Then   $\left|E(T_{x_{i}x_{j}}\cup T_{x_{j}x_{k}}\cup T_{x_{i}x_{k}},T_{x_{i}x_{j}x_{z}})\right| =0$; otherwise, either $G[\{x_i,\allowbreak x_j,\allowbreak x_k,\allowbreak q_{ij,1},\allowbreak q_{jk,1},\allowbreak w_{ijz,1}\}]\cong R_{20}$ or  $G$ contains a subdivision of $K_{5}$, a contradiction.
Then   $\left|E(T_{x_{j}x_{z}},T_{x_{i}x_{j}x_{z}})\right| =0$; otherwise, $G[\{x_i,\allowbreak x_k,\allowbreak x_z,\allowbreak q_{jk,1},\allowbreak q_{jz,1},\allowbreak w_{ijz,1}\}]\cong R_{14}$, a contradiction.
If $\left|T_{x_{i}x_{k}}\right| =1$, then $\left|T_{x_{j}}\right| =0$; otherwise, either $G[\{x_i,\allowbreak x_j,\allowbreak x_k,\allowbreak q_{ij,1},\allowbreak q_{ij,2},\allowbreak q_{ij,3},\allowbreak q_{ik,1},\allowbreak p_{j,1}\}]\cong R_{42}$ ($E(\{q_{ij,1},\allowbreak q_{ij,2},\allowbreak q_{ij,3},\allowbreak q_{ik,1},\allowbreak p_{j,1}\})=\{q_{ik,1}p_{j,1}\}$) or $G[\{x_i,\allowbreak x_j,\allowbreak x_k,\allowbreak x_z,\allowbreak q_{ij,1},\allowbreak q_{ij,2},\allowbreak q_{ij,3},\allowbreak q_{ik,1},\allowbreak p_{j,1}\}]\cong R_{63}$ ($E(\{q_{ij,1},\allowbreak q_{ij,2},\allowbreak q_{ij,3},\allowbreak q_{ik,1},\allowbreak p_{j,1}\})=\varnothing$), a contradiction.
If $\left|T_{x_{i}x_{k}}\right| =1$, then $\left|T_{x_{i}x_{j}x_{z}}\right| =0$; otherwise, either $G[\{x_i,\allowbreak x_j,\allowbreak x_k,\allowbreak q_{ij,1},\allowbreak q_{ij,2},\allowbreak q_{ij,3},\allowbreak q_{jk,1},\allowbreak q_{ik,1},\allowbreak w_{ijz,1}\}]\cong R_{61}$ ($E(\{q_{ij,1},\allowbreak q_{ij,2},\allowbreak q_{ij,3},\allowbreak q_{jk,1},\allowbreak q_{ik,1},\allowbreak w_{ijz,1}\})=\varnothing$) or   $G$ contains a subdivision of $K_{5}$, a contradiction.
If $\left|T_{x_{i}x_{k}}\right| =1$, then $\left|T_{x_{i}x_{j}x_{k}}\right| =0$; otherwise, either $R_{65} \subset G$ (if  $\left|E(T_{x_{i}x_{j}},T_{x_{i}x_{j}x_{k}})\right| =0$) or $R_{66} \subset G$ (if  $\left|E(T_{x_{i}x_{j}},T_{x_{i}x_{j}x_{k}})\right| =1$), a contradiction.
Thus, $\left|T_{x_{i}x_{k}}\right| =0$ for $\left|T- T_{2}\right| \geq 1$.
By an argument similar to that used in case $3$ of Lemma \ref{lem-10}, $G[T_{x_{j}}\cup T_{x_{i}x_{j}x_{k}}(T_{x_{i}x_{j}x_{z}})]\subset 3K_{1}\cup K_{2}$ (see $H_{1}$ in Fig.~\ref{fig-4}).
If $\left|E(T_{x_{i}x_{j}},T_{x_{i}x_{j}x_{k}})\right| =1$, then $\left|T_{x_{j}x_{z}}\right| =0$; otherwise, $G[\{x_i,\allowbreak x_j,\allowbreak x_z,\allowbreak q_{ij,1},\allowbreak q_{jz,1},\allowbreak w_{ijk,1}\}]\cong R_{20}$, a contradiction.
If $\left|E(T_{x_{i}x_{j}},T_{x_{i}x_{j}x_{k}})\right| =1$, then $\left|T_{0}\right| =0$ by (F\ref{f-16}) and Lemma \ref{lem-9}\ref{lem-9:item-6}; otherwise, either $R_{14} \subset G$ (if  $\left|E(T_{0},T_{x_{i}x_{j}x_{z}})\right| =1$) or $R_{28} \subset G$ (if  $\left|E(T_{0},T_{x_{i}x_{j}x_{k}})\right| =1$), a contradiction.
Thus, if $\left|E(T_{x_{i}x_{j}},T_{x_{i}x_{j}x_{k}})\right| =1$, then  either $G\subset H_{65}$ (if $\left|T_{3}\right| =1$) or $G\subset H_{27}$ (if $\left|T_{3}\right| =2$) by Lemma \ref{lem-9}\ref{lem-9:item-1}\ref{lem-9:item-2} and \ref{lem-9:item-5}.
If $\left|E(T_{x_{i}x_{j}},T_{x_{i}x_{j}x_{k}})\right| =0$, then $\left|T_{0}\right| =0$ by (F\ref{f-16}) and Lemma \ref{lem-9}\ref{lem-9:item-6}; otherwise, either $R_{14} \subset G$ (if  $\left|E(T_{0},T_{x_{i}x_{j}x_{z}})\right| =1$) or $R_{59} \subset G$ (if  $\left|E(T_{0},T_{x_{i}x_{j}x_{k}})\right| =1$), a contradiction.
Thus, if $\left|E(T_{x_{i}x_{j}},T_{x_{i}x_{j}x_{k}})\right| =0$, then either $G\subset H_{68}$ (if $\left|T_{3}\right| \leq1$) or $G\subset H_{65}$ (if $\left|T_{3}\right| =2$)  by (F\ref{f-16}) and Lemma \ref{lem-9}\ref{lem-9:item-1}\ref{lem-9:item-2} and \ref{lem-9:item-5}.

\noindent$\textbf{Subsubcase 2.2.3.}$ $\left|T_{x_{i}x_{j}}\right|=2$.

We first establish the following facts.
By Lemma \ref{lem-9}\ref{lem-9:item-12} and (F\ref{f-13}), we obtain the following fact \ref{f-2231}.

\begin{fact}\label{f-2231}
$\left|T_{x_{i}x_{k}x_{z}}\cup T_{x_{j}x_{k}x_{z}}\right|=\left|T_{x_{i}x_{z}}\cup T_{x_{k}x_{z}}\right|=0$.
\end{fact}

\begin{fact}\label{f-2232}
If  $\left|T_{x_{j}x_{z}}\right|=1$ and $\left|E(T_{x_{i}x_{j}},T_{x_{j}x_{k}}\cup T_{x_{j}x_{z}})\right| =0$, then $\left|T_{x_{i}}\right|=0$.
\end{fact}

\begin{proof}
Otherwise,  $R_{26} \subset G$ (whether  $\left|E(T_{x_{i}x_{j}},T_{x_{j}x_{k}})\right| =0$ or $1$) by  Lemma \ref{lem-9}\ref{lem-9:item-2} and (F\ref{f-2.1}), a contradiction.
\end{proof}

By Lemma \ref{lem-9}\ref{lem-9:item-2} and $R_{25} \not\subset G$, we immediately get the following one.

\begin{fact}\label{f-2233}
At most one vertex in $T_{x_{i}}$ is not adjacent to  vertex of $T_{x_{j}x_{k}}$.
\end{fact}

\begin{fact}\label{f-2234}
$\left|T_{x_{i}}\right| \leq 2$.
\end{fact}
\begin{proof}
Otherwise,  $R_{26} \subset G$ (whether  $\left|E(T_{x_{i}x_{j}},T_{x_{j}x_{k}})\right| =0$ or $1$) by Lemma \ref{lem-9}\ref{lem-9:item-2} and (F\ref{f-2233}), a contradiction.
\end{proof}

\begin{fact}\label{f-2235}
$\left|E(T_{x_{i}x_{j}x_{z}},T_{2}\cup T_{0})\right| =0$.
\end{fact}
\begin{proof}
If $\left|E(T_{x_{i}x_{j}x_{z}},T_{x_{i}x_{j}}\cup T_{x_{j}x_{k}}\cup T_{x_{i}x_{k}})\right| \geq 1$, then either $G[\{x_i,\allowbreak x_j,\allowbreak x_k,\allowbreak q_{ij,1},\allowbreak q_{jk,1},\allowbreak w_{ijz,1}\}]\cong R_{20}$ or  $G$ contains a subdivision of $K_{5}$, a contradiction.
If $\left|E(T_{x_{i}x_{j}x_{z}},T_{x_{j}x_{z}})\right| \geq 1$, then either $G[\{x_i,\allowbreak x_k,\allowbreak x_z,\allowbreak q_{jk,1},\allowbreak q_{jz,1},\allowbreak w_{ijz,1}\}]\cong R_{14}$ or  $G$ contains a subdivision of $K_{5}$, a contradiction.
If $\left|E(T_{x_{i}x_{j}x_{z}},T_{0})\right| \geq 1$, then  $G[\{x_i,\allowbreak x_j,\allowbreak x_k,\allowbreak q_{jk,1},\allowbreak w_{ijz,1},\allowbreak t_{1}\}]\cong R_{14}$ by Lemma \ref{lem-9}\ref{lem-9:item-6}, a contradiction.
This completes the proof of (F\ref{f-2235}) by (F\ref{f-2231}).
\end{proof}

\begin{fact}\label{f-2236}
$\left|E(T_{x_{z}},T_{2}\cup T_{3})\right| =0$.
\end{fact}
\begin{proof}
If $\left|E(T_{x_{z}},T_{x_{i}x_{j}}\cup T_{x_{j}x_{k}}\cup T_{x_{i}x_{k}})\right| \geq 1$, then  either $G[\{x_i,\allowbreak x_j,\allowbreak x_k,\allowbreak q_{ij,1},\allowbreak q_{jk,1},\allowbreak p_{z,1}\}]\cong R_{14}$  or   $G$ contains a subdivision of $K_{3,3}$, a contradiction.
By (F\ref{f-2231}) and Lemma \ref{lem-9}\ref{lem-9:item-2}, $\left|E(T_{x_{z}},T_{2})\right| =0$.
If $\left|E(T_{x_{z}},T_{x_{i}x_{j}x_{z}})\right| \geq 1$, then  $G[\{x_i,\allowbreak x_j,\allowbreak x_k,\allowbreak q_{jk,1},\allowbreak w_{ijz,1},\allowbreak p_{z,1}\}]\cong R_{14}$ by (F\ref{f-2235}), a contradiction.
If $\left|E(T_{x_{i}x_{j}x_{z}},T_{0})\right| \geq 1$, then  $G[\{x_i,\allowbreak x_j,\allowbreak x_k,\allowbreak q_{jk,1},\allowbreak w_{ijz,1},\allowbreak t_{1}\}]\cong R_{14}$ by Lemma \ref{lem-9}\ref{lem-9:item-6}, a contradiction.
This completes the proof of (F\ref{f-2236}) by Lemma \ref{lem-9}\ref{lem-9:item-5} and (F\ref{f-2231}).
\end{proof}

\begin{fact}\label{f-2237}
$\left|T_{x_{k}}\right|\leq 1$ and $\left|T_{x_{z}}\right|\leq 1$.
\end{fact}
\begin{proof}
If $\left|T_{x_{k}}\right|\geq 2$, then either $R_{25} \subset G$ (if some vertex $q_{ij,1}\in T_{x_{i}x_{j}}$ is not adjacent to every vertex of $T_{x_{k}}\cup T_{x_{j}x_{k}}$, after deleting $x_{z}$) or $G$ contains a subdivision of $K_{3,3}$ (in all other cases) by Lemma \ref{lem-9}\ref{lem-9:item-1} and \ref{lem-9:item-2}, a contradiction.
If $\left|T_{x_{z}}\right|\geq 2$, then either $R_{14} \subset G$ (if $\left|E(T_{x_{i}x_{j}},T_{x_{j}x_{k}})\right| =1$) or   $R_{22} \subset G$ (the other cases) by (F\ref{f-2236}), a contradiction.
\end{proof}

By symmetry, (F\ref{f-2233}) and (F\ref{f-2234}) yield 
\setcounter{equation}{2}
\begin{equation}\label{eq4}
If \left|T_{x_{j}}\right| \geq 3, then \left|T_{x_{i}x_{k}}\right| =0; if \left|T_{x_{j}}\right| =2, then \left|E(T_{x_{j}},T_{x_{i}x_{k}})\right| =1.
\tag{3}
\end{equation}
Without loss of generality, we  distinguish the following three subsubsubcases  by Lemma \ref{lem-9}\ref{lem-9:item-1}.

~\

\noindent$\textbf{Subsubsubcase 2.2.3.1.}$ $\left|T_{x_{i}x_{j}x_{z}}\right|=1$.

Then  $\left|E(T_{2})\right| =0$ by (F\ref{f-2.1}) and (F\ref{f-14}).
Then $\left|T_{x_{j}x_{k}}\right| =1$  by  Lemma \ref{lem-9}\ref{lem-9:item-12}.
Then  $\left|T_{x_{k}}\right| =0$ by Lemma \ref{lem-9}\ref{lem-9:item-2} and \ref{lem-9:item-5}; otherwise, either $R_{12} \subset G$ (if  $\left|E(T_{x_{k}},T_{x_{i}x_{j}})\right| =1$) or $R_{28} \subset G$ (if  $\left|E(T_{x_{k}},T_{x_{i}x_{j}})\right| =0$), a contradiction.
Then  $\left|E(T_{x_{i}x_{j}x_{z}}\cup T_{x_{j}x_{k}},T_{x_{i}})\right| =0$ by Lemma \ref{lem-9}\ref{lem-9:item-2}; otherwise, either $R_{14} \subset G$ (if  $\left|E(T_{x_{i}x_{j}x_{z}}\cup T_{x_{j}x_{k}},T_{x_{i}})\right| =1$) or  $G$ contains a subdivision of $K_{5}$ (the other cases), a contradiction.
By symmetry, $\left|E(T_{x_{i}x_{j}x_{k}}\cup T_{x_{i}x_{k}},T_{x_{j}})\right| =0$.
Thus, $\left|T_{x_{i}}\right| \leq 1$ by (F\ref{f-2233}).

Suppose that $\left|T_{x_{i}x_{j}x_{k}}\right| =1$.
Then  $\left|E(T_{x_{i}x_{j}x_{k}},T_{x_{j}x_{k}})\right| =0$ by Lemma \ref{lem-9}\ref{lem-9:item-1} and (F\ref{f-2235}); otherwise, $R_{20} \subset G$ (whether $\left|E(T_{x_{i}x_{j}x_{k}},T_{x_{i}x_{j}})\right| =0$ or $1$), a contradiction.
Then  $\left| T_{x_{z}}\right| =0$ by  (F\ref{f-2235}) and  Lemma \ref{lem-9}\ref{lem-9:item-1} and \ref{lem-9:item-2}; otherwise, either $R_{14} \subset G$ (if $\left|E(T_{x_{i}x_{j}x_{k}},T_{x_{i}x_{j}})\right| =1$) or $R_{28} \subset G$ (if $\left|E(T_{x_{i}x_{j}x_{k}},T_{x_{i}x_{j}})\right| =0$), a contradiction.
By  Lemma \ref{lem-9}\ref{lem-9:item-1} and  \ref{lem-9:item-5},    $\left| T_{x_{j}}\right| \leq 3$.
By  Lemma \ref{lem-9}\ref{lem-9:item-15}, $\left|E(T_{x_{i}x_{j}x_{k}},T_{x_{i}x_{j}})\right|\leq 1$.
If  $\left|E(T_{x_{i}x_{j}x_{k}},T_{x_{i}x_{j}})\right|=0$, then $\left| T_{x_{i}x_{k}}\right| =0$ by  (F\ref{f-2235}) and  Lemma \ref{lem-9}\ref{lem-9:item-1}; otherwise, either $R_{20} \subset G$ (if $\left|E(T_{x_{i}x_{j}x_{k}},T_{x_{i}x_{k}})\right| =1$) or $R_{66} \subset G$ (the other cases), a contradiction.
If   $\left|E(T_{x_{i}x_{j}x_{k}},T_{x_{i}x_{j}})\right|=0$, then $\left| T_{x_{i}}\right| =0$ by  (F\ref{f-2235}) and  Lemma \ref{lem-9}\ref{lem-9:item-1}; otherwise, either $R_{20} \subset G$ (if $\left|E(T_{x_{i}x_{j}x_{k}},T_{x_{i}})\right| =1$) or $R_{63} \subset G$ (the other cases), a contradiction.
Thus, if $\left|E(T_{x_{i}x_{j}x_{k}},T_{x_{i}x_{j}})\right|=0$, then $G\subset H_{65}$ by (F\ref{f-2231}), (F\ref{f-2235})   and  Lemma \ref{lem-9}\ref{lem-9:item-1}\ref{lem-9:item-2} and \ref{lem-9:item-5}.
If $\left|E(T_{x_{i}x_{j}x_{k}},T_{x_{i}x_{j}})\right|=1$, then $\left| T_{x_{j}x_{z}}\right| =0$ by  (F\ref{f-2235}) and  Lemma \ref{lem-9}\ref{lem-9:item-1}; otherwise, either $G[\{x_i,\allowbreak x_j,\allowbreak x_z,\allowbreak q_{ij,1},\allowbreak q_{jz,1},\allowbreak w_{ijk,1}\}]\cong R_{20}$ or $G$ contains a subdivision of $K_{5}$, a contradiction.
Note that $\left|E(T_{x_{j}x_{k}},T_{x_{i}})\right| =\left|E(T_{x_{i}x_{k}},T_{x_{j}})\right| =0$.
Thus, if  $\left|E(T_{x_{i}x_{j}x_{k}},T_{x_{i}x_{j}})\right|=1$ and $\left| T_{x_{j}}\right| \leq 1$, then $G\subset H_{26}$ by (F\ref{f-2231}), (F\ref{f-2235})   and  Lemma \ref{lem-9}\ref{lem-9:item-1}\ref{lem-9:item-2} and \ref{lem-9:item-5}.
If  $\left|E(T_{x_{i}x_{j}x_{k}},T_{x_{i}x_{j}})\right|=1$ and $\left| T_{x_{j}}\right| \geq 2$, then $\left| T_{x_{i}}\right| =0$ by  (F\ref{f-2235})   and  Lemma \ref{lem-9}\ref{lem-9:item-1}\ref{lem-9:item-2} and \ref{lem-9:item-5}; otherwise, $G[\{x_i,\allowbreak x_j,\allowbreak x_k,\allowbreak x_z,\allowbreak q_{jk,1},\allowbreak w_{ijz,1},\allowbreak p_{i,1},\allowbreak p_{j,1},\allowbreak p_{j,2}\}]\cong R_{57}$, a contradiction.
If  $\left|E(T_{x_{i}x_{j}x_{k}},T_{x_{i}x_{j}})\right|=1$ and $\left| T_{x_{j}}\right| \geq 2$, then $\left| T_{x_{i}x_{k}}\right| =0$ by \eqref{eq4}.
Thus, if  $\left|E(T_{x_{i}x_{j}x_{k}},T_{x_{i}x_{j}})\right|=1$ and $\left| T_{x_{j}}\right| \geq 2$, then $G\subset H_{27}$  by (F\ref{f-2231}), (F\ref{f-2235})   and  Lemma \ref{lem-9}\ref{lem-9:item-1}\ref{lem-9:item-2} and \ref{lem-9:item-5}.

Suppose that $\left|T_{x_{i}x_{j}x_{k}}\right| =0$.
By an argument similar to that used in case $2$ of Lemma \ref{lem-10},  if $\left|T_{x_{j}}\right| \geq 3$, then $\left|T_{x_{i}}\cup T_{x_{z}}\right| =0$.
By an argument similar to that used in case $3$ of Lemma \ref{lem-10}, $G[T_{x_{j}}\cup T_{x_{i}x_{j}x_{z}}]\subset 3K_{1}\cup K_{2}$ (see $H_{1}$ in Fig.~\ref{fig-4}).
Thus, if $\left|T_{x_{j}}\right| \geq 3$, then $G\subset H_{67}$ by \eqref{eq4},  (F\ref{f-2231}), (F\ref{f-2235})   and  Lemma \ref{lem-9}\ref{lem-9:item-1}\ref{lem-9:item-2}\ref{lem-9:item-5} and \ref{lem-9:item-6}.
Let $\left|T_{x_{j}}\right| \leq 2$.
If $\left|E(T_{x_{i}x_{j}x_{z}},T_{x_{j}})\right| =1$, then $\left|T_{x_{i}x_{k}}\right| =0$ for $R_{18} \not\subset G$.
Then   $\left|T_{x_{j}x_{z}}\cup T_{x_{z}}\right| \leq 1$ and $\left|T_{x_{j}x_{z}}\cup T_{x_{i}}\right| \leq 1$; otherwise, either $G[\{x_i,\allowbreak x_j,\allowbreak x_z,\allowbreak q_{ij,1},\allowbreak q_{ij,2},\allowbreak q_{jz,1},\allowbreak p_{z,1}\}]\cong R_{28}$ (if $E(\{q_{ij,1},\allowbreak q_{ij,2},\allowbreak q_{jz,1},\allowbreak p_{z,1}\})=\varnothing$) or $G[\{x_i,\allowbreak x_j,\allowbreak x_k,\allowbreak x_z,\allowbreak q_{ij,1},\allowbreak q_{ij,2},\allowbreak q_{jk,1},\allowbreak q_{jz,1},\allowbreak p_{i,1}\}]\cong R_{60}$ (if $E(\{q_{ij,1},\allowbreak q_{ij,2},\allowbreak q_{jk,1},\allowbreak q_{jz,1},\allowbreak p_{i,1}\})=\varnothing$), a contradiction.
Thus, if $\left|E(T_{x_{i}x_{j}x_{z}},T_{x_{j}})\right| =1$, then either $G\subset H_{59}$ (if $\left|T_{x_{j}x_{z}}\right| =0$) or $G\subset H_{67}$ (if $\left|T_{x_{j}x_{z}}\right| =1$)  by   (F\ref{f-2231}), (F\ref{f-2235})   and  Lemma \ref{lem-9}\ref{lem-9:item-1}\ref{lem-9:item-2}\ref{lem-9:item-5} and \ref{lem-9:item-6}.
Let $\left|E(T_{x_{i}x_{j}x_{z}},T_{x_{j}})\right| =0$.
If $\left|T_{x_{z}}\right|=1$ or $\left|T_{x_{i}}\right| =1$, then $\left|T_{x_{j}}\right| \leq 1$ for $R_{57} \not\subset G$.  
If $\left|T_{x_{i}x_{k}}\right| =1$, then $\left|T_{x_{j}}\right| \leq 1$ for $R_{18} \not\subset G$ and $R_{25} \not\subset G$.
Thus, either $G\subset H_{30}$ (if $\left|T_{x_{j}}\right| \leq 1$) or $G\subset H_{67}$ (if $\left|T_{x_{j}}\right| \geq 2$) by   (F\ref{f-2231}), (F\ref{f-2235})   and  Lemma \ref{lem-9}\ref{lem-9:item-1}\ref{lem-9:item-2}\ref{lem-9:item-5} and \ref{lem-9:item-6}.

~\

\noindent$\textbf{Subsubsubcase 2.2.3.2.}$ $\left|T_{x_{i}x_{j}x_{z}}\right|=0$ and $\left|T_{x_{i}x_{j}x_{k}}\right|=1$.

Then $\left|E(T_{2})\right| =0$ by   (F\ref{f-2.1}) and  $G$ contains no subdivision of $K_{3,3}$ or $K_{5}$.
Then $\left|T_{x_{i}}\right| \leq 1$ and $\left|E(T_{x_{j}x_{k}},T_{x_{i}})\right| =\left|E(T_{x_{i}x_{j}},T_{x_{k}})\right| =0$ by (F\ref{f-2233}); otherwise, $G$ contains a subdivision of $K_{3,3}$, a contradiction.
Then $\left|T_{x_{i}}\cup T_{x_{j}x_{k}}\right| \leq 2$ by Lemma \ref{lem-9}\ref{lem-9:item-2} and \ref{lem-9:item-15}; otherwise, $R_{28} \subset G$ (according to the adjacency relations between $T_{x_{i}x_{j}x_{k}}$ and the vertex sets $T_{x_{i}}\cup T_{x_{i}x_{j}}\cup T_{x_{j}x_{k}}$, there are eight possible cases to consider), a contradiction.
Then  $\left|T_{x_{k}}\right| =0$  by Lemma \ref{lem-9}\ref{lem-9:item-2} and \ref{lem-9:item-15}; otherwise, $R_{28} \subset G$ (according to the adjacency relations between $T_{x_{i}x_{j}x_{k}}$ and the vertex sets $T_{x_{k}}\cup T_{x_{i}x_{j}}\cup T_{x_{j}x_{k}}$, there are eight possible cases to consider), a contradiction.
If  $\left|T_{x_{j}x_{z}}\right|=1$, then $\left|E(T_{x_{i}x_{j}x_{k}},T_{2})\right| =0$ for $R_{20} \not\subset G$ and $G$ contains no subdivision of $K_{5}$.

Suppose that $\left|E(T_{x_{i}x_{j}x_{k}},T_{0})\right| =1$.
Then  $\left|T_{x_{j}x_{z}}\right| =0$ and $\left|E(T_{x_{i}x_{j}x_{k}},T_{x_{i}}\cup T_{x_{j}})\right| =0$ by Lemma \ref{lem-9}\ref{lem-9:item-17}. 
Then $\left|E(T_{x_{i}x_{j}x_{k}},T_{2})\right| =0$; otherwise, $G[\{x_j,\allowbreak x_k,\allowbreak x_z,\allowbreak q_{ij,1},\allowbreak q_{jk,1},\allowbreak w_{ijk,1},\allowbreak t_{1}\}]\cong R_{28}$ (if $E(\{q_{ij,1},\allowbreak q_{jk,1},\allowbreak w_{ijk,1},\allowbreak t_{1}\})=\{q_{ij,1}w_{ijk,1},\allowbreak w_{ijk,1}t_{1}\}$) or $G[\{x_i,\allowbreak x_j,\allowbreak x_k,\allowbreak q_{ij,1},\allowbreak q_{ij,2},\allowbreak q_{jk,1},\allowbreak w_{ijk,1}\}]\cong R_{36}$ (if $E(\{q_{ij,1},\allowbreak q_{ij,2},\allowbreak q_{jk,1},\allowbreak w_{ijk,1}\})=\{q_{jk,1}w_{ijk,1}\}$) or $G$ contains subdivision of $K_{3,3}$ or $K_{5}$.
Then  $\left|T_{x_{z}}\right| =0$ by $d(G) < 4$ and the subdivision of  $K_{5}$ is not a subgraph of $G$. 
Then $\left|T_{x_{j}}\right| \leq 2$ by Lemma \ref{lem-9}\ref{lem-9:item-16} (here we consider a different $K_{4}$).
If $\left|T_{x_{j}}\right| =2$, then $\left|T_{x_{i}x_{k}}\right| =0$  by \eqref{eq4}; otherwise, $G$ contains a subdivision of $K_{3,3}$, a contradiction.
Note that $\left|T_{x_{i}}\cup T_{x_{j}x_{k}}\right| \leq 2$.
Thus, if $\left|T_{x_{j}}\right| =2$, then either $G\subset A_{4}$  or $G\subset A_{6}$ by   (F\ref{f-2231})  and  Lemma \ref{lem-9}\ref{lem-9:item-1}\ref{lem-9:item-2}\ref{lem-9:item-5} and \ref{lem-9:item-6}.
Let $\left|T_{x_{j}}\right| \leq 1$.
If $\left|T_{x_{i}x_{k}}\right| =2$, then $\left|T_{x_{j}}\cup T_{x_{k}}\right| =0$ for $R_{28} \not\subset G$.
Thus, if $\left|T_{x_{i}x_{k}}\right| =2$, then  $G\subset A_{2}$ (note that $\left|T_{x_{j}x_{k}}\right| \geq \left|T_{x_{i}x_{k}}\right|$) by   (F\ref{f-2231})    and  Lemma \ref{lem-9}\ref{lem-9:item-1}\ref{lem-9:item-2}\ref{lem-9:item-5} and \ref{lem-9:item-6}.
If $\left|T_{x_{i}x_{k}}\right| \leq 1$, then either $G\subset A_{3}$  or  $G\subset A_{5}$ by   (F\ref{f-2231})    and  Lemma \ref{lem-9}\ref{lem-9:item-1}\ref{lem-9:item-2}\ref{lem-9:item-5} and \ref{lem-9:item-6}.
Let  $\left|E(T_{x_{i}x_{j}x_{k}},T_{0})\right| =0$. 

Suppose that $\left|T_{x_{z}}\right| \geq 1$.
Then, by symmetry, we obtain the same conclusion as in the previous paragraph.
Let $\left|T_{x_{z}}\right|=0$.

Suppose that $\left|T_{x_{j}x_{z}}\right| \geq  1$.
Then $\left|E(T_{2},T_{x_{i}x_{j}x_{k}})\right| =0$ by Lemma \ref{lem-9}\ref{lem-9:item-5} (here we consider a different $K_{4}$).
Thus, $\left|T_{x_{j}x_{z}}\right| = 1$ for $R_{36} \not\subset G$.
Then $\left|T_{x_{i}}\cup T_{x_{i}x_{k}}\right| =0$  by (F\ref{f-13}) and (F\ref{f-2232}).
Note that $G[T_{x_{k}}\cup T_{x_{i}x_{j}x_{k}}]\subset 3K_{1}\cup K_{2}$.
Thus,  $G\subset H_{67}$ by   (F\ref{f-2231})    and  Lemma \ref{lem-9}\ref{lem-9:item-1}\ref{lem-9:item-2}\ref{lem-9:item-5} and \ref{lem-9:item-6}.
Let $\left|T_{x_{j}x_{z}}\right| = 0$.

Suppose that $\left|T_{x_{i}x_{k}}\right| = 2$.
Then $\left|T_{1}\right| = 0$ by $\left|T_{x_{i}}\cup T_{x_{j}x_{k}}\right| \leq 2$ and $\left|T_{x_{j}}\cup T_{x_{i}x_{k}}\right| \leq 2$.
If $\left|E(T_{x_{i}x_{j}x_{k}},T_{2})\right| \geq 1$, then $\left|E(T_{x_{i}x_{j}x_{k}},T_{x_{i}x_{j}})\right| =\left|E(T_{x_{i}x_{j}x_{k}},T_{x_{j}x_{k}})\right| =\left|E(T_{x_{i}x_{j}x_{k}},T_{x_{i}x_{k}})\right| = 1$ by Lemma \ref{lem-9}\ref{lem-9:item-2}\ref{lem-9:item-15} and $R_{36} \not\subset G$.
Thus, if $\left|T_{x_{i}x_{k}}\right| = 2$, then either $G\subset H_{35}$  or $G\subset A_{2}$.
Let $\left|T_{x_{i}x_{k}}\right| \leq 1$.

Then $\left|E(T_{x_{i}}\cup T_{x_{j}},T_{2})\right| =0$ by Lemma \ref{lem-9}\ref{lem-9:item-2} and  $G$ contains no subdivision of $K_{3,3}$. 
If $\left|E(T_{x_{i}x_{j}x_{k}},T_{2})\right|\geq 1$, then $\left|E(T_{x_{i}x_{j}x_{k}},T_{x_{i}x_{j}})\right|=1$ by Lemma \ref{lem-9}\ref{lem-9:item-15} and  $R_{36} \not\subset G$.
If $\left|E(T_{x_{i}x_{j}x_{k}},T_{2})\right|\leq 1$, then $\left|E(T_{x_{i}},T_{3})\right|=0$ for $R_{35} \not\subset G$.
If $\left|E(T_{x_{i}x_{j}x_{k}},T_{2})\right|\geq 2$, then $\left|E(T_{x_{i}},T_{3})\right|=1$  by Lemma \ref{lem-9}\ref{lem-9:item-5} and   $R_{25} \not\subset G$.
If $\left|E(T_{x_{i}x_{j}x_{k}},T_{2})\right|=1$, then $\left|T_{x_{j}x_{k}}\right| \leq 1$ for $R_{36} \not\subset G$.

Suppose that $\left|T_{x_i x_k}\right|=1$.
Write $T_{x_i x_j}=\{q_{ij,1},q_{ij,2}\}$, $T_{x_j x_k}=\{q_{jk,1}\}$, $T_{x_i x_k}=\{q_{ik,1}\}$ and $T_{x_i x_j x_k}=\{w_{ijk,1}\}$.
Recall that $T_3=\{w_{ijk,1}\}$, $E(T_2)=\varnothing$, $E(T_{x_j},T_2)=\varnothing$ and $E(T_{x_j})=\varnothing$.
If $p_{j,1},p_{j,2}\in T_{x_j}$, then $G[\{x_i,x_j,x_k,q_{ij,1},q_{ik,1},p_{j,1},p_{j,2}\}]\cong R_{25}$, a contradiction. 
Thus, $|T_{x_j}|\leq1$.
If $\left|E(T_{x_{i}x_{j}x_{k}},T_{2})\right|\leq 1$, then $\left|E(T_{x_{j}},T_{3})\right|=0$; otherwise,  $G[\{x_i,x_j,x_k,x_z,w_{ijk,1},q_{ik,1},p_{j,1}\}]\cong R_{25}$ (if $|E(w_{ijk,1},T_2)|=0$ and $p_{j,1}\in T_{x_j}$) or $G[\{x_i,x_j,x_k,x_z,w_{ijk,1},q_{ik,1},p_{j,1}\}]\cong R_{35}$ (if $|E(w_{ijk,1},T_2)|\geq1$ and $|E(q_{ik,1},T_3)|=0$), a contradiction.
Thus, if $\left|E(T_{x_{i}x_{j}x_{k}},T_{2})\right|=0$, then either $G\subset A_{3}$  or $G\subset A_{5}$ by (F\ref{f-2231}), Lemma \ref{lem-9}\ref{lem-9:item-1}\ref{lem-9:item-5} and \ref{lem-9:item-6}.
Thus, if $\left|E(T_{x_{i}x_{j}x_{k}},T_{2})\right|=1$, then  $G\subset H_{26}$ by (F\ref{f-2231}), Lemma \ref{lem-9}\ref{lem-9:item-1}\ref{lem-9:item-5}\ref{lem-9:item-6} and \ref{lem-9:item-15}.
Note that $\left|T_{x_{i}}\cup T_{x_{j}x_{k}}\right| \leq 2$.
If $w_{ijk,1}\nsim q_{ik,1}$ and $w_{ijk,1}\sim q_{jk,1}$, then $|E(T_{x_i},w_{ijk,1})|=|T_{x_i}|$ and $|E(T_{x_j},w_{ijk,1})|=0$; otherwise, either $G[\{x_i,x_j,x_z,w_{ijk,1},q_{jk,1},q_{ik,1},p_{i,1}\}]\cong R_{25}$ (if $p_{i,1}\in T_{x_i}$ and $w_{ijk,1}\nsim p_{i,1}$) or $G[\{x_i,x_j,x_k,x_z,w_{ijk,1},q_{ik,1},p_{j,1}\}]\cong R_{35}$ (if $p_{j,1}\in T_{x_j}$ and $w_{ijk,1}\sim p_{j,1}$), a contradiction.
If $w_{ijk,1}\sim q_{ik,1}$, then $|E(T_{x_i},w_{ijk,1})|=0$ and $|E(T_{x_j},w_{ijk,1})|=|T_{x_j}|$; otherwise, either $G[\{x_i,x_j,x_k,x_z,w_{ijk,1},q_{jk,1},p_{i,1}\}]\cong R_{35}$ (if $p_{i,1}\in T_{x_i}$ and $w_{ijk,1}\sim p_{i,1}$) or $G[\{x_i,x_j,x_z,w_{ijk,1},q_{ik,1},q_{jk,1},p_{j,1}\}]\cong R_{25}$ (if $p_{j,1}\in T_{x_j}$ and $w_{ijk,1}\nsim p_{j,1}$), a contradiction.
Thus, if $\left|E(T_{x_{i}x_{j}x_{k}},T_{2})\right|=2$, then either $G\subset H_{37}$  or $G\subset H_{41}$ by (F\ref{f-2231}), Lemma \ref{lem-9}\ref{lem-9:item-1}\ref{lem-9:item-5}\ref{lem-9:item-6} and \ref{lem-9:item-15}.
If $w_{ijk,1}\nsim q_{ij,2}$ and $w_{ijk,1}\sim q_{jk,1}$, then  $|E(T_{x_i},w_{ijk,1})|=|T_{x_i}|$; otherwise, $G[\{x_i,x_k,x_z,w_{ijk,1},q_{jk,1},q_{ij,2},p_{i,1}\}]\cong R_{25}$ (if $p_{i,1}\in T_{x_i}$ and $w_{ijk,1}\nsim p_{i,1}$), a contradiction.
If $\left|E(T_{x_{i}x_{j}x_{k}},T_{2})\right|=3$, then $\left|T_{x_{j}}\right| =0$ by Lemma \ref{lem-9}\ref{lem-9:item-2}, $R_{32} \not\subset G$, and  $G$ contains no subdivision of $K_{3,3}$.
Thus, if $\left|E(T_{x_{i}x_{j}x_{k}},T_{2})\right|=3$, then  $G\subset H_{35}$ by (F\ref{f-2231}), Lemma \ref{lem-9}\ref{lem-9:item-1}\ref{lem-9:item-5}\ref{lem-9:item-6} and \ref{lem-9:item-15}.
Let $\left|T_{x_{i}x_{k}}\right| =0$.

If  $\left|T_{x_{i}}\right| = 1$, then $\left|T_{x_{j}}\right| \leq 2$  for $R_{51} \not\subset G$.
If $\left|T_{x_{i}}\right| = 1$ and $\left|E(T_{x_{i}x_{j}x_{k}},T_{2})\right|=0$, then $\left|E(T_{x_{j}},T_{3})\right|=0$ for $R_{60} \not\subset G$.
Thus, if  $\left|E(T_{x_{i}x_{j}x_{k}},T_{2})\right|=0$, then  either $G\subset H_{67}$  or $G\subset A_{6}$ by (F\ref{f-2231}), Lemma \ref{lem-9}\ref{lem-9:item-1}\ref{lem-9:item-5} and \ref{lem-9:item-6}.
Since $R_{36} \not\subset G$, $\left|E(T_{x_{j}},T_{3})\right|\leq 1$.
If $\left|T_{x_{i}}\right| = 1$,  $\left|T_{x_{j}}\right| = 2$ and $1\leq \left|E(T_{x_{i}x_{j}x_{k}},T_{2})\right|\leq 2$, then $\left|E(T_{x_{j}},T_{3})\right|= 1$; otherwise, either $G[\{x_i,\allowbreak x_j,\allowbreak x_z,\allowbreak w_{ijk,1},\allowbreak p_{i,1},\allowbreak p_{j,1},\allowbreak p_{j,2}\}]\cong R_{25}$ (if $E(\{w_{ijk,1},\allowbreak p_{i,1},\allowbreak p_{j,1},\allowbreak p_{j,2}\})=\{w_{ijk,1}p_{i,1}\}$) or $G[\{x_i,\allowbreak x_j,\allowbreak x_k,\allowbreak q_{ij,1},\allowbreak q_{jk,1},\allowbreak w_{ijk,1},\allowbreak p_{i,1},\allowbreak p_{j,1},\allowbreak p_{j,2}\}]\cong R_{57}$ (if $E(\{q_{ij,1},\allowbreak q_{jk,1},\allowbreak w_{ijk,1},\allowbreak p_{i,1},\allowbreak p_{j,1},\allowbreak p_{j,2}\})=\{q_{ij,1}w_{ijk,1}\}$), a contradiction.
Thus, if  $\left|E(T_{x_{i}x_{j}x_{k}},T_{2})\right|=1$, then  either $G\subset H_{49}$  or $G\subset H_{64}$ by (F\ref{f-2231}), Lemma \ref{lem-9}\ref{lem-9:item-1}\ref{lem-9:item-5}\ref{lem-9:item-6} and \ref{lem-9:item-15}.
Thus, if  $\left|E(T_{x_{i}x_{j}x_{k}},T_{2})\right|=2$, then  either $G\subset H_{37}$  or $G\subset H_{64}$ by (F\ref{f-2231}), Lemma \ref{lem-9}\ref{lem-9:item-1}\ref{lem-9:item-5}\ref{lem-9:item-6} and \ref{lem-9:item-15}.

~\

\noindent$\textbf{Subsubsubcase 2.2.3.3.}$ $\left|T_{x_{i}x_{j}x_{z}}\cup T_{x_{i}x_{j}x_{k}}\right|=0$.

Then $\left|T_{0}\right| =0$ by Lemma \ref{lem-9}\ref{lem-9:item-6} and (F\ref{f-16}).
Then $E(T_1)=T_3=\varnothing$  by (F\ref{f-15}) and (F\ref{f-2231}).
Since $|T|>|T_2|$, we have $|T_1|\geq1$.

Suppose that $\left|T_{x_{j}x_{z}}\right| \geq 1$.
Then $\left|T_{x_{i}x_{k}}\right| =0$ by (F\ref{f-13}).
Then $\left|E(T_{x_{j}x_{z}},T_{x_{j}x_{k}})\right| =0$; otherwise either $G[\{x_i,\allowbreak x_j,\allowbreak q_{ij,1},\allowbreak q_{ij,2},\allowbreak q_{jk,1},\allowbreak q_{jz,1}\}]\cong R_{16}$ or $G$ contains a subdivision of $K_{5}$, a contradiction.
If $\left|E(T_{2})\right| =0$, then $\left|T_{x_{i}}\right|=0$  by Lemma \ref{lem-9}\ref{lem-9:item-2}; otherwise, either
$G[\{x_i,\allowbreak x_j,\allowbreak x_z,\allowbreak q_{jk,1},\allowbreak q_{jz,1},\allowbreak p_{i,1}\}]\cong R_{12}$ (if $E(\{q_{jk,1},\allowbreak q_{jz,1},\allowbreak p_{i,1}\})=\{q_{jk,1}p_{i,1}\}$) or $G[\{x_i,\allowbreak x_j,\allowbreak x_k,\allowbreak x_z,\allowbreak q_{ij,1},\allowbreak q_{ij,2},\allowbreak q_{jk,1},\allowbreak q_{jz,1},\allowbreak p_{i,1}\}]\cong R_{60}$ (if $E(\{q_{ij,1},\allowbreak q_{ij,2},\allowbreak q_{jk,1},\allowbreak q_{jz,1},\allowbreak p_{i,1}\})=\varnothing$), a contradiction.
Note that $\left|E(T_{x_{z}},T_{2})\right| =0$ and $\left|T_{x_{z}}\right|\leq 1$.
Thus, if $\left|E(T_{2})\right| =0$, then $\left|T_{x_{z}}\right|=0$; otherwise $G[\{x_i,\allowbreak x_j,\allowbreak x_z,\allowbreak q_{ij,1},\allowbreak q_{ij,2},\allowbreak q_{jz,1},\allowbreak p_{z,1}\}]\cong R_{28}$, a contradiction.
By symmetry, if $\left|E(T_{2})\right| =0$, then $\left|T_{x_{k}}\right|=0$.
Thus, if $\left|E(T_{2})\right| =0$, then $G\subset H_{46}$ by (F\ref{f-15}), (F\ref{f-2231}), Lemma \ref{lem-9}\ref{lem-9:item-1} and \ref{lem-9:item-2}.
Let $\left|E(T_{2})\right| \geq 1$.
Then $\left|E(T_{2})\right| = 1$; otherwise either $G$ contains a subdivision of $K_{3,3}$ or $G[\{x_i,\allowbreak x_z,\allowbreak q_{ij,1},\allowbreak q_{ij,2},\allowbreak q_{jk,1},\allowbreak q_{jz,1}\}]\cong R_{2}$, a contradiction.
Without loss of generality, let $\left|E(T_{x_{i}x_{j}},T_{x_{j}x_{k}})\right| =1$.
Then  $\left|E(T_{x_{i}}\cup T_{x_{k}},T_{2})\right| =0$ and $\left|T_{x_{i}}\right|\leq 1$ by Lemma \ref{lem-9}\ref{lem-9:item-2}; otherwise $G[\{x_i,\allowbreak x_j,\allowbreak x_k,\allowbreak q_{jk,1},\allowbreak q_{jz,1},\allowbreak p_{i,1}\}]\cong R_{12}$ or $G[\{x_i,\allowbreak x_j,\allowbreak x_k,\allowbreak q_{ij,2},\allowbreak q_{jk,1},\allowbreak p_{i,1},\allowbreak p_{i,2}\}]\cong R_{25}$ or $G$ contains a subdivision of $K_{3,3}$, a contradiction.
Note that $\left|T_{x_{k}}\right|\leq 1$.
Then $\left|T_{x_{z}}\right|=0$ by Lemma \ref{lem-9}\ref{lem-9:item-2}; otherwise either $G[\{x_i,\allowbreak x_j,\allowbreak x_z,\allowbreak q_{ij,1},\allowbreak q_{jk,1},\allowbreak p_{z,1}\}]\cong R_{14}$ or $G$ contains a subdivision of $K_{5}$, a contradiction.
Then $\left|T_{x_{j}x_{z}}\right|=1$; otherwise $G[\{x_i,\allowbreak x_z,\allowbreak q_{ij,1},\allowbreak q_{jk,1},\allowbreak q_{jz,1},\allowbreak q_{jz,2}\}]\cong R_{1}$, a contradiction.
If  $\left|T_{x_{i}}\cup T_{x_{k}}\right| \geq 1$, then $\left|T_{x_{j}}\right| \leq 2$ ; otherwise, by symmetry we may assume
$p_{i,1}\in T_{x_i}$,  $G[\{x_i,\allowbreak x_j,\allowbreak x_k,\allowbreak x_z,\allowbreak q_{jk,1},\allowbreak p_{i,1},\allowbreak p_{j,1},\allowbreak p_{j,2},\allowbreak p_{j,3}\}]\cong R_{51}$, a contradiction.
If  $\left|T_{x_{i}}\cup T_{x_{k}}\right| \geq 1$, then $\left|T_{x_{j}x_{k}}\right| \leq 1$; otherwise  $G[\{x_i,\allowbreak x_j,\allowbreak x_k,\allowbreak q_{ij,2},\allowbreak q_{jk,1},\allowbreak q_{jk,2},\allowbreak p_{i,1}\}]\cong R_{28}$, a contradiction.
Thus, either $G\subset H_{48}$  or $G\subset H_{54}$ by (F\ref{f-15}), (F\ref{f-2231}), Lemma \ref{lem-9}\ref{lem-9:item-1} and \ref{lem-9:item-2}.
Let $\left|T_{x_{j}x_{z}}\right| =0$.

Suppose that $\left|T_{x_{i}x_{k}}\right| \geq 1$.
If $\left|E(T_{2})\right| = 1$, then $\left|E(T_{x_{i}x_{k}},T_{x_{j}x_{k}})\right|=0$; otherwise $G[\{x_i,\allowbreak x_j,\allowbreak x_k,\allowbreak q_{ij,1},\allowbreak q_{ij,2},\allowbreak q_{ik,1},\allowbreak q_{jk,1}\}]\cong R_{33}$, a contradiction.
If $\left|E(T_{2})\right| = 1$, then $\left|T_{x_{i}}\right| = 0$; otherwise either $G[\{x_i,\allowbreak x_j,\allowbreak x_k,\allowbreak x_z,\allowbreak q_{ij,1},\allowbreak q_{ik,1},\allowbreak q_{jk,1},\allowbreak p_{i,1}\}]\cong R_{44}$ or $G$ contains a subdivision of $K_{3,3}$, a contradiction.
If $\left|E(T_{2})\right| = 1$, then $\left|T_{x_{j}}\right| = 0$; otherwise either $G[\{x_j,\allowbreak x_k,\allowbreak x_z,\allowbreak q_{ij,1},\allowbreak q_{ik,1},\allowbreak q_{jk,1},\allowbreak p_{j,1}\}]\cong R_{26}$ or $G[\{x_i,\allowbreak x_j,\allowbreak x_k,\allowbreak q_{ij,1},\allowbreak q_{ik,1},\allowbreak q_{jk,1},\allowbreak p_{j,1}\}]\cong R_{30}$, a contradiction.
If $\left|E(T_{2})\right| = 1$, then $\left|T_{x_{k}}\right| = 0$; otherwise $G[\{x_i,\allowbreak x_j,\allowbreak x_k,\allowbreak q_{ij,1},\allowbreak q_{ij,2},\allowbreak q_{jk,1} (q_{ik,1}),\allowbreak p_{k,1}\}]\cong R_{28}$, a contradiction.
If $\left|E(T_{2})\right| = 1$, then $\left|T_{x_{z}}\right| = 0$; otherwise $G[\{x_i,\allowbreak x_j,\allowbreak x_z,\allowbreak q_{ij,1},\allowbreak q_{ik,1} (q_{jk,1}),\allowbreak p_{z,1}\}]\cong R_{14}$, a contradiction for $\left|T_{1}\right| \geq 1$.
Thus, $\left|E(T_{2})\right| \neq 1$.
Then $\left|E(T_{2})\right| \leq 2$; otherwise either $G[\{x_i,\allowbreak x_j,\allowbreak x_k,\allowbreak q_{ij,1},\allowbreak q_{ik,1},\allowbreak q_{jk,1},\allowbreak q_{jk,2}\}]\cong R_{33}$ or $G$ contains a subdivision of $K_{5}$, a contradiction.
Thus, $\left|T_{x_{i}x_{k}}\right| \leq 1$ and $\left|T_{x_{j}x_{k}}\right| \leq 1$; otherwise, $\left|T_{1}\right| =0$, a contradiction.
If $\left|E(T_{2})\right| =2$ and $\left|E(T_{x_{j}x_{k}},T_{x_{i}x_{k}})\right| =0$, then $\left|E(T_{x_{i}}\cup T_{x_{j}},T_{2})\right| =0$ and $\left|T_{x_{i}}\right|  \leq 1$; otherwise either $G[\{x_i,\allowbreak x_j,\allowbreak x_k,\allowbreak q_{ij,1},\allowbreak q_{jk,1},\allowbreak p_{i,1},\allowbreak p_{i,2}\}]\cong R_{25}$ or $G$ contains a subdivision of $K_{3,3}$, a contradiction.
If $\left|E(T_{2})\right| =2$ and $\left|E(T_{x_{j}x_{k}},T_{x_{i}x_{k}})\right| =0$, then $\left|E(T_{x_{k}},T_{2})\right| =0$; otherwise $G[\{x_i,\allowbreak x_j,\allowbreak x_z,\allowbreak q_{ij,1},\allowbreak q_{ij,2},\allowbreak q_{ik,1} (q_{jk,1}),\allowbreak p_{k,1}\}]\cong R_{28}$, a contradiction.
If $\left|E(T_{2})\right| =2$ and $\left|E(T_{x_{j}x_{k}},T_{x_{i}x_{k}})\right| =0$, then $\left|T_{x_{z}}\right| =0$; otherwise either $G[\{x_i,\allowbreak x_j,\allowbreak x_z,\allowbreak q_{ij,1},\allowbreak q_{ik,1},\allowbreak p_{z,1}\}]\cong R_{14}$ or $G$ contains a subdivision of $K_{5}$, a contradiction.
If $\left|E(T_{2})\right| =2$, then $\left|E(T_{1},T_{2})\right| =0$; otherwise $G$ contains a subdivision of $K_{3,3}$, a contradiction.
Thus, if $\left|E(T_{2})\right| =2$ and $\left|E(T_{x_{j}x_{k}},T_{x_{i}x_{k}})\right| =0$, then $G\subset H_{44}$ by (F\ref{f-15}), (F\ref{f-2231}), (F\ref{f-2237}), Lemma \ref{lem-9}\ref{lem-9:item-1} and \ref{lem-9:item-2}.
If $\left|E(T_{2})\right| =2$ and $\left|E(T_{x_{j}x_{k}},T_{x_{i}x_{k}})\right| =1$, then $\left|T_{x_{i}}\right| = 0$; otherwise either $G[\{x_i,\allowbreak x_j,\allowbreak x_k,\allowbreak q_{ij,2},\allowbreak q_{ik,1},\allowbreak q_{jk,1},\allowbreak p_{i,1}\}]\cong R_{30}$ or $G$ contains a subdivision of $K_{3,3}$, a contradiction.
If $\left|E(T_{2})\right| =2$ and $\left|E(T_{x_{j}x_{k}},T_{x_{i}x_{k}})\right| =1$, then $\left|T_{x_{j}}\right| = 0$; otherwise either $G[\{x_i,\allowbreak x_j,\allowbreak x_z,\allowbreak q_{ij,2},\allowbreak q_{ik,1},\allowbreak q_{jk,1},\allowbreak p_{j,1}\}]\cong R_{26}$ or $G[\{x_i,\allowbreak x_j,\allowbreak x_k,\allowbreak q_{ij,2},\allowbreak q_{ik,1},\allowbreak q_{jk,1},\allowbreak p_{j,1}\}]\cong R_{30}$, a contradiction.
If $\left|E(T_{2})\right| =2$ and $\left|E(T_{x_{j}x_{k}},T_{x_{i}x_{k}})\right| =1$, then $\left|T_{x_{k}}\right| = 0$; otherwise either $G[\{x_i,\allowbreak x_j,\allowbreak x_k,\allowbreak q_{ij,1},\allowbreak q_{ij,2},\allowbreak q_{jk,1} (q_{ik,1}),\allowbreak p_{k,1}\}]\cong R_{28}$ or $G$ contains a subdivision of $K_{3,3}$, a contradiction.
If $\left|E(T_{2})\right| =2$ and $\left|E(T_{x_{j}x_{k}},T_{x_{i}x_{k}})\right| =1$, then $\left|T_{x_{i}}\right| = 0$; otherwise either $G[\{x_j,\allowbreak x_k,\allowbreak x_z,\allowbreak q_{ik,1},\allowbreak q_{jk,1},\allowbreak p_{i,1}\}]\cong R_{14}$ or $G$ contains a subdivision of $K_{5}$, a contradiction for $\left|T_{1}\right| \geq 1$.
Let $\left|E(T_{2})\right| =0$.
Then $\left|E(T_{x_{k}},T_{x_{i}x_{j}})\right| =1$ (if $\left|T_{x_{k}}\right| =1$), $\left|T_{x_{j}x_{k}}\cup T_{x_{i}}\right| \leq 3$ and $\left|E(T_{x_{j}x_{k}},T_{x_{i}})\right|=1$ (if $\left|T_{x_{j}x_{k}}\cup T_{x_{i}}\right| =3$); otherwise $G[\{x_i,\allowbreak x_j,\allowbreak x_k,\allowbreak q_{ij,1},\allowbreak q_{jk,1},\allowbreak p_{i,1},\allowbreak p_{i,2}\}]\cong R_{25}$ or $G[\{x_i,\allowbreak x_j,\allowbreak x_k,\allowbreak q_{ij,1},\allowbreak q_{ij,2},\allowbreak q_{ik,1},\allowbreak p_{k,1}\}]\cong R_{28}$ or $G[\{x_i,\allowbreak x_j,\allowbreak x_k,\allowbreak q_{ij,1},\allowbreak q_{jk,1},\allowbreak q_{jk,2},\allowbreak p_{i,1}\}]\cong R_{28}$ or $G$ contains a subdivision of $K_{3,3}$, a contradiction.
Similarly, $\left|T_{x_{i}x_{k}}\cup T_{x_{j}}\right| \leq 3$ and $\left|E(T_{x_{i}x_{k}},T_{x_{j}})\right|=1$ (if $\left|T_{x_{i}x_{k}}\cup T_{x_{j}}\right| =3$).
If $\left|T_{x_{i}x_{k}}\right| =2$, then $G\subset H_{39}$ by (F\ref{f-15}), (F\ref{f-2231}), (F\ref{f-2236}), (F\ref{f-2237}) and Lemma \ref{lem-9}\ref{lem-9:item-1}\ref{lem-9:item-2}; otherwise $G$ contains a subdivision of $K_{3,3}$, a contradiction.
Let $\left|T_{x_{i}x_{k}}\right| =1$.
If $\left|T_{x_{j}x_{k}}\right| =2$, then either $G\subset H_{42}$ (if $\left|T_{x_{i}}\right| =0$) or $G\subset H_{43}$ (if $\left|T_{x_{i}}\right| =1$) by (F\ref{f-15}), (F\ref{f-2231}), (F\ref{f-2236}), (F\ref{f-2237}) and Lemma \ref{lem-9}\ref{lem-9:item-1}\ref{lem-9:item-2}; otherwise $G$ contains a subdivision of $K_{3,3}$, a contradiction.
Let $\left|T_{x_{j}x_{k}}\right| =1$.
If $\left|E(T_{x_{j}x_{k}},T_{x_{i}})\right|=1$, then $G\subset H_{57}$ by (F\ref{f-15}), (F\ref{f-2231}), (F\ref{f-2236}), (F\ref{f-2237}) and Lemma \ref{lem-9}\ref{lem-9:item-1}\ref{lem-9:item-2}; otherwise $G$ contains a subdivision of $K_{3,3}$, a contradiction.
Without loss of generality, let $\left|E(T_{x_{j}x_{k}},T_{x_{i}})\right|=0$ and $\left|E(T_{x_{i}x_{k}},T_{x_{j}})\right|=0$.
Then $\left|T_{x_{i}}\right| \leq 1$ and $\left|T_{x_{j}}\right| \leq 1$.
Thus, $G\subset H_{58}$ by (F\ref{f-15}), (F\ref{f-2231}), (F\ref{f-2236}), (F\ref{f-2237}), Lemma \ref{lem-9}\ref{lem-9:item-1} and \ref{lem-9:item-2}.
Let $\left|T_{x_{i}x_{k}}\right| =0$.

Suppose that $\left|E(T_{2})\right|=0$.
Then $\left|E(T_{x_{k}},T_{x_{i}x_{j}})\right| =1$ (if $\left|T_{x_{k}}\right| =1$), $\left|T_{x_{j}x_{k}}\cup T_{x_{i}}\right| \leq 3$ and $\left|E(T_{x_{j}x_{k}},T_{x_{i}})\right|=1$ (if $\left|T_{x_{j}x_{k}}\cup T_{x_{i}}\right| =3$); otherwise $G[\{x_i,\allowbreak x_j,\allowbreak x_k,\allowbreak q_{ij,1},\allowbreak q_{jk,1},\allowbreak p_{i,1},\allowbreak p_{i,2}\}]\cong R_{25}$ or $G[\{x_i,\allowbreak x_j,\allowbreak x_k,\allowbreak q_{ij,1},\allowbreak q_{ij,2},\allowbreak q_{jk,1},\allowbreak p_{k,1}\}]\cong R_{28}$ or $G[\{x_i,\allowbreak x_j,\allowbreak x_k,\allowbreak q_{ij,1},\allowbreak q_{jk,1},\allowbreak q_{jk,2},\allowbreak p_{i,1}\}]\cong R_{28}$ or $G$ contains a subdivision of $K_{3,3}$, a contradiction.
If $\left|E(T_{x_{j}x_{k}},T_{x_{i}})\right|=1$, then $\left|T_{x_{k}}\right| =0$; otherwise $G$ contains a subdivision of $K_{3,3}$, a contradiction.
Thus, if $\left|E(T_{x_{j}x_{k}},T_{x_{i}})\right|=1$, then  either $G\subset H_{56}$ (if $\left|T_{x_{i}}\right|=1$) or $G\subset H_{69}$ (if $\left|T_{x_{i}}\right|=2$) by (F\ref{f-15}), (F\ref{f-2231}), (F\ref{f-2236}), (F\ref{f-2237}), Lemma \ref{lem-9}\ref{lem-9:item-1} and \ref{lem-9:item-2}; otherwise $G$ contains a subdivision of $K_{3,3}$, a contradiction.
Thus, if $\left|E(T_{x_{j}x_{k}},T_{x_{i}})\right|=0$, then  either $G\subset H_{56}$ (if $\left|T_{x_{j}x_{k}}\right|=2$) or $G\subset H_{70}$ (if $\left|T_{x_{j}x_{k}}\right|=1$) by (F\ref{f-15}), (F\ref{f-2231}), (F\ref{f-2236}), (F\ref{f-2237}), Lemma \ref{lem-9}\ref{lem-9:item-1}\ref{lem-9:item-2}; otherwise $G$ contains a subdivision of $K_{3,3}$, a contradiction.

Suppose that $\left|E(T_{2})\right|=1$.
Then $\left|T_{x_{z}}\right|=0$; otherwise $G[\{x_i,\allowbreak x_j,\allowbreak x_z,\allowbreak q_{ij,1},\allowbreak q_{jk,1},\allowbreak p_{z,1}\}]\cong R_{14}$, a contradiction.
Then $\left|E(T_{x_{k}},T_{x_{i}x_{j}})\right| =0$; otherwise either $G[\{x_i,\allowbreak x_j,\allowbreak x_z,\allowbreak q_{ij,1},\allowbreak q_{ij,2},\allowbreak q_{jk,1},\allowbreak p_{k,1}\}]\cong R_{28}$ or $G$ contains a subdivision of $K_{3,3}$, a contradiction.
Thus, if $\left|T_{x_{j}x_{k}}\right|=1$, then $G\subset H_{60}$ by (F\ref{f-15}), (F\ref{f-2231}), (F\ref{f-2237}), Lemma \ref{lem-9}\ref{lem-9:item-1} and \ref{lem-9:item-2}.
If $\left|T_{x_{j}x_{k}}\right|=2$, then $\left|T_{x_{i}}\cup T_{x_{k}}\right|=0$; otherwise either $G[\{x_i,\allowbreak x_j,\allowbreak x_k,\allowbreak q_{ij,2},\allowbreak q_{jk,1},\allowbreak q_{jk,2},\allowbreak p_{i,1}\}]\cong R_{28}$ or $G[\{x_i,\allowbreak x_j,\allowbreak x_k,\allowbreak q_{ij,1},\allowbreak q_{ij,2},\allowbreak q_{jk,2},\allowbreak p_{k,1}\}]\cong R_{28}$, a contradiction.
Thus, if $\left|T_{x_{j}x_{k}}\right|=2$, then $G\subset H_{48}$  by (F\ref{f-15}), (F\ref{f-2231}), (F\ref{f-2237}), Lemma \ref{lem-9}\ref{lem-9:item-1} and \ref{lem-9:item-2}.

\noindent$\textbf{Subsubcase 2.2.4.}$ $\left|T_{x_{i}x_{j}}\right|=1$.

By the choice of $T_{x_ix_j}$ and (F\ref{f-13}), $|T_{x_jx_k}|=1$, $T_{x_ix_z}=T_{x_kx_z}=\varnothing$ and $|T_{x_ix_k}\cup T_{x_jx_z}|\leq1$.
Moreover, $E(T_1)=\varnothing$ and $E(T_2,T_0)=\varnothing$ by (F\ref{f-15}) and (F\ref{f-16}), respectively.
Then $\left|E(T_{x_{z}},T_{2})\right| =0$; otherwise either $G[\{x_i,\allowbreak x_j,\allowbreak x_k,\allowbreak q_{ij,1},\allowbreak q_{jk,1} (q_{ik,1}),\allowbreak p_{z,1}\}]\cong R_{14}$  or $G$ contains a subdivision of $K_{5}$, a contradiction.
Thus, $\left|T_{x_{z}}\right|\leq 1$; otherwise either $G[\{x_i,\allowbreak x_j,\allowbreak x_z,\allowbreak q_{ij,1},\allowbreak q_{jk,1},\allowbreak p_{z,1}\}]\cong R_{14}$ (if $\left|E(T_{x_{i}x_{j}},T_{x_{j}x_{k}})\right| =1$) or $G[\{x_i,\allowbreak x_j,\allowbreak x_z,\allowbreak q_{ij,1},\allowbreak q_{jk,1},\allowbreak p_{z,1},\allowbreak p_{z,2}\}]\cong R_{22}$ (if $\left|E(T_{x_{i}x_{j}},T_{x_{j}x_{k}})\right| =0$), a contradiction.
Then $\left|T_{x_i}\right|\leq 2$ and $\left|T_{x_k}\right|\leq 2$ by Lemma
\ref{lem-9}\ref{lem-9:item-1} and \ref{lem-9:item-2}; otherwise, by symmetry we may assume $\left|T_{x_i}\right|\geq 3$, and $G[\{x_i,x_j,x_z,q_{ij,1},q_{jk,1},p_{i,1},p_{i,2}\}]\cong R_{26}$ (if $\left|E(T_{x_ix_j},T_{x_jx_k})\right|=0$ and $\left|E(T_{x_i},T_{x_jx_k})\right|\geq 2$)
or $G[\{x_i,x_j,x_k,q_{ij,1},q_{jk,1},p_{i,1} (p_{i,3}),p_{i,2}\}]\cong R_{25}$ (the remaining cases), a contradiction.
If $\left|T_{x_{j}}\right| \geq 2$, then  $\left|E(T_{x_{j}},T_{x_{i}x_{k}})\right| =1$; otherwise $G[\{x_i,\allowbreak x_j,\allowbreak x_k,\allowbreak q_{ij,1},\allowbreak q_{ik,1},\allowbreak p_{j,1},\allowbreak p_{j,2}\}]\cong R_{25}$, a contradiction.

Suppose that $\left|T_{x_{i}}\right|= 2$.
Then $\left|E(T_{x_{i}},T_{x_{j}x_{k}})\right| =1$; otherwise either $G[\{x_i,\allowbreak x_j,\allowbreak x_k,\allowbreak q_{ij,1},\allowbreak q_{jk,1},\allowbreak p_{i,1},\allowbreak p_{i,2}\}]\cong R_{25}$ or $G[\{x_i,\allowbreak x_j,\allowbreak x_z,\allowbreak q_{ij,1},\allowbreak q_{jk,1},\allowbreak p_{i,1},\allowbreak p_{i,2}\}]\cong R_{26}$, a contradiction.
Then $\left|T_{3}\right|= 0$; otherwise either $R_{12}\subset G$ or $G$ contains a subdivision of $K_{3,3}$, a contradiction.
Then $T_0=\varnothing$ by (F\ref{f-16}) and Lemma \ref{lem-9}\ref{lem-9:item-6}.
Then $\left|T_{x_{j}x_{z}}\right| =0$; otherwise $G[\{x_i,\allowbreak x_j,\allowbreak x_z,\allowbreak q_{jk,1},\allowbreak q_{jz,1},\allowbreak p_{i,1}\}]\cong R_{12}$, a contradiction.
Then $\left|E(T_{x_{k}},T_{x_{i}x_{j}})\right| =0$; otherwise $G$ contains a subdivision of $K_{3,3}$, a contradiction.
Thus, $\left|T_{x_{k}}\right|\leq 1$; otherwise, $\left|E(T_{x_{k}},T_{x_{i}x_{j}})\right| =1$, a contradiction.
By symmetry,  $\left|T_{x_{j}}\cup T_{x_{i}x_{k}}\right|\leq 2$ and  $\left|E(T_{x_{j}},T_{x_{i}x_{k}})\right| =0$.
If  $\left|E(T_{x_{i}x_{j}},T_{x_{j}x_{k}})\right| =0$, then $\left|E(T_{x_{i}x_{k}},T_{2})\right| =0$; otherwise $G[\{x_i,\allowbreak x_j,\allowbreak q_{ij,1},\allowbreak q_{ik,1},\allowbreak q_{jk,1},\allowbreak p_{i,1}\}]\cong R_{12}$ or $G[\{x_i,\allowbreak x_j,\allowbreak x_k,\allowbreak q_{ij,1},\allowbreak q_{ik,1},\allowbreak q_{jk,1},\allowbreak p_{i,2}\}]\cong R_{30}$ or $G$ contains a subdivision of $K_{3,3}$, a contradiction.
Thus, if $\left|E(T_{x_{i}x_{j}},T_{x_{j}x_{k}})\right| =0$, then  either $G\subset H_{2}$ or $G\subset H_{61}$ by (F\ref{f-15}), Lemma \ref{lem-9}\ref{lem-9:item-1} and \ref{lem-9:item-2}.
Let $\left|E(T_{x_{i}x_{j}},T_{x_{j}x_{k}})\right| =1$.
Then $\left|T_{x_{z}}\right|= 0$ by (F\ref{f-14}); otherwise $G[\{x_i,\allowbreak x_j,\allowbreak x_z,\allowbreak q_{ij,1},\allowbreak q_{jk,1},\allowbreak p_{z,1}\}]\cong R_{14}$, a contradiction.
Then  $\left|E(T_{x_{i}x_{k}},T_{x_{i}x_{j}})\right| =0$ and $\left|E(T_{x_{i}x_{k}},T_{x_{j}x_{k}})\right| =1$; otherwise either $G[\{x_i,\allowbreak x_j,\allowbreak x_k,\allowbreak q_{ij,1},\allowbreak q_{ik,1},\allowbreak q_{jk,1},\allowbreak p_{i,2}\}]\cong R_{30}$ or $G$ contains a subdivision of $K_{3,3}$, a contradiction.
Thus, either $G\subset H_{53}$ or $G\subset H_{60}$ by (F\ref{f-15}), Lemma \ref{lem-9}\ref{lem-9:item-1} and \ref{lem-9:item-2}.
By symmetry, if $\left|T_{x_{k}}\right|= 2$, then $G$ is an induced subgraph of $\{H_{2},H_{53},H_{60},H_{61}\}$.
Let $\left|T_{x_{i}}\right|\leq 1$ and $\left|T_{x_{k}}\right|\leq 1$.

Suppose that $\left|T_{x_{j}}\right| \geq 3$.
By Lemma \ref{lem-9}\ref{lem-9:item-1} and \ref{lem-9:item-5}, $|T_{x_j}|\leq4$ and $\left|T_{x_{i}x_{k}x_{z}}\right| =0$.
Then $\left|E(T_{x_{i}x_{j}x_{k}}\cup T_{x_{j}x_{k}x_{z}}, T_{2})\right| =0$; otherwise either $R_{20}\subset G$  or $G$ contains a subdivision of $K_{5}$, a contradiction.
Then $\left|T_{x_{i}}\cup T_{x_{k}}\cup T_{x_{z}}\right| =0$; otherwise either $R_{38}\subset G$ or $R_{51}\subset G$, a contradiction.
Then $\left|T_{x_{i}x_{k}}\right| =0$; otherwise either $G[\{x_i,\allowbreak x_j,\allowbreak x_k,\allowbreak q_{ij,1},\allowbreak q_{ik,1},\allowbreak p_{j,2},\allowbreak p_{j,3}\}]\cong R_{25}$ or $G[\{x_i,\allowbreak x_j,\allowbreak x_k,\allowbreak q_{ij,1},\allowbreak q_{ik,1},\allowbreak p_{j,2},\allowbreak p_{j,3}\}]\cong R_{27}$, a contradiction.
If $\left|E(T_{x_{i}x_{j}},T_{x_{j}x_{k}})\right| =1$, then $\left|T_{3}\right| =0$ by (F\ref{f-14}) and; otherwise $G[\{x_i,\allowbreak x_j,\allowbreak x_z,\allowbreak q_{ij,1},\allowbreak q_{jk,1},\allowbreak w_{ikz,1}\}]\cong R_{18}$ or $G[\{x_i,\allowbreak x_j,\allowbreak x_z,\allowbreak q_{ij,1},\allowbreak q_{jk,1},\allowbreak w_{ijz,1}\}]\cong R_{20}$ or $G[\{x_j,\allowbreak x_k,\allowbreak x_z,\allowbreak q_{ij,1},\allowbreak q_{jk,1},\allowbreak w_{jkz,1}\}]\cong R_{20}$ or $G$ contains a subdivision of $K_{5}$ or $G$ contains a subdivision of $K_{3,3}$, a contradiction.
If $\left|E(T_{x_{i}x_{j}},T_{x_{j}x_{k}})\right| =1$, then $\left|E(T_{2})\right| =1$; otherwise $G$ contains a subdivision of $K_{5}$, a contradiction.
Thus, if $\left|E(T_{x_{i}x_{j}},T_{x_{j}x_{k}})\right| =1$, then $G\subset H_{47}$ by (F\ref{f-15}), (F\ref{f-16}), Lemma \ref{lem-9}\ref{lem-9:item-1}\ref{lem-9:item-2} and \ref{lem-9:item-6}.
Let $\left|E(T_{2})\right| =0$.
Note that, if $\left|T_{x_{j}}\right| =4$, then $\left|E(T_{x_{j}},T_{x_{i}x_{j}x_{k}})\right| =1$, $\left|E(T_{x_{j}},T_{x_{i}x_{j}x_{z}})\right| =1$ and $\left|E(T_{x_{j}},T_{x_{j}x_{k}x_{z}})\right| =1$.
Thus, if $\left|T_{x_{j}}\right| =4$, then $\left|T_{3}\right| =1$; otherwise $G$ contains a subdivision of $K_{5}$, a contradiction.
Thus, if $\left|E(T_{3},T_{2})\right| =0$, then either $G\subset H_{64}$ (if $\left|T_{x_{j}}\right| =3$) or $G\subset H_{68}$ (if $\left|T_{x_{j}}\right| =4$) by (F\ref{f-15}), (F\ref{f-16}), Lemma \ref{lem-9}\ref{lem-9:item-1}\ref{lem-9:item-5}\ref{lem-9:item-6} and \ref{lem-9:item-15}.
If $\left|E(T_{3},T_{2})\right| \geq 1$, then $\left|T_{x_{j}x_{z}}\right| =0$; otherwise $R_{20}\subset G$, a contradiction.
If $\left|E(T_{3},T_{2})\right| \geq 1$, then $\left|T_{3}\right| \leq 2$; otherwise $G[\{x_i,\allowbreak x_j,\allowbreak x_z,\allowbreak q_{ij,1},\allowbreak w_{ijk,1},\allowbreak w_{jkz,1}\}]\cong R_{20}$, a contradiction.
If $\left|E(T_{3},T_{2})\right| \geq 2$, then  $\left|T_{x_{i}x_{j}x_{z}}\cup T_{x_{j}x_{k}x_{z}}\right| =0$; otherwise  $R_{20}\subset G$, a contradiction.
Thus, either $G\subset H_{27}$ (if $\left|E(T_{3},T_{2})\right| =1$) or $G\subset H_{64}$ (if $\left|E(T_{3},T_{2})\right| =2$) by (F\ref{f-15}), (F\ref{f-16}), Lemma \ref{lem-9}\ref{lem-9:item-1}\ref{lem-9:item-5}\ref{lem-9:item-6} and \ref{lem-9:item-15}.
Let $\left|T_{x_{j}}\right|  \leq 2$.

Suppose that $\left|T_{x_{j}x_{z}}\right| =1$.
Then $\left|E(T_{2})\right| \leq 1$; otherwise $G$ contains a subdivision of $K_{5}$, a contradiction.
Then $\left|E(T_{1},T_{2})\right| =0$ by Lemma \ref{lem-9}\ref{lem-9:item-2}; otherwise $G[\{x_j,\allowbreak x_k,\allowbreak x_z,\allowbreak q_{jk,1},\allowbreak q_{jz,1},\allowbreak p_{i,1}\}]\cong R_{14}$ or $G[\{x_i,\allowbreak x_j,\allowbreak x_z,\allowbreak q_{ij,1},\allowbreak q_{jz,1},\allowbreak p_{k,1}\}]\cong R_{14}$ or $G$ contains a subdivision of $K_{5}$, a contradiction.
If $\left|E(T_{2})\right| =1$, then $\left|T_{x_{z}}\right| =0$ by (F\ref{f-14}); otherwise $G[\{x_i,\allowbreak x_j,\allowbreak x_z,\allowbreak q_{ij,1},\allowbreak q_{jk,1},\allowbreak p_{z,1}\}]\cong R_{14}$, a contradiction.
If $\left|E(T_{2})\right| =1$, then $\left|T_{3}\right| =0$; otherwise $G[\{x_i,\allowbreak x_j,\allowbreak x_z,\allowbreak q_{ij,1},\allowbreak q_{jk,1},\allowbreak w_{ikz,1}\}]\cong R_{18}$ or $R_{20}\subset G$  or $G$ contains a subdivision of $K_{5}$ or $G$ contains a subdivision of $K_{3,3}$, a contradiction.
Thus, if $\left|E(T_{2})\right| =1$, then $G\subset H_{54}$ by (F\ref{f-15}).
Let $\left|E(T_{2})\right| =0$.
If $\left|T_{3}\right| =0$, then $G\subset H_{63}$ by (F\ref{f-15}).
Then $\left|E(T_{3},T_{2})\right| =0$ by Lemma \ref{lem-9}\ref{lem-9:item-15}; otherwise either $R_{20}\subset G$ or $G$ contains a subdivision of $K_{5}$, a contradiction.
Thus, if $\left|T_{1}\right| =0$, then $G\subset H_{38}$ by Lemma \ref{lem-9}\ref{lem-9:item-1}.
Let $\left|T_{1}\right| \geq 1$ and $\left|T_{3}\right| \geq 1$.
Then $\left|E(T_{x_{i}x_{j}x_{k}}\cup T_{x_{i}x_{j}x_{z}}\cup T_{x_{j}x_{k}x_{z}},T_{x_i}\cup T_{x_k}\cup T_{x_z})\right|=0$ and $\left|E(T_{x_{i}x_{k}x_{z}},T_{1})\right|=0$ by Lemma \ref{lem-9}\ref{lem-9:item-5}; otherwise $R_{14}\subset G$ or $G$ contains a subdivision of $K_{5}$, a contradiction.
If $\left|T_{x_{i}x_{j}x_{k}}\right| =1$, then $\left|T_{x_{z}}\right| =0$; otherwise either $G[\{x_i,\allowbreak x_j,\allowbreak x_z,\allowbreak q_{ij,1},\allowbreak q_{jz,1},\allowbreak p_{z,1},\allowbreak w_{ijk,1}\}]\cong R_{28}$ or $G$ contains a subdivision of $K_{5}$, a contradiction.
If $\left|T_{x_{i}x_{j}x_{k}}\right| =1$  and  $\left|E(T_{x_{j}},T_{x_{i}x_{j}x_{k}})\right| =1$, then $\left|T_{3}- T_{x_{i}x_{j}x_{k}}\right| =0$ by Lemma \ref{lem-9}\ref{lem-9:item-5}; otherwise $R_{3}\subset G$, a contradiction.
Thus, if $\left|T_{x_{i}x_{j}x_{k}}\right| =1$ and  $\left|E(T_{x_{j}},T_{x_{i}x_{j}x_{k}})\right| =1$, then $G\subset H_{62}$ by (F\ref{f-15}), Lemma \ref{lem-9}\ref{lem-9:item-1} and \ref{lem-9:item-5}.
By symmetry and a similar argument, if $\left|E(T_{1},T_{3})\right|\geq1$, then $G\subset H_{62}$.
Let $\left|E(T_{1},T_{3})\right| =0$.
Then $\left|T_{x_{i}}\cup T_{x_{j}x_{k}x_{z}}\right| \leq 1$; otherwise $G[\{x_i,\allowbreak x_j,\allowbreak x_k,\allowbreak q_{ij,1},\allowbreak q_{jk,1},\allowbreak p_{i,1},\allowbreak w_{jkz,1}\}]\cong R_{28}$, a contradiction.
Similarly, if $\left|T_{x_{i}x_{k}x_{z}}\right| =1$, then $\left|T_{x_{k}}\cup T_{x_{i}x_{j}x_{z}}\right| \leq 1$ and $\left|T_{x_{z}}\cup T_{x_{i}x_{j}x_{k}}\right| \leq 1$.
If $\left|T_{x_{i}x_{k}x_{z}}\right| =1$, then  $\left|T_{x_{j}}\right| =0$; otherwise either $G[\{x_i,\allowbreak x_j,\allowbreak x_k,\allowbreak q_{ij,1},\allowbreak q_{jz,1},\allowbreak p_{j,1},\allowbreak w_{ikz,1}\}]\cong R_{25}$ or $G$ contains a subdivision of $K_{5}$, a contradiction.
Thus, if $\left|T_{x_{i}x_{k}x_{z}}\right| =1$, then  $G\subset H_{3}$ (if $\left|T_{1}\right| =3$),  $G\subset H_{32}$ (if $\left|T_{1}\right| =2$), or $G\subset H_{41}$ (if $\left|T_{1}\right| =1$) by (F\ref{f-15}) and Lemma \ref{lem-9}\ref{lem-9:item-1}.
Let $\left|T_{x_{i}x_{k}x_{z}}\right| =0$.
Without loss of generality, let $\left|T_{x_{i}x_{j}x_{k}}\right| =1$.
Then $\left|T_{x_i}\cup T_{x_j}\right|\leq2$ and $\left|T_{x_k}\cup T_{x_j}\right|\leq2$; otherwise $R_{57}\subset G$, a contradiction.
Thus, if $\left|T_{x_{j}}\right| =2$, then $G\subset H_{64}$ by (F\ref{f-15}) and Lemma \ref{lem-9}\ref{lem-9:item-1}.
Let $\left|T_{x_{j}}\right| \leq 1$.
Note that $\left|T_{x_{i}}\cup T_{x_{j}x_{k}x_{z}}\right| \leq 1$ and $\left|T_{x_{k}}\cup T_{x_{i}x_{j}x_{z}}\right| \leq 1$.
Thus, either $G\subset H_{62}$ (if $\left|T_{x_{i}}\cup T_{x_{k}}\right| =2$),  $G\subset H_{49}$ (if $\left|T_{x_{i}}\cup T_{x_{k}}\right| =1$), or $G\subset H_{64}$ (if $\left|T_{x_{i}}\cup T_{x_{k}}\right| =0$) by (F\ref{f-15}) and Lemma \ref{lem-9}\ref{lem-9:item-1}.
Let $\left|T_{x_{j}x_{z}}\right| =0$.

Suppose that $\left|E(T_{2})\right| \geq 2$.
If $\left|E(T_{2})\right| = 3$, then $\left|T_{1}\cup T_{3}\right| = 0$; otherwise $R_{6}\subset G$ or $R_{30}\subset G$ or $G$ contains a subdivision of $K_{5}$, a contradiction.
Thus, $\left|E(T_{2})\right| = 2$.
Then $\left|T_{x_{z}}\right| =0$; otherwise either $G[\{x_i,\allowbreak x_j,\allowbreak x_z,\allowbreak q_{ij,1},\allowbreak q_{ik,1},\allowbreak p_{z,1}\}]\cong R_{14}$ or $G$ contains a subdivision of $K_{5}$, a contradiction.
If $\left|E(T_{x_{j}},T_{x_{i}x_{k}})\right| =1$, then  $\left|E(T_{x_{i}x_{j}},T_{x_{j}x_{k}})\right| =0$; otherwise $G$ contains a subdivision of $K_{3,3}$, a contradiction.
If $\left|E(T_{x_{j}},T_{x_{i}x_{k}})\right| =1$, then  $\left|E(T_{x_{i}}\cup T_{x_{k}},T_{2})\right| = 0$; otherwise $G$ contains a subdivision of $K_{3,3}$, a contradiction.
Thus, if  $\left|E(T_{x_{j}},T_{x_{i}x_{k}})\right| =1$, then $G\subset H_{53}$ by (F\ref{f-15}) and Lemma \ref{lem-9}\ref{lem-9:item-2}.
By symmetry, if $\left|E(T_{x_{i}},T_{x_{j}x_{k}})\right| =1$ or $\left|E(T_{x_{k}},T_{x_{i}x_{j}})\right| =1$, then $G\subset H_{53}$ by (F\ref{f-15}) and Lemma \ref{lem-9}\ref{lem-9:item-2}.
Let $\left|E(T_{1},T_{2})\right| = 0$.
If $\left|T_{x_{j}}\right| =2$, then  $\left|E(T_{x_{j}},T_{x_{i}x_{k}})\right| =1$; otherwise $G[\{x_i,\allowbreak x_j,\allowbreak x_z,\allowbreak q_{ij,1},\allowbreak q_{ik,1},\allowbreak p_{j,1},\allowbreak p_{j,2}\}]\cong R_{25}$, a contradiction.
Thus, $\left|T_{x_{j}}\right| \leq 1$.
Thus, $G\subset H_{53}$ by (F\ref{f-15}) and Lemma \ref{lem-9}\ref{lem-9:item-2}.
Let $\left|E(T_{2})\right| \leq 1$.

Suppose that $\left|E(T_{2})\right| = 1$.
Then $\left|T_{x_{z}}\right| =0$ for $R_{14}\not\subset G$ and $\left|E(T_{x_{z}},T_{2})\right| =0$.
Then $\left|T_{3}\right| =0$; otherwise $G[\{x_j,\allowbreak x_k,\allowbreak x_z,\allowbreak q_{ij,1},\allowbreak q_{jk,1},\allowbreak w_{ikz,1}\}]\cong R_{18}$ or $R_{20}\subset G$ or $G$ contains a subdivision of $K_{5}$, a contradiction.
Then $T_0=\varnothing$ by (F\ref{f-16}), Lemma \ref{lem-9}\ref{lem-9:item-6} and the connectedness of $G$.
Since $|T|>|T_2|$, we have $|T_1|\geq1$.
If $\left|T_{x_{i}x_{k}}\right| = 1$, then $\left|T_{1}-T_{x_{z}}\right| =0$ by Lemma \ref{lem-9}\ref{lem-9:item-2}; otherwise $G[\{x_j,\allowbreak x_k,\allowbreak x_z,\allowbreak q_{ij,1},\allowbreak q_{ik,1},\allowbreak q_{jk,1},\allowbreak p_{j,1}\}]\cong R_{26}$ or $G[\{x_i,\allowbreak x_j,\allowbreak x_k,\allowbreak q_{ij,1},\allowbreak q_{ik,1},\allowbreak q_{jk,1},\allowbreak p_{j,1}\}]\cong R_{30}$ or $G[\{x_i,\allowbreak x_j,\allowbreak x_k,\allowbreak x_z,\allowbreak q_{ij,1},\allowbreak q_{ik,1},\allowbreak q_{jk,1},\allowbreak p_{i,1}\}]\cong R_{44}$ or $G$ contains a subdivision of $K_{3,3}$, a contradiction for $\left|T_{1}\cup T_{3}\right| \geq 1$.
Thus, $\left|T_{x_{i}x_{k}}\right| = 0$.
If $\left|E(T_{x_{i}},T_{x_{j}x_{k}})\right| =1$, then $\left|E(T_{x_{k}},T_{x_{i}x_{j}})\right| =0$; otherwise $G$ contains a subdivision of $K_{3,3}$, a contradiction.
Thus, either $G\subset H_{50}$ (if $\left|E(T_{1},T_{2})\right| =1$) or $G\subset H_{54}$ (if $\left|E(T_{1},T_{2})\right| =0$) by (F\ref{f-15}) and Lemma \ref{lem-9}\ref{lem-9:item-2}.
Let $\left|E(T_{2})\right| = 0$.

If $\left|T_{x_{j}}\right| =2$, then $\left|E(T_{x_{j}},T_{x_{i}x_{k}})\right| =1$; otherwise $G[\{x_i,\allowbreak x_j,\allowbreak x_k,\allowbreak q_{ij,1},\allowbreak q_{ik,1},\allowbreak p_{j,1},\allowbreak p_{j,2}\}]\cong R_{25}$, a contradiction.
Thus, if $\left|T_{3}\right| =0$, then $G\subset H_{61}$ by (F\ref{f-15}), (F\ref{f-16}), Lemma \ref{lem-9}\ref{lem-9:item-2} and \ref{lem-9:item-6}.
Let $\left|T_{3}\right| \geq 1$.
If $|T_{x_ix_jx_k}|=1$, then $E(T_1,T_2)=\varnothing$ by Lemma \ref{lem-9}\ref{lem-9:item-2}; otherwise $G$ contains a subdivision of $K_{3,3}$, a contradiction.
Then $\left|E(T_{3}-T_{x_{i}x_{j}x_{k}},T_{2})\right| =0$; otherwise either $R_{20}\subset G$ or $G$ contains a subdivision of $K_{5}$, a contradiction.
Then $\left|E(T_{3}-T_{x_{i}x_{j}x_{k}},T_{1})\right| =0$; otherwise either $R_{18}\subset G$ or $G$ contains a subdivision of $K_{5}$, a contradiction.
Then $\left|T_{x_ix_jx_z}\cup T_{x_k}\right|\leq 1$, $\left|T_{x_jx_kx_z}\cup T_{x_i}\right|\leq 1$ and $\left|T_{x_ix_kx_z}\cup T_{x_j}\right|\leq 1$; otherwise, by symmetry we may assume that $\left|T_{x_ix_jx_z}\cup T_{x_k}\right|\geq 2$ or
$\left|T_{x_ix_kx_z}\cup T_{x_j}\right|\geq 2$, and $G[\{x_i,x_k,x_z,q_{ij,1},p_{k,1},w_{ijz,1}\}]\cong R_{12}$ or $G[\{x_i,x_j,x_k,q_{ij,1},q_{jk,1},w_{ijz,1}\}]\cong R_{20}$ or $G[\{x_i,x_j,x_k,q_{ij,1},q_{jk,1},p_{k,1},w_{ijz,1}\}]\cong R_{28}$ or $G$ contains a subdivision of $K_{5}$, a contradiction.
Then  $\left|E(T_{x_{j}},T_{x_{i}x_{k}})\right| =0$ and $\left|T_{x_{j}}\right| \leq 1$; otherwise either $R_{14}\subset G$  or $G$ contains a subdivision of $K_{5}$, a contradiction.
If $\left|T_{x_{i}x_{j}x_{k}}\cup T_{x_{z}}\right| =2$, then $\left|E(T_{x_{i}x_{j}x_{k}},T_{2})\right| =0$; otherwise $G[\{x_i,\allowbreak x_k,\allowbreak x_z,\allowbreak q_{ij,1},\allowbreak p_{z,1},\allowbreak w_{ijk,1}\}]\cong R_{14}$, a contradiction.
If $\left|T_{x_{i}x_{j}x_{k}}\cup T_{x_{z}}\right| =2$, then $\left|T_{3}-T_{x_{i}x_{j}x_{k}}\right| =0$; otherwise either $R_{14}\subset G$ or $R_{28}\subset G$, a contradiction.
If $\left|T_{x_{i}x_{j}x_{k}}\cup T_{x_{z}}\right| =2$, then $T_0=\varnothing$ by (F\ref{f-16}), Lemma \ref{lem-9}\ref{lem-9:item-5}\ref{lem-9:item-6} and $d(G)<4$.
Thus, if $\left|T_{x_{i}x_{j}x_{k}}\cup T_{x_{z}}\right| =2$, then $G\subset A_{9}$ by (F\ref{f-15}), Lemma \ref{lem-9}\ref{lem-9:item-1}\ref{lem-9:item-2} and \ref{lem-9:item-5}.
Let $\left|T_{x_{i}x_{j}x_{k}}\cup T_{x_{z}}\right| \leq 1$.
If $T_{x_ix_jx_k}=\varnothing$, then $T_0=\varnothing$ by (F\ref{f-16}), Lemma \ref{lem-9}\ref{lem-9:item-6} and \ref{lem-9:item-17}.
Thus, if  $\left|T_{x_{i}x_{j}x_{k}}\right| =0$, then $G\subset H_{51}$ (if $\left|T_{3}\right| =3$) or $G\subset H_{52}$ (if $\left|T_{3}\right| =2$) or  $G\subset H_{30}$ (if $\left|T_{3}\right| =1$) by (F\ref{f-15}), Lemma \ref{lem-9}\ref{lem-9:item-1}\ref{lem-9:item-2} and \ref{lem-9:item-5}.
Let $\left|T_{x_{i}x_{j}x_{k}}\right| =1$.
If $\left|E(T_{x_{i}x_{j}x_{k}},T_{2})\right| =0$, then $\left|E(T_{x_{i}x_{j}x_{k}},T_{1})\right| =0$; otherwise $R_{35}\subset G$, a contradiction.
Thus, if $\left|E(T_{x_{i}x_{j}x_{k}},T_{2})\right| =0$, then $G\subset H_{35}$ (if $\left|T_{3}\right| =4$) or $G\subset H_{37}$ (if $\left|T_{3}\right| =3$) or  $G\subset H_{26}$ (if $\left|T_{3}\right| =2$) or  $G\subset A_{9}$ (if $\left|T_{3}\right| =1$) by (F\ref{f-15}), Lemma \ref{lem-9}\ref{lem-9:item-1}\ref{lem-9:item-2}\ref{lem-9:item-5} and \ref{lem-9:item-6}.
Let $\left|E(T_{x_{i}x_{j}x_{k}},T_{2})\right| \geq 1$.

Suppose that $\left|E(T_{x_{i}x_{j}x_{k}},T_{2})\right| = 1$.
Without loss of generality, let $\left|E(T_{x_{i}x_{j}x_{k}},T_{x_{i}x_{j}})\right| = 1$.
Then $\left|E(T_{x_{i}x_{j}x_{k}},T_{1}-T_{x_{k}})\right| = 0$; otherwise either $R_{35}\subset G$ or $G$ contains a subdivision of $K_{5}$, a contradiction.
Then $\left|E(T_{x_{i}x_{j}x_{k}},T_{x_{k}})\right| = 1$; otherwise $G[\{x_i,\allowbreak x_k,\allowbreak x_z,\allowbreak q_{ij,1},\allowbreak q_{jk,1},\allowbreak p_{k,1},\allowbreak w_{ijk,1}\}]\cong R_{25}$, a contradiction.
Then $\left|T_{x_{i}x_{k}x_{z}}\cup T_{x_{j}x_{k}x_{z}}\right| = 0$; otherwise $R_{20}\subset G$, a contradiction.
Thus,  $G\subset H_{26}$ (if $\left|T_{3}\right| =2$) or $G\subset H_{32}$ (if $\left|T_{3}\right| =1$).

Suppose that $\left|E(T_{x_{i}x_{j}x_{k}},T_{2})\right| = 2$.
Without loss of generality, let $\left|E(T_{x_{i}x_{j}x_{k}},T_{x_{i}x_{k}})\right| = 0$.
Then $\left|E(T_{x_{i}x_{j}x_{k}}, T_{x_{i}})\right| = 1$ and $\left|E(T_{x_{i}x_{j}x_{k}}, T_{x_{k}})\right| = 1$; otherwise $R_{25}\subset G$, a contradiction.
Then $\left|E(T_{x_{i}x_{j}x_{k}},T_{x_{j}})\right| = 0$; otherwise either $G[\{x_i,\allowbreak x_j,\allowbreak x_k,\allowbreak q_{ik,1},\allowbreak q_{jk,1},\allowbreak p_{j,1},\allowbreak w_{ijk,1}\}]\cong R_{35}$ or $G$ contains a subdivision of $K_{5}$, a contradiction.
Then $\left|T_{3}-T_{x_{i}x_{j}x_{k}}\right| = 0$; otherwise $R_{20}\subset G$, a contradiction.
Thus,  $G\subset H_{51}$.

Then $\left|T_{3}-T_{x_{i}x_{j}x_{k}}\right| = 0$; otherwise $R_{20}\subset G$, a contradiction.
If $\left|E(T_{x_{i}x_{j}x_{k}}, T_{x_{i}})\right| = 1$, then $\left|E(T_{x_{i}x_{j}x_{k}}, T_{x_{j}})\right| = 1$ and $\left|E(T_{x_{i}x_{j}x_{k}}, T_{x_{k}})\right| = 1$; otherwise $R_{28}\subset G$, a contradiction.
If $\left|E(T_{x_{i}x_{j}x_{k}}, T_{x_{i}})\right| = 1$, then $\left|T_{0}\right| = 0$; otherwise $G[\{x_j,\allowbreak x_k,\allowbreak x_z,\allowbreak q_{ij,1},\allowbreak p_{i,1},\allowbreak w_{ijk,1},\allowbreak t_{1}\}]\cong R_{25}$, a contradiction.
Thus,  $G\subset H_{6}$ (if $\left|E(T_{x_{i}x_{j}x_{k}}, T_{1})\right| =0$) or $G\subset H_{38}$ (if $\left|E(T_{x_{i}x_{j}x_{k}}, T_{1})\right| \geq 1$) by Lemma \ref{lem-9}\ref{lem-9:item-1}\ref{lem-9:item-5} and \ref{lem-9:item-6}.
\end{proof}

The following theorem follows immediately from Lemmas \ref{lem-10}, \ref{lem-11}, and \ref{lem-12}.
\begin{thm}\label{thm-31}
If $G \in \mathcal{G}_4$, then $G$ is an induced subgraph of a graph in
$\{H_{1},\allowbreak H_{2},\allowbreak H_{3},\allowbreak H_{5},\allowbreak \ldots,\allowbreak  H_{70},\allowbreak A_{2},\allowbreak \ldots,\allowbreak A_{6},\allowbreak A_{9},\allowbreak A_{11},\allowbreak A_{14},\allowbreak A_{91},\allowbreak D_2\}.$
\end{thm}

\section{Planar graphs with $\omega(G)=3$ and $\lambda_{2}(G)\leq 1$}\label{s-4}

Throughout this section, we assume that $\omega(G)=3$.
Let $V(K_{4}-e)=\{x_{0},x_{1},x_{2},x_{3}\}\triangleq \mathcal{K}$, where $x_0\nsim x_3$.
Denote by $T$ the set $V(G)\backslash \mathcal{K}$.
Thus, $T = T_{0} \cup T_{1} \cup  T_{2}\cup  T_{3}$, where $T_{i}$ $(i=0,1,2,3)$ is the set of vertices of $T$ which are adjacent to exactly $i$ vertices of $K_{4}-e$.
Also, we may suppose that $T_{k}=\bigcup\limits_{0\leq i_{1}< \cdots < i_{k}\leq 3} T_{x_{i_{1}}\ldots x_{i_{k}}}$ $(1\leq k \leq 3)$,
where $T_{x_{i_{1}}\ldots x_{i_{k}}}$ is the set of vertices adjacent precisely to $x_{i_{1}},\ldots,x_{i_{k}}$.

For the notation used below, we label the vertices of each class as follows:
$T_0=\{t_r\mid 1\le r\le |T_0|\}$,
$T_{x_0}=\{p_{0,r}\mid 1\le r\le |T_{x_0}|\}$,
$T_{x_1}=\{p_{1,r}\mid 1\le r\le |T_{x_1}|\}$,
$T_{x_2}=\{p_{2,r}\mid 1\le r\le |T_{x_2}|\}$,
$T_{x_3}=\{p_{3,r}\mid 1\le r\le |T_{x_3}|\}$,
$T_{x_0x_1}=\{q_{01,r}\mid 1\le r\le |T_{x_0x_1}|\}$,
$T_{x_0x_2}=\{q_{02,r}\mid 1\le r\le |T_{x_0x_2}|\}$,
$T_{x_0x_3}=\{q_{03,r}\mid 1\le r\le |T_{x_0x_3}|\}$,
$T_{x_1x_2}=\{q_{12,r}\mid 1\le r\le |T_{x_1x_2}|\}$,
$T_{x_1x_3}=\{q_{13,r}\mid 1\le r\le |T_{x_1x_3}|\}$,
$T_{x_2x_3}=\{q_{23,r}\mid 1\le r\le |T_{x_2x_3}|\}$,
$T_{x_0x_1x_3}=\{w_{013,r}\mid 1\le r\le |T_{x_0x_1x_3}|\}$,
$T_{x_0x_2x_3}=\{w_{023,r}\mid 1\le r\le |T_{x_0x_2x_3}|\}$,
$T_{x_0x_1x_2}=\{w_{012,r}\mid 1\le r\le |T_{x_0x_1x_2}|\}$, and
$T_{x_1x_2x_3}=\{w_{123,r}\mid 1\le r\le |T_{x_1x_2x_3}|\}$.
If the upper bound is zero, the corresponding vertex class is understood to be empty.

~\

\begin{lem}\label{lem-13}
Let $K_{4}-e$ be a fixed induced subgraph of $G \in \mathcal{G}_3$, with $\mathcal{K}=\{x_{0},x_{1},x_{2},x_{3}\}$ and $x_{0}\nsim x_{3}$.
Let $i,j,k,z$ be distinct integers in $\{0,1,2,3\}$.
Then the following statements hold:
\begin{enumerate}[label=(\arabic*),leftmargin=2em,itemsep=2pt]

\item\label{lem-13:item-1} $\left|E(T_{x_{i}})\right|=0$,  $\left|T_{x_{0}}\right| \leq 2$, $\left|T_{x_{3}}\right| \leq 2$ and $G[T_{x_{1}}\cup T_{x_{2}}]=aK_{1}\cup bK_{2}$ ($a+2b\leq 8$);

\item\label{lem-13:item-2} $\left|E(T_{x_{i}}, T_{x_{i}x_{j}}\cup T_{x_{1}}\cup T_{x_{2}})\right|=0$ ($i=0,3$);

\item\label{lem-13:item-3} If   $\left|T_{x_{1}}\cup T_{x_{2}}\right| \geq1$, then  $\left|T_{x_{0}}\right|\leq 1$ and $\left|T_{x_{3}}\right|\leq 1$;

\item\label{lem-13:item-4} $G[T_{x_{0}x_{3}}]=sK_{2}\cup (\left|T_{x_{0}x_{3}}\right|-2s)K_{1}$,  $\left|E(T_{x_{1}x_{2}})\right|=\left|E(T_{x_{1}x_{2}},T_{x_{0}x_{3}}\cup T_{3})\right|=0$; if $x_{i}x_{j}\neq x_{0}x_{3}$, $x_{1}x_{2}$, then $\left|E(T_{x_{i}x_{j}})\right|=0$ and $\left|T_{x_{i}x_{j}}\right| \leq 8$;

\item\label{lem-13:item-6} If  $\left|T_{x_{0}x_{1}}\right|\geq 1$ and $\left|T_{x_{2}x_{3}}\right|\geq 1$, then $\left|E(T_{x_{0}x_{1}},T_{x_{2}x_{3}})\right|=1$,   $\left|T_{x_{0}x_{1}}\right|=1$, $\left|T_{x_{2}x_{3}}\right|= 1$ and   $\left|T_{3}\cup T_{x_{0}x_{3}}\right|=0$;

\item\label{lem-13:item-7} $\left|T_{x_{0}x_{1}x_{2}}\cup T_{x_{1}x_{2}x_{3}}\right|=0$,  $\left|T_{x_{0}x_{1}x_{3}}\right|\leq 1$ and $\left|T_{x_{0}x_{2}x_{3}}\right|\leq 1$;

\item\label{lem-13:item-8} $\left|E(T_{0},T_{1})\right|=\left|E(T_{0},T_{2}-(T_{x_{1}x_{2}}\cup T_{x_{0}x_{3}}))\right|=\left|E(T_{0})\right|=0$; if  $u\in T_{x_{1}x_{2}}\cup T_{3}$, then $\left|E(u,T_{0})\right|\leq 1$;  $G[T_{x_{0}x_{3}}\cup T_{0}]=aK_{1}\cup bK_{2}$ ($a,b\geq 0$).
\end{enumerate}
\end{lem}

\begin{proof}
If $\left|E(T_{x_{0}})\right|\geq 1$, then $G[\mathcal{K}\cup\{p_{0,1},\allowbreak p_{0,2}\}]\cong R_{13}$, a contradiction.
If $\left|E(T_{x_{1}})\right|\geq 1$, then $G[\mathcal{K}\cup\{p_{1,1},\allowbreak p_{1,2}\}]\cong R_{16}$, a contradiction.
Thus, $\left|E(T_{x_{i}})\right|=0$.
Then $\left|T_{x_{0}}\right| \leq 2$ and $\left|T_{x_{3}}\right| \leq 2$; otherwise, $G[\mathcal{K}\cup\{p_{0,1},\allowbreak p_{0,2},\allowbreak p_{0,3}\}]\cong R_{23}$, a contradiction.
If $\left|a\right| \geq 9$, then either $R_{97}\subset G$ (if $\left|T_{x_{1}}\right|\geq 9$ or $\left|T_{x_{2}}\right|\geq 9$) or $R_{47}\subset G$ (the remaining cases), a contradiction.
If $P_{3}\subset G[T_{x_{1}}\cup T_{x_{2}}]$ ($K_{3}\not\subset G[T_{x_{1}}\cup T_{x_{2}}]$), then $G[\mathcal{K}\cup\{p_{1,1},\allowbreak p_{2,1},\allowbreak p_{2,2}\}]\cong R_{26}$, a contradiction.
If $\left|b\right| \geq 5$, then $G[\{x_{1},\allowbreak x_{2},\allowbreak p_{1,1},\allowbreak p_{1,2},\allowbreak p_{1,3},\allowbreak p_{2,4},\allowbreak p_{2,5}\}]\cong R_{21}$ ($p_{1,i}\sim p_{2,i}$), a contradiction.
If $\left|a+2b\right| \geq 9$, then $G[\mathcal{K}\cup\{p_{1,2},\ldots,p_{1,5},\allowbreak p_{2,1}\}]\cong R_{47}$ ($p_{1,i}\sim p_{2,i}$), a contradiction.
Thus, \ref{lem-13:item-1} holds.

If $\left|E(T_{x_{0}}, T_{x_{0}x_{j}})\right|\geq 1$, then either $R_{17}\subset G$ (if $j=3$) or $R_{18}\subset G$ (the remaining cases), a contradiction.
If $\left|E(T_{x_{0}}, T_{x_{1}}\cup T_{x_{2}})\right|\geq 1$, then $G[\mathcal{K}\cup\{p_{0,1},\allowbreak p_{1,1}\}]\cong R_{12}$, a contradiction.
Thus, \ref{lem-13:item-2} holds.

If $\left|T_{x_{1}}\cup T_{x_{2}}\right| \geq 1$ and   $\left|T_{x_{0}}\right|\geq 2$, then $G[\mathcal{K}\cup\{p_{0,1},\allowbreak p_{0,2},\allowbreak p_{1,1}\}]\cong R_{22}$ by Lemma \ref{lem-13}\ref{lem-13:item-2}, a contradiction.
Thus,  \ref{lem-13:item-3} holds.

If $P_{3}\subset G[T_{x_{0}x_{3}}]$ or $K_{3}\subset G[T_{x_{0}x_{3}}]$, then either $G[\{x_{0},\allowbreak x_{1},\allowbreak x_{2},\allowbreak q_{03,1},\allowbreak q_{03,2},\allowbreak q_{03,3}\}]\cong R_{16}$ or $\omega(G)=4$, a contradiction.
If $\left|E(T_{x_{1}x_{2}})\right|\geq 1$, then $\omega(G)=4$, a contradiction.
If $\left|E(T_{x_{1}x_{2}}, T_{x_{0}x_{3}}\cup T_{3})\right|\geq 1$, then $G$ contains a subdivision of $K_{3,3}$ by Lemma \ref{lem-13}\ref{lem-13:item-7}, a contradiction.
If $\left|E(T_{x_{i}x_{j}})\right|\geq 1$, then $\omega(G)=4$, a contradiction.
If $\left|T_{x_{i}x_{j}}\right|\geq 9$, then  $G[\mathcal{K}\cup\{q_{ij,1},\ldots,q_{ij,9}\}]\cong R_{70}$, a contradiction.
Thus,  \ref{lem-13:item-4} holds.

If  $\left|E(T_{x_{0}x_{1}},T_{x_{2}x_{3}})\right|\neq 1$, then either $G[\mathcal{K}\cup\{q_{01,1},\allowbreak q_{23,1}\}]\cong R_{18}$ ($\left|E(T_{x_{0}x_{1}},T_{x_{2}x_{3}})\right|=0$) or $G$ contains a subdivision of $K_{3,3}$, a contradiction.
If   $\left|T_{x_{0}x_{1}}\right|\geq 2$ or $\left|T_{x_{2}x_{3}}\right|\geq 2$, then $G[\mathcal{K}\cup\{q_{01,1},\allowbreak q_{23,1}\}]\cong R_{18}$, a contradiction.
If   $\left|T_{3}\cup T_{x_{0}x_{3}}\right|\geq 1$, then either $\omega(G)=4$ or  $G$ contains a subdivision of $K_{3,3}$, a contradiction.
Thus,  \ref{lem-13:item-6} holds.

If $\left|T_{x_{0}x_{1}x_{2}}\cup T_{x_{1}x_{2}x_{3}}\right|\geq 1$, then  $\omega(G)=4$, a contradiction.
If  $\left|T_{x_{0}x_{1}x_{3}}\right|\geq 2$ or $\left|T_{x_{0}x_{2}x_{3}}\right|\geq 2$, then  $G$ contains a subdivision of $K_{3,3}$, a contradiction.
Thus,  \ref{lem-13:item-7} holds.

If $\left|E(T_{0},T_{1})\right|\geq 1$, then either $R_{6}\subset G$ (if $\left|E(T_{0},T_{x_{0}})\right|\geq 1$) or $R_{9}\subset G$ (if $\left|E(T_{0},T_{x_{1}})\right|\geq 1$), a contradiction.
If $\left|E(T_{0},T_{2}-(T_{x_{1}x_{2}}\cup T_{x_{0}x_{3}}))\right|\geq 1$, then $G[\mathcal{K}\cup\{q_{01,1},\allowbreak t_{1}\}]\cong R_{14}$, a contradiction.
If $t_{1},t_{2}\in T_{0}$, $t_{1}\sim t_{2}$, $u\in T_{x_{0}x_{3}}\cup T_{x_{1}x_{2}}\cup T_{x_{0}x_{1}x_{3}}$ and  $\left|E(u,\{t_{1},t_{2}\})\right|=1$, then $d(G)=4$, a contradiction.
If $u\in T_{x_{0}x_{3}}$ and $\left|E(u,\{t_{1},t_{2}\})\right|=2$, then $G[\{x_{0},\allowbreak x_{1},\allowbreak x_{2},\allowbreak u,\allowbreak t_1,\allowbreak t_2\}]\cong R_{5}$, a contradiction.
If $u\in T_{x_{1}x_{2}}\cup T_{x_{0}x_{1}x_{3}}$ and $\left|E(u,\{t_{1},t_{2}\})\right|=2$, then $G[\{x_{0},\allowbreak x_{1},\allowbreak x_{2},\allowbreak u,\allowbreak t_1,\allowbreak t_2\}]\cong R_{13}$, a contradiction.
Thus,  $\left|E(T_{0})\right|=0$.
If $u\in T_{x_{1}x_{2}}$ and $\left|E(u,T_{0})\right|\geq 2$, then  $G[\mathcal{K}\cup\{u,\allowbreak t_{1},\allowbreak t_{2}\}]\cong R_{24}$, a contradiction.
If $u\in T_{x_{0}x_{1}x_{3}}$ and $\left|E(u,T_{0})\right|\geq 2$, then  $G[\mathcal{K}\cup\{u,\allowbreak t_{1},\allowbreak t_{2}\}]\cong R_{27}$, a contradiction.
Thus,  $\left|E(u,T_{0})\right|\leq 1$.
If  $t_{1}\in T_{0}$, $x,y\in T_{x_{0}x_{3}}$, $x\sim y$ and $\left|E(t_{1},\{x,y\})\right|\geq 1$, then either  $R_{10}\subset G$ (if $\left|E(t_{1},\{x,y\})\right|=1$) or   $R_{13}\subset G$ (if $\left|E(t_{1},\{x,y\})\right|=2$), a contradiction.
If $u\in T_{x_{0}x_{3}}$ and $\left|E(u,T_{0})\right|\geq 2$, then  $G[\{x_{0},\allowbreak x_{1},\allowbreak x_{2},\allowbreak u,\allowbreak t_{1},\allowbreak t_{2}\}]\cong R_{4}$, a contradiction.
If $t\in T_{0}$ and $\left|E(t,T_{x_{0}x_{3}})\right|\geq 2$, then $G[\{x_{0},\allowbreak x_{1},\allowbreak x_{2},\allowbreak q_{03,1},\allowbreak q_{03,2},\allowbreak t\}]\cong R_{8}$, a contradiction.
Thus,  \ref{lem-13:item-8} holds.
\end{proof}

\begin{lem}\label{lem-42}
If $G \in \mathcal{G}_3$, $E(T_{1})\neq \varnothing$, and $G$ contains $K_{4}-e$ as an induced subgraph, then $G$ is an induced subgraph of a graph in $\{H_{71},\ldots,H_{101},H_{111},H_{115},H_{131},D_1,D_2\}$.
\end{lem}

\begin{proof}
Let $K_{4}-e$ be an induced subgraph of $G$ with $\mathcal{K}=\{x_{0},x_{1},x_{2},x_{3}\}$ and $x_{0}\nsim x_{3}$.
The vertices listed in each induced subgraph below are distinct.
By  Lemma \ref{lem-13}\ref{lem-13:item-1} and \ref{lem-13:item-2}, we distinguish the following two cases.

\noindent$\textbf{Case 1.}$ Suppose that $\left|E(T_{x_{0}},T_{x_{3}})\right|\geq 1$.

Then $\left|T_{x_{1}}\cup T_{x_{2}}\right|=0$ by  Lemma \ref{lem-13}\ref{lem-13:item-2}; otherwise, $G[\{x_{1},\allowbreak x_{2},\allowbreak x_{3},\allowbreak p_{0,1},\allowbreak p_{3,1},\allowbreak p_{1,1}\}]\cong R_{3}$, a contradiction.
Then $\left|T_{2}-(T_{x_{0}x_{3}}\cup T_{x_{1}x_{2}})\right|=0$ by  Lemma \ref{lem-13}\ref{lem-13:item-2}; otherwise, either $R_{12}\subset G$ (if $\left|E(T_{x_{0}x_{1}},T_{x_{3}})\right|\geq 1$) or $R_{14}\subset G$ (the remaining cases), a contradiction.
Then $\left|T_{x_{1}x_{2}}\right|=0$; otherwise, either $R_{6}\subset G$ (if $\left|E(T_{x_{1}x_{2}},T_{x_{0}}\cup T_{x_{3}})\right|=0$) or  $G$ contains a subdivision of $K_{3,3}$ (the remaining cases), a contradiction.
Then $\left|T_{3}\right|=0$ by  Lemma \ref{lem-13}\ref{lem-13:item-7}; otherwise, $R_{9}\subset G$ (if $\left|E(T_{x_{0}x_{1}x_{3}},T_{x_{0}}\cup T_{x_{3}})\right|=0$) or  $R_{14}\subset G$ (if $\left|E(T_{x_{0}x_{1}x_{3}},T_{x_{0}}\cup T_{x_{3}})\right|=1$) or $R_{18}\subset G$ (if $\left|E(T_{x_{0}x_{1}x_{3}},T_{x_{0}}\cup T_{x_{3}})\right|=2$), a contradiction.
If  $p\in T_{x_{0}}$, then $\left|E(p,T_{x_{3}})\right|\leq 1$ by  Lemma \ref{lem-13}\ref{lem-13:item-1}; otherwise, $G[\{x_{1},\allowbreak x_{2},\allowbreak x_{3},\allowbreak p,\allowbreak p_{3,1},\allowbreak p_{3,2}\}]\cong R_{8}$, a contradiction.
Thus, if  $\left|T_{x_{0}}\cup T_{x_{3}}\right| =4$, then $\left|E(T_{x_{0}},T_{x_{3}})\right|=2$ by  Lemma \ref{lem-13}\ref{lem-13:item-1}; otherwise, $G[\{x_{0},\allowbreak x_{1},\allowbreak x_{3},\allowbreak p_{0,1},\allowbreak p_{0,2},\allowbreak p_{3,2}\}]\cong R_{1}$ ($E(\{p_{0,1},p_{0,2}\},\{p_{3,1},p_{3,2}\})=\{p_{0,1}p_{3,1}\}$), a contradiction.
Thus, $G\subset D_{1}$ by  Lemma \ref{lem-13}\ref{lem-13:item-1}\ref{lem-13:item-2} and \ref{lem-13:item-8}.

In the subsequent discussion, we assume that $\left|E(T_{x_{0}},T_{x_{3}})\right|=0$.
Thus, since $R_{1}\not\subset G$,  we obtain the following fact.

\begin{fact}\label{case2-f1}
If $\left|T_{x_{0}}\right| \geq 1$, then $\left|T_{x_{3}}\right| \leq 1$.
\end{fact}

\noindent$\textbf{Case 2.}$ Suppose that $\left|E(T_{x_{1}},T_{x_{2}})\right|\geq 1$.

\begin{fact}\label{case2-f2}
$\left|T_{3}\cup T_{x_{0}x_{3}}\right|=0$.
\end{fact}
\begin{proof}
If  $\left|T_{3}\right|\geq 1$, then  either $R_{12}\subset G$ (if $\left|E(T_{3},T_{x_{1}}\cup T_{x_{2}})\right|=0$) or $G$ contains a subdivision of $K_{3,3}$ (the remaining cases), a contradiction.
If  $\left|T_{x_{0}x_{3}}\right|\geq 1$, then either $R_{2}\subset G$ (if $\left|E(T_{x_{0}x_{3}},T_{x_{1}}\cup T_{x_{2}})\right|=0$) or $G$ contains a subdivision of $K_{3,3}$ (the remaining cases), a contradiction.
\end{proof}

\begin{fact}\label{case2-f3}
For  $q\in T_{x_{1}}$, $y\in T_{x_{2}}$, $q\sim y$  and $s\in T_{2}-(T_{x_{0}x_{3}}\cup T_{x_{1}x_{2}})$, we have $\left|E(s,\{q,y\})\right|=1$.
\end{fact}
\begin{proof}
If  $\left|E(s,\{q,y\})\right|\neq 1$, then either $R_{7}\subset G$ (if $\left|E(s,\{q,y\})\right|=2$) or $R_{12}\subset G$ (if $\left|E(s,\{q,y\})\right|=0$), a contradiction.
\end{proof}

Let  $q\in T_{x_{1}}$, $y\in T_{x_{2}}$ and $q\sim y$.
We distinguish the following two subcases according to  the edge configuration between $T_{x_{1}x_{2}}$ and $T_{x_{1}}\cup T_{x_{2}}$.

\noindent$\textbf{Subcase 2.1.}$ Suppose that $\left|E(T_{x_{1}x_{2}},\{q,y\})\right| \geq 1$.

\begin{fact}\label{case2-f4}
$\left|T_{x_{0}x_{i}}\cup T_{x_{j}x_{3}}\right| \leq 1$ for all $i,j\in \{1,2\}$.
\end{fact}
\begin{proof}
By (F\ref{case2-f3}) and  $G$ contains no subdivision of $K_{3,3}$, the (F\ref{case2-f4}) holds.
\end{proof}

If  $q_{12,1}\in T_{x_{1}x_{2}}$, then $\left|E(q_{12,1},\{q,y\})\right|\leq 1$; otherwise, $G[\mathcal{K}\cup\{q,\allowbreak y,\allowbreak q_{12,1}\}]\cong R_{33}$, a contradiction.
If $q_{12,1},q_{12,2}\in T_{x_{1}x_{2}}$ and $q_{12,1}\sim q$, then $\left|E(q_{12,2},\{q,y\})\right|=1$  by Lemma \ref{lem-13}\ref{lem-13:item-4}; otherwise, $G[\mathcal{K}\cup\{q,\allowbreak y,\allowbreak q_{12,1},\allowbreak q_{12,2}\}]\cong R_{42}$, a contradiction.
Thus, $\left|T_{0}\right| =0$ by Lemma \ref{lem-13}\ref{lem-13:item-8} and (F\ref{case2-f2}); otherwise, $G[\mathcal{K}\cup\{q,\allowbreak q_{12,1},\allowbreak t_{1}\}]\cong R_{28}$ ($q_{12,1}\sim q,t_1$), a contradiction.
Moreover, $\left|T_{x_{1}x_{2}}\right| \leq 2$ for $G$ contains no subdivision of $K_{3,3}$.
We distinguish the following two subsubcases.

\noindent$\textbf{Subsubcase 2.1.1.}$ Suppose that $\left|T_{x_{1}x_{2}}\right|= 2$.

Then $\left|T_{x_{0}}\cup T_{x_{3}}\right| =0$  by Lemma \ref{lem-13}\ref{lem-13:item-2}; otherwise, either $R_{14}\subset G$ (if $\left|E(T_{x_{0}},T_{x_{1}x_{2}})\right|\leq 1$) or $G$ contains a subdivision of $K_{3,3}$ (the remaining cases), a contradiction.
Then $T_{2}=T_{x_{1}x_{2}}$ by (F\ref{case2-f2}), (F\ref{case2-f3}) and $G$ contains no subdivision of $K_{3,3}$.
Then $\left|E(T_{x_{1}},T_{x_{2}})\right|= 1$ by Lemma \ref{lem-13}\ref{lem-13:item-1}; otherwise, either $R_{12}\subset G$ (if $\left|E(T_{x_{1}}\cup T_{x_{2}},T_{x_{1}x_{2}})\right|\leq 2$) or $G$ contains a subdivision of $K_{3,3}$ (the remaining cases), a contradiction.
Let $q_{12,1},q_{12,2}\in T_{x_{1}x_{2}}$.

Suppose that $q_{12,1}\sim q$ and $q_{12,2}\sim y$.
Then $\left|E(T_{x_{1}}-\{q\},q_{12,1})\right| =\left|E(T_{x_{2}}-\{y\},q_{12,2})\right| =0$ by Lemma \ref{lem-13}\ref{lem-13:item-1} and \ref{lem-13:item-4}; otherwise, either $R_{18}\subset G$ (if $\left|E(T_{x_{1}}-\{q\},q_{12,2})\right| =0$) or $G$ contains a subdivision of $K_{3,3}$ (the remaining cases), a contradiction.
Then $\left|E(T_{x_{1}},q_{12,2})\right|\leq 1$ and $\left|E(T_{x_{2}},q_{12,1})\right| \leq 1$ by Lemma \ref{lem-13}\ref{lem-13:item-1} and \ref{lem-13:item-4}; otherwise, $G[\{x_{1},\allowbreak q,\allowbreak p_{1,2},\allowbreak p_{1,3},\allowbreak q_{12,1},\allowbreak q_{12,2}\}]\cong R_{16}$ ($p_{1,2},p_{1,3}\sim q_{12,2}$), a contradiction.
Then at most two vertices in $T_{x_{1}}$ ($T_{x_{2}}$) are not adjacent to $q_{12,2}$ ($q_{12,1}$)  by Lemma \ref{lem-13}\ref{lem-13:item-1} and \ref{lem-13:item-4}; otherwise, $G[\{x_{0},\allowbreak x_{1},\allowbreak x_{2},\allowbreak p_{1,2},\allowbreak p_{1,3},\allowbreak y,\allowbreak q_{12,2}\}]\cong R_{25}$, a contradiction.
Thus,  $\left|T_{x_{1}}\right|  \leq 3$ and $\left|T_{x_{2}}\right|  \leq 3$.
If $\left|E(T_{x_{1}},q_{12,2})\right| =1$, then $\left|E(T_{x_{2}},q_{12,1})\right| =0$; otherwise, $G[\{x_{1},\allowbreak x_{2},\allowbreak p_{1,2},\allowbreak p_{2,2},\allowbreak q_{12,1},\allowbreak q_{12,2}\}]\cong R_{18}$ ($p_{1,2}\sim q_{12,2},\ p_{2,2}\sim q_{12,1}$), a contradiction.
Thus,   $G\subset H_{71}$ by (F\ref{case2-f2}).

Suppose that $q_{12,1}\sim q$ and $q_{12,2}\sim q$.
If $s\in T_{x_{1}}\cup T_{x_{2}}-\{q\}$, then  $\left|E(s,\{q_{12,1},q_{12,2}\})\right|\leq 1$ for $G$ contains no subdivision of $K_{3,3}$.
Then  $\left|E(T_{x_{1}}-\{q\},q_{12,1})\right|\leq 1$ and $\left|E(T_{x_{1}}-\{q\},q_{12,2})\right|\leq 1$ by Lemma \ref{lem-13}\ref{lem-13:item-1}  and \ref{lem-13:item-4}; otherwise, $G[\{x_{0},\allowbreak x_{1},\allowbreak x_{2},\allowbreak p_{1,2},\allowbreak p_{1,3},\allowbreak y,\allowbreak q_{12,1}\}]\cong R_{28}$, a contradiction.
Then at most one vertex in $T_{x_{1}}$ is not adjacent to $q_{12,1}$ and $q_{12,2}$ by Lemma \ref{lem-13}\ref{lem-13:item-1} and \ref{lem-13:item-4}; otherwise, $R_{47}\subset G$, a contradiction.
Thus,  $\left|T_{x_{1}}-\{q\}\right|  \leq 3$.
If $p_{2,1}\in T_{x_{2}}-\{y\}$, then  $\left|E(p_{2,1},\{q_{12,1},q_{12,2}\})\right|= 1$; otherwise, $G[\mathcal{K}\cup\{q,\allowbreak p_{2,1},\allowbreak q_{12,1},\allowbreak q_{12,2}\}]\cong R_{43}$ ($p_{2,1}\nsim q_{12,1},q_{12,2}$), a contradiction.
Then  $\left|E(T_{x_{2}}-\{y\},q_{12,1})\right|\leq 1$  and $\left|E(T_{x_{2}}-\{y\},q_{12,2})\right|\leq 1$ by Lemma \ref{lem-13}\ref{lem-13:item-1}  and \ref{lem-13:item-4}; otherwise, $G[\{x_{0},\allowbreak x_{1},\allowbreak x_{2},\allowbreak q,\allowbreak p_{2,2},\allowbreak p_{2,3},\allowbreak q_{12,2}\}]\cong R_{25}$, a contradiction.
Thus,  $\left|T_{x_{2}}-\{y\}\right|  \leq 2$.
If  $p_{2,1}\in T_{x_{2}}-\{y\}$, $p_{1,1}\in T_{x_{1}}-\{q\}$ and $p_{2,1}\sim q_{12,1}$, then $p_{1,1}\nsim q_{12,2}$; otherwise, $G[\{x_{1},\allowbreak x_{2},\allowbreak p_{1,1},\allowbreak p_{2,1},\allowbreak q_{12,1},\allowbreak q_{12,2}\}]\cong R_{18}$, a contradiction.
Thus, $G\subset H_{111}$ (if $\left|T_{x_{2}}-\{y\}\right|=2$) or $G\subset H_{115}$ (if $\left|T_{x_{2}}-\{y\}\right|=1$) or $G\subset H_{131}$ (if $\left|T_{x_{2}}-\{y\}\right|=0$) by (F\ref{case2-f2}).

\noindent$\textbf{Subsubcase 2.1.2.}$ Suppose that $\left|T_{x_{1}x_{2}}\right| = 1$.

Without loss of generality, let $\left|E(T_{x_{1}x_{2}},q)\right|=1$.

\begin{fact}\label{case2-f5}
$\left|E(T_{x_{0}x_{i}}\cup T_{x_{i}x_{3}},T_{x_{1}x_{2}})\right| = \left|E(T_{x_{j}x_{2}},T_{x_{2}})\right| = \left|E(T_{x_{j}x_{1}},T_{x_{j}x_{2}})\right| =0$ ($i=1,2$, $j=0,3$).
\end{fact}
\begin{proof}
If $\left|E(T_{x_{0}x_{i}},T_{x_{1}x_{2}})\right|\geq 1$, then either $G$ contains a subdivision of $K_{3,3}$ (if $\left|E(T_{x_{0}x_{1}},q)\right|=0$) or $\omega=4$ (the remaining cases)  by (F\ref{case2-f3}) and (F\ref{case2-f4}), a contradiction.
If $\left|E(T_{x_{i}x_{3}},T_{x_{1}x_{2}})\right|\geq 1$, then  $G$ contains a subdivision of $K_{3,3}$  by (F\ref{case2-f3}), a contradiction.
If $\left|E(T_{x_{j}x_{2}},y)\right|\geq 1$, then $G[\{x_{0},\allowbreak x_{1},\allowbreak x_{3},\allowbreak y,\allowbreak q_{12,1},\allowbreak q_{02,1}\}]\cong R_{1}$ ($j=0,\ q_{02,1}\sim y$)  by (F\ref{case2-f3}), a contradiction.
If $\left|E(T_{x_{j}x_{2}},T_{x_{2}}-\{y\})\right|\geq 1$, then either $R_{7}\subset G$ (if $\left|E(T_{x_{2}}-\{y\},T_{x_{1}x_{2}})\right|=0$) or $G$ contains a subdivision of $K_{3,3}$ (the remaining cases) by Lemma \ref{lem-13}\ref{lem-13:item-1}, a contradiction.
If $\left|E(T_{x_{j}x_{1}},T_{x_{j}x_{2}})\right|\geq 1$, then either $R_{14}\subset G$ (if $\left|E(T_{x_{0}x_{1}},q)\right|=1$) or $G$ contains a subdivision of $K_{3,3}$ (the remaining cases)  by (F\ref{case2-f3}), a contradiction.
\end{proof}

\begin{fact}\label{case2-f6}
$\left|E(T_{x_{1}x_{2}},T_{x_{1}}-\{q\})\right| \leq 1$, $\left|E(T_{x_{1}x_{2}},T_{x_{2}}-\{y\})\right| \leq 2$; if $\left|E(T_{x_{1}x_{2}},T_{x_{1}}-\{q\})\right| =0$, then $\left|T_{x_{1}}-\{q\}\right| \leq 2$; if $\left|E(T_{x_{1}x_{2}},T_{x_{2}}-\{y\})\right| =0$, then $\left|T_{x_{2}}-\{y\}\right| \leq 1$; $\left|T_{x_{1}}\right| \leq 4$ and $\left|T_{x_{2}}\right| \leq 4$.
\end{fact}
\begin{proof}
If $\left|E(T_{x_{1}x_{2}},T_{x_{1}}-\{q\})\right| \geq 2$, then $G[\{x_{0},\allowbreak x_{1},\allowbreak x_{2},\allowbreak p_{1,2},\allowbreak p_{1,3},\allowbreak y,\allowbreak q_{12,1}\}]\cong R_{28}$ by Lemma \ref{lem-13}\ref{lem-13:item-1}, a contradiction.
If $\left|E(T_{x_{1}x_{2}},T_{x_{2}}-\{y\})\right| \geq 3$, then $G[\{x_{0},\allowbreak x_{1},\allowbreak x_{3},\allowbreak p_{2,2},\allowbreak p_{2,3},\allowbreak p_{2,4},\allowbreak q_{12,1}\}]\cong R_{21}$ by Lemma \ref{lem-13}\ref{lem-13:item-1}, a contradiction.
Then at most two vertices in $T_{x_{1}}-\{q\}$ are not adjacent to $T_{x_{1}x_{2}}$ by Lemma \ref{lem-13}\ref{lem-13:item-1}; otherwise, $G[\{x_{1},\allowbreak x_{2},\allowbreak q,\allowbreak p_{1,2},\allowbreak p_{1,3},\allowbreak p_{1,4},\allowbreak y,\allowbreak q_{12,1}\}]\cong R_{38}$, a contradiction.
Then at most one vertex in $T_{x_{2}}-\{y\}$ is not adjacent to $T_{x_{1}x_{2}}$ by Lemma \ref{lem-13}\ref{lem-13:item-1}; otherwise, $G[\{x_{0},\allowbreak x_{1},\allowbreak x_{2},\allowbreak q,\allowbreak p_{2,2},\allowbreak p_{2,3},\allowbreak q_{12,1}\}]\cong R_{25}$, a contradiction.
Thus,  $\left|T_{x_{1}}\right| \leq 4$ and $\left|T_{x_{2}}\right| \leq 4$.
\end{proof}

Suppose that $\left|T_{x_{0}}\right| \geq 1$ (or $\left|T_{x_{3}}\right| \geq 1$).
Then  $\left|E(T_{x_{0}},T_{x_{1}x_{2}})\right| =\left|T_{x_{0}}\right|$ ($\left|E(T_{x_{3}},T_{x_{1}x_{2}})\right| =\left|T_{x_{3}}\right|$) by Lemma \ref{lem-13}\ref{lem-13:item-2}; otherwise, $G[\{x_{0},\allowbreak x_{1},\allowbreak x_{2},\allowbreak p_{0,1},\allowbreak q,\allowbreak q_{12,1}\}]\cong R_{14}$ ($p_{0,1}\nsim q_{12,1}$), a contradiction.
Then $\left|T_{x_{0}}\right| =1$; otherwise, $G[\{x_{0},\allowbreak x_{1},\allowbreak p_{0,1},\allowbreak p_{0,2},\allowbreak q,\allowbreak y\}]\cong R_{1}$, a contradiction.
Then $\left|T_{x_{3}}\right| =0$ for $G$ contains no subdivision of $K_{3,3}$.
Then $\left|T_{2}-T_{x_{1}x_{2}}\right| =0$ by (F\ref{case2-f2}), (F\ref{case2-f3}) and  $G$ contains no subdivision of $K_{3,3}$.
Then $\left|E(T_{x_{1}},T_{x_{2}})\right| =1$ by Lemma \ref{lem-13}\ref{lem-13:item-1}; otherwise, either $R_{2}\subset G$ (if $\left|E(T_{x_{1}}\cup T_{x_{2}},T_{x_{1}x_{2}})\right| =1$) or $G$ contains a subdivision of $K_{3,3}$ (the remaining cases), a contradiction.
By (F\ref{case2-f6}), $\left|T_{x_{1}}\right| \leq 4$ and $1\leq \left|E(T_{x_{1}x_{2}},T_{x_{1}})\right| \leq 2$.
If  $\left|T_{x_{1}}\right| =1$, then $G\subset H_{73}$ by (F\ref{case2-f2}) and (F\ref{case2-f6}).
If  $\left|T_{x_{1}}\right| =2$ and  $\left|E(T_{x_{1}x_{2}},T_{x_{1}})\right| =1$, then  $\left|E(T_{x_{1}x_{2}},T_{x_{2}}-\{y\})\right| \leq 1$ by Lemma \ref{lem-13}\ref{lem-13:item-1} and \ref{lem-13:item-2}; otherwise, $G[\{x_{0},\allowbreak x_{1},\allowbreak x_{2},\allowbreak p_{1,2},\allowbreak p_{2,2},\allowbreak p_{2,3},\allowbreak q_{12,1}\}]\cong R_{28}$ ($p_{2,2},p_{2,3}\sim q_{12,1},\ p_{1,2}\nsim q_{12,1}$), a contradiction.
Thus, if  $\left|T_{x_{1}}\right| =2$ and  $\left|E(T_{x_{1}x_{2}},T_{x_{1}})\right| =1$, then $G\subset H_{74}$  by (F\ref{case2-f2}) and (F\ref{case2-f6}).
If   $\left|E(T_{x_{1}x_{2}},T_{x_{1}})\right| =2$, then no vertex in $T_{x_{2}}-\{y\}$ is not adjacent to $T_{x_{1}x_{2}}$  by Lemma \ref{lem-13}\ref{lem-13:item-2}; otherwise, $G[\{x_{0},\allowbreak x_{1},\allowbreak x_{2},\allowbreak p_{1,2},\allowbreak y,\allowbreak p_{2,2},\allowbreak q_{12,1}\}]\cong R_{25}$ ($p_{1,2}\sim q_{12,1},\ p_{2,2}\nsim q_{12,1}$), a contradiction.
Thus, if  $\left|T_{x_{1}}\right| =2$ and  $\left|E(T_{x_{1}x_{2}},T_{x_{1}})\right| =2$, then $G\subset H_{76}$ by (F\ref{case2-f2}) and (F\ref{case2-f6}).
Thus, if  $\left|T_{x_{1}}\right| =3$ and  $\left|E(T_{x_{1}x_{2}},T_{x_{1}})\right| =2$, then $G\subset H_{77}$ by (F\ref{case2-f2}) and (F\ref{case2-f6}).
If  $\left|T_{x_{1}}\right| \geq 3$ and  $\left|E(T_{x_{1}x_{2}},T_{x_{1}})\right| =1$, then no vertex in $T_{x_{2}}-\{y\}$ is adjacent to $T_{x_{1}x_{2}}$  by Lemma \ref{lem-13}\ref{lem-13:item-2}; otherwise, $G[\{x_{0},\allowbreak x_{1},\allowbreak x_{2},\allowbreak p_{1,2},\allowbreak p_{1,3},\allowbreak p_{2,2},\allowbreak q_{12,1}\}]\cong R_{25}$ ($p_{1,2},p_{1,3}\nsim q_{12,1},\ p_{2,2}\sim q_{12,1}$), a contradiction.
Thus, if  $\left|T_{x_{1}}\right| =3$ and  $\left|E(T_{x_{1}x_{2}},T_{x_{1}})\right| =1$, then $G\subset H_{75}$ by (F\ref{case2-f2}) and (F\ref{case2-f6}).
Thus, if  $\left|T_{x_{1}}\right| =4$, then $G\subset H_{78}$  by (F\ref{case2-f2}) and (F\ref{case2-f6}).
Let $\left|T_{x_{0}}\cup T_{x_{3}}\right| = 0$.

Suppose that $\left|E(T_{x_{1}},T_{x_{2}})\right| \geq 2$.
Let $p_{1,1}\in T_{x_{1}}-\{q\}$, $p_{2,1}\in T_{x_{2}}-\{y\}$ and $p_{1,1}\sim p_{2,1}$.
Note that  $\left|E(T_{x_{1}x_{2}},\{p_{1,1},p_{2,1}\})\right|\leq 1$.
Then $\left|E(T_{x_{1}x_{2}},\{p_{1,1},p_{2,1}\})\right|=1$  by Lemma \ref{lem-13}\ref{lem-13:item-1}; otherwise, $G[\{x_{1},\allowbreak x_{2},\allowbreak q,\allowbreak p_{1,1},\allowbreak p_{2,1},\allowbreak q_{12,1}\}]\cong R_{12}$ ($q_{12,1}\nsim p_{1,1},p_{2,1}$), a contradiction.
Thus, $\left|E(T_{x_{1}},T_{x_{2}})\right| = 2$ for $G$ contains no subdivision of $K_{3,3}$.
Then $\left|T_{2}-T_{x_{1}x_{2}}\right| =0$ by (F\ref{case2-f2}), (F\ref{case2-f3}) and  $G$ contains no subdivision of $K_{3,3}$.
Then  $\left|E(T_{x_{1}x_{2}},T_{x_{1}})\right|\leq 2$ by Lemma \ref{lem-13}\ref{lem-13:item-1}; otherwise, $R_{28}\subset G$ (whether $\left|E(T_{x_{1}x_{2}},p_{1,1})\right|=0$ or $1$), a contradiction.
Then  $\left|E(T_{x_{1}x_{2}},T_{x_{2}})\right|\leq 2$ by Lemma \ref{lem-13}\ref{lem-13:item-1}; otherwise, $R_{21}\subset G$ (whether $\left|E(T_{x_{1}x_{2}},p_{1,1})\right|=0$ or $1$), a contradiction.
If $\left|E(T_{x_{1}x_{2}},p_{1,1})\right|= 1$, then at most two vertices in $T_{x_{2}}$ are not adjacent to $T_{x_{1}x_{2}}$; otherwise, $G[\{x_{0},\allowbreak x_{1},\allowbreak x_{2},\allowbreak q,\allowbreak p_{2,1},\allowbreak p_{2,3},\allowbreak q_{12,1}\}]\cong R_{25}$ ($p_{2,3}\nsim q_{12,1}$), a contradiction.
Thus, if $\left|E(T_{x_{1}x_{2}},p_{1,1})\right|= 1$ and $\left|T_{x_{2}}\right|= 2$, then $G\subset H_{84}$ by (F\ref{case2-f2}) and (F\ref{case2-f6}).
If $\left|E(T_{x_{1}x_{2}},p_{1,1})\right|= 1$ and $\left|T_{x_{2}}\right|= 3$, then $\left|T_{x_{1}}\right|\leq 3$; otherwise, $G[\{x_{0},\allowbreak x_{1},\allowbreak x_{2},\allowbreak p_{1,3},\allowbreak p_{1,4},\allowbreak p_{2,3},\allowbreak q_{12,1}\}]\cong R_{25}$ ($p_{1,3},p_{1,4}\nsim q_{12,1},\ p_{2,3}\sim q_{12,1}$), a contradiction.
Thus, if $\left|E(T_{x_{1}x_{2}},p_{1,1})\right|= 1$ and $\left|T_{x_{2}}\right|= 3$, then $G\subset H_{82}$  by (F\ref{case2-f2}) and (F\ref{case2-f6}).
If $\left|E(T_{x_{1}x_{2}},p_{1,1})\right|= 1$ and $\left|T_{x_{2}}\right|= 4$, then $\left|T_{x_{1}}\right|=2$; otherwise, $G[\{x_{0},\allowbreak x_{1},\allowbreak x_{2},\allowbreak p_{1,3},\allowbreak p_{2,3},\allowbreak p_{2,4},\allowbreak q_{12,1}\}]\cong R_{28}$ ($p_{1,3}\nsim q_{12,1},\ p_{2,3},p_{2,4}\sim q_{12,1}$), a contradiction.
Thus, if $\left|E(T_{x_{1}x_{2}},p_{1,1})\right|= 1$ and $\left|T_{x_{2}}\right|= 4$, then $G= H_{79}$  by (F\ref{case2-f2}) and (F\ref{case2-f6}).
Let $\left|E(T_{x_{1}x_{2}},p_{2,1})\right|= 1$.
If $\left|E(T_{x_{1}x_{2}},T_{x_{1}})\right|= 2$ (or $\left|E(T_{x_{1}x_{2}},T_{x_{2}})\right|= 2$), then at most one vertex in $T_{x_{2}}$ (or $T_{x_{1}}$) is not adjacent to $T_{x_{1}x_{2}}$; otherwise, $G[\{x_{0},\allowbreak x_{1},\allowbreak x_{2},\allowbreak p_{1,3},\allowbreak y,\allowbreak p_{2,3},\allowbreak q_{12,1}\}]\cong R_{25}$ ($p_{1,3}\sim q_{12,1},\ p_{2,3}\nsim q_{12,1}$), a contradiction.
Thus, if  $\left|E(T_{x_{1}x_{2}},T_{x_{1}})\right|= 2$ (or $\left|E(T_{x_{1}x_{2}},T_{x_{2}})\right|= 2$), then either $G\subset H_{73}$ (if $\left|E(T_{x_{1}x_{2}},T_{x_{2}})\right|= 1$) or $G\subset H_{76}$ (if $\left|E(T_{x_{1}x_{2}},T_{x_{2}})\right|= 2$) by (F\ref{case2-f2}) and (F\ref{case2-f6}).
If $\left|E(T_{x_{1}x_{2}},T_{x_{1}})\right|= \left|E(T_{x_{1}x_{2}},T_{x_{2}})\right|=1$, then $\left|T_{x_{1}}\right|\leq 3$ and $\left|T_{x_{2}}\right|\leq 3$; otherwise, $G[\{x_{0},\allowbreak x_{1},\allowbreak x_{2},\allowbreak p_{1,3},\allowbreak p_{1,4},\allowbreak p_{2,1},\allowbreak q_{12,1}\}]\cong R_{25}$ ($p_{1,3},p_{1,4}\nsim q_{12,1}$), a contradiction.
Thus, if $\left|E(T_{x_{1}x_{2}},T_{x_{1}})\right|= \left|E(T_{x_{1}x_{2}},T_{x_{2}})\right|=1$, then $G\subset H_{81}$ by (F\ref{case2-f2}).
Let $\left|E(T_{x_{1}},T_{x_{2}})\right| =1$.

Suppose that $\left|T_{x_{0}x_{2}}\right| \geq 1$.
Then  $\left|T_{x_{0}x_{2}}\right| = 1$ and $\left|T_{x_{1}x_{3}}\cup T_{x_{2}x_{3}}\right| = 0$ by (F\ref{case2-f4}).
Then $\left|E(T_{x_{0}x_{2}},q)\right| = 1$ by (F\ref{case2-f3}) and (F\ref{case2-f5}); otherwise, $G[\{x_{0},\allowbreak x_{1},\allowbreak x_{2},\allowbreak q,\allowbreak q_{12,1},\allowbreak q_{02,1}\}]\cong R_{18}$, a contradiction.
If $\left|T_{x_{0}x_{1}}\right| = 1$, then  $\left|E(T_{x_{0}x_{1}},q)\right| =1$ by (F\ref{case2-f3}), (F\ref{case2-f5}) and  $G$ contains no subdivision of $K_{3,3}$.
If $\left|T_{x_{0}x_{1}}\right| = 1$, then  $\left|E(T_{x_{1}}\cup T_{x_{2}}-\{q\},T_{x_{1}x_{2}})\right| =0$; otherwise, either $G$ contains a subdivision of $K_{3,3}$ (if  $\left|E(T_{x_{1}}\cup T_{x_{2}}-\{q\},T_{x_{0}x_{1}}\cup T_{x_{0}x_{2}})\right| \geq 1$) or $R_{18}\subset G$ (the remaining cases), a contradiction.
If $\left|T_{x_{0}x_{1}}\right| = 1$, then  $\left|E(T_{x_{1}}\cup T_{x_{2}}-\{q\},T_{x_{0}x_{1}}\cup T_{x_{0}x_{2}})\right| =0$; otherwise, either $G$ contains a subdivision of $K_{3,3}$ (if  $\left|E(T_{x_{1}}-\{q\},T_{x_{0}x_{2}})\right| \geq 1$ or $\left|E(T_{x_{2}},T_{x_{0}x_{1}})\right| \geq 1$) or $R_{14}\subset G$ (the remaining cases), a contradiction.
If $\left|T_{x_{0}x_{1}}\right| = 1$, then  $\left|T_{x_{2}}\right|=1$; otherwise, $G[\{x_{1},\allowbreak x_{2},\allowbreak q,\allowbreak p_{2,2},\allowbreak q_{12,1},\allowbreak q_{01,1}\}]\cong R_{14}$ ($p_{2,2}\nsim q_{01,1},q_{02,1},q_{12,1}$), a contradiction.
If $\left|T_{x_{0}x_{1}}\right| = 1$, then  $\left|T_{x_{1}}\right|\leq 2$; otherwise, $G[\{x_{0},\allowbreak x_{1},\allowbreak x_{2},\allowbreak p_{1,2},\allowbreak p_{1,3},\allowbreak q_{01,1},\allowbreak q_{02,1}\}]\cong R_{25}$ ($p_{1,2},p_{1,3}\nsim q_{01,1},q_{02,1},q_{12,1}$), a contradiction.
Thus, if $\left|T_{x_{0}x_{1}}\right| = 1$, then $G\subset H_{72}$  by (F\ref{case2-f2}).
Let $\left|T_{x_{0}x_{1}}\right| = 0$ by (F\ref{case2-f4}).
Then $\left|E(T_{x_{1}}-\{q\},T_{x_{1}x_{2}})\right| = 0$ for $R_{18}\not\subset G$ and $G$ contains no subdivision of $K_{3,3}$.
Then $\left|E(T_{x_{1}},T_{x_{0}x_{2}})\right| \leq 2$; otherwise, $G[\{x_{2},\allowbreak x_{3},\allowbreak p_{1,2},\allowbreak p_{1,3},\allowbreak y,\allowbreak q_{12,1},\allowbreak q_{02,1}\}]\cong R_{21}$ ($p_{1,2},p_{1,3}\sim q_{02,1}$), a contradiction.
Then at most one vertex in $T_{x_{1}}$ is not adjacent to $T_{x_{0}x_{2}}$; otherwise, $G[\mathcal{K}\cup\{p_{1,2},\allowbreak p_{1,3},\allowbreak q_{02,1}\}]\cong R_{25}$, a contradiction.
Then $\left|E(T_{x_{2}},T_{x_{0}x_{2}})\right| = 0$; otherwise, either $R_{1}\subset G$ (if  $\left|E(T_{x_{2}},T_{x_{1}x_{2}})\right| =0$) or $G$ contains a subdivision of $K_{3,3}$  (the remaining cases), a contradiction.
Then $\left|E(T_{x_{2}},T_{x_{1}x_{2}})\right| \leq 1$; otherwise, $G[\{x_{0},\allowbreak x_{1},\allowbreak p_{2,2},\allowbreak p_{2,3},\allowbreak q_{12,1},\allowbreak q_{02,1}\}]\cong R_{1}$, a contradiction.
Then $\left|E(T_{x_{1}},T_{x_{0}x_{2}})\right| + \left|E(T_{x_{2}},T_{x_{1}x_{2}})\right| \leq 2$; otherwise, $G[\{x_{0},\allowbreak x_{1},\allowbreak p_{1,2},\allowbreak p_{2,2},\allowbreak q_{12,1},\allowbreak q_{02,1}\}]\cong R_{2}$ ($p_{1,2}\sim q_{02,1},\ p_{2,2}\sim q_{12,1}$), a contradiction.
Thus, either $G\subset H_{80}$ (if $\left|E(T_{x_{2}},T_{x_{1}x_{2}})\right| = 1$) or $G\subset H_{83}$ (if $\left|E(T_{x_{2}},T_{x_{1}x_{2}})\right| = 0$) by (F\ref{case2-f2}) and (F\ref{case2-f6}).
By symmetry, let $\left|T_{x_{0}x_{2}}\cup T_{x_{2}x_{3}}\right| = 0$.

Suppose that $\left|T_{x_{0}x_{1}}\right|\geq 1$.
Then  $\left|T_{x_{0}x_{1}}\right| = 1$ and $\left|T_{x_{1}x_{3}}\right| = \left|E(T_{x_{0}x_{1}},T_{x_{1}x_{2}})\right| =0$ by (F\ref{case2-f4})  and (F\ref{case2-f5}).
If $\left|E(T_{x_{0}x_{1}},y)\right| = 1$, then $\left|E(T_{x_{1}}-\{q\},T_{x_{1}x_{2}}\cup T_{x_{0}x_{1}})\right| = 0$ for $R_{16}\not\subset G$ and $G$ contains no subdivision of $K_{3,3}$. 
If $\left|E(T_{x_{0}x_{1}},y)\right| = 1$, then $\left|E(T_{x_{2}}-\{y\},T_{x_{1}x_{2}}\cup T_{x_{0}x_{1}})\right| = 0$ for $R_{10}\not\subset G$ and $G$ contains no subdivision of $K_{3,3}$.
If $\left|E(T_{x_{0}x_{1}},y)\right| = 1$, then $\left|T_{x_{1}}\right| \leq 3$; otherwise, $G[\{x_{0},\allowbreak x_{1},\allowbreak x_{2},\allowbreak p_{1,2},\allowbreak p_{1,3},\allowbreak p_{1,4},\allowbreak y,\allowbreak q_{01,1}\}]\cong R_{38}$, a contradiction.
If $\left|E(T_{x_{0}x_{1}},y)\right| = 1$, then $\left|T_{x_{2}}\right| \leq 2$; otherwise, $G[\mathcal{K}\cup\{p_{2,2},\allowbreak p_{2,3},\allowbreak q_{01,1}\}]\cong R_{25}$, a contradiction.
Thus, if $\left|E(T_{x_{0}x_{1}},y)\right| = 1$, then  $G\subset H_{86}$  by (F\ref{case2-f2}).
Let $\left|E(T_{x_{0}x_{1}},q)\right| = 1$.
Then  $\left|E(T_{x_{1}}-\{q\},T_{x_{0}x_{1}})\right| = 0$ for $R_{14}\not\subset G$ and $G$ contains no subdivision of $K_{3,3}$.
Then  $\left|E(T_{x_{2}},T_{x_{1}x_{2}})\right| = 0$; otherwise, either $G[\{x_{0},\allowbreak x_{1},\allowbreak x_{2},\allowbreak q,\allowbreak p_{2,2},\allowbreak q_{01,1}\}]\cong R_{14}$ ($p_{2,2}\nsim q_{01,1}$) or $G$ contains a subdivision of $K_{3,3}$, a contradiction.
Then $\left|T_{x_{2}}-\{y\}\right| \leq 1$ by (F\ref{case2-f6}).
Then $\left|E(T_{x_{2}}-\{y\},T_{x_{0}x_{1}})\right| = 1$ for $R_{14}\not\subset G$.
Then $\left|E(T_{x_{1}}-\{q\},T_{x_{1}x_{2}})\right| + \left|T_{x_{2}}-\{y\}\right| \leq 1$ by (F\ref{case2-f6}) and $R_{10}\not\subset G$.
Thus, either $G\subset H_{85}$ (if $\left|T_{x_{2}}-\{y\}\right| = 1$) or $G\subset H_{87}$ (if $\left|T_{x_{2}}-\{y\}\right| = 0$).
By symmetry, let $\left|T_{x_{0}x_{1}}\cup T_{x_{1}x_{3}}\right| = 0$.

If  $\left|E(T_{x_{1}x_{2}},T_{x_{2}})\right| = 2$, then $\left|E(T_{x_{1}x_{2}},T_{x_{1}})\right| = \left|T_{x_{1}}\right|$ by (F\ref{case2-f6}); otherwise, $G[\{x_{0},\allowbreak x_{1},\allowbreak x_{2},\allowbreak p_{1,2},\allowbreak p_{2,2},\allowbreak p_{2,3},\allowbreak q_{12,1}\}]\cong R_{28}$ ($p_{1,2}\nsim q_{12,1},\ p_{2,2},p_{2,3}\sim q_{12,1}$), a contradiction.
If  at least one vertex in $T_{x_{2}}-\{y\}$ is not adjacent to $T_{x_{1}x_{2}}$, then  $\left|E(T_{x_{1}x_{2}},T_{x_{1}}-\{q\})\right| = 0$; otherwise, $G[\{x_{0},\allowbreak x_{1},\allowbreak x_{2},\allowbreak p_{1,2},\allowbreak y,\allowbreak p_{2,2},\allowbreak q_{12,1}\}]\cong R_{25}$ ($p_{1,2}\sim q_{12,1},\ p_{2,2}\nsim q_{12,1}$), a contradiction.
Thus, if $\left|E(T_{x_{1}x_{2}},T_{x_{2}})\right| = 2$, then either $G\subset H_{76}$ (if $\left|T_{x_{2}}\right| = 3$) or $G\subset H_{79}$ (if $\left|T_{x_{2}}\right| = 4$) by (F\ref{case2-f2}) and (F\ref{case2-f6}).
If  $\left|E(T_{x_{1}x_{2}},T_{x_{2}})\right| = 1$, then at most one vertex in $T_{x_{1}}$ is not adjacent to $T_{x_{1}x_{2}}$; otherwise, $G[\{x_{0},\allowbreak x_{1},\allowbreak x_{2},\allowbreak p_{1,2},\allowbreak p_{1,3},\allowbreak p_{2,2},\allowbreak q_{12,1}\}]\cong R_{25}$ ($p_{1,2},p_{1,3}\nsim q_{12,1},\ p_{2,2}\sim q_{12,1}$), a contradiction.
Thus, if $\left|E(T_{x_{1}x_{2}},T_{x_{2}})\right| = 1$, then either $G\subset H_{111}$ (if $\left|T_{x_{2}}\right| = 3$) or $G\subset H_{115}$ (if $\left|T_{x_{2}}\right| = 2$) by (F\ref{case2-f2}) and (F\ref{case2-f6}).
Thus,  if $\left|E(T_{x_{1}x_{2}},T_{x_{2}})\right| = 0$, then either $G\subset H_{78}$ (if $\left|T_{x_{2}}\right| = 1$) or $G\subset H_{81}$ (if $\left|T_{x_{2}}\right| = 2$) by (F\ref{case2-f2}) and (F\ref{case2-f6}).

\noindent$\textbf{Subcase 2.2.}$ $\left|E(T_{x_{1}x_{2}},\{q,y\})\right| =0$.

If  $\left|T_{x_{1}x_{2}}\right| \geq 1$, then $\left|T_{2}-T_{x_{1}x_{2}}\right|=0$ by (F\ref{case2-f2}) and (F\ref{case2-f3}); otherwise, $G$ contains a subdivision of $K_{3,3}$  (if $\left|E(T_{x_{0}x_{1}},T_{x_{1}x_{2}})\right| =1$) or $R_{1}\subset G$ (if $\left|E(T_{x_{0}x_{1}},T_{x_{1}x_{2}})\right| =0$ and $\left|E(T_{x_{0}x_{1}},q)\right| =1$) or  $R_{26}\subset G$ (if $\left|E(T_{x_{0}x_{1}},T_{x_{1}x_{2}})\right| =0$ and $\left|E(T_{x_{0}x_{1}},q)\right| =0$), a contradiction.
If  $\left|T_{x_{1}x_{2}}\right| \geq 1$, then  $\left|E(T_{x_{1}x_{2}},T_{x_{0}}\cup T_{x_{3}})\right| =0$; otherwise, $G[\{x_{0},\allowbreak x_{1},\allowbreak p_{0,1},\allowbreak q,\allowbreak y,\allowbreak q_{12,1}\}]\cong R_{2}$ ($p_{0,1}\sim q_{12,1}$), up to symmetry, a contradiction.
If  $\left|T_{x_{1}x_{2}}\right| \geq 1$, then  $\left|E(T_{x_{1}x_{2}},T_{x_{1}}\cup T_{x_{2}}-\{q,y\})\right| =0$ for $R_{12}\not\subset G$.
If  $\left|T_{x_{1}x_{2}}\right| \geq 1$, then  $\left|T_{x_{1}}\right| \leq 4$ and $\left|T_{x_{2}}\right| \leq 4$; otherwise, $G[\{x_{1},\allowbreak x_{2},\allowbreak q,\allowbreak p_{1,2},\ldots,p_{1,5},\allowbreak y\}]\cong R_{37}$, a contradiction.
If $\left|T_{x_{i}}\right| \geq 1$, then  $G[E(T_{x_{1}x_{2}}\cup T_{0})]=aK_{2}$ ($a\geq 0$); otherwise, $G[\mathcal{K}\cup\{y,\allowbreak q_{12,1},\allowbreak q_{12,2},\allowbreak t_{1}\}]\cong R_{40}$ ($t_1\sim q_{12,1},q_{12,2}$), a contradiction.
Thus, if $\left|T_{x_{1}x_{2}}\right| \geq 1$, then $G\subset D_{2}$  by (F\ref{case2-f2}), Lemma \ref{lem-13}\ref{lem-13:item-1} and \ref{lem-13:item-8}.
Let $\left|T_{x_{1}x_{2}}\right| =0$.
Then $\left|T_{0}\right| =0$ by (F\ref{case2-f2}) and Lemma \ref{lem-13}\ref{lem-13:item-8}.

By a simple argument, if $\left|T_{2}\right| =0$, then $G\subset D_{2}$.
Let $\left|T_{2}\right| \geq 1$.
Then $\left|T_{x_{i}x_{j}}\right| \leq 1$ ($i=0,3$ and $j=1,2$) by  (F\ref{case2-f3}) and  $G$ contains no subdivision of $K_{3,3}$.
Without loss of generality, let $\left|T_{x_{0}x_{1}}\right| =1$ by (F\ref{case2-f3}).
Then $\left|T_{x_{2}x_{3}}\right| =0$ by (F\ref{case2-f3}), Lemma \ref{lem-13}\ref{lem-13:item-6} and  $G$ contains no subdivision of $K_{3,3}$.
Then $\left|T_{x_{3}}\right| =0$; otherwise, either $G[\mathcal{K}\cup\{q_{01,1},\allowbreak p_{3,1}\}]\cong R_{14}$ ($p_{3,1}\nsim q_{01,1}$)  or $G$ contains a subdivision of $K_{3,3}$, a contradiction.
Then  $\left|T_{x_{0}}\right| \leq 1$ by Lemma \ref{lem-13}\ref{lem-13:item-1} and \ref{lem-13:item-2}; otherwise, $G[\{x_{0},\allowbreak x_{1},\allowbreak p_{0,1},\allowbreak p_{0,2},\allowbreak q,\allowbreak y\}]\cong R_{1}$, a contradiction.
Then $\left|E(T_{x_{1}},T_{x_{2}})\right| =1$ and $\left|E(T_{x_{0}x_{1}},T_{x_{1}x_{3}})\right|=0$ by (F\ref{case2-f3}) and  $G$ contains no subdivision of $K_{3,3}$.
If $\left|E(T_{x_{0}x_{1}},q)\right|=1$, then $\left|E(T_{x_{0}x_{1}},T_{x_{1}}-\{q\})\right|=0$ for $R_{14}\not\subset G$.

Suppose that $\left|T_{x_{0}x_{2}}\right| =1$.
Then $\left|T_{x_{1}x_{3}}\right| =0$ by (F\ref{case2-f3}), Lemma \ref{lem-13}\ref{lem-13:item-6} and  $G$ contains no subdivision of $K_{3,3}$.
If $\left|E(T_{x_{0}x_{1}},q)\right|=\left|E(T_{x_{0}x_{2}},q)\right|=1$, then  $\left|E(T_{x_{0}x_{1}},T_{x_{0}x_{2}})\right|=0$ by (F\ref{case2-f3}); otherwise, $G[\{x_{0},\allowbreak x_{2},\allowbreak x_{3},\allowbreak q,\allowbreak q_{01,1},\allowbreak q_{02,1}\}]\cong R_{14}$, a contradiction.
If $\left|E(T_{x_{0}x_{1}},q)\right|=\left|E(T_{x_{0}x_{2}},q)\right|=1$, then  $\left|E(T_{x_{0}x_{1}},T_{x_{2}})\right|=\left|E(T_{x_{0}x_{2}},T_{x_{1}}-\{q\})\right|=0$ for $G$ contains no subdivision of $K_{3,3}$.
Thus, if $\left|E(T_{x_{0}x_{1}},q)\right|=\left|E(T_{x_{0}x_{2}},q)\right|=1$, then  $\left|T_{x_{2}}-\{y\}\right|=0$; otherwise, $G[\{x_{0},\allowbreak x_{1},\allowbreak x_{2},\allowbreak q,\allowbreak p_{2,2},\allowbreak q_{01,1}\}]\cong R_{14}$ ($p_{2,2}\nsim q_{01,1},q_{02,1}$), a contradiction.
If $\left|E(T_{x_{0}x_{1}},q)\right|=\left|E(T_{x_{0}x_{2}},q)\right|=1$, then $\left|E(T_{x_{0}x_{1}},T_{x_{1}}-\{q\})\right|=0$; otherwise, $G[\{x_{0},\allowbreak x_{1},\allowbreak x_{2},\allowbreak p_{1,2},\allowbreak y,\allowbreak q_{01,1}\}]\cong R_{14}$ ($p_{1,2}\sim q_{01,1}$), a contradiction.
If $\left|E(T_{x_{0}x_{1}},q)\right|=\left|E(T_{x_{0}x_{2}},q)\right|=1$, then $\left|T_{x_{1}}\right|\leq 2$; otherwise, $G[\mathcal{K}\cup\{p_{1,2},\allowbreak p_{1,3},\allowbreak q_{02,1}\}]\cong R_{25}$ ($p_{1,2},p_{1,3}\nsim q_{01,1},q_{02,1}$), a contradiction.
Thus, if $\left|E(T_{x_{0}x_{1}},q)\right|=\left|E(T_{x_{0}x_{2}},q)\right|=1$, then $G\subset H_{89}$  by  Lemma \ref{lem-13}\ref{lem-13:item-2}.
If $\left|E(T_{x_{0}x_{1}},y)\right|=\left|E(T_{x_{0}x_{2}},q)\right|=1$, then $G$ contains a subdivision of $K_{3,3}$, a contradiction.
If  $\left|E(T_{x_{0}x_{1}},q)\right|=\left|E(T_{x_{0}x_{2}},y)\right|=1$, then $\left|E(T_{x_{0}x_{1}},T_{x_{0}x_{2}})\right|=1$; otherwise, $G[\{x_{0},\allowbreak x_{1},\allowbreak x_{2},\allowbreak q,\allowbreak q_{01,1},\allowbreak q_{02,1}\}]\cong R_{18}$, a contradiction.
If  $\left|E(T_{x_{0}x_{1}},q)\right|=\left|E(T_{x_{0}x_{2}},y)\right|=1$, then  $\left|T_{x_{1}}\cup T_{x_{2}}-\{q,y\}\right|=0$ for $R_{14}\not\subset G$.
Thus, if  $\left|E(T_{x_{0}x_{1}},q)\right|=\left|E(T_{x_{0}x_{2}},y)\right|=1$, then $G\subset H_{88}$ by  Lemma \ref{lem-13}\ref{lem-13:item-2}.
Let $\left|T_{x_{0}x_{2}}\right| =0$.

Suppose that $\left|T_{x_{1}x_{3}}\right| =1$.
Note that $\left|E(T_{x_{0}x_{1}},T_{x_{1}x_{3}})\right|=0$.
By (F\ref{case2-f3}), either $\left|E(T_{x_{0}x_{1}}\cup T_{x_{1}x_{3}},q)\right|=2$ or $\left|E(T_{x_{0}x_{1}}\cup T_{x_{1}x_{3}},y)\right|=2$; otherwise, $G[\{x_{0},\allowbreak x_{1},\allowbreak x_{2},\allowbreak y,\allowbreak q_{01,1},\allowbreak q_{13,1}\}]\cong R_{12}$ ($q_{01,1}\sim q,\ q_{13,1}\sim y$), a contradiction.
Then  $\left|T_{x_{0}}\right| =0$; otherwise, either $G[\mathcal{K}\cup\{q_{13,1},\allowbreak p_{0,1}\}]\cong R_{14}$ ($q_{01,1},q_{13,1}\sim q,\ p_{0,1}\nsim q_{13,1}$) or $G$ contains a subdivision of $K_{3,3}$, a contradiction.
If $\left|E(T_{x_{0}x_{1}},q)\right|=1$, then $\left|E(T_{x_{0}x_{1}}\cup T_{x_{1}x_{3}},T_{x_{1}}-\{q\})\right|=0$ for $R_{14}\not\subset G$.
Thus, if $\left|E(T_{x_{0}x_{1}}\cup T_{x_{1}x_{3}},q)\right|=2$, then $\left|T_{x_{1}}\right|\leq 4$; otherwise, $G[\{x_{1},\allowbreak x_{2},\allowbreak q,\allowbreak p_{1,2},\ldots,p_{1,5},\allowbreak y\}]\cong R_{37}$, a contradiction.
If $\left|E(T_{x_{0}x_{1}}\cup T_{x_{1}x_{3}},q)\right|=2$, then $\left|T_{x_{2}}\right|=1$; otherwise, $G[\{x_{0},\allowbreak x_{1},\allowbreak x_{2},\allowbreak q,\allowbreak p_{2,2},\allowbreak q_{13,1}\}]\cong R_{10}$ ($p_{2,2}\nsim q_{13,1}$) or $G[\{x_{1},\allowbreak x_{2},\allowbreak x_{3},\allowbreak q,\allowbreak p_{2,2},\allowbreak q_{01,1}\}]\cong R_{10}$ ($p_{2,2}\sim q_{13,1},\ p_{2,2}\nsim q_{01,1}$) or $G$ contains a subdivision of $K_{3,3}$, a contradiction.
Thus, if $\left|E(T_{x_{0}x_{1}}\cup T_{x_{1}x_{3}},q)\right|=2$, then  $G\subset H_{100}$ by (F\ref{case2-f2}).
Let $\left|E(T_{x_{0}x_{1}}\cup T_{x_{1}x_{3}},y)\right|=2$.
Then  $\left|E(T_{x_{0}x_{1}}\cup T_{x_{1}x_{3}},T_{x_{1}}-\{q\})\right|=0$; otherwise either $G[x_{1},x_{2},x_{3},p_{1,1},q_{01,1},q_{13,1}]\cong R_{16}$ ($p_{1,1}\sim q_{01,1}$ and $p_{1,1}\nsim q_{13,1}$) or $G$ contains a subdivision of $K_{3,3}$ ($p_{1,1}\sim q_{01,1},q_{13,1}$), a contradiction.
Thus, $\left|T_{x_{1}}\right|\leq 3$; otherwise, $G[\{x_{0},\allowbreak x_{1},\allowbreak x_{2},\allowbreak p_{1,2},\allowbreak p_{1,3},\allowbreak p_{1,4},\allowbreak y,\allowbreak q_{01,1}\}]\cong R_{38}$ ($p_{1,2},p_{1,3},p_{1,4}\nsim q_{01,1},q_{13,1}$), a contradiction.
Then  $\left|E(T_{x_{0}x_{1}}\cup T_{x_{1}x_{3}},T_{x_{2}})\right|=2$; otherwise, either $G[\{x_{0},\allowbreak x_{1},\allowbreak x_{3},\allowbreak p_{2,2},\allowbreak q_{01,1},\allowbreak q_{13,1}\}]\cong R_{10}$ ($\left|E(p_{2,2},\{q_{01,1},q_{13,1}\}\right|=1)$)  or $G$ contains a subdivision of $K_{3,3}$, a contradiction.
Thus, $\left|T_{x_{2}}\right|\leq 2$; otherwise, $G[\mathcal{K}\cup\{p_{2,2},\allowbreak p_{2,3},\allowbreak q_{13,1}\}]\cong R_{25}$ ($p_{2,2},p_{2,3}\nsim q_{01,1},q_{13,1}$), a contradiction.
Thus, $G\subset H_{101}$  by (F\ref{case2-f2}).
Let $\left|T_{x_{1}x_{3}}\right| =0$.

Suppose that $\left|E(T_{x_{0}x_{1}},q)\right|=1$.
Then $\left|E(T_{x_{0}x_{1}},T_{x_{2}})\right|=\left|T_{x_{2}}-\{y\}\right|$; otherwise, $G[\{x_{0},\allowbreak x_{1},\allowbreak x_{2},\allowbreak q,\allowbreak p_{2,2},\allowbreak q_{01,1}\}]\cong R_{14}$ ($p_{2,2}\nsim q_{01,1}$), a contradiction.
Then  $\left|E(T_{x_{0}x_{1}},T_{x_{1}})\right|=1$; otherwise, $G[\{x_{0},\allowbreak x_{1},\allowbreak x_{2},\allowbreak p_{1,2},\allowbreak y,\allowbreak q_{01,1}\}]\cong R_{14}$ ($p_{1,2}\sim q_{01,1}$), a contradiction.
Then $\left|T_{x_{1}}\cup T_{x_{2}}\right|\leq 5$ by Lemma \ref{lem-13}\ref{lem-13:item-1}; otherwise, either $R_{37}\subset G$ (if $\left|T_{x_{1}}\right|\geq 4$ or $\left|T_{x_{2}}\right|\geq 4$) or $R_{21}\subset G$ (the remaining cases), a contradiction.
Thus, $G\subset H_{91}$ (if  $\left|T_{x_{1}}\right|=4$) or $G\subset H_{96}$ (if  $\left|T_{x_{1}}\right|=3$) or $G\subset H_{97}$ (if  $\left|T_{x_{1}}\right|=2$)  or $G\subset H_{92}$ (if  $\left|T_{x_{1}}\right|=1$) by (F\ref{case2-f2}).

Suppose that $\left|E(T_{x_{0}x_{1}},y)\right|=1$.
Then $\left|E(T_{x_{0}x_{1}},T_{x_{1}})\right|\leq 1$; otherwise, $G[\{x_{0},\allowbreak x_{2},\allowbreak x_{3},\allowbreak p_{1,2},\allowbreak p_{1,3},\allowbreak q_{01,1}\}]\cong R_{1}$, a contradiction.
Then at most three vertices in $T_{x_{1}}$ are not adjacent to $T_{x_{0}x_{1}}$; otherwise, $G[\{x_{0},\allowbreak x_{1},\allowbreak x_{2},\allowbreak p_{1,2},\allowbreak p_{1,3},\allowbreak p_{1,4},\allowbreak y,\allowbreak q_{01,1}\}]\cong R_{38}$, a contradiction.
Thus, $\left|T_{x_{1}}\right|\leq 4$.
If $\left|T_{x_{1}}\right|=4$, then $\left|T_{x_{2}}\right|=1$; otherwise, $G[\{x_{0},\allowbreak x_{1},\allowbreak x_{2},\allowbreak p_{1,2},\allowbreak p_{2,2},\allowbreak q_{01,1}\}]\cong R_{14}$ ($p_{1,2}\sim q_{01,1},\ p_{2,2}\nsim q_{01,1}$) or $G[\{x_{1},\allowbreak x_{2},\allowbreak p_{1,2},\allowbreak p_{1,3},\allowbreak y,\allowbreak p_{2,2},\allowbreak q_{01,1}\}]\cong R_{77}$ ($p_{1,2},p_{1,3},p_{1,4}\nsim q_{01,1},\ p_{2,2}\sim q_{01,1}$) or $G[\{x_{1},\allowbreak x_{2},\allowbreak p_{1,3},\allowbreak p_{1,4},\allowbreak y,\allowbreak p_{2,2},\allowbreak q_{01,1}\}]\cong R_{77}$ ($p_{1,2},p_{2,2}\sim q_{01,1}$), a contradiction.
Thus, if $\left|T_{x_{1}}\right|=4$, then $G\subset H_{90}$  by (F\ref{case2-f2}).
Then $\left|E(T_{x_{0}x_{1}},T_{x_{2}})\right|\leq 3$; otherwise, $G[\{x_{1},\allowbreak x_{3},\allowbreak q,\allowbreak p_{2,2},\allowbreak p_{2,3},\allowbreak p_{2,4},\allowbreak q_{01,1}\}]\cong R_{21}$, a contradiction.
Then at most one vertex in $T_{x_{2}}$ is not adjacent to $T_{x_{0}x_{1}}$; otherwise, $G[\mathcal{K}\cup\{p_{2,2},\allowbreak p_{2,3},\allowbreak q_{01,1}\}]\cong R_{25}$, a contradiction.
Thus, if $\left|T_{x_{1}}\right|=1$, then $G\subset H_{93}$  by (F\ref{case2-f2}).
If $\left|T_{x_{1}}\right|=3$ and $\left|E(T_{x_{0}x_{1}},T_{x_{1}})\right|=0$, then $\left|E(T_{x_{0}x_{1}},T_{x_{2}})\right|=1$ and  $\left|T_{x_{2}}\right|\leq 2$; otherwise, $G[\{x_{1},\allowbreak x_{2},\allowbreak p_{1,2},\allowbreak p_{1,3},\allowbreak y,\allowbreak p_{2,2},\allowbreak q_{01,1}\}]\cong R_{77}$ ($p_{2,2}\sim q_{01,1}$), a contradiction.
Thus, if $\left|T_{x_{1}}\right|=3$ and $\left|E(T_{x_{0}x_{1}},T_{x_{1}})\right|=0$, then $G\subset H_{94}$  by (F\ref{case2-f2}).
If $\left|T_{x_{1}}\right|=3$ and $\left|E(T_{x_{0}x_{1}},T_{x_{1}})\right|=1$, then $\left|E(T_{x_{0}x_{1}},T_{x_{2}})\right|=\left|T_{x_{2}}\right|$ and  $\left|T_{x_{2}}\right|\leq 2$; otherwise, either $G[\{x_{0},\allowbreak x_{1},\allowbreak x_{2},\allowbreak p_{1,2},\allowbreak p_{2,2},\allowbreak q_{01,1}\}]\cong R_{14}$ ($p_{1,2}\sim q_{01,1},\ p_{2,2}\nsim q_{01,1}$) or $G[\{x_{1},\allowbreak x_{2},\allowbreak q,\allowbreak p_{1,2},\allowbreak p_{1,3},\allowbreak p_{2,2},\allowbreak p_{2,3}\}]\cong R_{21}$ ($p_{1,2},p_{2,2},p_{2,3}\sim q_{01,1}$), a contradiction.
Thus, if $\left|T_{x_{1}}\right|=3$ and $\left|E(T_{x_{0}x_{1}},T_{x_{1}})\right|=1$, then $G\subset H_{95}$ by (F\ref{case2-f2}).
If $\left|T_{x_{1}}\right|=2$ and $\left|E(T_{x_{0}x_{1}},T_{x_{1}})\right|=0$, then $\left|E(T_{x_{0}x_{1}},T_{x_{2}})\right|\leq 2$; otherwise, $G[\{x_{1},\allowbreak x_{2},\allowbreak q,\allowbreak p_{1,2},\allowbreak p_{2,2},\allowbreak p_{2,3},\allowbreak q_{01,1}\}]\cong R_{77}$, a contradiction.
Thus, if $\left|T_{x_{1}}\right|=2$ and $\left|E(T_{x_{0}x_{1}},T_{x_{1}})\right|=0$, then $G\subset H_{98}$ by (F\ref{case2-f2}).
If $\left|T_{x_{1}}\right|=2$ and $\left|E(T_{x_{0}x_{1}},T_{x_{1}})\right|=1$, then  $\left|E(T_{x_{0}x_{1}},T_{x_{2}})\right|=\left|T_{x_{2}}\right|$; otherwise, $G[\{x_{0},\allowbreak x_{1},\allowbreak x_{2},\allowbreak p_{1,2},\allowbreak p_{2,2},\allowbreak q_{01,1}\}]\cong R_{14}$ ($p_{1,2}\sim q_{01,1},\ p_{2,2}\nsim q_{01,1}$), a contradiction.
Thus, if $\left|T_{x_{1}}\right|=2$ and $\left|E(T_{x_{0}x_{1}},T_{x_{1}})\right|=1$, then $G\subset H_{99}$ by (F\ref{case2-f2}).
\end{proof}

\begin{lem}\label{lem-43}
If $G \in \mathcal{G}_3$, $E(T_{1})=\varnothing$, $E(T_{2})-E(T_{x_{0}x_{3}})\neq\varnothing$, and $G$ contains $K_{4}-e$ as an induced subgraph, then $G$ is an induced subgraph of a graph in $\{H_ {11},\allowbreak H_{71},\allowbreak H_{72},\allowbreak H_{80},\allowbreak H_{83},\allowbreak H_{85},\allowbreak H_{87},\allowbreak H_{88},\allowbreak H_{89},\allowbreak H_{102},\allowbreak \ldots,\allowbreak H_{136},\allowbreak  A_{16},\allowbreak A_{17},\allowbreak A_{21},\allowbreak A_{24},\allowbreak A_{25},\allowbreak A_{27},\allowbreak A_{28},\allowbreak  A_{30},\allowbreak A_{31},\allowbreak A_{34},\allowbreak D_1\}.$
\end{lem}
\begin{proof}
Let $K_{4}-e$ be an induced subgraph of $G$ with $\mathcal{K}=\{x_{0},x_{1},x_{2},x_{3}\}$ and $x_{0}\nsim x_{3}$.
By  Lemma \ref{lem-13}\ref{lem-13:item-4}, we distinguish the following four cases.

\noindent$\textbf{Case 1.}$ Suppose that $\left|E(T_{x_{0}x_{1}},T_{x_{0}x_{2}})\right|\geq 1$.

Let $a\in T_{x_{0}x_{1}}$, $b\in T_{x_{0}x_{2}}$ and $a\sim b$.
Since $R_{12}\not\subset G$ and  $G$ contains no subdivision of $K_{3,3}$, we have the following fact.

\begin{fact}\label{43-f7}
$\left|T_{3}\right|=\left|T_{x_{0}x_{3}}\right|=0$.
\end{fact}

\begin{fact}\label{43-f8}
$\left|E(T_{x_{1}x_{2}},\{a,b\})\right|\geq 1$ (if $\left|T_{x_{1}x_{2}}\right|\geq 1$), $\left|T_{x_{1}x_{2}}\right|\leq 1$ and $\left|E(T_{x_{0}x_{1}}\cup T_{x_{0}x_{2}}-\{a,b\},T_{x_{1}x_{2}})\right|=0$.
\end{fact}
\begin{proof}
If $\left|T_{x_{1}x_{2}}\right|\geq 1$ and $\left|E(T_{x_{1}x_{2}},\{a,b\})\right|=0$, then $R_{33}\subset G$, a contradiction.
If $\left|T_{x_{1}x_{2}}\right|\geq 2$, then $G$ contains a subdivision of $K_{3,3}$, a contradiction.
If $\left|E(T_{x_{0}x_{1}}\cup T_{x_{0}x_{2}}-\{a,b\},T_{x_{1}x_{2}})\right|\geq 1$, then $G$ contains a subdivision of $K_{3,3}$, a contradiction.
\end{proof}

\begin{fact}\label{43-f9}
$\left|E(T_{x_{0}x_{1}}\cup T_{x_{0}x_{2}}-\{a,b\},\{a,b\})\right|=0$,  $\left|T_{x_{0}x_{1}}\right|\leq 4$ and $\left|T_{x_{0}x_{2}}\right|\leq 4$.
\end{fact}
\begin{proof}
If $\left|E(T_{x_{0}x_{1}}\cup T_{x_{0}x_{2}}-\{a,b\},\{a,b\})\right|\geq 1$, then  $G$ contains a subdivision of $K_{3,3}$ by Lemma \ref{lem-13}\ref{lem-13:item-4}, a contradiction.
If $\left|T_{x_{0}x_{1}}\right|\geq 5$ or $\left|T_{x_{0}x_{2}}\right|\geq 5$, then $G[\{x_{1},\allowbreak x_{2},\allowbreak a,\allowbreak q_{01,2},\ldots,q_{01,5},\allowbreak b\}]\cong R_{37}$, a contradiction.
\end{proof}
\begin{fact}\label{43-f10}
If $s\in T_{x_{1}}\cup T_{x_{2}}$, then $\left|E(s,\{a,b\})\right|=1$.
\end{fact}
\begin{proof}
Otherwise either $G[\{x_{0},\allowbreak x_{2},\allowbreak x_{3},\allowbreak a,\allowbreak b,\allowbreak s\}]\cong R_{14}$ ($E(\{a,\allowbreak b,\allowbreak s\})=\{ab,\allowbreak as,\allowbreak bs\}$) or  $G[\mathcal{K}\cup\{a,\allowbreak b,\allowbreak s\}]\cong R_{30}$ ($E(\{a,\allowbreak b,\allowbreak s\})=\{ab\}$), a contradiction.
\end{proof}

\begin{fact}\label{43-f11}
$\left|E(T_{x_{1}},a)\right|\leq 1$, $\left|E(T_{x_{2}},b)\right|\leq 1$, $\left|E(T_{x_{1}},b)\right|\leq 3$, $\left|E(T_{x_{2}},a)\right|\leq 3$  and  we  have $\left|E(T_{x_{1}}\cup T_{x_{2}},T_{x_{0}x_{1}}\cup T_{x_{0}x_{2}}-\{a,b\})\right|=0$.
\end{fact}
\begin{proof}
If $\left|E(T_{x_{1}},a)\right|\geq 2$ or  $\left|E(T_{x_{2}},b)\right|\geq 2$, then  $G[\{x_{0},\allowbreak x_{2},\allowbreak x_{3},\allowbreak a,\allowbreak p_{1,1},\allowbreak p_{1,2}\}]\cong R_{1}$ ($E(\{a,\allowbreak p_{1,1},\allowbreak p_{1,2}\})=\{ap_{1,1},\allowbreak ap_{1,2}\}$) by (F\ref{43-f10}), a contradiction.
If $\left|E(T_{x_{1}},b)\right|\geq 4$ or  $\left|E(T_{x_{2}},a)\right|\geq 4$, then  $G[\mathcal{K}\cup\{b,\allowbreak p_{1,1},\ldots,p_{1,4}\}]\cong R_{79}$ by (F\ref{43-f10}), a contradiction.
We get $\left|E(T_{x_{0}}\cup T_{x_{3}},T_{x_{1}x_{2}})\right|=0$ for $G$ contains no subdivision of $K_{3,3}$.
Thus, if $\left|E(T_{x_{1}}\cup T_{x_{2}},T_{x_{0}x_{1}}\cup T_{x_{0}x_{2}}-\{a,b\})\right|\geq 1$, then either $G[\{x_{0},\allowbreak x_{1},\allowbreak x_{2},\allowbreak q_{01,2},\allowbreak b,\allowbreak p_{1,1}\}]\cong R_{18}$ ($E(\{q_{01,2},\allowbreak b,\allowbreak p_{1,1}\})=\{q_{01,2}p_{1,1}\}$) or $G$ contains a subdivision of $K_{3,3}$ by (F\ref{43-f9}) and by (F\ref{43-f10}), a contradiction.
\end{proof}

\begin{fact}\label{43-f12}
$\left|T_{0}\cup T_{x_{0}}\cup T_{x_{3}}\right|\leq 1$ and if  $\left|E(T_{x_{1}x_{2}},T_{0})\right|= 1$, then $\left|E(T_{x_{1}x_{2}},\{a,b\})\right|=2$.
\end{fact}
\begin{proof}
If $\left|T_{x_{0}}\cup T_{x_{3}}\right|\geq 2$, then either $R_{14}\subset G$ (if $\left|T_{x_{0}}\right|\geq2$) or $d(G)=4$ (if $\left|T_{x_{0}}\right|=\left|T_{x_{3}}\right|=1$) by Lemma \ref{lem-13}\ref{lem-13:item-2}, a contradiction.
By Lemma \ref{lem-13}\ref{lem-13:item-8} and (F\ref{43-f7}),  $\left|E(T-T_{x_{1}x_{2}},T_{0})\right|= 0$.
Thus, by Lemma \ref{lem-13}\ref{lem-13:item-8} and (F\ref{43-f8}), $\left|T_{0}\right|\leq 1$.
If $\left|E(T_{x_{1}x_{2}},T_{0})\right|= 1$ and $\left|E(T_{x_{1}x_{2}},\{a,b\})\right|=1$, then  $G[\{x_{0},\allowbreak x_{1},\allowbreak x_{2},\allowbreak b,\allowbreak q_{12,1},\allowbreak t_{1}\}]\cong R_{14}$ by (F\ref{43-f8}), a contradiction.
If $\left|T_{0}\cup T_{x_{0}}\right|\geq 2$, then either $R_{25}\subset G$ (if  $\left|E(T_{x_{1}x_{2}},T_{x_{0}})\right|= 1$) or $d(G)=4$ (if  $\left|E(T_{x_{1}x_{2}},T_{x_{0}})\right|= 0$) by Lemma \ref{lem-13}\ref{lem-13:item-2}, a contradiction.
\end{proof}

\begin{fact}\label{43-f13}
$\left|E(T_{x_{3}},\{a,b\})\right|=2$ (if $\left|T_{x_{3}}\right|\geq 1$), $\left|T_{x_{3}}\right|\leq 1$, $\left|E(T_{x_{3}},T_{x_{0}x_{1}}\cup T_{x_{0}x_{2}}\cup T_{x_{1}x_{2}}-\{a,b\})\right|=0$; $\left|E(T_{x_{0}},T_{x_{1}x_{2}})\right|=0$ and   $\left|T_{x_{0}}\right|\leq 1$.
\end{fact}
\begin{proof}
If $\left|T_{x_{3}}\right|\geq 1$ and $\left|E(T_{x_{3}},\{a,b\})\right|\leq 1$, then  $R_{14}\subset G$, a contradiction.
If $\left|T_{x_{3}}\right|\geq 2$ or  $\left|E(T_{x_{3}},T_{x_{0}x_{1}}\cup T_{x_{0}x_{2}}\cup T_{x_{1}x_{2}}-\{a,b\})\right|\geq 1$ or $\left|E(T_{x_{0}},T_{x_{1}x_{2}})\right|\geq 1$, then $G$ contains a subdivision of $K_{3,3}$, a contradiction.
If $\left|T_{x_{0}}\right|\geq 2$, then  $G[\mathcal{K}\cup\{a,\allowbreak p_{0,1},\allowbreak p_{0,2}\}]\cong R_{25}$ by Lemma \ref{lem-13}\ref{lem-13:item-2}, a contradiction.
\end{proof}

By (F\ref{43-f9}) and  $G$ contains no subdivision of $K_{3,3}$, $\left|E(T_{x_{0}x_{1}},T_{x_{0}x_{2}})\right|\leq 2$.
Then $\left|E(\{a,b\},T_{x_{1}x_{3}}\cup T_{x_{2}x_{3}})\right|=2\left|T_{x_{1}x_{3}}\cup T_{x_{2}x_{3}}\right|$; otherwise, $G[\{x_{1},\allowbreak x_{2},\allowbreak x_{3},\allowbreak a,\allowbreak b,\allowbreak q_{13,1}\}]\cong R_{12}$ ($E(\{a,\allowbreak b,\allowbreak q_{13,1}\})=\{ab\}$) or $G[\{x_{0},\allowbreak x_{2},\allowbreak x_{3},\allowbreak a,\allowbreak b,\allowbreak q_{13,1}\}]\cong R_{12}$ ($E(\{a,\allowbreak b,\allowbreak q_{13,1}\})=\{ab,\allowbreak bq_{13,1}\}$) or $G[\{x_{0},\allowbreak x_{1},\allowbreak x_{3},\allowbreak a,\allowbreak b,\allowbreak q_{13,1}\}]\cong R_{18}$ ($E(\{a,\allowbreak b,\allowbreak q_{13,1}\})=\{ab,\allowbreak aq_{13,1}\}$), a contradiction.
Thus, $\left|T_{x_{1}x_{3}}\cup T_{x_{2}x_{3}}\right|\leq 1$  for $G$ contains no subdivision of $K_{3,3}$.
Then $\left|E(T_{x_{1}},T_{x_{1}x_{3}})\right|=0$ and $\left|E(T_{x_{2}},T_{x_{2}x_{3}})\right|=0$; otherwise, either $G[\{x_{0},\allowbreak x_{1},\allowbreak x_{2},\allowbreak a,\allowbreak q_{13,1},\allowbreak p_{1,1}\}]\cong R_{20}$ ($E(\{a,\allowbreak q_{13,1},\allowbreak p_{1,1}\})=\{aq_{13,1},\allowbreak ap_{1,1},\allowbreak q_{13,1}p_{1,1}\}$) or $G$ contains a subdivision of $K_{3,3}$, a contradiction.

Suppose that $\left|T_{x_{1}x_{3}}\cup T_{x_{2}x_{3}}\right|=1$.
Without loss of generality, let  $\left|T_{x_{1}x_{3}}\right|=1$.
Then $\left|E(T_{x_{0}x_{1}}\cup T_{x_{0}x_{2}}-\{a,b\},T_{x_{1}x_{3}})\right|=0$ and  $\left|T_{x_{1}x_{2}}\right| =0$ for $G$ contains no subdivision of $K_{3,3}$.
Thus, $\left|T_{0}\right| =0$ by Lemma \ref{lem-13}\ref{lem-13:item-8}, (F\ref{43-f7}) and (F\ref{43-f12}).
Then $\left|T_{x_{0}x_{1}}\right| \leq 2$  by (F\ref{43-f9}); otherwise, $G[\{x_{0},\allowbreak x_{1},\allowbreak x_{3},\allowbreak q_{01,2},\allowbreak q_{01,3},\allowbreak q_{13,1}\}]\cong R_{16}$, a contradiction.
Then $\left|T_{x_{0}x_{2}}\right| =1$  by (F\ref{43-f9}); otherwise, $G[\{x_{0},\allowbreak x_{1},\allowbreak x_{2},\allowbreak a,\allowbreak q_{02,2},\allowbreak q_{13,1}\}]\cong R_{18}$, a contradiction.
Then $\left|T_{x_{0}}\right| =0$  by Lemma \ref{lem-13}\ref{lem-13:item-2}; otherwise, either $G[\mathcal{K}\cup\{p_{0,1},\allowbreak q_{13,1}\}]\cong R_{14}$ or $G$ contains a subdivision of $K_{3,3}$, a contradiction.
Then $\left|T_{x_{3}}\right| =0$  by (F\ref{43-f13}); otherwise, $G$ contains a subdivision of $K_{3,3}$, a contradiction.
Then $\left|T_{x_{2}}\right|=\left|E(T_{x_{2}},b)\right|=\left|E(T_{x_{2}},T_{x_{1}x_{3}})\right|$ by (F\ref{43-f10}); otherwise, either $G$ contains a subdivision of $K_{3,3}$ (if $\left|E(T_{x_{2}},a)\right|\geq 1$) or $R_{14}\subset G$ (the remaining cases), a contradiction.
Thus, $\left|T_{x_{2}}\right|\leq 1$ by (F\ref{43-f11}).
Then $\left|T_{x_{1}}\right|=\left|E(T_{x_{1}},a)\right|$  and $\left|T_{x_{1}}\right|\leq 1$ by (F\ref{43-f10}), (F\ref{43-f11}) and  $G$ contains no subdivision of $K_{3,3}$.
Then $\left|E(T_{x_{1}},T_{x_{1}x_{3}})\right|=0$ for $K_{4}\not\subset G$.
Thus, $\left|T_{x_{1}}\cup T_{x_{2}}\right|\leq 1$; otherwise, $G[\{x_{0},\allowbreak x_{1},\allowbreak x_{2},\allowbreak b,\allowbreak p_{1,1},\allowbreak p_{2,1}\}]\cong R_{14}$, a contradiction.
Thus, $\left|T_{x_{2}}\cup T_{x_{0}x_{1}}\right|\leq 2$; otherwise, $G[\{x_{0},\allowbreak x_{1},\allowbreak x_{2},\allowbreak q_{01,2},\allowbreak b,\allowbreak p_{2,1}\}]\cong R_{18}$ ($E(\{q_{01,2},\allowbreak b,\allowbreak p_{2,1}\})=\{bp_{2,1}\}$), a contradiction.
Thus, either $G\subset H_{102}$ (if $\left|T_{x_{2}}\right|=0$) or $G\subset H_{103}$ (if $\left|T_{x_{2}}\right|=1$) by (F\ref{43-f7}).
Let $\left|T_{x_{1}x_{3}}\cup T_{x_{2}x_{3}}\right|=0$.

Suppose that $\left|E(T_{x_{1}x_{2}},\{a,b\})\right|\geq 2$.
By (F\ref{43-f8}), $\left|E(T_{x_{1}x_{2}},\{a,b\})\right|=2$ and $\left|T_{x_{1}x_{2}}\right|=1$.
Then $\left|T_{x_{0}x_{1}}\cup T_{x_{0}x_{2}}-\{a,b\}\right|=0$ by (F\ref{43-f8}) and (F\ref{43-f9}); otherwise, $G[\{x_{1},\allowbreak x_{2},\allowbreak x_{3},\allowbreak a,\allowbreak q_{01,2},\allowbreak b,\allowbreak q_{12,1}\}]\cong R_{30}$, a contradiction.
Then $\left|T_{x_{3}}\right|=0$ by (F\ref{43-f13}) and  $G$ contains no subdivision of $K_{3,3}$.
Then $\left|E(T_{x_{1}x_{2}},T_{x_{1}}\cup T_{x_{2}})\right|=0$ by (F\ref{43-f10}); otherwise, either $G$ contains a subdivision of $K_{3,3}$ or $G[\{x_{0},\allowbreak x_{2},\allowbreak x_{3},\allowbreak a,\allowbreak q_{12,1},\allowbreak p_{1,1}\}]\cong R_{7}$ ($E(\{a,\allowbreak q_{12,1},\allowbreak p_{1,1}\})=\{aq_{12,1},\allowbreak ap_{1,1},\allowbreak q_{12,1}p_{1,1}\}$), a contradiction.
Thus, $\left|E(T_{x_{1}},b)\right|=\left|E(T_{x_{2}},a)\right|=0$ by (F\ref{43-f10}); otherwise, $G$ contains a subdivision of $K_{3,3}$, a contradiction.
Thus, $\left|T_{x_{1}}\cup T_{x_{2}}\right|\leq 1$ by (F\ref{43-f11}); otherwise, $G[\{x_{0},\allowbreak x_{2},\allowbreak x_{3},\allowbreak a,\allowbreak p_{1,1},\allowbreak p_{2,1}\}]\cong R_{1}$ ($E(\{a,\allowbreak p_{1,1},\allowbreak p_{2,1}\})=\{ap_{1,1}\}$), a contradiction.
Thus, $G\subset A_{16}$ by (F\ref{43-f7}), (F\ref{43-f12}) and (F\ref{43-f13}).
Let $\left|E(T_{x_{1}x_{2}},\{a,b\})\right|=\left|T_{x_{1}x_{2}}\right|$ by (F\ref{43-f8}).
Then $\left|T_{0}\right|=0$ by (F\ref{43-f12}).

Suppose that $\left|E(T_{x_{0}x_{1}},T_{x_{0}x_{2}})\right|=2$.
By (F\ref{43-f9}), $G[E(T_{x_{0}x_{1}},T_{x_{0}x_{2}})]=2K_{2}$.
Then $\left|T_{x_{1}}\cup T_{x_{2}}\cup T_{x_{3}}\cup T_{x_{1}x_{2}}\right|=0$ by (F\ref{43-f8}), (F\ref{43-f10}) and (F\ref{43-f13}); otherwise, $G$ contains a subdivision of $K_{3,3}$, a contradiction.
Then $\left|T_{x_{0}x_{1}}\cup T_{x_{0}x_{2}}\right|\leq 6$ by (F\ref{43-f9}); otherwise, $G[\{x_{1},\allowbreak x_{2},\allowbreak q_{01,2},\allowbreak q_{01,3},\allowbreak q_{01,4},\allowbreak b,\allowbreak q_{02,3}\}]\cong R_{21}$, a contradiction.
Thus, either $G\subset H_{116}$ (if $\left|T_{x_{0}x_{1}}\right|=4$) or $G\subset H_{117}$ (if $\left|T_{x_{0}x_{1}}\right|=3$)  by (F\ref{43-f7}) and (F\ref{43-f13}).
Let $\left|E(T_{x_{0}x_{1}},T_{x_{0}x_{2}})\right|=1$.

Suppose that $\left|T_{x_{3}}\right|\geq 1$.
Then $\left|T_{x_{3}}\right|= 1$ by (F\ref{43-f13}).
Then,  by (F\ref{43-f8}) and (F\ref{43-f13}), $\left|T_{x_{0}x_{1}}\cup T_{x_{0}x_{2}}\cup T_{x_{1}x_{2}}-\{a,b\}\right|=0$; otherwise, $R_{14}\subset G$, a contradiction.
Since  $G$ contains no subdivision of $K_{3,3}$, $\left|E(T_{x_{1}},b)\right|=\left|E(T_{x_{2}},a)\right|=0$.
Thus, $\left|T_{x_{1}}\cup T_{x_{2}}\right|\leq 1$ by (F\ref{43-f10}) and (F\ref{43-f11}); otherwise, $G[\{x_{0},\allowbreak x_{1},\allowbreak x_{2},\allowbreak a,\allowbreak p_{1,1},\allowbreak p_{2,1}\}]\cong R_{14}$ ($E(\{a,\allowbreak p_{1,1},\allowbreak p_{2,1}\})=\{ap_{1,1}\}$), a contradiction.
Thus, $G\subset H_{88}$ by (F\ref{43-f7}) and (F\ref{43-f13}).
Let  $\left|T_{x_{3}}\right|=0$.
By (F\ref{43-f9}),  $\left|T_{x_{0}x_{1}}\right|\leq 4$.

Suppose that $\left|T_{x_{0}x_{1}}\right|=4$.
Then $\left|T_{x_{0}x_{2}}-\{b\}\right|=\left|T_{x_{1}x_{2}}\right|=0$ by (F\ref{43-f8}) and (F\ref{43-f9}); otherwise, $G[\{x_{1},\allowbreak x_{2},\allowbreak q_{01,2},\allowbreak q_{01,3},\allowbreak q_{01,4},\allowbreak b,\allowbreak q_{02,2}\}]\cong R_{21}$ or $G[\{x_{0},\allowbreak x_{2},\allowbreak x_{3},\allowbreak q_{01,2},\allowbreak q_{01,3},\allowbreak q_{01,4},\allowbreak q_{12,1}\}]\cong R_{21}$, a contradiction.
Then $\left|T_{x_{2}}\right|=0$ by (F\ref{43-f10}) and (F\ref{43-f11}); otherwise, either $G[\{x_{0},\allowbreak x_{1},\allowbreak x_{2},\allowbreak q_{01,2},\allowbreak b,\allowbreak p_{2,1}\}]\cong R_{18}$ ($b\sim p_{2,1}$) or $G[\{x_{0},\allowbreak x_{1},\allowbreak x_{2},\allowbreak q_{01,2},\allowbreak q_{01,3},\allowbreak b,\allowbreak p_{2,1}\}]\cong R_{28}$ ($b\nsim p_{2,1}$), a contradiction.
Then $\left|E(T_{x_{1}},b)\right|=0$ by (F\ref{43-f10}) and (F\ref{43-f11}); otherwise, $G[\{x_{0},\allowbreak x_{1},\allowbreak x_{2},\allowbreak q_{01,2},\allowbreak q_{01,3},\allowbreak q_{01,4},\allowbreak b,\allowbreak p_{1,1}\}]\cong R_{42}$, a contradiction.
Thus, $\left|T_{x_{1}}\right|\leq 1$  by (F\ref{43-f10}) and (F\ref{43-f11}).
Thus, $G\subset H_{114}$ by Lemma \ref{lem-13}\ref{lem-13:item-2}, (F\ref{43-f7}) and (F\ref{43-f13}).
By symmetry, if $\left|T_{x_{0}x_{2}}\right|=4$, we obtain a similar conclusion.
Let $\left|T_{x_{0}x_{1}}\right|\leq 3$ and $\left|T_{x_{0}x_{2}}\right|\leq 3$.

Suppose that $\left|T_{x_{0}x_{1}}\right|=3$.
Then $\left|T_{x_{0}x_{2}}\right|\leq 2$ by (F\ref{43-f9}); otherwise, $G[\{x_{1},\allowbreak x_{2},\allowbreak a,\allowbreak q_{01,2},\allowbreak q_{01,3},\allowbreak q_{02,2},\allowbreak q_{02,3}\}]\cong R_{21}$, a contradiction.
Then $\left|E(T_{x_{1}},b)\right|\leq 1$ by (F\ref{43-f10}) and (F\ref{43-f11}); otherwise, $G[\{x_{0},\allowbreak x_{1},\allowbreak q_{01,2},\allowbreak q_{01,3},\allowbreak b,\allowbreak p_{1,1},\allowbreak p_{1,2}\}]\cong R_{26}$, a contradiction.
Then $\left|T_{x_{2}}\right|=0$ by (F\ref{43-f10}) and (F\ref{43-f11}); otherwise, either $G[\{x_{0},\allowbreak x_{1},\allowbreak x_{2},\allowbreak q_{01,2},\allowbreak b,\allowbreak p_{2,1}\}]\cong R_{18}$ ($b\sim p_{2,1}$) or $G[\{x_{0},\allowbreak x_{1},\allowbreak x_{2},\allowbreak q_{01,2},\allowbreak q_{01,3},\allowbreak b,\allowbreak p_{2,1}\}]\cong R_{28}$ ($b\nsim p_{2,1}$), a contradiction.
Then $\left|E(T_{x_{1}x_{2}},b)\right|=0$ by (F\ref{43-f8}); otherwise, $G[\{x_{0},\allowbreak x_{1},\allowbreak x_{2},\allowbreak q_{01,2},\allowbreak q_{01,3},\allowbreak b,\allowbreak q_{12,1}\}]\cong R_{33}$, a contradiction.
Then $\left|E(T_{x_{1}x_{2}},T_{x_{1}})\right|=0$; otherwise, $G[\{x_{0},\allowbreak x_{1},\allowbreak q_{01,2},\allowbreak q_{01,3},\allowbreak q_{12,1},\allowbreak p_{1,1}\}]\cong R_{16}$, a contradiction.
Then $\left|T_{x_{1}x_{2}}\cup T_{x_{0}}\right|\leq 1$ by Lemma \ref{lem-13}\ref{lem-13:item-2} and (F\ref{43-f13}); otherwise, $G[\mathcal{K}\cup\{q_{01,2},\allowbreak q_{12,1},\allowbreak p_{0,1}\}]\cong R_{28}$, a contradiction.
If $\left|E(T_{x_{1}},b)\right|= 1$, then $\left|T_{x_{1}x_{2}}\right|=0$ for $G$ contains no subdivision of $K_{3,3}$.
If  $\left|E(T_{x_{1}},a)\right|= 1$, then  $\left|T_{x_{0}x_{2}}\right|=1$; otherwise, $G[\{x_{0},\allowbreak x_{1},\allowbreak x_{2},\allowbreak a,\allowbreak q_{02,2},\allowbreak p_{1,1}\}]\cong R_{18}$, a contradiction.
Thus,  
$G\subset H_{106}$ (if $\left|E(T_{x_{1}},a)\right|= 1$ and $\left|T_{x_{1}x_{2}}\right|=1$)  or 
$G\subset H_{11}$ (if $\left|E(T_{x_{1}},a)\right|= 0$ and $\left|T_{x_{1}x_{2}}\right|=1$)  or 
$G\subset H_{111}$ (if $\left|E(T_{x_{1}},a)\right|= 1$ and $\left|T_{x_{1}x_{2}}\right|=0$)  or 
$G\subset H_{115}$ (if $\left|E(T_{x_{1}},a)\right|= 0$ and $\left|T_{x_{1}x_{2}}\right|=0$) 
by Lemma \ref{lem-13}\ref{lem-13:item-2}, (F\ref{43-f7}) and (F\ref{43-f13}).
By symmetry, let $\left|T_{x_{0}x_{1}}\right|\leq 2$ and $\left|T_{x_{0}x_{2}}\right|\leq 2$.

Suppose that $\left|T_{x_{0}x_{1}}\right|=2$.
Then $\left|E(T_{x_{1}},b)\right|\leq 2$ by (F\ref{43-f10}) and (F\ref{43-f11}); otherwise, $G[\mathcal{K}\cup\{q_{01,2},\allowbreak b,\allowbreak p_{1,1},\allowbreak p_{1,2},\allowbreak p_{1,3}\}]\cong R_{78}$, a contradiction.
Then $\left|E(T_{x_{2}},b)\right|=0$ by (F\ref{43-f10}) and (F\ref{43-f11}); otherwise, $G[\{x_{0},\allowbreak x_{1},\allowbreak x_{2},\allowbreak q_{01,2},\allowbreak b,\allowbreak p_{2,1}\}]\cong R_{18}$, a contradiction.
Thus, $\left|T_{x_{2}}\right|\leq 1$ by (F\ref{43-f10}) and (F\ref{43-f11}); otherwise, $G[\{x_{0},\allowbreak x_{1},\allowbreak x_{2},\allowbreak q_{01,2},\allowbreak b,\allowbreak p_{2,1},\allowbreak p_{2,2}\}]\cong R_{25}$, a contradiction.
By symmetry, if $\left|T_{x_{0}x_{2}}\right|=2$, then $\left|E(T_{x_{1}},a)\right|=0$ and $\left|T_{x_{1}}\right|\leq 1$.
If $\left|E(T_{x_{1}x_{2}},a)\right|=1$, then $\left|E(T_{x_{1}x_{2}},T_{x_{1}})\right|=0$; otherwise, $G[\{x_{0},\allowbreak x_{1},\allowbreak x_{2},\allowbreak b,\allowbreak q_{12,1},\allowbreak p_{1,1}\}]\cong R_{18}$, a contradiction.
Then $\left|E(T_{x_{2}},T_{x_{1}x_{2}})\right|=0$; otherwise, $G[\{x_{0},\allowbreak x_{1},\allowbreak x_{2},\allowbreak q_{01,2},\allowbreak q_{12,1},\allowbreak p_{2,1}\}]\cong R_{18}$, a contradiction.
If $\left|E(T_{x_{1}},T_{x_{1}x_{2}})\right|\geq 1$, then $\left|E(T_{x_{1}},b)\right|=0$ for $G$ contains no subdivision of $K_{3,3}$.
If $\left|E(T_{x_{1}},T_{x_{1}x_{2}})\right|\geq 1$, then $\left|T_{x_{1}}\right|\leq 1$ by (F\ref{43-f10}) and (F\ref{43-f11}).
If $\left|E(T_{x_{1}},T_{x_{1}x_{2}})\right|\geq 1$, then $\left|T_{x_{0}x_{2}}\right|=1$ and  $\left|T_{x_{2}}\right|=0$ by (F\ref{43-f8}) and (F\ref{43-f11}); otherwise, either $G[\{x_{0},\allowbreak a,\allowbreak q_{01,2},\allowbreak q_{02,2},\allowbreak q_{12,1},\allowbreak p_{1,1}\}]\cong R_{1}$ ($E(\{a,\allowbreak q_{01,2},\allowbreak q_{02,2},\allowbreak q_{12,1},\allowbreak p_{1,1}\})=\{ap_{1,1},\allowbreak q_{12,1}p_{1,1}\}$) or $G$ contains a subdivision of $K_{3,3}$, a contradiction.
Then $\left|T_{x_{1}x_{2}}\cup T_{x_{0}}\right|\leq 1$; otherwise, $G[\mathcal{K}\cup\{a,\allowbreak q_{12,1},\allowbreak p_{0,1}\}]\cong R_{28}$, a contradiction.
Thus, if $\left|E(T_{x_{1}},T_{x_{1}x_{2}})\right|\geq 1$, then $G\subset H_{102}$ by (F\ref{43-f7}).
Let $\left|E(T_{x_{1}},T_{x_{1}x_{2}})\right|=0$.
If $\left|T_{x_{0}x_{2}}\right|=2$, then $\left|T_{x_{1}} \cup T_{x_{2}}\right|\leq 1$ for $G$ contains no subdivision of $K_{3,3}$.
If $\left|E(T_{x_{1}},b)\right|=1$ (or $\left|E(T_{x_{2}},a)\right|=1$), then   $\left|E(T_{x_{1}x_{2}},a)\right|=0$ (or $\left|E(T_{x_{1}x_{2}},b)\right|=0$) for $G$ contains no subdivision of $K_{3,3}$.
Thus, if $\left|T_{x_{0}x_{2}}\right|=2$, then either $G\subset H_{108}$ (if $\left|T_{x_{1}x_{2}}\right|=1$) or $G\subset H_{115}$ (if $\left|T_{x_{1}x_{2}}\right|=0$) by Lemma \ref{lem-13}\ref{lem-13:item-2}, (F\ref{43-f7}), (F\ref{43-f8}) and (F\ref{43-f13}).
Let $\left|T_{x_{0}x_{2}}\right|=1$.
If $\left|E(T_{x_{1}x_{2}},b)\right|=1$, then $\left|T_{x_{2}}\right|=0$ and $\left|E(T_{x_{1}},a)\right|=0$ by (F\ref{43-f10}); otherwise, either $G[\{x_{0},\allowbreak x_{2},\allowbreak x_{3},\allowbreak a,\allowbreak q_{12,1},\allowbreak p_{1,1}\}]\cong R_{1}$ ($E(\{a,\allowbreak q_{12,1},\allowbreak p_{1,1}\})=\{ap_{1,1}\}$) or  $G$ contains a subdivision of $K_{3,3}$, a contradiction.
If $\left|E(T_{x_{1}x_{2}},a)\right|=1$, then $\left|E(T_{x_{1}},b)\right|=0$ by (F\ref{43-f10}) and $G$ contains no subdivision of $K_{3,3}$.
Thus, if $\left|T_{x_{1}x_{2}}\right|=1$, then either $G\subset H_{104}$ (if $\left|E(T_{x_{1}x_{2}},a)\right|=1$) or $G\subset H_{107}$ (if $\left|E(T_{x_{1}x_{2}},a)\right|=0$) by  (F\ref{43-f7}), (F\ref{43-f8}) and (F\ref{43-f13}).
Let $\left|T_{x_{1}x_{2}}\right|=0$.
If $\left|T_{x_{2}}\right|=1$, then $\left|E(T_{x_{1}},b)\right|=0$ for $G$ contains no subdivision of $K_{3,3}$.
Thus, either $G\subset H_{105}$ (if $\left|T_{x_{2}}\right|=1$) or $G\subset H_{110}$ (if $\left|T_{x_{2}}\right|=0$) by Lemma \ref{lem-13}\ref{lem-13:item-2}, (F\ref{43-f7}), (F\ref{43-f8}) and (F\ref{43-f13}).
By symmetry, let  $\left|T_{x_{0}x_{1}}\right|=\left|T_{x_{0}x_{2}}\right|=1$.

Suppose that $\left|T_{x_{1}x_{2}}\right|= 1$.
Without loss of generality, let $\left|E(T_{x_{1}x_{2}},b)\right|=1$ by (F\ref{43-f8}).
Then $\left|T_{x_{0}}\right|=0$; otherwise, $G[\mathcal{K}\cup\{a,\allowbreak q_{12,1},\allowbreak p_{0,1}\}]\cong R_{28}$, a contradiction.
Then $\left|E(T_{x_{2}},a)\right|=0$ for $G$ contains no subdivision of $K_{3,3}$.
Thus,  $\left|E(T_{x_{2}},T_{x_{1}x_{2}})\right|=0$ by (F\ref{43-f10}); otherwise, $G[\{x_{0},\allowbreak x_{1},\allowbreak x_{2},\allowbreak a,\allowbreak q_{12,1},\allowbreak p_{2,1}\}]\cong R_{18}$, a contradiction.
Then $\left|T_{x_{2}}\right|\leq 1$ by (F\ref{43-f11}).
If $q\in T_{x_{1}}$ and $\left|E(q,a)\right|=1$, then $\left|E(q,T_{x_{1}x_{2}})\right|=1$; otherwise, $G[\{x_{0},\allowbreak x_{2},\allowbreak x_{3},\allowbreak a,\allowbreak q_{12,1},\allowbreak q\}]\cong R_{1}$, a contradiction.
If $\left|E(T_{x_{1}},a)\right|=1$, then $\left|T_{x_{2}}\right|=0$; otherwise, $G[\{x_{0},\allowbreak x_{1},\allowbreak x_{2},\allowbreak a,\allowbreak p_{1,1},\allowbreak p_{2,1}\}]\cong R_{14}$, a contradiction.
If $\left|E(T_{x_{1}},a)\right|=1$, then $\left|T_{x_{1}}\right|=1$ for $G$ contains no subdivision of $K_{3,3}$.
Thus, if  $\left|E(T_{x_{1}},a)\right|=1$, then $G\subset H_{102}$ by (F\ref{43-f7}) and (F\ref{43-f8}).
Let $\left|E(T_{x_{1}},a)\right|=0$.
Then at most two vertices in $T_{x_{1}}$ are not adjacent to $T_{x_{1}x_{2}}$ by (F\ref{43-f10}); otherwise, $G[\{x_{1},\allowbreak x_{2},\allowbreak x_{3},\allowbreak a,\allowbreak b,\allowbreak q_{12,1},\allowbreak p_{1,1},\allowbreak p_{1,2},\allowbreak p_{1,3}\}]\cong R_{79}$, a contradiction.
Moreover, at most one vertex in $T_{x_{1}}$ is adjacent to $T_{x_{1}x_{2}}$ by (F\ref{43-f10}) and  $G$ contains no subdivision of $K_{3,3}$.
Thus, $G\subset H_{112}$ by (F\ref{43-f7}) and (F\ref{43-f8}).
Let $\left|T_{x_{1}x_{2}}\right|= 0$.

If  $\left|E(T_{x_{1}},a)\right|=1$, then $\left|E(T_{x_{2}},b)\right|=0$ by (F\ref{43-f10}); otherwise, $G[\{x_{0},\allowbreak x_{1},\allowbreak x_{2},\allowbreak b,\allowbreak p_{1,1},\allowbreak p_{2,1}\}]\cong R_{14}$, a contradiction.
Thus, if $\left|E(T_{x_{1}},a)\right|=1$, then either $G\subset H_{109}$ or $G\subset H_{113}$ by (F\ref{43-f7}), (F\ref{43-f10}) and (F\ref{43-f11}).
Let $\left|E(T_{x_{1}},a)\right|=0$.
If  $\left|E(T_{x_{1}},b)\right|\geq 1$, then $\left|E(T_{x_{2}},a)\right|=0$  for $G$ contains no subdivision of $K_{3,3}$.
Thus, $G\subset H_{113}$  by (F\ref{43-f7}), (F\ref{43-f10}) and (F\ref{43-f11}).
Let $\left|E(T_{x_{1}},a)\right|=0$.

In the following cases, let $\left|E(T_{x_{0}x_{1}},T_{x_{0}x_{2}})\right|=0$ and $\left|E(T_{x_{1}x_{3}},T_{x_{2}x_{3}})\right|=0$.

\noindent$\textbf{Case 2.}$ Suppose that $\left|E(T_{x_{0}x_{1}},T_{x_{1}x_{3}})\right|\geq 1$.

Let  $a\in T_{x_{0}x_{1}}$, $e\in T_{x_{1}x_{3}}$ and $a\sim e$.
Similarly, we first establish the following facts.

\begin{fact}\label{43-f14}
$\left|E(T_{x_{0}},a)\right|=\left|E(T_{x_{3}},e)\right|=0$, $\left|E(T_{x_{0}},e)\right|=\left|T_{x_{0}}\right|$, $\left|E(T_{x_{3}},a)\right|=\left|T_{x_{3}}\right|$ and $\left|T_{x_{0}}\cup T_{x_{3}}\right|\leq 1$.
\end{fact}
\begin{proof}
By Lemma \ref{lem-13}\ref{lem-13:item-2}, $\left|E(T_{x_{0}},a)\right|=\left|E(T_{x_{3}},e)\right|=0$.
If $\left|E(T_{x_{0}},e)\right|\neq\left|T_{x_{0}}\right|$, then $G[\mathcal{K}\cup\{e,\allowbreak p_{0,1}\}]\cong R_{14}$, a contradiction.
If $\left|T_{x_{0}}\right|\geq 2$, then $R_{25}\subset G$, a contradiction.
If $\left|T_{x_{0}}\right|\geq 1$ and $\left|T_{x_{3}}\right|\geq 1$, then $G$ contains a subdivision of $K_{3,3}$, a contradiction.
\end{proof}

\begin{fact}\label{43-f15}
If $q\in T_{x_{1}}$, then  $\left|E(q,\{a,e\})\right|\leq 1$; $\left|E(T_{x_{1}},a)\right|\leq 1$ and $\left|E(T_{x_{1}},e)\right|\leq 1$; $\left|T_{x_{1}}\right|\leq 8$.
\end{fact}
\begin{proof}
If $\left|E(q,\{a,e\})\right|=2$, then $\omega=4$, a contradiction. 
If $\left|E(T_{x_{1}},a)\right|\geq 2$, then $G[\{x_{2},\allowbreak x_{3},\allowbreak a,\allowbreak e,\allowbreak p_{1,1},\allowbreak p_{1,2}\}]\cong R_{1}$, a contradiction.
If $\left|T_{x_{1}}\right|\geq 9$, then $G[\mathcal{K}\cup\{a,\allowbreak e,\allowbreak p_{1,1},\ldots,p_{1,7}\}]\cong R_{80}$, a contradiction. 
\end{proof}

\begin{fact}\label{43-f16}
If $y\in T_{x_{2}}$, then  $\left|E(y,\{a,e\})\right|\geq  1$; at most one vertex in $T_{x_{2}}$ is not adjacent to $a$ (or $e$); at most one vertex in $T_{x_{2}}$ is adjacent to all vertices of $\{a,e\}$; $\left|T_{x_{2}}\right|\leq 3$; if $\left|T_{x_{0}}\right|= 1$, then $\left|E(T_{x_{2}},a)\right|=0$.
\end{fact}
\begin{proof}
If $\left|E(y,\{a,e\})\right|=0$, then $G[\{x_{0},\allowbreak x_{1},\allowbreak x_{2},\allowbreak a,\allowbreak e,\allowbreak y\}]\cong R_{14}$, a contradiction.
If two vertices in $T_{x_{2}}$ are not adjacent to $a$, then $G[\mathcal{K}\cup\{a,\allowbreak p_{2,1},\allowbreak p_{2,2}\}]\cong R_{25}$, a contradiction.
If two vertices in $T_{x_{2}}$ are adjacent to both $a$ and $e$, then $G$ contains a subdivision of $K_{3,3}$, a contradiction.
Thus, $\left|T_{x_{2}}\right|\leq 3$.
If  $\left|E(T_{x_{2}},a)\right|\geq1$, then  $G$ contains a subdivision of $K_{3,3}$, a contradiction.
\end{proof}

\begin{fact}\label{43-f17}
$\left|E(T_{x_{0}x_{1}},T_{x_{1}x_{3}})\right|= 1$, $\left|T_{x_{0}x_{1}}\cup T_{x_{1}x_{3}}\right|\leq 3$  and $\left|E(T_{x_{0}x_{1}}\cup T_{x_{1}x_{3}}-\{a,e\},T_{1})\right|=0$.
\end{fact}
\begin{proof}
If $\left|E(T_{x_{0}x_{1}},T_{x_{1}x_{3}})\right| \geq 2$, then $G$ contains a subdivision of $K_{3,3}$, a contradiction.
If $\left|T_{x_{0}x_{1}}\cup T_{x_{1}x_{3}}\right|\geq 4$, then  $G[\{x_{0},\allowbreak x_{1},\allowbreak x_{3},\allowbreak q_{01,2},\allowbreak q_{01,3},\allowbreak e\}]\cong R_{16}$, a contradiction.
By Lemma \ref{lem-13}\ref{lem-13:item-2}, $\left|E(T_{x_{0}x_{1}},T_{x_{0}})\right|=0$.
If $\left|E(T_{x_{0}},T_{x_{1}x_{3}}-\{e\})\right| \geq 1$, then $G$ contains a subdivision of $K_{3,3}$, a contradiction.
If $\left|E(T_{x_{1}},T_{x_{0}x_{1}}-\{a\})\right| \geq 1$, then either $G$ contains a subdivision of $K_{3,3}$ (if $\left|E(T_{x_{1}},e)\right|=1$) or  $R_{16}\subset G$ (if $\left|E(T_{x_{1}},e)\right|=0$), a contradiction.
If $\left|E(T_{x_{0}x_{1}}-\{a\},T_{x_{2}})\right|\geq 1$, then  either $G$ contains a subdivision of $K_{3,3}$ (if $\left|E(T_{x_{2}}, \{a,e\})\right|\geq 1$) or $R_{14}\subset G$ by (F\ref{43-f16}), a contradiction.
\end{proof}

\begin{fact}\label{43-f18}
$\left|E(T_{x_{1}x_{2}},T_{x_{0}x_{1}}\cup T_{x_{1}x_{3}}\cup T_{1})\right|=0$ and $\left|T_{x_{1}x_{2}}\right|\leq 1$.
\end{fact}
\begin{proof}
If $\left|E(T_{x_{1}x_{2}},\{a,e\})\right|\geq 1$, then $G$ contains a subdivision of $K_{3,3}$, a contradiction.
If $\left|E(T_{x_{1}x_{2}},T_{x_{0}x_{1}}\cup T_{x_{1}x_{3}}-\{a,e\})\right|\geq 1$, then $G[\{x_{0},\allowbreak x_{1},\allowbreak x_{3},\allowbreak q_{01,2},\allowbreak e,\allowbreak q_{12,1}\}]\cong R_{16}$, a contradiction.
If $\left|E(T_{x_{1}x_{2}},T_{x_{0}}\cup T_{x_{3}})\right|\geq 1$, then  $G$ contains a subdivision of $K_{3,3}$ by (F\ref{43-f14}), a contradiction.
If $\left|E(T_{x_{1}x_{2}},T_{x_{1}})\right|\geq 1$, then  either $G$ contains a subdivision of $K_{3,3}$ (if $\left|E(T_{x_{1}}, \{a,e\})\right|\geq 1$) or $G[\{x_{0},\allowbreak x_{1},\allowbreak a,\allowbreak e,\allowbreak p_{1,1},\allowbreak q_{12,1}\}]\cong R_{16}$ ($E(\{a,\allowbreak e,\allowbreak p_{1,1},\allowbreak q_{12,1}\})=\{ae,\allowbreak p_{1,1}q_{12,1}\}$), a contradiction.
If $\left|E(T_{x_{1}x_{2}},T_{x_{2}})\right|\geq 1$, then  either $G$ contains a subdivision of $K_{3,3}$ (if $\left|E(T_{x_{2}}, \{a,e\})\right|\geq 1$) or $G[\{x_{0},\allowbreak x_{1},\allowbreak x_{2},\allowbreak a,\allowbreak e,\allowbreak p_{2,1}\}]\cong R_{14}$ ($E(\{a,\allowbreak e,\allowbreak p_{2,1}\})=\{ae\}$) by (F\ref{43-f16}), a contradiction.
If $\left|T_{x_{1}x_{2}}\right|\geq 2$, then  $G[\{x_{2},\allowbreak x_{3},\allowbreak a,\allowbreak e,\allowbreak q_{12,1},\allowbreak q_{12,2}\}]\cong R_{1}$, a contradiction.
\end{proof}

\begin{fact}\label{43-f19}
$\left|E(T_{0},T_{x_{1}x_{2}}\cup T_{x_{0}x_{3}})\right|=0$; if $\left|E(T_{0},T_{x_{0}x_{2}x_{3}})\right|=1$, then $\left|E(T_{x_{0}x_{2}x_{3}},\{a,e\})\right|=2$.
\end{fact}
\begin{proof}
If $\left|E(T_{0},T_{x_{1}x_{2}})\right|\geq 1$, then $G[\{x_{0},\allowbreak x_{1},\allowbreak x_{2},\allowbreak a,\allowbreak q_{12,1},\allowbreak t_{1}\}]\cong R_{14}$ by (F\ref{43-f18}) and Lemma \ref{lem-13}\ref{lem-13:item-8}, a contradiction.
If $\left|E(T_{0},T_{x_{0}x_{3}})\right|\geq 1$, then $G[\{x_{0},\allowbreak x_{1},\allowbreak x_{3},\allowbreak a,\allowbreak q_{03,1},\allowbreak t_{1}\}]\cong R_{7}$ by Lemma \ref{lem-13}\ref{lem-13:item-8}, a contradiction.
If $\left|E(T_{x_{0}x_{2}x_{3}},\{a,e\})\right|\leq 1$, then $G[\{x_{0},\allowbreak x_{1},\allowbreak x_{3},\allowbreak a,\allowbreak w_{023,1},\allowbreak t_{1}\}]\cong R_{7}$ ($E(\{a,\allowbreak w_{023,1},\allowbreak t_{1}\})=\{w_{023,1}t_{1}\}$) by Lemma \ref{lem-13}\ref{lem-13:item-8}, a contradiction.
\end{proof}

Suppose that $\left|T_{3}\right|\geq 1$.
Since $G$ contains no subdivision of $K_{3,3}$, $\left|T_{x_{0}x_{1}x_{3}}\right|=0$.
By Lemma \ref{lem-13}\ref{lem-13:item-7},  $\left|T_{x_{0}x_{2}x_{3}}\right|=1$.
Then $\left|E(T_{x_{0}x_{2}x_{3}},\{a,e\})\right|\geq 1$; otherwise, $G[\{x_{1},\allowbreak x_{2},\allowbreak x_{3},\allowbreak a,\allowbreak e,\allowbreak w_{023,1}\}]\cong R_{18}$, a contradiction.
By (F\ref{43-f16}), Lemma \ref{lem-13}\ref{lem-13:item-6}, and  $G$ contains no subdivision of $K_{3,3}$, $\left|E(T_{x_{0}x_{2}x_{3}},T_{x_{1}}\cup T_{x_{1}x_{2}}\cup T_{x_{0}x_{1}}\cup T_{x_{1}x_{3}}-\{a,e\})\right|=0$ and $\left|T_{x_{2}}\cup T_{x_{0}x_{2}}\cup T_{x_{2}x_{3}}\cup T_{x_{0}x_{3}}\right|=0$.
If $\left|E(T_{x_{0}x_{2}x_{3}},\{a,e\})\right|=2$, then  $\left|T_{x_{0}}\cup T_{x_{3}}\right|=0$ by (F\ref{43-f14}) and  $G$ contains no subdivision of $K_{3,3}$.
If $\left|E(T_{x_{0}x_{2}x_{3}},\{a,e\})\right|=2$, then $\left|T_{x_{1}x_{2}}\right|=0$; otherwise, $G[\{x_{1},\allowbreak x_{2},\allowbreak a,\allowbreak e,\allowbreak w_{023,1},\allowbreak q_{12,1}\}]\cong R_{17}$, a contradiction.
If $\left|E(T_{x_{0}x_{2}x_{3}},\{a,e\})\right|=2$, then $\left|T_{x_{0}x_{1}}\cup T_{x_{1}x_{3}}\right|=2$ by (F\ref{43-f17}); otherwise, $G[\{x_{0},\allowbreak x_{1},\allowbreak x_{3},\allowbreak q_{01,2},\allowbreak e,\allowbreak w_{023,1}\}]\cong R_{17}$, a contradiction.
If $\left|E(T_{x_{0}x_{2}x_{3}},\{a,e\})\right|=2$, then $\left|E(T_{x_{1}},\{a,e\})\right|=0$; otherwise, $G[\{x_{1},\allowbreak x_{2},\allowbreak x_{3},\allowbreak a,\allowbreak w_{023,1},\allowbreak p_{1,1}\}]\cong R_{17}$, a contradiction.
Thus, if $\left|E(T_{x_{0}x_{2}x_{3}},\{a,e\})\right|=2$, then $\left|T_{x_{1}}\right|\leq 1$; otherwise, $G[\{x_{0},\allowbreak x_{1},\allowbreak x_{2},\allowbreak a,\allowbreak w_{023,1},\allowbreak p_{1,1},\allowbreak p_{1,2}\}]\cong R_{27}$, a contradiction.
By (F\ref{43-f19}), Lemma \ref{lem-13}\ref{lem-13:item-8} and $d(G)\leq 3$,  $\left|T_{0}\cup T_{x_{1}}\right|\leq 1$.
Thus, if $\left|E(T_{x_{0}x_{2}x_{3}},\{a,e\})\right|=2$, then $G\subset A_{17}$.
Let $\left|E(T_{x_{0}x_{2}x_{3}},\{a,e\})\right|=1$.
Without loss of generality, let $\left|E(T_{x_{0}x_{2}x_{3}},a)\right|=1$.
Then  $\left|T_{0}\right|=0$ by (F\ref{43-f19}).
Then  $\left|T_{x_{1}}\right|=0$; otherwise, $G[\mathcal{K}\cup\{e,\allowbreak w_{023,1},\allowbreak p_{1,1}\}]\cong R_{30}$ ($E(\{e,\allowbreak w_{023,1},\allowbreak p_{1,1}\})=\varnothing$) or $G[\{x_{2},\allowbreak x_{3},\allowbreak a,\allowbreak e,\allowbreak w_{023,1},\allowbreak p_{1,1}\}]\cong R_{7}$ ($E(\{a,\allowbreak e,\allowbreak w_{023,1},\allowbreak p_{1,1}\})=\{ae,\allowbreak aw_{023,1},\allowbreak ep_{1,1}\}$), a contradiction.
Then $\left|T_{x_{0}}\right|=0$ by (F\ref{43-f14}) and  $G$ contains no subdivision of $K_{3,3}$.
Then  $\left|T_{x_{1}x_{3}}-\{e\}\right|=0$ by (F\ref{43-f17}) and $R_{17}\not\subset G$.
Then $\left|E(T_{x_{0}x_{2}x_{3}},T_{x_{3}})\right|=\left|T_{x_{3}}\right|$ by (F\ref{43-f14}); otherwise, $G[\{x_{0},\allowbreak x_{1},\allowbreak x_{2},\allowbreak a,\allowbreak e,\allowbreak w_{023,1},\allowbreak p_{3,1}\}]\cong R_{30}$, a contradiction.
Then $\left|T_{x_{3}}\cup T_{x_{0}x_{1}}-\{a\}\right|\leq 1$ and $\left|T_{x_{3}}\cup T_{x_{1}x_{2}}\right|\leq 1$ by (F\ref{43-f17}), (F\ref{43-f18}), $R_{14}\not\subset G$ and  $G$ contains no subdivision of $K_{3,3}$.
Thus, either $G\subset H_{102}$ (if $\left|T_{x_{3}}\right|=0$) or $G\subset H_{103}$ (if $\left|T_{x_{3}}\right|=1$) by (F\ref{43-f17}) and (F\ref{43-f18}).
Let $\left|T_{3}\right|= 0$. 
Then $\left|T_{0}\right|=0$ by (F\ref{43-f19}) and Lemma \ref{lem-13}\ref{lem-13:item-8}.

Suppose that $\left|T_{x_{0}x_{2}}\right|\geq 1$.
Then $\left|T_{x_{0}x_{2}}\right|=\left|T_{x_{1}x_{3}}\right|= \left|E(T_{x_{0}x_{2}},e)\right|=1$ and $\left|T_{x_{2}x_{3}}\right|=0$ by Lemma \ref{lem-13}\ref{lem-13:item-6} and  $G$ contains no subdivision of $K_{3,3}$.
Then $\left|T_{x_{0}x_{3}}\cup T_{x_{3}} \right|=\left|E(T_{x_{0}x_{2}},T_{x_{1}x_{2}}\cup T_{x_{1}})\right|=\left|E(T_{x_{2}},a)\right|=0$ by (F\ref{43-f14}) and  $G$ contains no subdivision of $K_{3,3}$.
Thus, $\left|E(T_{x_{1}},e)\right|=0$. Moreover, if $\left|E(T_{x_{1}},a)\right|\geq1$, then $\left|E(T_{x_{0}x_{2}},\{a,e\})\right|=2$; otherwise, $G[\{x_{0},\allowbreak x_{2},\allowbreak a,\allowbreak e,\allowbreak q_{02,1},\allowbreak p_{1,1}\}]\cong R_{7}$, a contradiction.
Then $\left|E(T_{x_{0}x_{1}}-\{a\},T_{x_{1}x_{2}})\right|=0$; otherwise, $G[\{x_{0},\allowbreak x_{1},\allowbreak x_{3},\allowbreak q_{01,2},\allowbreak e,\allowbreak q_{12,1}\}]\cong R_{16}$, a contradiction.
If $\left|E(T_{x_{0}x_{2}},a)\right|=1$, then $\left|T_{x_{0}}\right|=0$ by (F\ref{43-f14}) and  $G$ contains no subdivision of $K_{3,3}$.
If $\left|E(T_{x_{0}x_{2}},a)\right|=1$, then $\left|T_{x_{1}x_{2}}\right|=0$ by (F\ref{43-f18}); otherwise, either $G[\{x_{0},\allowbreak x_{2},\allowbreak a,\allowbreak e,\allowbreak q_{02,1},\allowbreak q_{12,1}\}]\cong R_{14}$ ($E(\{a,\allowbreak e,\allowbreak q_{02,1},\allowbreak q_{12,1}\})=\{ae,\allowbreak aq_{02,1},\allowbreak eq_{02,1}\}$) or  $G$ contains a subdivision of $K_{3,3}$, a contradiction.
If $\left|E(T_{x_{0}x_{2}},a)\right|=1$, then $\left|T_{x_{2}}\right|=0$ by (F\ref{43-f16}) and  $G$ contains no subdivision of $K_{3,3}$.
If $\left|E(T_{x_{0}x_{2}},a)\right|=1$, then $\left|E(T_{x_{1}},a)\right|=\left|T_{x_{1}}\right|$ and $\left|T_{x_{1}}\right|\leq 1$; otherwise, either $G[\{x_{2},\allowbreak x_{3},\allowbreak a,\allowbreak e,\allowbreak p_{1,1},\allowbreak p_{1,2}\}]\cong R_{1}$ ($E(\{a,\allowbreak e,\allowbreak p_{1,1},\allowbreak p_{1,2}\})=\{ae,\allowbreak ap_{1,1},\allowbreak ap_{1,2}\}$) or $G[\mathcal{K}\cup\{a,\allowbreak q_{02,1},\allowbreak p_{1,1}\}]\cong R_{30}$ ($E(\{a,\allowbreak q_{02,1},\allowbreak p_{1,1}\})=\{aq_{02,1}\}$), a contradiction.
Thus, if $\left|E(T_{x_{0}x_{2}},a)\right|=1$, then $G\subset H_{102}$.
Let $\left|E(T_{x_{0}x_{2}},a)\right|=0$.
Thus, $\left|E(T_{x_{1}},\{a,e\})\right|=0$.
Then $\left|T_{x_{1}}\right|\leq 1$; otherwise, $G[\mathcal{K}\cup\{q_{02,1},\allowbreak p_{1,1},\allowbreak p_{1,2}\}]\cong R_{25}$, a contradiction.
Then $\left|T_{x_{0}}\cup T_{x_{1}x_{2}}\right|\leq 1$ and $\left|T_{x_{2}}\cup T_{x_{0}x_{1}}-\{a\}\right|\leq 1$ by (F\ref{43-f17}), (F\ref{43-f18}) and $R_{12}\not\subset G$.
Thus,   $G\subset H_{119}$ (if $\left|T_{x_{0}x_{1}}\right|= 2$ and $\left|T_{x_{1}x_{2}}\right|=1$) or   $G\subset H_{72}$ (if $\left|T_{x_{0}x_{1}}\right|= 2$ and $\left|T_{x_{1}x_{2}}\right|=0$) or  $G\subset H_{72}$ (if $\left|T_{x_{0}x_{1}}\right|= 1$ and $\left|T_{x_{1}x_{2}}\right|=1$) or  $G\subset H_{118}$ (if $\left|T_{x_{0}x_{1}}\right|= 1$ and $\left|T_{x_{1}x_{2}}\right|=0$).
By symmetry, let   $\left|T_{x_{0}x_{2}}\cup T_{x_{2}x_{3}}\right|= 0$.

Suppose that $\left|T_{x_{0}x_{3}}\right|\geq 1$.
Since $R_{25}\not\subset G$,  at most one vertex in $T_{x_{0}x_{3}}$ is not adjacent to $\{a,e\}$.
Since $G$ contains no subdivision of $K_{3,3}$, $\left|E(T_{x_{0}x_{3}},a)\right|\leq 1$ and $\left|E(T_{x_{0}x_{3}},e)\right|\leq 1$.
Thus,   $\left|T_{x_{0}x_{3}}\right|\leq 2$ for  $G$ contains no subdivision of $K_{3,3}$.
Then $\left|T_{x_{2}}\right|=\left|E(T_{x_{0}x_{3}},T_{x_{1}})\right|=0$ by (F\ref{43-f16}),  and  $G$ contains no subdivision of $K_{3,3}$.
If $\left|E(T_{x_{0}x_{3}},\{a,e\})\right|=2$, then $\left|T_{x_{0}} \cup T_{x_{3}}\right|=0$ by (F\ref{43-f14}) and   $G$ contains no subdivision of $K_{3,3}$.
If $c,q_{03,1}\in T_{x_{0}x_{3}}$ and  $c\sim a,e$, then $c\sim q_{03,1}$; otherwise $G[\{x_{0},x_{1},x_{2},a, e, c,q_{03,1}\}]\cong R_{30}$, a contradiction.
If $c\sim a,e$, then $\left|T_{x_{1}x_{2}}\right|=0$ by (F\ref{43-f18}) and Lemma \ref{lem-13}\ref{lem-13:item-4}; otherwise $G[\{x_{0},x_{1},x_{2},a,c,d\}]\cong R_{18}$, a contradiction.
If $c\sim a,e$, then $\left|T_{x_{0}x_{1}}\cup T_{x_{0}x_{3}}\right|=2$ by (F\ref{43-f17}) and Lemma \ref{lem-13}\ref{lem-13:item-4}; otherwise either $G[\{x_{0},x_{1},x_{3},e,c,q_{01,1}\}]\cong R_{17}$ (if $q_{01,1}\nsim c$) or  $G$ contains a subdivision of $K_{5}$ (if $q_{01,1}\sim c$), a contradiction.
If $c\sim a,e$, then every vertex of $T_{x_{1}}$ is  adjacent to $a$ or $e$ by (F\ref{43-f15}); otherwise $G[\{x_{0},x_{1},x_{2},e,c,a,p_{1,1}\}]\cong R_{30}$, a contradiction.
If $c\sim a,e$ and $\left|T_{x_{1}}\right|\geq 1$, then $\left|T_{x_{0}x_{3}}\right|=1$; otherwise either $G[\{x_{0},x_{2},x_{3},a,p_{1,1}\}]\cong R_{2}$ (if $p_{1,1}\nsim q_{03,1}$ and $p_{1,1}\sim a$) or  $G$ contains a subdivision of $K_{3,3}$, a contradiction.
Thus, if $c\sim a,e$, then either $G\subset H_{102}$ (if $\left|T_{x_{1}}\right|\geq 1$) or $G\subset H_{103}$ (if $\left|T_{x_{1}}\right|=0$) by Lemma \ref{lem-13}\ref{lem-13:item-8}, (F\ref{43-f15}) and (F\ref{43-f19}).
Let  $\left|E(c,\{a,e\})\right|\leq 1$ for any vertex of $T_{x_{0}x_{3}}$.
Then $\left|E(T_{x_{0}x_{3}},\{a,e\})\right|\leq 1$ for $G$ contains no subdivision of $K_{3,3}$.
If $\left|T_{x_{1}x_{2}}\right|\geq1$, then  $\left|E(T_{x_{0}x_{3}},\{a,e\})\right|=0$ by (F\ref{43-f18}); otherwise $G[\{x_{0},x_{2},x_{3},a,q_{03,1},d\}]\cong R_{7}$, a contradiction.
If $\left|T_{x_{0}}\right|\geq 1$, then  $\left|E(T_{x_{0}x_{3}},a)\right|=0$ for  $G$ contains no subdivision of $K_{3,3}$.
We have $\left|T_{x_{0}}\cup T_{x_{3}}\cup T_{x_{1}x_{2}}\right|\leq 1$ by (F\ref{43-f14}); otherwise $G[\{x_{0},x_{2},a,e,q_{12,1},t_1\}]\cong R_{2}$ ( $\left|T_{x_{0}}\right|=\left|T_{x_{1}x_{2}}\right|=1$), a contradiction.
If $\left|T_{x_{0}x_{1}}\right|=2$, then  $\left|E(T_{x_{0}x_{3}},e)\right|=0$ by (F\ref{43-f17}); otherwise $G[\{x_{0},x_{2},x_{3},e,q_{01,1},q_{03,1}\}]\cong R_{7}$, a contradiction.
We have $\left|T_{x_{3}}\cup T_{x_{0}x_{1}}\right|\leq 1$ by (F\ref{43-f14}) and (F\ref{43-f17}); otherwise $G[\{x_{0},x_{1},x_{3},a,q_{02,1},p_{0,1}\}]\cong R_{10}$, a contradiction.
If $p_{1,1}\in T_{x_{1}}$ and  $p_{1,1}\sim a$, then $\left|E(T_{x_{0}x_{3}},e)\right|=0$, $\left|E(T_{x_{0}x_{3}},a)\right|=1$  and $\left|T_{x_{0}x_{3}}\right|=1$ by (F\ref{43-f15}); otherwise $G[\{x_{0},x_{2},x_{3},a,q_{03,1},p_{1,1}\}]\cong R_{2}$, a contradiction.
If $p_{1,1},p_{1,2}\in T_{x_{1}}$ and  $p_{1,1},p_{1,2}\nsim a,e$, then  $\left|E(T_{x_{0}x_{3}},\{a,e\})\right|=0$; otherwise $G[\{x_{0},x_{1},a,e,q_{03,1},p_{1,1},p_{1,2}\}]\cong R_{25}$, a contradiction.
Then $\left|T_{x_{1}}\right|\leq 2$ for $R_{38}\not\subset G$.
Thus, if $\left|T_{x_{1}}\right|=2$, then  $\left|T_{x_{0}x_{3}}\right| \leq1$.
Thus, 
$G\subset H_{118}$ (if $\left|T_{x_{0}x_{3}}\right|=2$ and $\left|T_{x_{0}}\cup T_{x_{3}}\right|=1$) or 
$G\subset H_{72}$ (if $\left|T_{x_{0}x_{3}}\right|=2$ and $\left|T_{x_{0}}\cup T_{x_{3}}\right|=0$) or 
$G\subset H_{87}$ (if $\left|E(T_{x_{0}x_{3}},\{a,e\})\right|=0$ and $\left|T_{x_{1}x_{2}}\right|\geq1$) or 
$G\subset H_{85}$ (if $\left|E(T_{x_{0}x_{3}},\{a,e\})\right|=0$ and $\left|T_{x_{1}x_{2}}\right|=0$) or 
$G\subset H_{119}$ (if $\left|T_{x_{0}x_{3}}\right|=1$, $\left|E(T_{x_{0}x_{3}},a)\right|=1$ and $\left|T_{x_{0}x_{1}}\cup T_{x_{1}x_{3}}\right|=3$) or 
$G\subset H_{72}$ (if $\left|T_{x_{0}x_{3}}\right|=1$, $\left|E(T_{x_{0}x_{3}},a)\right|=1$ and $\left|T_{x_{0}x_{1}}\cup T_{x_{1}x_{3}}\right|=2$) 
by Lemma \ref{lem-13}\ref{lem-13:item-8}, (F\ref{43-f14}), (F\ref{43-f15}), (F\ref{43-f17}), (F\ref{43-f18}), (F\ref{43-f19}).

Suppose that $\left|T_{x_{1}x_{2}}\right|= 1$.
Then $\left|T_{x_{0}}\cup T_{x_{3}}\right|=0$ by (F\ref{43-f14}) and (F\ref{43-f18}); otherwise, $G[\{x_{0},\allowbreak x_{1},\allowbreak x_{2},\allowbreak e,\allowbreak q_{12,1},\allowbreak p_{0,1}\}]\cong R_{12}$, a contradiction.
Then $\left|E(T_{x_{1}},\{a,e\})\right|=0$  by (F\ref{43-f15}), (F\ref{43-f18}); otherwise, $G[\{x_{0},\allowbreak x_{2},\allowbreak x_{3},\allowbreak a,\allowbreak q_{12,1},\allowbreak p_{1,1}\}]\cong R_{1}$, a contradiction.
Then $\left|E(T_{x_{2}},\{a,e\})\right|=\left|T_{x_{2}}\right|$   by (F\ref{43-f16}) and (F\ref{43-f18}); otherwise, $G[\{x_{0},\allowbreak x_{2},\allowbreak a,\allowbreak e,\allowbreak q_{12,1},\allowbreak p_{2,1}\}]\cong R_{7}$, a contradiction.
If $\left|T_{x_{2}}\right|\geq 1$, then $\left|T_{x_{1}x_{3}}-\{e\}\right|=0$ by (F\ref{43-f17}), (F\ref{43-f18}) and  $R_{10}\not\subset G$.
If $\left|T_{x_{2}}\right|\geq 1$, then $\left|T_{x_{1}}\right|\leq 2$ by (F\ref{43-f18}); otherwise, $G[\{x_{0},\allowbreak x_{1},\allowbreak x_{2},\allowbreak a,\allowbreak p_{2,1},\allowbreak p_{1,1},\allowbreak p_{1,2},\allowbreak p_{1,3}\}]\cong R_{38}$, a contradiction.
Thus, either $G\subset H_{85}$ (if  $\left|T_{x_{2}}\right|\geq 1$) or $G\subset H_{120}$ (if  $\left|T_{x_{2}}\right|=0$) by (F\ref{43-f15}), (F\ref{43-f16}), (F\ref{43-f17}) and (F\ref{43-f18}).
Let $\left|T_{x_{1}x_{2}}\right|= 0$.

Suppose that $\left|T_{x_{0}x_{1}}-\{a\}\right|= 1$.
Then $\left|T_{x_{3}}\right|= 0$  by (F\ref{43-f14}) and (F\ref{43-f17}); otherwise, $G[\mathcal{K}\cup\{q_{01,2},\allowbreak p_{3,1}\}]\cong R_{14}$, a contradiction.
Then  $\left|E(T_{x_{2}},e)\right|=0$ by (F\ref{43-f16}) and (F\ref{43-f17}); otherwise, $G[\{x_{0},\allowbreak x_{1},\allowbreak x_{2},\allowbreak q_{01,2},\allowbreak e,\allowbreak p_{2,1}\}]\cong R_{12}$, a contradiction.
Thus, $\left|T_{x_{2}}\right|\leq  1$ by (F\ref{43-f16}).
Then $\left|E(T_{x_{1}},e)\right|=0$ by (F\ref{43-f15}) and (F\ref{43-f17}); otherwise, $G[\{x_{0},\allowbreak x_{1},\allowbreak x_{2},\allowbreak q_{01,2},\allowbreak e,\allowbreak p_{1,1}\}]\cong R_{16}$, a contradiction.
Then $\left|T_{x_{0}}\cup T_{x_{2}}\right|\leq 1$ by (F\ref{43-f14}), (F\ref{43-f16}) and   $G$ contains no subdivision of $K_{3,3}$.
If $\left|T_{x_{2}}\cup T_{x_{0}}\right|= 1$, then  at most two vertices in $T_{x_{1}}$ are not adjacent to $\{a,e\}$  by (F\ref{43-f14}), (F\ref{43-f16}) and (F\ref{43-f17}); otherwise, $G[\{x_{0},\allowbreak x_{1},\allowbreak a,\allowbreak e,\allowbreak p_{1,1},\allowbreak p_{1,2},\allowbreak p_{1,3},\allowbreak p_{0,1}\}]\cong R_{38}$, a contradiction.
Thus, either $G\subset H_{87}$ (if $\left|T_{x_{2}}\cup T_{x_{0}}\right|= 1$) or $G\subset H_{120}$ (if $\left|T_{x_{2}}\cup T_{x_{0}}\right|= 0$) by (F\ref{43-f14}), (F\ref{43-f15}), (F\ref{43-f16}) and (F\ref{43-f17}).
Let $\left|T_{x_{0}x_{1}}\cup T_{x_{1}x_{3}}-\{a,e\}\right|= 0$.

If $\left|T_{x_{0}}\right|\geq 1$, then  at most two vertices in $T_{x_{1}}$ are not adjacent to $\{a,e\}$ by (F\ref{43-f14}) and (F\ref{43-f16}); otherwise, $G[\{x_{0},\allowbreak x_{1},\allowbreak a,\allowbreak e,\allowbreak p_{0,1},\allowbreak p_{1,1},\allowbreak p_{1,2},\allowbreak p_{1,3}\}]\cong R_{38}$, a contradiction.
If $\left|T_{x_{2}}\right|= 1$ and $\left|E(T_{x_{2}},a)\right|=0$, then $\left|E(T_{x_{1}},a)\right|=0$ by (F\ref{43-f15}) and (F\ref{43-f16}); otherwise, $G[\{x_{0},\allowbreak x_{1},\allowbreak x_{2},\allowbreak a,\allowbreak p_{1,1},\allowbreak p_{2,1}\}]\cong R_{14}$, a contradiction.
Thus, if $\left|T_{x_{0}}\right|\geq 1$, then either $G\subset H_{85}$ (if $\left|T_{x_{2}}\right|\geq 1$) or $G\subset H_{87}$ (if $\left|T_{x_{2}}\right|= 0$) by (F\ref{43-f14}), (F\ref{43-f15}) and (F\ref{43-f16}).
Let $\left|T_{x_{0}}\cup T_{x_{3}}\right|= 0$.
If $\left|T_{x_{2}}\right|=0$, then $G\subset H_{120}$  by (F\ref{43-f15}).
Let $\left|T_{x_{2}}\right| \geq1$.
If $p_{2,1}\in T_{x_{2}}$ and $\left|E(p_{2,1},\{a,e\})\right|=2$, then at most one vertex in $T_{x_{1}}$ is not adjacent to $\{a,e\}$  by (F\ref{43-f15}) and  (F\ref{43-f16}); otherwise, $G[\{x_{0},\allowbreak x_{1},\allowbreak a,\allowbreak e,\allowbreak p_{2,1},\allowbreak p_{1,1},\allowbreak p_{1,2}\}]\cong R_{25}$, a contradiction.
If  $p_{2,2}\in T_{x_{2}}$ and $\left|E(p_{2,2},e)\right|=0$ (or $\left|E(p_{2,2},a)\right|=0$), then at most two vertices in $T_{x_{1}}$ are not adjacent to $\{a,e\}$; otherwise, $G[\{x_{0},\allowbreak x_{1},\allowbreak x_{2},\allowbreak a,\allowbreak p_{2,2},\allowbreak p_{1,1},\allowbreak p_{1,2},\allowbreak p_{1,3}\}]\cong R_{38}$, a contradiction.
If $p_{2,2}\in T_{x_{2}}$ and $\left|E(p_{2,2},e)\right|=0$ (or $\left|E(p_{2,2},a)\right|=0$), then $\left|E(T_{x_{1}},e)\right|=0$(or $\left|E(T_{x_{1}},a)\right|=0$); otherwise, $G[\{x_{0},\allowbreak x_{1},\allowbreak x_{2},\allowbreak e,\allowbreak p_{2,1},\allowbreak p_{1,1}\}]\cong R_{10}$, a contradiction.
Thus, $G\subset H_{119}$ (if $\left|T_{x_{2}}\right|=1$ and $\left|E(T_{x_{2}},\{a,e\})\right|=2$) or  $G\subset H_{85}$ (if $\left|T_{x_{2}}\right|=1$ and $\left|E(T_{x_{2}},\{a,e\})\right|=1$)  or  $G\subset H_{72}$ (if $\left|T_{x_{2}}\right|=2$ and $\left|E(T_{x_{2}},\{a,e\})\right|=3$)  or  $G\subset H_{85}$ (if $\left|T_{x_{2}}\right|=2$ and $\left|E(T_{x_{2}},\{a,e\})\right|=2$)  or   $G\subset H_{118}$ (if $\left|T_{x_{2}}\right|=3$)  by (F\ref{43-f15}) and (F\ref{43-f16}).

In the following cases, let $\left|E(T_{x_{0}x_{1}},T_{x_{1}x_{3}})\right|=\left|E(T_{x_{0}x_{2}},T_{x_{2}x_{3}})\right|=0$.

\noindent$\textbf{Case 3.}$ Suppose that $\left|E(T_{x_{0}x_{1}},T_{x_{2}x_{3}})\right|\geq 1$ (or $\left|E(T_{x_{0}x_{1}},T_{x_{0}x_{3}})\right|\geq 1$).

The conditions $\left|E(T_{x_{0}x_{1}},T_{x_{2}x_{3}})\right|\geq 1$ and $\left|E(T_{x_{0}x_{1}},T_{x_{0}x_{3}})\right|\geq 1$ correspond to the same configuration up to isomorphism, obtained by choosing different copies of $K_{4}-e$ in the same graph.
Since the two choices lead to the same conclusion, it suffices to consider $\left|E(T_{x_{0}x_{1}},T_{x_{2}x_{3}})\right|\geq 1$.
By Lemma \ref{lem-13}\ref{lem-13:item-6},  $\left|E(T_{x_{0}x_{1}},T_{x_{2}x_{3}})\right|=1$, $\left|T_{x_{0}x_{1}}\right|=1$, $\left|T_{x_{2}x_{3}}\right|= 1$ and   $\left|T_{3}\cup T_{x_{0}x_{3}}\right|=0$.
Let  $a\in T_{x_{0}x_{1}}$, $f\in T_{x_{2}x_{3}}$ and $a\sim f$.

\begin{fact}\label{43-f20}
At most one vertex in $T_{x_{1}}$ (or $T_{x_{2}}$) is not adjacent to $\{a,f\}$; at most two vertices in $T_{x_{1}}$ (or $T_{x_{2}}$) are adjacent to $f$ (or $a$) and are not adjacent to $a$ (or $f$); at most one vertex in $T_{x_{1}}$ (or $T_{x_{2}}$) is adjacent to $\{a,f\}$; if $q\in T_{x_{1}}$ and  $q\sim a$, then $q\sim f$; $\left|T_{x_{1}}\right|\leq 4$ and $\left|T_{x_{2}}\right|\leq 4$.
\end{fact}
\begin{proof}
If  at least two vertices in $T_{x_{1}}$ (or $T_{x_{2}}$) are not adjacent to $\{a,f\}$, then $G[\mathcal{K}\cup\{f,\allowbreak p_{1,1},\allowbreak p_{1,2}\}]\cong R_{25}$, a contradiction.
If  at least three vertices in $T_{x_{1}}$ (or $T_{x_{2}}$) are adjacent to $f$ (or $a$) and are not adjacent to $a$ (or $f$), then $G[\{x_{0},\allowbreak x_{1},\allowbreak x_{2},\allowbreak a,\allowbreak f,\allowbreak p_{1,1},\allowbreak p_{1,2},\allowbreak p_{1,3}\}]\cong R_{81}$, a contradiction.
If  at least two vertices in $T_{x_{1}}$ (or $T_{x_{2}}$) are adjacent to $\{a,f\}$, then $G$ contains a subdivision of $K_{3,3}$, a contradiction.
If $q\nsim f$, then $G[\{x_{0},\allowbreak x_{2},\allowbreak x_{3},\allowbreak a,\allowbreak f,\allowbreak q\}]\cong R_{7}$, a contradiction.
Thus, $\left|T_{x_{1}}\right|\leq 4$ and $\left|T_{x_{2}}\right|\leq 4$.
\end{proof}

\begin{fact}\label{43-f21}
$\left|E(T_{x_{0}},a)\right|=\left|E(T_{x_{3}},f)\right|=0$, $\left|E(T_{x_{0}},f)\right|=\left|T_{x_{0}}\right|$, $\left|E(T_{x_{3}},a)\right|=\left|T_{x_{3}}\right|$, $\left|T_{x_{0}}\right|\leq 1$ and  $\left|T_{x_{3}}\right|\leq 1$.
\end{fact}
\begin{proof}
By Lemma \ref{lem-13}\ref{lem-13:item-2}, $\left|E(T_{x_{0}},a)\right|=\left|E(T_{x_{3}},f)\right|=0$.
If  $\left|E(T_{x_{0}},f)\right|< \left|T_{x_{0}}\right|$, then $G[\mathcal{K}\cup\{f,\allowbreak p_{0,1}\}]\cong R_{14}$, a contradiction.
Thus, $\left|T_{x_{1}}\right|\leq 4$ and $\left|T_{x_{2}}\right|\leq 4$.
If $\left|T_{x_{0}}\right|\geq 2$, then $G[\mathcal{K}\cup\{a,\allowbreak p_{0,1},\allowbreak p_{0,2}\}]\cong R_{25}$, a contradiction.
\end{proof}

\begin{fact}\label{43-f22}
$\left|T_{0}\right|=0$; $\left|T_{x_{0}x_{2}}\cup T_{x_{1}x_{3}}\right|\leq 1$; $\left|T_{x_{1}x_{3}}\cup T_{x_{0}}\right|\leq 1$; $\left|T_{x_{0}x_{2}}\cup T_{x_{3}}\right|\leq 1$.
\end{fact}
\begin{proof}
If $\left|T_{0}\right|\geq 1$, then either $G[\{x_{0},\allowbreak x_{1},\allowbreak x_{2},\allowbreak a,\allowbreak q_{12,1},\allowbreak t_{1}\}]\cong R_{14}$ ($E(\{a,\allowbreak q_{12,1},\allowbreak t_{1}\})=\{q_{12,1}t_{1}\}$) or $G$ contains a subdivision of $K_{3,3}$ by Lemma \ref{lem-13}\ref{lem-13:item-8}, a contradiction.
If $\left|T_{x_{0}x_{2}}\cup T_{x_{1}x_{3}}\right|\geq 2$, then either $R_{1}\subset G$ (if $\left|T_{x_{0}x_{2}}\right|\geq 2$ or $=0$) or $R_{18}\subset G$ (if $\left|T_{x_{0}x_{2}}\right|=1$), a contradiction.
If $\left|T_{x_{1}x_{3}}\cup T_{x_{0}}\right|\geq 2$, then  $G[\mathcal{K}\cup\{q_{13,1},\allowbreak p_{0,1}\}]\cong R_{14}$ by (F\ref{43-f21}), a contradiction.
\end{proof}

Suppose that $\left|T_{x_{1}x_{2}}\right|\geq 1$.
Since $G$ contains no subdivision of $K_{3,3}$, $\left|E(T_{x_{1}x_{2}},\{a,f\})\right|=0$.
Thus,  $\left|T_{x_{1}x_{2}}\right|\leq 2$; otherwise, $G[\{x_{1},\allowbreak x_{2},\allowbreak x_{3},\allowbreak a,\allowbreak f,\allowbreak q_{12,1},\allowbreak q_{12,2},\allowbreak q_{12,3}\}]\cong R_{42}$, a contradiction.
Then $\left|E(T_{x_{1}x_{2}},T_{x_{0}x_{2}}\cup T_{x_{1}x_{3}})\right|=0$; otherwise, $G[\{x_{0},\allowbreak x_{1},\allowbreak x_{3},\allowbreak f,\allowbreak q_{12,1},\allowbreak q_{02,1}\}]\cong R_{2}$, a contradiction.
Then $\left|E(T_{x_{1}x_{2}},T_{x_{0}}\cup T_{x_{3}})\right|=0$ by (F\ref{43-f21}); otherwise, $G$ contains a subdivision of $K_{3,3}$, a contradiction.
Thus,  $\left|T_{x_{0}}\cup T_{x_{3}}\right|=0$ by (F\ref{43-f21}); otherwise, $G[\{x_{0},\allowbreak x_{1},\allowbreak x_{2},\allowbreak f,\allowbreak q_{12,1},\allowbreak p_{0,1}\}]\cong R_{12}$, a contradiction.
Then $\left|E(T_{x_{1}x_{2}},T_{x_{1}}\cup T_{x_{2}})\right|=0$  by (F\ref{43-f20}); otherwise, either $G$ contains a subdivision of $K_{3,3}$ (if $\left|E(\{a,f\},T_{x_{1}}\cup T_{x_{2}})\right|\geq 1$) or $R_{10}\subset G$ (if $\left|E(\{a,f\},T_{x_{1}}\cup T_{x_{2}})\right|=0$), a contradiction.
Then no vertex in $T_{x_{1}}\cup T_{x_{2}}$ is adjacent to all vertices of $\{a,f\}$; otherwise, $G[\{x_{1},\allowbreak x_{2},\allowbreak x_{3},\allowbreak a,\allowbreak q_{12,1},\allowbreak p_{1,1}\}]\cong R_{16}$, a contradiction.
If $\left|T_{x_{1}x_{2}}\right|=2$, then $\left|E(T_{x_{1}}\cup T_{x_{2}},\{a,f\})\right|=0$ by (F\ref{43-f20}); otherwise, $G[\{x_{1},\allowbreak x_{2},\allowbreak a,\allowbreak f,\allowbreak q_{12,1},\allowbreak q_{12,2},\allowbreak p_{1,1}\}]\cong R_{26}$, a contradiction.
Thus, if $\left|T_{x_{1}x_{2}}\right|=2$, then $G\subset H_{71}$ by (F\ref{43-f20}) and (F\ref{43-f22}).
If $\left|T_{x_{1}x_{2}}\right|=1$, then $\left|E(T_{x_{1}},f)\right|\leq 1$ and $\left|E(T_{x_{2}},a)\right|\leq 1$  by (F\ref{43-f20}); otherwise, either $G$ contains a subdivision of $K_{3,3}$ or $G[\{x_{0},\allowbreak x_{1},\allowbreak x_{2},\allowbreak f,\allowbreak q_{12,1},\allowbreak p_{1,1},\allowbreak p_{1,2}\}]\cong R_{26}$ ($E(\{f,\allowbreak q_{12,1},\allowbreak p_{1,1},\allowbreak p_{1,2}\})=\{fp_{1,1},\allowbreak fp_{1,2}\}$), a contradiction.
If $\left|T_{x_{1}x_{2}}\right|=1$ and $\left|T_{x_{0}x_{2}}\right|=1$, then $\left|E(T_{x_{1}},f)\right|=0$ by (F\ref{43-f20}); otherwise, either $G$ contains a subdivision of $K_{3,3}$ or $G[\{x_{0},\allowbreak x_{1},\allowbreak x_{2},\allowbreak f,\allowbreak q_{02,1},\allowbreak p_{1,1}\}]\cong R_{12}$ ($E(\{f,\allowbreak q_{02,1},\allowbreak p_{1,1}\})=\{fp_{1,1}\}$), a contradiction.
Thus, if $\left|T_{x_{1}x_{2}}\right|=1$, then either $G\subset H_{80}$ (if $\left|T_{x_{0}x_{2}}\cup T_{x_{1}x_{3}}\right|= 1$) or $G\subset H_{83}$ (if $\left|T_{x_{0}x_{2}}\cup T_{x_{1}x_{3}}\right|= 0$) by (F\ref{43-f20}) and (F\ref{43-f22}).
Let $\left|T_{x_{1}x_{2}}\right|=0$.

Suppose that $\left|T_{x_{0}x_{2}}\cup T_{x_{1}x_{3}}\right|\geq 1$.
By (F\ref{43-f22}), $\left|T_{x_{0}x_{2}}\cup T_{x_{1}x_{3}}\right|=1$.
Without loss of generality, let $\left|T_{x_{0}x_{2}}\right|=1$.
Then $\left|T_{x_{3}}\right|=0$ by (F\ref{43-f22}).
Then $\left|E(T_{x_{1}},T_{x_{0}x_{2}}\cup\{a,f\})\right|=0$ by (F\ref{43-f20}); otherwise, either $G$ contains a subdivision of $K_{3,3}$ or $G[\{x_{0},\allowbreak x_{1},\allowbreak x_{2},\allowbreak f,\allowbreak q_{02,1},\allowbreak p_{1,1}\}]\cong R_{12}$ ($E(\{f,\allowbreak q_{02,1},\allowbreak p_{1,1}\})=\{fp_{1,1}\}$) or $G[\{x_{1},\allowbreak x_{2},\allowbreak x_{3},\allowbreak f,\allowbreak q_{02,1},\allowbreak p_{1,1}\}]\cong R_{12}$ ($E(\{f,\allowbreak q_{02,1},\allowbreak p_{1,1}\})=\{q_{02,1}p_{1,1}\}$), a contradiction.
If $p_{2,1}\in T_{x_{2}}$ and $p_{2,1}\sim f$, then  $\left|E(p_{2,1},T_{x_{0}x_{2}})\right|=1$  by (F\ref{43-f20}); otherwise, $G[\{x_{0},\allowbreak x_{2},\allowbreak a,\allowbreak f,\allowbreak q_{02,1},\allowbreak p_{2,1}\}]\cong R_{17}$, a contradiction.
If $p_{2,2}\in T_{x_{2}}$ and $\left|E(p_{2,2},f)\right|=0$, then $\left|E(p_{2,2},T_{x_{0}x_{2}})\right|=0$ by (F\ref{43-f20}) and  $P_{6}\not\subset G$.
If  $\left|E(T_{x_{2}},\{a,f\})\right|\geq 1$, then $\left|T_{x_{0}}\right|=0$ by (F\ref{43-f20}) and (F\ref{43-f21}); otherwise, $G$ contains a subdivision of $K_{3,3}$, a contradiction.
If  $\left|E(T_{x_{2}},f)\right|\geq 1$, then  no vertex in $T_{x_{1}}$ is not adjacent to $\{a,f\}$; otherwise, $G[\{x_{0},\allowbreak x_{1},\allowbreak x_{2},\allowbreak f,\allowbreak p_{2,1},\allowbreak p_{1,1}\}]\cong R_{10}$, a contradiction.
Thus, $G\subset H_{121}$ (if  $\left|E(T_{x_{2}},f)\right|\geq 1$)  or $G\subset H_{89}$ (if $\left|E(T_{x_{2}},f)\right|=0$ and $\left|T_{x_{0}}\right|= 1$) or $G\subset H_{129}$ (if $\left|E(T_{x_{2}},f)\right|=0$ and if $\left|T_{x_{0}}\right|= 0$) by (F\ref{43-f20}), (F\ref{43-f21}) and (F\ref{43-f22}).
Let $\left|T_{x_{0}x_{2}}\cup T_{x_{1}x_{3}}\right|=0$.

If  $\left|E(T_{x_{2}},\{a,f\})\right|\geq 1$ (or $\left|E(T_{x_{1}},\{a,f\})\right|\geq 1$), then $\left|T_{x_{0}}\right|=0$ (or $\left|T_{x_{3}}\right|=0$) by (F\ref{43-f20}) and (F\ref{43-f21}); otherwise, $G$ contains a subdivision of $K_{3,3}$, a contradiction.
Thus, if $\left|T_{x_{0}}\right|=\left|T_{x_{3}}\right|=1$, then $G\subset H_{122}$ by (F\ref{43-f20}), (F\ref{43-f21}) and (F\ref{43-f22}).
Let $\left|T_{x_{0}}\cup T_{x_{3}}\right|\leq 1$.
If $\left|E(T_{x_{2}},f)\right|\geq 1$, then  no vertex in $T_{x_{1}}$ is not adjacent to $\{a,f\}$ by (F\ref{43-f20}).
Thus, if $\left|T_{x_{0}}\cup T_{x_{3}}\right|=1$, then either $G\subset H_{127}$ (if $\left|E(T_{x_{2}},f)\right|\geq 1$) or $G\subset H_{128}$ (if $\left|E(T_{x_{2}},f)\right|=0$) by (F\ref{43-f20}) and (F\ref{43-f21}).
Let $\left|T_{x_{0}}\cup T_{x_{3}}\right|=0$.
If  $\left|T_{x_{1}}\right|=0$ or $\left|T_{x_{2}}\right|=0$, then  $G\subset H_{127}$  by (F\ref{43-f20}).
Let $\left|T_{x_{1}}\right|\geq 1$ and $\left|T_{x_{2}}\right|\geq 1$.
If $S_{1}\subset T_{x_{1}}$, $S_{2}\subset T_{x_{2}}$, $\left|E(S_{1},a)\right|=\left|E(S_{2},f)\right|=0$, $\left|E(S_{1},f)\right|=\left|S_{1}\right|$ and $\left|E(S_{2},a)\right|=\left|S_{2}\right|$, then  $\left|S_{1}\cup S_{2}\right|\leq 2$; otherwise, $G[\{x_{0},\allowbreak x_{3},\allowbreak a,\allowbreak f,\allowbreak p_{1,1},\allowbreak p_{1,2},\allowbreak p_{2,1}\}]\cong R_{21}$, a contradiction.
Thus, $\left|T_{x_{1}}\cup T_{x_{2}}\right|\leq 4$ by (F\ref{43-f20}).
Note that if $\left|E(T_{x_{2}},f)\right|\geq 1$, then  no vertex in $T_{x_{1}}$ is not adjacent to $\{a,f\}$.
Then $G\subset H_{125}$ (if $\left|T_{x_{2}}\right|=3$ and $\left|E(T_{x_{2}},\{a,f\})\right|=4$) or  
$G\subset H_{123}$ (if  $\left|T_{x_{2}}\right|=3$, $\left|E(T_{x_{2}},a)\right|=1$ and $\left|E(T_{x_{2}},f)\right|=2$) or  
$G\subset H_{128}$ (if $\left|T_{x_{2}}\right|=3$, $\left|E(T_{x_{2}},a)\right|=0$ and $\left|E(T_{x_{2}},f)\right|=2$) or  
$G\subset H_{126}$ (if $\left|T_{x_{2}}\right|=2$, $\left|E(T_{x_{2}},a)\right|=1$ and $\left|E(T_{x_{2}},f)\right|=2$) or  
$G\subset H_{83}$ (if $\left|T_{x_{2}}\right|=2$, $\left|E(T_{x_{2}},a)\right|=0$ and $\left|E(T_{x_{2}},f)\right|=1$) or  
$G\subset H_{124}$ (if $\left|T_{x_{2}}\right|=2$, $\left|E(T_{x_{2}},a)\right|=1$ and $\left|E(T_{x_{2}},f)\right|=1$).

In the following cases, let $E(T_{2})\subseteq E(T_{x_{0}x_{3}})\cup E(T_{x_{1}x_{2}},T_{2}-T_{x_{1}x_{2}})$.

\noindent$\textbf{Case 4.}$ $\left|E(T_{x_{0}x_{1}},T_{x_{1}x_{2}})\right|\geq 1$.

Let  $a\in T_{x_{0}x_{1}}$, $d\in T_{x_{1}x_{2}}$ and $a\sim d$.

\begin{fact}\label{case7-f1}
$\left|T_{x_{3}}\cup T_{x_{2}x_{3}}\cup T_{x_{0}x_{3}}\cup T_{x_{0}x_{2}x_{3}}\right|=0$.
\end{fact}
\begin{proof}
If $\left|T_{x_{3}}\right|\geq 1$, then  either $G$ contains a subdivision of $K_{3,3}$ or $G[\mathcal{K}\cup\{a,\allowbreak p_{3,1}\}]\cong R_{14}$ ($E(\{a,\allowbreak p_{3,1}\})=\varnothing$), a contradiction.
By Lemma \ref{lem-13}\ref{lem-13:item-6} and $\left|E(T_{x_{0}x_{1}},T_{x_{2}x_{3}})\right|=0$, $\left|T_{x_{2}x_{3}}\right|=0$.
If $\left|T_{x_{0}x_{3}}\right|\geq 1$, then   $G[\{x_{0},\allowbreak x_{1},\allowbreak x_{3},\allowbreak a,\allowbreak d,\allowbreak q_{03,1}\}]\cong R_{12}$ by  Lemma \ref{lem-13}\ref{lem-13:item-4}, a contradiction.
If $\left|T_{x_{0}x_{2}x_{3}}\right|\geq 1$, then  either $G$ contains a subdivision of $K_{3,3}$ or $G[\{x_{0},\allowbreak x_{1},\allowbreak x_{3},\allowbreak a,\allowbreak d,\allowbreak w_{023,1}\}]\cong R_{12}$ ($E(\{a,\allowbreak d,\allowbreak w_{023,1}\})=\{ad\}$) by  Lemma \ref{lem-13}\ref{lem-13:item-4}, a contradiction.
\end{proof}

\begin{fact}\label{case7-f2}
$\left|T_{x_{0}}\right|\leq 1$ and $\left|E(T_{x_{0}},T_{x_{1}x_{2}})\right|\leq 1$.
\end{fact}
\begin{proof}
If $\left|T_{x_{0}}\right|\geq 2$, then $G[\mathcal{K}\cup\{a,\allowbreak p_{0,1},\allowbreak p_{0,2}\}]\cong R_{25}$ by Lemma \ref{lem-13}\ref{lem-13:item-2}, a contradiction.
If $\left|E(T_{x_{0}},T_{x_{1}x_{2}})\right|\geq2$, then $G$ contains a subdivision of $K_{3,3}$, a contradiction.
\end{proof}

\begin{fact}\label{case7-f4}
No vertex in $T_{x_{1}}$ is adjacent to $\{a,d\}$; at most five vertices in $T_{x_{1}}$ are not adjacent to $\{a,d\}$; at most one vertex in $T_{x_{1}}$ is not adjacent to $d$ and is adjacent to $a$; at most two vertices in $T_{x_{1}}$ are not adjacent to $a$ and are adjacent to $d$; $\left|E(T_{x_{1}}, T_{x_{0}x_{1}}\cup T_{x_{0}x_{2}}- \{a\})\right|=0$;  $\left|E(T_{x_{1}},T_{x_{1}x_{2}}-\{d\})\right|\leq  1$; if $q\in T_{x_{1}}$ and $\left|E(q,T_{x_{1}x_{2}}-\{d\})\right|= 1$, then $\left|E(q,d)\right|=1$ and $\left|E(q,a)\right|=0$.
\end{fact}
\begin{proof}
If a vertex in $T_{x_{1}}$ is adjacent to $\{a,d\}$, then $K_{4}\subset G$, a contradiction.
If at least six vertices in $T_{x_{1}}$ are not adjacent to $\{a,d\}$, then $G[\{x_{0},\allowbreak x_{1},\allowbreak x_{2},\allowbreak a,\allowbreak d,\allowbreak p_{1,1},\ldots,p_{1,6}\}]\cong R_{84}$, a contradiction.
If at least two vertices in $T_{x_{1}}$ are not adjacent to $d$ and are adjacent to $a$, then $G[\{x_{0},\allowbreak x_{2},\allowbreak x_{3},\allowbreak a,\allowbreak p_{1,1},\allowbreak p_{1,2}\}]\cong R_{1}$, a contradiction.
If at least three vertices in $T_{x_{1}}$ are not adjacent to $a$ and are adjacent to $d$, then $G[\{x_{0},\allowbreak x_{2},\allowbreak x_{3},\allowbreak d,\allowbreak p_{1,1},\allowbreak p_{1,2},\allowbreak p_{1,3}\}]\cong R_{21}$, a contradiction.
Since $G$ contains no subdivision of $K_{3,3}$, $\left|E(\{a,d\},T_{x_{1}x_{2}}\cup T_{x_{0}x_{1}}-\{a,d\})\right|=0$.
If $\left|E(T_{x_{1}}, T_{x_{0}x_{1}}- \{a\})\right|\geq 1$, then either $G$ contains a subdivision of $K_{3,3}$ or  $G[\{x_{1},\allowbreak x_{2},\allowbreak x_{3},\allowbreak q_{01,2},\allowbreak d,\allowbreak p_{1,1}\}]\cong R_{16}$ ($E(\{q_{01,2},\allowbreak d,\allowbreak p_{1,1}\})=\{q_{01,2}p_{1,1}\}$) by Lemma \ref{lem-13}\ref{lem-13:item-4}, a contradiction.
If $\left|E(T_{x_{1}}, T_{x_{0}x_{2}})\right|\geq 1$, then either $G$ contains a subdivision of $K_{3,3}$ or  $G[\{x_{0},\allowbreak x_{1},\allowbreak a,\allowbreak d,\allowbreak q_{02,1},\allowbreak p_{1,1}\}]\cong R_{12}$ by Lemma \ref{lem-13}\ref{lem-13:item-4}, a contradiction.
If $\left|E(T_{x_{1}},T_{x_{1}x_{2}}-\{d\})\right|\geq 2$, then either $R_{42}\subset G$ (if $\left|E(T_{x_{1}},\{a,d\})\right|=0$)  or $G$ contains a subdivision of $K_{3,3}$ (if $\left|E(T_{x_{1}},\{a,d\})\right|\geq 1$), a contradiction.
If $q\sim a$, then $G$ contains a subdivision of $K_{3,3}$, a contradiction.
If $q\nsim d$, then $G[\{x_{0},\allowbreak x_{1},\allowbreak a,\allowbreak d,\allowbreak q_{12,2},\allowbreak q\}]\cong R_{16}$, a contradiction.
\end{proof}

\begin{fact}\label{case7-f8}
At most one vertex in $T_{x_{2}}$ is not adjacent to $a$; at most one vertex in $T_{x_{2}}$ is adjacent to $\{a,d\}$; at most three vertices in $T_{x_{2}}$ are not adjacent to $d$ and are adjacent to $a$.
\end{fact}
\begin{proof}
If at least two vertices in $T_{x_{2}}$ are not adjacent to $a$, then $G[\mathcal{K}\cup\{a,\allowbreak p_{2,1},\allowbreak p_{2,2}\}]\cong R_{25}$, a contradiction.
If at least two vertices in $T_{x_{2}}$ are adjacent to $\{a,d\}$, then $G$ contains a subdivision of $K_{3,3}$, a contradiction.
If at least four vertices in $T_{x_{2}}$ are not adjacent to $d$ and are adjacent to $a$, then $G[\mathcal{K}\cup\{a,\allowbreak p_{2,1},\ldots,p_{2,4}\}]\cong R_{79}$, a contradiction.
\end{proof}

\begin{fact}\label{case7-f13}
$\left|T_{x_{1}x_{2}}-\{d\}\right|\leq 2$; $\left|T_{x_{0}x_{1}}-\{a\}\right|\leq 2$; $\left|E(T_{x_{2}},T_{x_{1}x_{2}}\cup T_{x_{0}x_{1}}-\{a,d\})\right|=\left|E(T_{x_{0}},T_{x_{0}x_{1}})\right|\\=0$.
\end{fact}
\begin{proof}
If $\left|T_{x_{1}x_{2}}-\{d\}\right| \geq 3$, then either $G$ contains a subdivision of $K_{3,3}$ or  $G[\{x_{0},\allowbreak x_{2},\allowbreak x_{3},\allowbreak a,\allowbreak d,\allowbreak q_{12,2},\allowbreak q_{12,3},\allowbreak q_{12,4}\}]\cong R_{37}$ ($E(\{a,\allowbreak d,\allowbreak q_{12,2},\allowbreak q_{12,3},\allowbreak q_{12,4}\})=\{ad\}$) by Lemma \ref{lem-13}\ref{lem-13:item-4}, a contradiction.
If $\left|T_{x_{0}x_{1}}-\{a\}\right| \geq 3$, then  either $G$ contains a subdivision of $K_{3,3}$ or $G[\{x_{0},\allowbreak x_{2},\allowbreak x_{3},\allowbreak q_{01,2},\allowbreak q_{01,3},\allowbreak q_{01,4},\allowbreak d\}]\cong R_{21}$ ($E(\{q_{01,2},\allowbreak q_{01,3},\allowbreak q_{01,4},\allowbreak d\})=\varnothing$) by Lemma \ref{lem-13}\ref{lem-13:item-4}, a contradiction.
If $\left|E(T_{x_{2}},T_{x_{1}x_{2}}-\{d\})\right| \geq 1$ or $\left|E(T_{x_{0}},T_{x_{0}x_{1}})\right|\geq 1$, then  either $G$ contains a subdivision of $K_{3,3}$ or $G[\{x_{0},\allowbreak x_{1},\allowbreak x_{2},\allowbreak a,\allowbreak q_{12,2},\allowbreak p_{2,1}\}]\cong R_{18}$ ($E(\{a,\allowbreak q_{12,2},\allowbreak p_{2,1}\})=\{q_{12,2}p_{2,1}\}$) by Lemma \ref{lem-13}\ref{lem-13:item-4}, a contradiction.
If $\left|E(T_{x_{2}},T_{x_{0}x_{1}}-\{a\})\right| \geq 1$, then  either $G$ contains a subdivision of $K_{3,3}$ or $G[\{x_{0},\allowbreak x_{1},\allowbreak a,\allowbreak q_{01,2},\allowbreak d,\allowbreak p_{2,1}\}]\cong R_{14}$ ($E(\{a,\allowbreak q_{01,2},\allowbreak d,\allowbreak p_{2,1}\})=\{ad,\allowbreak q_{01,2}p_{2,1}\}$) by Lemma \ref{lem-13}\ref{lem-13:item-4}, a contradiction.
\end{proof}

\begin{fact}\label{case7-f15}
$\left|E(T_{x_{1}x_{3}},d)\right|=1$ and $\left|T_{x_{1}x_{3}}\right|\leq 1$.
\end{fact}
\begin{proof}
If $\left|T_{x_{1}x_{3}}\right| \geq1$ and $\left|E(T_{x_{1}x_{3}},d)\right|=0$, then $R_{37}\subset G$, a contradiction.
If $\left|T_{x_{1}x_{3}}\right|\geq 2$, then  $G$ contains a subdivision of $K_{3,3}$, a contradiction.
\end{proof}

\begin{fact}\label{case7-f16}
At most one vertex in $T_{x_{0}x_{2}}$ is not adjacent to $d$; at most one vertex in $T_{x_{0}x_{2}}$ is adjacent to $d$; $\left|E(T_{x_{0}x_{2}},T_{x_{1}}\cup T_{x_{1}x_{2}}-\{d\})\right|=0$.
\end{fact}
\begin{proof}
If at least two vertices in $T_{x_{0}x_{2}}$ are not adjacent to $d$, then $G[\{x_{0},\allowbreak x_{1},\allowbreak x_{2},\allowbreak a,\allowbreak d,\allowbreak q_{02,1},\allowbreak q_{02,2}\}]\cong R_{33}$, a contradiction.
If $\left|E(T_{x_{0}x_{2}},d)\right|\geq 2$, then  $G$ contains a subdivision of $K_{3,3}$, a contradiction.
If $\left|E(T_{x_{0}x_{2}},T_{x_{1}}\cup T_{x_{1}x_{2}}-\{d\})\right|\geq 1$, then either $G$ contains a subdivision of $K_{3,3}$ or $G[\{x_{0},\allowbreak x_{1},\allowbreak a,\allowbreak d,\allowbreak q_{02,1},\allowbreak p_{1,1}(q_{12,2})\}]\cong R_{12}$ by (F\ref{case7-f4}) and Lemma \ref{lem-13}\ref{lem-13:item-4}, a contradiction.
\end{proof}

\begin{fact}\label{case7-f22}
$\left|E(T_{0},T-\{d\})\right|=0$; if $\left|E(T_{x_{1}}\cup T_{x_{2}},d)\right|+\left|T_{x_{0}x_{1}}\cup T_{x_{1}x_{2}}-\{a,d\}\right|\geq 1$, then  $\left|T_{0}\right|=0$.
\end{fact}
\begin{proof}
If $\left|E(T_{0},T_{x_{0}x_{1}x_{3}}\cup T_{x_{1}x_{2}}-\{d\})\right|\geq 1$, then either $G$ contains a subdivision of $K_{3,3}$ or $G[\{x_{0},\allowbreak x_{1},\allowbreak x_{2},\allowbreak d,\allowbreak w_{013,1}(q_{12,2}),\allowbreak t_{1}\}]\cong R_{14}$ by Lemma \ref{lem-13}\ref{lem-13:item-8}, a contradiction.
Thus, $\left|E(T_{0},T-\{d\})\right|=0$ by Lemma \ref{lem-13}\ref{lem-13:item-8} and (F\ref{case7-f1}).
If $\left|E(T_{x_{1}}\cup T_{x_{2}},d)\right|\geq 1$ and $\left|T_{0}\right|\geq 1$, then $G[\mathcal{K}\cup\{d,\allowbreak p_{1,1},\allowbreak t_{1}\}]\cong R_{28}$ by (F\ref{case7-f4}), a contradiction.
If $\left|T_{x_{0}x_{1}}\cup T_{x_{1}x_{2}}-\{a,d\}\right|\geq 1$ and $\left|T_{0}\right|\geq 1$, then either $R_{14}\subset G$ (if $\left|T_{x_{0}x_{1}}-\{a\}\right|\geq 1$) or $R_{28}\subset G$ (if $\left|T_{x_{1}x_{2}}-\{d\}\right|\geq 1$), a contradiction.
\end{proof}

\begin{fact}\label{case7-f27}
If $\left|E(T_{x_{1}},T_{x_{1}x_{2}}-\{d\})\right|=0$, then $\left|E(T_{x_{1}},d)\right|+\left|T_{x_{1}x_{2}}-\{d\}\right|\leq 2$.
\end{fact}
\begin{proof}
If $\left|E(T_{x_{1}},d)\right|+\left|T_{x_{1}x_{2}}-\{d\}\right|\geq 3$, then either $G$ contains a subdivision of $K_{3,3}$ or $G[\{x_{0},\allowbreak x_{2},\allowbreak x_{3},\allowbreak d,\allowbreak p_{1,1},\allowbreak p_{1,2},\allowbreak p_{1,3}\}]\cong R_{21}$ by Lemma \ref{lem-13}\ref{lem-13:item-8}, a contradiction.
\end{proof}

\begin{fact}\label{case7-f29}
Suppose that $\left|T_{x_{1}x_{2}}-\{d\}\right|\geq 1$. 
Then we have:   $\left|T_{x_{0}}\right|=0$;  no vertex in $T_{x_{2}}$ is not adjacent to $\{a,d\}$;  at most one vertex in $T_{x_{2}}$ is adjacent to $a$ and  is not adjacent to $d$;   if $\left|T_{x_{1}x_{2}}-\{d\}\right|=2$,  then $\left|E(T_{x_{2}},a)\right|=0$.
\end{fact}
\begin{proof}
Since $G$ contains no subdivision of $K_{3,3}$, we have $\left|E(a,T_{x_{1}x_{2}}-\{d\})\right|=0$.
If $\left|T_{x_{0}}\right|\geq 1$, then either $R_{12}\subset G$ (if $\left|E(T_{x_{0}},T_{x_{1}x_{2}}-\{d\})\right|\geq 1$) or $R_{28}\subset G$ (if $\left|E(T_{x_{0}},T_{x_{1}x_{2}}-\{d\})\right|=0$) by Lemma \ref{lem-13}\ref{lem-13:item-4}, (F\ref{case7-f2}) and (F\ref{case7-f13}), a contradiction.
If at least one vertex in $T_{x_{2}}$ is not adjacent to $\{a,d\}$, then  $G[\mathcal{K}\cup\{a,\allowbreak d,\allowbreak q_{12,2},\allowbreak p_{2,1}\}]\cong R_{43}$ by Lemma \ref{lem-13}\ref{lem-13:item-4} and (F\ref{case7-f2}), a contradiction.
If at least two vertices in $T_{x_{2}}$ are adjacent to $a$ and  are not adjacent to $d$, then  $G[\{x_{1},\allowbreak x_{2},\allowbreak x_{3},\allowbreak a,\allowbreak q_{12,2},\allowbreak p_{2,1},\allowbreak p_{2,2}\}]\cong R_{26}$ by Lemma \ref{lem-13}\ref{lem-13:item-4} and (F\ref{case7-f13}), a contradiction.
If $\left|E(T_{x_{2}},a)\right|\geq 1$, then  either $G[\{x_{0},\allowbreak x_{1},\allowbreak x_{3},\allowbreak a,\allowbreak d,\allowbreak q_{12,2},\allowbreak p_{2,1}\}]\cong R_{25}$ ($E(\{a,\allowbreak d,\allowbreak q_{12,2},\allowbreak p_{2,1}\})=\{ad,\allowbreak ap_{2,1},\allowbreak dp_{2,1}\}$) or $G[\mathcal{K}\cup\{a,\allowbreak q_{12,2},\allowbreak q_{12,3},\allowbreak p_{2,1}\}]\cong R_{42}$ ($E(\{a,\allowbreak q_{12,2},\allowbreak q_{12,3},\allowbreak p_{2,1}\})=\{ap_{2,1}\}$) by Lemma \ref{lem-13}\ref{lem-13:item-4} and (F\ref{case7-f13}), a contradiction.
\end{proof}

\begin{fact}\label{case7-f31}
If $\left|T_{x_{0}x_{1}}-\{a\}\right|\geq 1$, then   $\left|E(T_{x_{0}},d)\right|=1$, $\left|E(T_{x_{2}},d)\right|=0$ and $\left|E(T_{x_{2}},a)\right|+\left|T_{x_{0}x_{1}}-\{a\}\right|\leq 2$.
\end{fact}
\begin{proof}
Since $G$ contains no subdivision of $K_{3,3}$, $\left|E(d,T_{x_{0}x_{1}}-\{a\})\right|=0$.
If $\left|E(T_{x_{0}},d)\right|=0$, then  $R_{28}\subset G$ by Lemma \ref{lem-13}\ref{lem-13:item-2} and \ref{lem-13:item-4}, a contradiction.
If $\left|E(T_{x_{2}},d)\right|\geq 1$, then $G[\{x_{0},\allowbreak x_{1},\allowbreak x_{2},\allowbreak q_{01,2},\allowbreak d,\allowbreak p_{2,1}\}]\cong R_{18}$ by Lemma \ref{lem-13}\ref{lem-13:item-4} and (F\ref{case7-f13}), a contradiction.
If $\left|E(T_{x_{2}},a)\right|+\left|T_{x_{0}x_{1}}-\{a\}\right|\geq 3$, then  $R_{25}\subset G$ (if $\left|E(T_{x_{2}},d)\right|=0$ and $\left|T_{x_{2}}\right|=2$) or $R_{28}\subset G$ (if $\left|E(T_{x_{2}},d)\right|=0$ and $\left|T_{x_{2}}\right|=1$) or  $R_{18}\subset G$ (the remaining case)  by Lemma \ref{lem-13}\ref{lem-13:item-4}, (F\ref{case7-f8}) and (F\ref{case7-f13}), a contradiction.
\end{proof}

\begin{fact}\label{case7-f34} 
If $\left|E(T_{x_{0}x_{1}},T_{x_{1}x_{2}})\right|=1$, then $\left|T_{x_{0}x_{1}}\cup T_{x_{1}x_{2}}-\{a,d\}\right|\leq 2$.
\end{fact}
\begin{proof}
If $\left|T_{x_{0}x_{1}}\cup T_{x_{1}x_{2}}-\{a,d\}\right|\geq 3$, then  $G[\{x_{0},\allowbreak x_{2},\allowbreak x_{3},\allowbreak a,\allowbreak q_{01,2},\allowbreak q_{01,3},\allowbreak q_{12,2}\}]\cong R_{21}$  or $G[\{x_{0},\allowbreak x_{2},\allowbreak x_{3},\allowbreak a,\allowbreak q_{01,2},\allowbreak q_{12,2},\allowbreak q_{12,3}\}]\cong R_{21}$, a contradiction.
\end{proof}

\begin{fact}\label{case7-f37} 
If $\left|E(T_{x_{2}}\cup T_{x_{0}},d)\right|\geq 1$, then at most one vertex in $T_{x_{1}}$ is not adjacent to $\{a,d\}$.
\end{fact}
\begin{proof}
If at least two vertices in $T_{x_{1}}$ are not adjacent to $\{a,d\}$, then either $R_{25}\subset G$ (if $\left|E(T_{x_{2}},d)\right|\geq 1$) or $R_{38}\subset G$ (if $\left|E(T_{x_{0}},d)\right|\geq 1$), a contradiction.
\end{proof}

\begin{fact}\label{case7-f39} 
If $\left|E(T_{x_{1}},a)\right|\geq1$, then  $\left|E(T_{x_{2}},a)\right|=\left|T_{x_{2}}\right|$.
\end{fact}
\begin{proof}
If $\left|E(T_{x_{2}},a)\right|<\left|T_{x_{2}}\right|$, then  $G[\{x_{0},\allowbreak x_{1},\allowbreak x_{2},\allowbreak a,\allowbreak p_{1,1},\allowbreak p_{2,1}\}]\cong R_{14}$, a contradiction.
\end{proof}

\begin{fact}\label{case7-f41} 
If  $S_{1}\subset T_{x_{1}}$, $S_{2}\subset T_{x_{2}}$ ($\left|S_{2}\right| \geq1$), $\left|E(S_{1},\{a,d\})\right|=0$, $\left|E(S_{2},d)\right|=0$, and $\left|E(T_{x_{2}},a)\right|=\left|T_{x_{2}}\right|$, then $\left|S_{1}\cup S_{2}\right| \leq 3$.
\end{fact}
\begin{proof}
If $\left|S_{1}\cup S_{2}\right| \geq 4$, then either $G[\{x_{1},\allowbreak x_{3},\allowbreak a,\allowbreak p_{1,1},\allowbreak p_{2,1},\allowbreak p_{2,2},\allowbreak p_{2,3}\}]\cong R_{21}$ ($E(\{a,\allowbreak p_{1,1},\allowbreak p_{2,1},\allowbreak p_{2,2},\allowbreak p_{2,3}\})=\{ap_{2,1},\allowbreak ap_{2,2},\allowbreak ap_{2,3}\}$) or $G[\{x_{1},\allowbreak x_{3},\allowbreak a,\allowbreak p_{1,1},\allowbreak p_{1,2},\allowbreak p_{2,1},\allowbreak p_{2,2}\}]\cong R_{21}$ ($E(\{a,\allowbreak p_{1,1},\allowbreak p_{1,2},\allowbreak p_{2,1},\allowbreak p_{2,2}\})=\{ap_{2,1},\allowbreak ap_{2,2}\}$) or $G[\{x_{0},\allowbreak x_{1},\allowbreak x_{2},\allowbreak a,\allowbreak p_{1,1},\allowbreak p_{1,2},\allowbreak p_{1,3},\allowbreak p_{2,1}\}]\cong R_{38}$ ($E(\{a,\allowbreak p_{1,1},\allowbreak p_{1,2},\allowbreak p_{1,3},\allowbreak p_{2,1}\})=\{ap_{2,1}\}$) by (F\ref{case7-f8}), a contradiction.
\end{proof}

\begin{fact}\label{case7-f42} 
If   $S_{3}\subset T_{x_{1}}$ ($\left|S_{3}\right| \geq 1$), $S_{4}\subset T_{x_{2}}$, $\left|E(S_{3},a)\right|=\left|E(S_{4},d)\right|=0$, and all vertices of $S_{3}$ ($S_{4}$) are adjacent to $d$ ($a$), then $\left|S_{3}\cup S_{4}\right| \leq 2$. 
\end{fact}
\begin{proof}
If $\left|S_{3}\cup S_{4}\right| \geq 3$, then either $G[\{x_{0},\allowbreak x_{1},\allowbreak x_{2},\allowbreak d,\allowbreak p_{1,1},\allowbreak p_{2,1},\allowbreak p_{2,2}\}]\cong R_{25}$ ($E(\{d,\allowbreak p_{1,1},\allowbreak p_{2,1},\allowbreak p_{2,2}\})=\{dp_{1,1}\}$) or $G[\{x_{0},\allowbreak x_{1},\allowbreak x_{2},\allowbreak d,\allowbreak p_{1,1},\allowbreak p_{1,2},\allowbreak p_{2,1}\}]\cong R_{28}$ ($E(\{d,\allowbreak p_{1,1},\allowbreak p_{1,2},\allowbreak p_{2,1}\})=\{dp_{1,1},\allowbreak dp_{1,2}\}$) by (F\ref{case7-f4}), a contradiction.
\end{proof}

Suppose that $\left|E(T_{x_{0}x_{1}},T_{x_{1}x_{2}})\right|\geq 2$.
By Lemma \ref{lem-13}\ref{lem-13:item-4},  $G[E(T_{x_{0}x_{1}},T_{x_{1}x_{2}})]\neq aK_{3}$ ($a\geq 1$).
If $G[E(T_{x_{0}x_{1}},T_{x_{1}x_{2}})]= aP_{3}$, then $G$ contains a subdivision of $K_{3,3}$, a contradiction.
If $\left|E(T_{x_{0}x_{1}},T_{x_{1}x_{2}})\right|\geq 3$, then $G$ contains a subdivision of $K_{3,3}$, a contradiction.
Thus, $G[E(T_{x_{0}x_{1}},T_{x_{1}x_{2}})]=2K_{2}$.
Let  $q_{01,1}\in T_{x_{0}x_{1}}-\{a\}$, $q_{12,1}\in T_{x_{1}x_{2}}-\{d\}$ and $q_{01,1}\sim q_{12,1}$.
By (F\ref{case7-f15}) and   $G$ contains no subdivision of $K_{3,3}$, $\left|T_{x_{1}x_{3}}\right|=0$.
Since $G$ contains no subdivision of $K_{3,3}$, $\left|T_{x_{0}x_{1}x_{3}}\right|=0$.
Then $\left|T_{x_{2}}\right|=0$; otherwise, $G$ contains a subdivision of $K_{3,3}$ by (F\ref{case7-f29}), a contradiction.
Then $\left|E(T_{x_{1}},\{a,q_{01,1}\})\right|=0$; otherwise, either $G[\{x_{0},\allowbreak x_{2},\allowbreak x_{3},\allowbreak a,\allowbreak q_{12,1},\allowbreak p_{1,1}\}]\cong R_{1}$ ($E(\{a,\allowbreak q_{12,1},\allowbreak p_{1,1}\})=\{ap_{1,1}\}$) or $G$ contains a subdivision of $K_{3,3}$, a contradiction.
Thus, $\left|E(T_{x_{1}}, T_{x_{0}x_{1}})\right|=0$ by (F\ref{case7-f4}).
Then $\left|E(T_{x_{1}},T_{x_{1}x_{2}})\right|=0$; otherwise, either $G[\{x_{0},\allowbreak x_{1},\allowbreak q_{01,1},\allowbreak d,\allowbreak q_{12,1},\allowbreak p_{1,1}\}]\cong R_{16}$ ($E(\{q_{01,1},\allowbreak d,\allowbreak q_{12,1},\allowbreak p_{1,1}\})=\{q_{01,1}q_{12,1},\allowbreak dp_{1,1}\}$) or $G$ contains a subdivision of $K_{3,3}$, a contradiction.
Then $\left|E(T_{x_{0}x_{2}},\{d,q_{12,1}\})\right|=0$; otherwise, either $G[\{x_{0},\allowbreak x_{1},\allowbreak q_{01,1},\allowbreak d,\allowbreak q_{12,1},\allowbreak q_{02,1}\}]\cong R_{12}$ ($E(\{q_{01,1},\allowbreak d,\allowbreak q_{12,1},\allowbreak q_{02,1}\})=\{q_{01,1}q_{12,1},\allowbreak dq_{02,1}\}$) or $G$ contains a subdivision of $K_{3,3}$, a contradiction.
Thus,  $\left|E(T_{x_{0}x_{2}},T_{x_{1}}\cup T_{x_{1}x_{2}})\right|=0$ by (F\ref{case7-f16}).
Then $\left|T_{x_{0}x_{1}} \cup T_{x_{1}x_{2}}-\{a,q_{01,1},d,q_{12,1}\}\right|\leq 1$; otherwise, $R_{21}\subset G$, a contradiction.
Then  $\left|E(T_{x_{0}x_{2}},T_{x_{1}x_{2}})\right|=0$; otherwise,  $G[\{x_{1},x_{2},x_{3},a,q_{01,1},q_{12,2},q_{02,1}\}]\cong R_{14}$ ($q_{12,2}\sim q_{02,1}$), a contradiction.
If $\left|T_{x_{0}x_{2}}\right|\geq 1$, then $\left|T_{x_{1}}\right|\leq 1$; otherwise, $G[\{x_{0},\allowbreak x_{1},\allowbreak x_{2},\allowbreak a,\allowbreak q_{02,1},\allowbreak p_{1,1},\allowbreak p_{1,2}\}]\cong R_{25}$, a contradiction.
Thus,  
$G\subset H_{116}$ (if $\left|T_{x_{0}x_{2}}\right|\geq1$ and $\left|T_{x_{1}x_{2}}\right|=3$) or  
$G\subset H_{117}$ (if $\left|T_{x_{0}x_{2}}\right|\geq1$ and $\left|T_{x_{1}x_{2}}\right|=2$) or  
$G\subset H_{134}$ (if $\left|T_{x_{0}x_{2}}\right|=0$ and $\left|T_{x_{1}x_{2}}\right|=3$) or  
$G\subset H_{135}$ (if $\left|T_{x_{0}x_{2}}\right|=0$ and $\left|T_{x_{1}x_{2}}\right|=2$)  
by (F\ref{case7-f1}), (F\ref{case7-f4}),  (F\ref{case7-f13}), (F\ref{case7-f16}), (F\ref{case7-f22}) and (F\ref{case7-f29}).
Let $\left|E(T_{x_{0}x_{1}},T_{x_{1}x_{2}})\right|=1$.

Suppose that $\left|T_{x_{1}x_{3}}\right|\geq 1$.
By (F\ref{case7-f15}), $\left|E(T_{x_{1}x_{3}},d)\right|=1$ and $\left|T_{x_{1}x_{3}}\right|=1$.
Then $\left|T_{x_{0}x_{2}}\right|=0$ by Lemma \ref{lem-13}\ref{lem-13:item-6} and   $G$ contains no subdivision of $K_{3,3}$.
Then $\left|T_{x_{0}x_{1}}-\{a\}\right|=0$; otherwise, $G[\{x_{0},\allowbreak x_{1},\allowbreak x_{3},\allowbreak a,\allowbreak q_{01,2},\allowbreak q_{13,1}\}]\cong R_{16}$, a contradiction.
Then $\left|T_{x_{0}}\right|=0$; otherwise, either $G$ contains a subdivision of $K_{3,3}$ (if $\left|E(T_{x_{0}},T_{x_{1}x_{3}})\right|\geq 1$) or $R_{14}\subset G$ (the remaining case), a contradiction.
Then $\left|T_{x_{0}x_{1}x_{3}}\right|=0$ for $G$ contains no subdivision of $K_{3,3}$.
Since  $G$ contains no subdivision of $K_{3,3}$, $\left|E(T_{x_{1}x_{2}} -\{d\},T_{x_{1}x_{3}})\right|=0$.
Then $\left|E(T_{x_{1}},a)\right|=0$ by (F\ref{case7-f4}); otherwise, either $G[\{x_{1},\allowbreak x_{2},\allowbreak x_{3},\allowbreak a,\allowbreak q_{13,1},\allowbreak p_{1,1}\}]\cong R_{16}$ ($E(\{a,\allowbreak q_{13,1},\allowbreak p_{1,1}\})=\{ap_{1,1}\}$) or  $G$ contains a subdivision of $K_{3,3}$, a contradiction.
Then $\left|E(T_{x_{1}},T_{x_{1}x_{3}})\right|=0$; otherwise, $G[\{x_{0},\allowbreak x_{1},\allowbreak x_{2},\allowbreak a,\allowbreak q_{13,1},\allowbreak p_{1,1}\}]\cong R_{16}$, a contradiction.
Then $\left|E(T_{x_{2}},T_{x_{1}x_{3}})\right|=0$; otherwise, either $G[\{x_{0},\allowbreak x_{1},\allowbreak x_{2},\allowbreak a,\allowbreak q_{13,1},\allowbreak p_{2,1}\}]\cong R_{12}$ ($E(\{a,\allowbreak q_{13,1},\allowbreak p_{2,1}\})=\{q_{13,1}p_{2,1}\}$) or $G$ contains a subdivision of $K_{3,3}$, a contradiction.
Then $\left|E(T_{x_{2}},a)\right|=0$; otherwise, $G[\{x_{0},\allowbreak x_{1},\allowbreak x_{3},\allowbreak a,\allowbreak q_{13,1},\allowbreak p_{2,1}\}]\cong R_{10}$, a contradiction.
Thus, $\left|T_{x_{2}}\right|\leq 1$ by (F\ref{case7-f8}).
If $\left|T_{x_{2}}\right|=1$ and $\left|E(T_{x_{2}},d)\right|=0$, then $\left|T_{x_{1}x_{2}}-\{d\}\right|=0$ by (F\ref{case7-f13}); otherwise $G[\{x_1,x_2,a,d,q_{12,2},q_{13,1},p_{2,1}\}]\cong R_{28}$, a contradiction.
If $\left|T_{x_{2}}\right|=1$ and $\left|E(T_{x_{2}},d)\right|=0$, then $\left|E(T_{x_{1}},d)\right|=0$; otherwise, $G[\{x_{0},\allowbreak x_{1},\allowbreak x_{2},\allowbreak d,\allowbreak q_{13,1},\allowbreak p_{2,1},\allowbreak p_{1,1}\}]\cong R_{28}$, a contradiction.
If $\left|T_{x_{2}}\right|=1$ and $\left|E(T_{x_{2}},d)\right|=0$, then $\left|T_{x_{1}}\right|\leq 3$; otherwise, $G[\mathcal{K}\cup\{p_{2,1},\allowbreak p_{1,1},\ldots,p_{1,4}\}]\cong R_{47}$, a contradiction.
Thus, if $\left|T_{x_{2}}\right|=1$ and $\left|E(T_{x_{2}},d)\right|=0$, then $G\subset H_{132}$ by (F\ref{case7-f1}), (F\ref{case7-f22}) and (F\ref{case7-f29}).
By (F\ref{case7-f4}) and   $G$ contains no subdivision of $K_{3,3}$,  $\left|E(T_{x_{1}},T_{x_{1}x_{2}}-\{d\})\right|=0$.
Thus, if  $\left|T_{x_{2}}\right|=1$ and $\left|E(T_{x_{2}},d)\right|=1$, then either $G\subset H_{116}$ (if $\left|E(T_{x_{1}},d)\right|= 2$ or $0$) or $G\subset H_{117}$ (if $\left|E(T_{x_{1}},d)\right|=1$) by (F\ref{case7-f1}), (F\ref{case7-f13}), (F\ref{case7-f22}), (F\ref{case7-f27}) and (F\ref{case7-f37}).
If  $\left|T_{0}\right|=1$, then $\left|T_{x_{1}}\right|\leq 3$  by  (F\ref{case7-f4}) and (F\ref{case7-f22}); otherwise, $G[\{x_{1},\allowbreak x_{3},\allowbreak a,\allowbreak d,\allowbreak t_{1},\allowbreak p_{1,1},\ldots,p_{1,4}\}]\cong R_{46}$, a contradiction.
If $\left|T_{x_{2}}\right|=0$, then $G\subset H_{132}$ (if $\left|T_{0}\right|= 1$) or $G\subset H_{134}$ (if $\left|T_{0}\right|= 0$ and $\left|E(T_{x_{1}},d)\right|= 2$ or $0$) or $G\subset H_{135}$ (if $\left|T_{0}\right|= 0$ and $\left|E(T_{x_{1}},d)\right|=1$) by (F\ref{case7-f1}), (F\ref{case7-f4}), (F\ref{case7-f22}) and (F\ref{case7-f27}).
Let  $\left|T_{x_{1}x_{3}}\right|=0$.

Suppose that $\left|T_{x_{0}x_{1}x_{3}}\right|\geq1$.
By  Lemma \ref{lem-13}\ref{lem-13:item-7}, $\left|T_{x_{0}x_{1}x_{3}}\right|=1$.
Then $\left|E(T_{x_{0}x_{1}x_{3}},\{a,d\})\right|=0$ for $G$ contains no subdivision of $K_{3,3}$.
Then $\left|T_{0}\right|=0$ by (F\ref{case7-f22}); otherwise, $G[\{x_{0},\allowbreak x_{1},\allowbreak x_{2},\allowbreak d,\allowbreak w_{013,1},\allowbreak t_{1}\}]\cong R_{14}$, a contradiction.
By  Lemma \ref{lem-13}\ref{lem-13:item-2}, $\left|E(T_{x_{0}},T_{x_{0}x_{1}})\right|= 0$.
Then $\left|E(T_{x_{0}},d)\right|= 0$; otherwise, either $G[\{x_{0},\allowbreak x_{1},\allowbreak x_{3},\allowbreak d,\allowbreak w_{013,1},\allowbreak p_{0,1}\}]\cong R_{12}$ or  $G$ contains a subdivision of $K_{3,3}$, a contradiction.
Then $\left|E(T_{x_{0}x_{1}x_{3}},T_{x_{0}})\right|= 0$; otherwise, $G[\{x_{0},\allowbreak x_{1},\allowbreak x_{2},\allowbreak d,\allowbreak w_{013,1},\allowbreak p_{0,1}\}]\cong R_{18}$ ($E(\{d,\allowbreak w_{013,1},\allowbreak p_{0,1}\})=\{w_{013,1}p_{0,1}\}$), a contradiction.
Then $\left|E(T_{x_{0}x_{1}x_{3}},T_{x_{1}})\right|= 0$; otherwise, either $G[\{x_{0},\allowbreak x_{2},\allowbreak a,\allowbreak d,\allowbreak w_{013,1},\allowbreak p_{1,1}\}]\cong R_{2}$  or $G$ contains a subdivision of $K_{3,3}$, a contradiction.
Then $\left|E(T_{x_{1}},\{a,d\})\right|= 0$; otherwise, $G[\{x_{1},\allowbreak x_{2},\allowbreak x_{3},\allowbreak a,\allowbreak w_{013,1},\allowbreak p_{1,1}\}]\cong R_{16}$ ($E(\{a,\allowbreak w_{013,1},\allowbreak p_{1,1}\})=\{ap_{1,1}\}$) or $G[\{x_{1},\allowbreak x_{3},\allowbreak a,\allowbreak d,\allowbreak w_{013,1},\allowbreak p_{1,1}\}]\cong R_{16}$ ($E(\{a,\allowbreak d,\allowbreak w_{013,1},\allowbreak p_{1,1}\})=\{ad,\allowbreak dp_{1,1}\}$), a contradiction.
Since $G$ contains no subdivision of $K_{3,3}$,  $\left|E(T_{x_{1}},T_{x_{1}x_{2}})\right|=0$.
Then $\left|E(T_{x_{0}x_{1}x_{3}},T_{x_{2}})\right|= 0$; otherwise, either $G[\{x_{0},\allowbreak x_{1},\allowbreak a,\allowbreak d,\allowbreak w_{013,1},\allowbreak p_{2,1}\}]\cong R_{14}$ or  $G$ contains a subdivision of $K_{3,3}$, a contradiction.
Then $\left|E(T_{x_{2}},\{a,d\})\right|= 0$; otherwise, $G[\{x_{0},\allowbreak x_{1},\allowbreak x_{3},\allowbreak a,\allowbreak w_{013,1},\allowbreak p_{2,1}\}]\cong R_{14}$, a contradiction.
Since  $G$ contains no subdivision of $K_{3,3}$,   $\left|E(T_{x_{2}},T_{x_{0}x_{1}})\right|= 0$.
Then $\left|E(T_{x_{0}x_{1}x_{3}},T_{x_{0}x_{1}})\right|= 0$; otherwise, $G[\{x_{0},\allowbreak x_{1},\allowbreak x_{2},\allowbreak q_{01,2},\allowbreak d,\allowbreak w_{013,1}\}]\cong R_{20}$, a contradiction.
Since $G$ contains no subdivision of $K_{3,3}$,  $\left|E(T_{x_{0}x_{1}x_{3}},T_{x_{1}x_{2}})\right|= 0$.
Then $\left|E(T_{x_{0}x_{1}x_{3}},T_{x_{0}x_{2}})\right|= 0$; otherwise, either $G[\{x_{0},\allowbreak x_{1},\allowbreak a,\allowbreak d,\allowbreak w_{013,1},\allowbreak q_{02,1}\}]\cong R_{18}$ or  $G$ contains a subdivision of $K_{3,3}$, a contradiction.
Then $\left|E(T_{x_{0}x_{2}},d)\right|= 0$; otherwise, $G[\{x_{0},\allowbreak x_{1},\allowbreak x_{3},\allowbreak d,\allowbreak w_{013,1},\allowbreak q_{02,1}\}]\cong R_{12}$, a contradiction.
Thus, $\left|T_{x_{0}x_{2}}\right|\leq  1$ by (F\ref{case7-f16}).
If $\left|T_{x_{1}x_{2}}-\{d\}\right|\geq 1$, then $\left|T_{x_{0}}\cup T_{x_{2}}\right|=0$ by (F\ref{case7-f29}), $R_{18}\not\subset G$ and $R_{28}\not\subset G$.
If $\left|T_{x_{0}x_{1}}-\{a\}\right|\geq 1$, then $\left|T_{x_{0}}\cup T_{x_{2}}\right|=0$; otherwise, $G[\mathcal{K}\cup\{q_{01,2},\allowbreak d,\allowbreak p_{0,1}\}]\cong R_{28}$, a contradiction.
If $\left|T_{x_{0}x_{2}}\right|= 1$, then $\left|T_{x_{0}}\cup T_{x_{2}}\right|=0$ for $R_{18}\not\subset G$ and $R_{30}\not\subset G$.
Since  $G$ contains no subdivision of $K_{3,3}$,   $\left|E(T_{x_{1}},T_{x_{0}x_{2}})\right|= 0$.
If  $\left|T_{x_{0}x_{2}}\right|= 1$, then $\left|T_{x_{1}}\right|\leq 1$; otherwise, $G[\{x_{0},\allowbreak x_{1},\allowbreak x_{2},\allowbreak a,\allowbreak q_{02,1},\allowbreak p_{1,1},\allowbreak p_{1,2}\}]\cong R_{25}$, a contradiction.
Thus, if $\left|T_{x_{1}x_{2}}\cup T_{x_{0}x_{1}}\cup T_{x_{0}x_{2}}-\{a,d\}\right|\geq 1$, then   $G\subset H_{116}$ (if $\left|T_{x_{1}x_{2}}-\{d\}\right|=2$ or $0$, and $\left|T_{x_{0}x_{2}}\right|=1$) or  $G\subset H_{117}$ (if $\left|T_{x_{1}x_{2}}-\{d\}\right|=1$ and $\left|T_{x_{0}x_{2}}\right|=1$) or  $G\subset H_{134}$ (if $\left|T_{x_{1}x_{2}}-\{d\}\right|=2$ or $0$, and $\left|T_{x_{0}x_{2}}\right|=0$) or  $G\subset H_{135}$ (if $\left|T_{x_{1}x_{2}}-\{d\}\right|=1$ and $\left|T_{x_{0}x_{2}}\right|=0$) by (F\ref{case7-f1}), (F\ref{case7-f4}) and (F\ref{case7-f34}).
Let $\left|T_{x_{1}x_{2}}\cup T_{x_{0}x_{1}}\cup T_{x_{0}x_{2}}-\{a,d\}\right|=0$.
If  $\left|T_{x_{2}}\cup T_{x_{0}}\right|\geq 1$, then $\left|T_{x_{1}}\right|\leq 3$; otherwise, $G[\{x_{0},\allowbreak x_{1},\allowbreak x_{2},\allowbreak a,\allowbreak p_{1,1},\ldots,p_{1,4},\allowbreak p_{0,1}\}]\cong R_{47}$, a contradiction.
Thus, either $G\subset H_{132}$ (if $\left|T_{x_{2}}\cup T_{x_{0}}\right|\geq 1$) or $G\subset H_{134}$ (if $\left|T_{x_{2}}\cup T_{x_{0}}\right|=0$) by (F\ref{case7-f1}), (F\ref{case7-f2}) and (F\ref{case7-f8}).
Let $\left|T_{x_{0}x_{1}x_{3}}\right|=0$.

Suppose that $\left|T_{x_{0}x_{2}}\right|\geq 1$.
If  $b\in T_{x_{0}x_{2}}$ and $b\nsim d$, then its conclusion is covered by Case $1$, by taking $x_{0},x_{1},x_{2},b$ as the vertices of a new $K_{4}-e$.
Thus,  $\left|T_{x_{0}x_{2}}\right|=1$ and $\left|E(T_{x_{0}x_{2}},d)\right|= 1$ by (F\ref{case7-f16}).
Then $\left|E(T_{x_{1}},a)\right|= 0$ by (F\ref{case7-f4}); otherwise, $G[\{x_{0},\allowbreak x_{1},\allowbreak x_{2},\allowbreak a,\allowbreak q_{02,1},\allowbreak p_{1,1}\}]\cong R_{18}$, a contradiction.
By (F\ref{case7-f4}), $\left|E(T_{x_{1}},T_{x_{0}x_{1}}\cup T_{x_{0}x_{2}})\right|=0$.
Then $\left|E(T_{x_{1}},T_{x_{1}x_{2}}- \{d\})\right|=0$ by (F\ref{case7-f4}); otherwise, $G[\{x_{0},\allowbreak x_{1},\allowbreak x_{2},\allowbreak q_{12,2},\allowbreak q_{02,1},\allowbreak p_{1,1}\}]\cong R_{18}$, a contradiction.
Thus, by symmetry, $\left|E(T_{x_{2}},T_{x_{0}x_{1}}\cup T_{x_{0}x_{2}}\cup T_{x_{1}x_{2}}- \{d\})\right|= 0$.
By (F\ref{case7-f8}) and symmetry, at most one vertex in $T_{x_{2}}$ (or $T_{x_{1}}$) is not adjacent to $a$.
If  $\left|T_{0}\right|=1$, then  $\left|T_{x_{0}x_{1}}\cup T_{x_{0}x_{2}}\right|=2$ by (F\ref{case7-f22}); otherwise, $G[\{x_{0},\allowbreak x_{1},\allowbreak x_{2},\allowbreak q_{01,2},\allowbreak d,\allowbreak t_{1}\}]\cong R_{14}$, a contradiction.
If  $\left|T_{0}\right|=1$, then  $\left|T_{x_{1}x_{2}}-\{d\}\right|=\left|E(T_{x_{1}}\cup T_{x_{2}},d)\right|=0$ by (F\ref{case7-f22}).
If  $\left|T_{0}\right|=1$, then  $\left|T_{x_{0}}\right|=0$ by Lemma \ref{lem-13}\ref{lem-13:item-2} and $d(G)< 4$; otherwise, $G[\{x_{1},\allowbreak x_{2},\allowbreak x_{3},\allowbreak a,\allowbreak d,\allowbreak p_{0,1},\allowbreak t_{1}\}]\cong R_{25}$, a contradiction.
Thus, if  $\left|T_{0}\right|=1$, then $G\subset A_{21}$ by (F\ref{case7-f1}).
Let $\left|T_{0}\right|=0$.
Then $\left|T_{x_{0}x_{1}}\right|\leq 2$ and $\left|T_{x_{0}x_{2}}\right| \leq2$; otherwise, $G[\{x_{0},\allowbreak x_{1},\allowbreak x_{2},\allowbreak q_{01,2},\allowbreak q_{01,3},\allowbreak d,\allowbreak q_{02,1}\}]\cong R_{33}$, a contradiction.
Thus,  $\left|T_{x_{0}x_{1}}\right|+\left|T_{x_{2}}\right|\leq 2$ and $\left|T_{x_{0}x_{2}}\right|+\left|T_{x_{1}}\right|\leq 2$ for $R_{18}\not\subset G$ and $R_{28}\not\subset G$.
Then $\left|T_{x_{1}x_{2}}-\{d\}\right|\leq 1$; otherwise, $G[\{x_{0},\allowbreak x_{1},\allowbreak x_{3},\allowbreak a,\allowbreak d,\allowbreak q_{12,2},\allowbreak q_{12,3},\allowbreak q_{02,1}\}]\cong R_{38}$, a contradiction.
Thus,   $\left|T_{x_{1}x_{2}}\cup T_{x_{0}}-\{d\}\right|\leq1$ by (F\ref{case7-f29}).
If $\left|T_{x_{1}}\cup T_{x_{2}}\right|= 2$, then $\left|E(T_{x_{1}},d)\right|+\left|E(T_{x_{2}},d)\right|= 0$ or $2$; otherwise, $G[\{x_{1},\allowbreak x_{2},\allowbreak x_{3},\allowbreak a,\allowbreak d,\allowbreak p_{1,1},\allowbreak p_{2,1}\}]\cong R_{28}$, a contradiction.
If $\left|E(T_{x_{1}}\cup T_{x_{2}},d)\right|\geq 1$, then $\left|E(T_{x_{0}},d)\right|=1$ (if $\left|T_{x_{0}}\right|\geq1$) for $R_{14}\not\subset G$.
Thus, if $\left|E(T_{x_{1}}\cup T_{x_{2}},d)\right|\geq 1$, then  $G\subset H_{104}$ (if $\left|T_{x_{0}x_{1}}\cup T_{x_{0}x_{2}}\right|=3$ and $\left|T_{x_{0}}\right|=1$) or  $G\subset H_{106}$ (if $\left|T_{x_{0}x_{1}}\cup T_{x_{0}x_{2}}\right|=3$ and $\left|T_{x_{0}}\right|=0$) or  $G\subset H_{108}$ (if $\left|T_{x_{0}x_{1}}\cup T_{x_{0}x_{2}}\right|=2$ and $\left|T_{x_{0}}\right|=1$) or  $G\subset H_{11}$ (if $\left|T_{x_{0}x_{1}}\cup T_{x_{0}x_{2}}\right|=2$ and $\left|T_{x_{0}}\right|=0$) by (F\ref{case7-f1}).
Let $\left|E(T_{x_{1}}\cup T_{x_{2}},d)\right|=0$.
If $\left|T_{x_{1}x_{2}}-\{d\}\right|\geq 1$, then $\left|T_{x_{1}}\cup T_{x_{2}}\right|=0$ by (F\ref{case7-f29}).
Thus, if $\left|T_{x_{1}x_{2}}-\{d\}\right|\geq 1$, then $G\subset H_{11}$ by (F\ref{case7-f1}).
Let  $\left|T_{x_{1}x_{2}}-\{d\}\right|=0$.
If $\left|T_{x_{1}}\cup T_{x_{2}}\right|\geq 1$, then $\left|T_{x_{0}x_{1}}\cup T_{x_{0}x_{2}}\right| =2$; otherwise, $G[\{x_{0},\allowbreak x_{1},\allowbreak x_{2},\allowbreak q_{01,2},\allowbreak d,\allowbreak q_{02,1},\allowbreak p_{1,1}\}]\cong R_{30}$, a contradiction.
If  $\left|T_{x_{1}}\cup T_{x_{2}}\right|\geq 1$, then $\left|E(T_{x_{0}},d)\right|=0$; otherwise, $G[\{x_{0},\allowbreak x_{1},\allowbreak x_{3},\allowbreak d,\allowbreak q_{02,1},\allowbreak p_{1,1},\allowbreak p_{0,1}\}]\cong R_{77}$, a contradiction.
Thus, if $\left|T_{x_{1}}\cup T_{x_{2}}\right| \geq 1$, then $G\subset A_{21}$ by (F\ref{case7-f1}).
Let  $\left|T_{x_{1}}\cup T_{x_{2}}\right| =0$.
If $\left|E(T_{x_{0}},d)\right|=0$, then $\left|T_{x_{0}x_{1}}\cup T_{x_{0}x_{2}}\right| =2$ by (F\ref{case7-f31}).
Thus either  $G\subset H_{108}$ (if $\left|E(T_{x_{0}},d)\right|=1$) or  $G\subset A_{21}$ (if $\left|E(T_{x_{0}},d)\right|=0$) by (F\ref{case7-f1}).
Let $\left|T_{x_{0}x_{2}}\right|=0$.

Suppose that $\left|T_{x_{1}x_{2}}-\{d\}\right|\geq 1$.
By (F\ref{case7-f22}) and (F\ref{case7-f29}), $\left|T_{x_{0}}\cup T_{0}\right|=0$.
Then $\left|E(T_{x_{1}},T_{x_{0}x_{1}})\right|=0$ by (F\ref{case7-f4}); otherwise, either $G[\{x_{0},\allowbreak x_{2},\allowbreak x_{3},\allowbreak a,\allowbreak q_{12,2},\allowbreak p_{1,1}\}]\cong R_{1}$ or $G$ contains a subdivision of $K_{3,3}$, a contradiction.
If $\left|E(T_{x_{1}},T_{x_{1}x_{2}}-\{d\})\right|= 1$, then $\left|E(T_{x_{2}},a)\right|=0$ by (F\ref{case7-f4}) and (F\ref{case7-f29}); otherwise, $G[\{x_{0},\allowbreak x_{1},\allowbreak a,\allowbreak q_{12,2},\allowbreak p_{1,1},\allowbreak p_{2,1}\}]\cong R_{10}$, a contradiction.
If $\left|E(T_{x_{1}},T_{x_{1}x_{2}}-\{d\})\right|= 1$, then $\left|T_{x_{0}x_{1}}-\{a\}\right|=0$ by (F\ref{case7-f4}); otherwise, $G[\{x_{0},\allowbreak x_{1},\allowbreak a,\allowbreak q_{01,2},\allowbreak q_{12,2},\allowbreak p_{1,1}\}]\cong R_{16}$, a contradiction.
If $\left|E(T_{x_{1}},T_{x_{1}x_{2}}-\{d\})\right|= 1$ and $\left|T_{x_{1}x_{2}}-\{d\}\right|=2$, then $\left|E(T_{x_{1}},d)\right|=1$ by (F\ref{case7-f4}) and (F\ref{case7-f27}).
Thus, if $\left|E(T_{x_{1}},T_{x_{1}x_{2}}-\{d\})\right|= 1$, then  $G\subset H_{116}$ (if  $\left|T_{x_{1}x_{2}}-\{d\}\right|=2$ and  $\left|T_{x_{2}}\right|=1$) or  $G\subset H_{117}$ (if  $\left|T_{x_{1}x_{2}}-\{d\}\right|=1$ and  $\left|T_{x_{2}}\right|=1$) or  $G\subset H_{134}$ (if  $\left|T_{x_{1}x_{2}}-\{d\}\right|=2$ and  $\left|T_{x_{2}}\right|=0$) or  $G\subset H_{135}$ (if  $\left|T_{x_{1}x_{2}}-\{d\}\right|=1$ and  $\left|T_{x_{2}}\right|=0$) by (F\ref{case7-f1}), (F\ref{case7-f4}), (F\ref{case7-f8}), (F\ref{case7-f13}), (F\ref{case7-f29}) and (F\ref{case7-f37}).
Let $\left|E(T_{x_{1}},T_{x_{1}x_{2}}-\{d\})\right|= 0$.
If $\left|T_{x_{1}x_{2}}-\{d\}\right|=2$, then either   $G\subset H_{114}$ (if $\left|T_{x_{2}}\right|=1$) or   $G\subset H_{134}$ (if $\left|T_{x_{2}}\right|=0$) by (F\ref{case7-f1}),  (F\ref{case7-f8}), (F\ref{case7-f13}), (F\ref{case7-f27}), (F\ref{case7-f29}), (F\ref{case7-f34}) and (F\ref{case7-f37}).
Let $\left|T_{x_{1}x_{2}}-\{d\}\right|=1$.
Then no vertex in $T_{x_{2}}$ is adjacent to all vertices of $\{a,d\}$ by  (F\ref{case7-f13}); otherwise, $G[\{x_{0},\allowbreak x_{1},\allowbreak x_{3},\allowbreak a,\allowbreak d,\allowbreak q_{12,2},\allowbreak p_{2,1}\}]\cong R_{25}$, a contradiction.
If $\left|T_{x_{2}}\right|\geq 1$, then at most one vertex in $T_{x_{1}}$ is not adjacent to $\{a,d\}$ by (F\ref{case7-f37}); otherwise, $G[\{x_{0},\allowbreak x_{1},\allowbreak x_{3},\allowbreak a,\allowbreak d,\allowbreak q_{12,2},\allowbreak p_{1,1},\allowbreak p_{1,2},\allowbreak p_{2,1}\}]\cong R_{47}$, a contradiction.
Thus, $G\subset H_{115}$ (if  $\left|E(T_{x_{2}},a)\right|\leq 1$ and  $\left|E(T_{x_{2}},d)\right|=1$) or  $G\subset H_{131}$ (if  $\left|E(T_{x_{2}},a)\right|=1$ and  $\left|E(T_{x_{2}},d)\right|=0$) or  $G\subset H_{133}$ (if  $\left|E(T_{x_{2}},a)\right|=0$ and  $\left|E(T_{x_{2}},d)\right|=0$) by (F\ref{case7-f1}),  (F\ref{case7-f8}), (F\ref{case7-f13}), (F\ref{case7-f27}), (F\ref{case7-f29}) and (F\ref{case7-f31}).
Let $\left|T_{x_{1}x_{2}}-\{d\}\right|=0$.

Suppose that $\left|T_{x_{0}x_{1}}-\{a\}\right|\geq 1$.
Then $\left|E(T_{x_{2}},T_{x_{0}x_{1}}\cup \{d\}-\{a\})\right|=0$ and $\left|E(T_{x_{2}},a)\right|=\left|T_{x_{2}}\right|$ by (F\ref{case7-f13}) and (F\ref{case7-f31}); otherwise, $G[\mathcal{K}\cup\{a,\allowbreak q_{01,2},\allowbreak p_{2,1}\}]\cong R_{28}$, a contradiction.
Thus, $\left|T_{x_{2}}\right|\leq 1$; otherwise, $G[\mathcal{K}\cup\{q_{01,2},\allowbreak p_{2,1},\allowbreak p_{2,2}\}]\cong R_{25}$, a contradiction.
If $\left|T_{x_{0}x_{1}}-\{a\}\right|=2$, then  $\left|E(T_{x_{1}},d)\right|=0$ by (F\ref{case7-f4}); otherwise, $G[\{x_{0},\allowbreak x_{1},\allowbreak q_{01,2},\allowbreak q_{01,3},\allowbreak d,\allowbreak p_{1,1}\}]\cong R_{16}$, a contradiction.
Thus, if $\left|T_{x_{0}x_{1}}-\{a\}\right|=2$, then either  $G\subset H_{131}$ (if  $\left|T_{x_{0}}\right|=1$) or  $G\subset H_{133}$ (if  $\left|T_{x_{0}}\right|=0$) by (F\ref{case7-f1}), (F\ref{case7-f4}), (F\ref{case7-f22}), (F\ref{case7-f31}) and (F\ref{case7-f37}).
Let $\left|T_{x_{0}x_{1}}-\{a\}\right|=1$.
By (F\ref{case7-f31}) and   $G$ contains no subdivision of $K_{3,3}$, $\left|T_{x_{0}}\cup T_{x_{2}}\right|\leq 1$.
If $\left|T_{x_{0}}\cup T_{x_{2}}\right|=1$, then at most one vertex in $T_{x_{1}}$ is  not adjacent to $\{a,d\}$ by (F\ref{case7-f37}); otherwise, $G[\{x_{0},\allowbreak x_{1},\allowbreak x_{3},\allowbreak a,\allowbreak d,\allowbreak p_{1,1},\allowbreak p_{1,2},\allowbreak p_{0,1}(p_{2,1})\}]\cong R_{38}$, a contradiction.
Thus, either $G\subset H_{130}$ (if  $\left|T_{x_{0}}\cup T_{x_{2}}\right|=1$) or $G\subset H_{136}$ (if  $\left|T_{x_{0}}\cup T_{x_{2}}\right|=0$) by (F\ref{case7-f1}), (F\ref{case7-f4}), (F\ref{case7-f22}) and (F\ref{case7-f31}).
Let $\left|T_{x_{0}x_{1}}-\{a\}\right|=0$.

Suppose that $\left|T_{x_{0}}\right|=1$.
Then $\left|T_{0}\right|=0$ by (F\ref{case7-f22}) and $d(G)< 4$; otherwise, $G[\{x_{1},\allowbreak x_{2},\allowbreak x_{3},\allowbreak a,\allowbreak d,\allowbreak p_{0,1},\allowbreak t_{1}\}]\cong R_{25}$ ($E(\{a,\allowbreak d,\allowbreak p_{0,1},\allowbreak t_{1}\})=\{ad,\allowbreak dp_{0,1},\allowbreak dt_{1}\}$), a contradiction.
If $\left|E(T_{x_{0}},d)\right|=1$, then  $\left|E(T_{x_{2}},a)\right|=0$ and $\left|E(T_{x_{2}},d)\right|=\left|T_{x_{2}}\right|$; otherwise, either $G[\{x_{0},\allowbreak x_{2},\allowbreak x_{3},\allowbreak a,\allowbreak d,\allowbreak p_{0,1},\allowbreak p_{2,1}\}]\cong R_{77}$ or $G$ contains a subdivision of $K_{3,3}$, a contradiction.
Thus, if $\left|E(T_{x_{0}},d)\right|=1$, then either  $G\subset H_{115}$ (if  $\left|T_{x_{2}}\right|=1$) or  $G\subset H_{131}$ (if  $\left|T_{x_{2}}\right|=0$) by  (F\ref{case7-f1}), (F\ref{case7-f4}), (F\ref{case7-f8}), (F\ref{case7-f37}) and (F\ref{case7-f39}).
Let  $\left|E(T_{x_{0}},d)\right|=0$.
Then $\left|E(T_{x_{1}}\cup T_{x_{2}},d)\right|=0$ by (F\ref{case7-f4}); otherwise, $G[\{x_{0},\allowbreak x_{1},\allowbreak x_{2},\allowbreak d,\allowbreak p_{0,1},\allowbreak p_{1,1}(p_{2,1})\}]\cong R_{14}$, a contradiction.
Then at most three vertices in $T_{x_{1}}$ are not adjacent to $\{a,d\}$; otherwise, $G[\{x_{0},\allowbreak x_{1},\allowbreak x_{2},\allowbreak a,\allowbreak p_{0,1},\allowbreak p_{1,1},\ldots,p_{1,4}\}]\cong R_{47}$, a contradiction.
Thus,  $G\subset A_{25}$ (if  $\left|E(T_{x_{2}},a)\right|=3$ and $\left|E(T_{x_{1}},a)\right|=1$)  or  $G\subset A_{28}$ (if  $\left|E(T_{x_{2}},a)\right|=2$ and $\left|E(T_{x_{1}},a)\right|=1$)  or  $G\subset A_{31}$ (if  $\left|E(T_{x_{2}},a)\right|=1$ and $\left|E(T_{x_{1}},a)\right|=1$) or  $G\subset A_{34}$ (if  $\left|E(T_{x_{2}},a)\right|=0$ and $\left|E(T_{x_{1}},a)\right|=1$) or  $G\subset A_{24}$ (if  $\left|E(T_{x_{2}},a)\right|=3$ and $\left|E(T_{x_{1}},a)\right|=0$) or  $G\subset A_{27}$ (if  $\left|E(T_{x_{2}},a)\right|=2$ and $\left|E(T_{x_{1}},a)\right|=0$) or  $G\subset A_{30}$ (if  $\left|E(T_{x_{2}},a)\right|=1$ and $\left|E(T_{x_{1}},a)\right|=0$) or  $G\subset H_{132}$ (if  $\left|E(T_{x_{2}},a)\right|=0$ and $\left|E(T_{x_{1}},a)\right|=0$) by  (F\ref{case7-f1}), (F\ref{case7-f4}), (F\ref{case7-f8}), (F\ref{case7-f39}) and (F\ref{case7-f41}).
Let $\left|T_{x_{0}}\right|=0$ and $\left|T_{0}\right|=0$ by (F\ref{case7-f22}) and symmetry.

If $p_{2,1}\in T_{x_{2}}$ and $\left|E(p_{2,1},\{a,d\})\right|=2$, then its conclusion is covered by Case $1$, by taking $x_{1},x_{2},x_{3},d$ as the vertices of a new $K_{4}-e$.
If $p_{2,2}\in T_{x_{2}}$ and $\left|E(p_{2,2},d)\right|=1$, then its conclusion is covered by Case $1$, by taking $x_{1},x_{2},p_{2,2},d$ as the vertices of a new $K_{4}-e$.
If $\left|E(T_{x_{1}},d)\right|\geq 1$, then no vertex in $T_{x_{2}}$ is not adjacent to $\{a,d\}$ by (F\ref{case7-f4}); otherwise, $G[\{x_{1},\allowbreak x_{2},\allowbreak x_{3},\allowbreak a,\allowbreak d,\allowbreak p_{1,1},\allowbreak p_{2,1}\}]\cong R_{28}$, a contradiction.
If  $\left|E(T_{x_{1}},d)\right|\geq 1$ and  $\left|E(T_{x_{2}},a)\right|= 1$, then at most one vertex in $T_{x_{1}}$ is not adjacent to $\{a,d\}$ by (F\ref{case7-f4}); otherwise, $G[\{x_{0},\allowbreak x_{1},\allowbreak x_{2},\allowbreak a,\allowbreak p_{1,1},\allowbreak p_{1,2},\allowbreak p_{1,3},\allowbreak p_{2,1}\}]\cong R_{38}$, a contradiction.
Thus, either  $G\subset H_{133}$ (if $\left|S_{4}\right|=0$) or  $G\subset H_{130}$ (if $\left|S_{4}\right|=1$) by  (F\ref{case7-f1}), (F\ref{case7-f4}), (F\ref{case7-f8}) and (F\ref{case7-f42}).
Let $\left|E(T_{x_{1}}\cup T_{x_{2}},d)\right|=0$.
If  at least four vertices in $T_{x_{1}}$ are not adjacent to $\{a,d\}$,  then $\left|T_{x_{2}}\right|=0$; otherwise, $G[\mathcal{K}\cup\{p_{1,1},\ldots,p_{1,4},\allowbreak p_{2,1}\}]\cong R_{47}$, a contradiction.
Thus,   $G\subset H_{133}$ (if  $\left|S_{1}\right|\geq 4$) or   $G\subset A_{34}$ (if  $\left|S_{1}\right|=3$ and $\left|E(T_{x_{1}},a)\right|=1$) or $G\subset H_{132}$ (if  $\left|S_{1}\right|=3$ and $\left|E(T_{x_{1}},a)\right|=0$) or  $G\subset A_{31}$ (if  $\left|S_{1}\right|=2$ and $\left|E(T_{x_{1}},a)\right|=1$) or  $G\subset A_{30}$ (if  $\left|S_{1}\right|=2$ and $\left|E(T_{x_{1}},a)\right|=0$) or  $G\subset A_{28}$ (if  $\left|S_{1}\right|=1$ and $\left|E(T_{x_{1}},a)\right|=1$) or  $G\subset A_{27}$ (if  $\left|S_{1}\right|=1$ and $\left|E(T_{x_{1}},a)\right|=0$) or  $G\subset A_{25}$ (if  $\left|S_{1}\right|=0$ and $\left|E(T_{x_{1}},a)\right|=1$) or  $G\subset A_{24}$ (if  $\left|S_{1}\right|=0$ and $\left|E(T_{x_{1}},a)\right|=0$) by  (F\ref{case7-f1}), (F\ref{case7-f4}), (F\ref{case7-f8}), (F\ref{case7-f39}) and (F\ref{case7-f41}).
\end{proof}

\begin{lem}\label{lem-44}
If $G \in \mathcal{G}_3$, $E(T_{1})=\varnothing$, $E(T_{2})-E(T_{x_{0}x_{3}})=\varnothing$, $E(T_{3})\neq\varnothing$, and $G$ contains $K_{4}-e$ as an induced subgraph, then $G$ is an induced subgraph of $A_{15}$ or $A_{16}$.
\end{lem}

\begin{proof}
Let $K_{4}-e$ be an induced subgraph of $G$ with $\mathcal{K}=\{x_{0},x_{1},x_{2},x_{3}\}$ and $x_{0}\nsim x_{3}$.
By  Lemma \ref{lem-13}\ref{lem-13:item-7}, let  $u\in T_{x_{0}x_{1}x_{3}}$, $v\in T_{x_{0}x_{2}x_{3}}$ and $u\sim v$.
Then no vertex in $T_{x_{0}}\cup T_{x_{3}}$ is adjacent to all vertices of $\{u,v\}$ for $\omega\leq 3$.
If $p_{0,1}\in T_{x_{0}}$ and $p_{0,1}\sim u$, then its conclusion is covered by Case $1$ of Lemma \ref{lem-43} by taking $\{x_{0},x_{1},u,p_{0,1}\}$ as the vertices of a new $K_{4}-e$.
Since $G$ contains no subdivision of $K_{3,3}$ and $K_{4}\not\subset G$, $\left|E(T_{x_{1}},v)\right|=\left|E(T_{x_{2}},u)\right|=\left|E(T_{2},\{u,v\})\right|=\left|T_{x_{0}x_{3}}\right|=0$.
If $p_{1,1}\in T_{x_{1}}$ and $p_{1,1}\sim u$, then its conclusion is covered by Case $4$ of Lemma \ref{lem-43} by taking $\{x_{1},x_{3},u,p_{1,1}\}$ as the vertices of a new $K_{4}-e$.
If $q_{01,1}\in T_{x_{0}x_{1}}$, then its conclusion is covered by Case $4$ of Lemma \ref{lem-43} by taking $\{x_{0},x_{1},x_{2},q_{01,1}\}$ as the vertices of a new $K_{4}-e$.
If $q_{12,1}\in T_{x_{1}x_{2}}$, then its conclusion is covered by Case $1$ of Lemma \ref{lem-43} by taking $\{x_{0},x_{1},x_{2},q_{12,1}\}$ as the vertices of a new $K_{4}-e$.
By  symmetry, let $\left|E(T_{1},\{u,v\})\right|=0$ and $\left|T_{2}\right|=0$.
Then $\left|T_{x_{i}}\right| \leq 1$ (i=0,1,2,3); otherwise, $G[\mathcal{K}\cup\{u,\allowbreak p_{0,1},\allowbreak p_{0,2}\}]\cong R_{27}$ (if $i=0$ or $3$) or $G[\{x_{0},\allowbreak x_{1},\allowbreak x_{2},\allowbreak u,\allowbreak v,\allowbreak p_{1,1},\allowbreak p_{1,2}\}]\cong R_{27}$ (if $i=1$ or $2$), a contradiction.
If $t_0\in T_{0}$ and $\left|E(t_0,\{u,v\})\right|=2$, then $\left|T_{x_{1}}\cup T_{x_{2}}\right| =0$ by Lemma \ref{lem-13}\ref{lem-13:item-8}; otherwise, $G[\{x_{0},\allowbreak x_{1},\allowbreak u,\allowbreak v,\allowbreak t_0,\allowbreak p_{1,1}\}]\cong R_{14}$, a contradiction.
Thus, either $G\subset A_{16}$ (if $t_0\in T_{0}$ and $\left|E(t_0,\{u,v\})\right|=2$) or $G\subset A_{15}$ (the remaining case)  by Lemma \ref{lem-13}\ref{lem-13:item-8}.
\end{proof}

\begin{lem}\label{lem-45}
If $G \in \mathcal{G}_3$, $E(T-T_{0})=E(T_{x_{0}x_{3}})\cup E(T_{1},T_{2})\cup E(T_{1},T_{3})\cup E(T_{2},T_{3})$, $\left|E(T_{1},T_{2})\cup E(T_{1},T_{3})\cup E(T_{2},T_{3})\right| \geq1$, and $G$ contains $K_{4}-e$ as an induced subgraph, then $G$ is an induced subgraph of a graph in $\{H_{4},\allowbreak H_{10},\allowbreak H_{11},\allowbreak H_{71},\allowbreak H_{73},\allowbreak \ldots,\allowbreak H_{82},\allowbreak H_{84},\allowbreak \ldots,\allowbreak H_{88},\allowbreak H_{90},\allowbreak \ldots,\allowbreak H_{102},\allowbreak H_{106},\allowbreak H_{109},\allowbreak H_{111},\allowbreak H_{112},\allowbreak  H_{114},\allowbreak \ldots,\allowbreak H_{117},\allowbreak H_{120},\allowbreak H_{131},\allowbreak H_{132},\allowbreak H_{134},\allowbreak H_{135},\allowbreak H_{137},\allowbreak \ldots,\allowbreak H_{153},\allowbreak A_{3},\allowbreak A_{6},\allowbreak A_{18},\allowbreak \ldots,\allowbreak A_{24} ,\allowbreak A_{26},\allowbreak A_{29},\allowbreak  A_{30},\allowbreak A_{32},\allowbreak A_{37}\ldots,\allowbreak A_{46}\}$.
\end{lem}
\begin{proof}
Let $K_{4}-e$ be an induced subgraph of $G$ with $\mathcal{K}=\{x_{0},x_{1},x_{2},x_{3}\}$ and $x_{0}\nsim x_{3}$.
We distinguish the following three cases.

\noindent$\textbf{Case 1.}$ Suppose that $\left|E(T_{1},T_{2})\right|\geq 1$.

By Lemma \ref{lem-13}\ref{lem-13:item-2} and symmetry, we distinguish the following six subcases.

\noindent$\textbf{Subcase 1.1}$ Suppose that $\left|E(T_{x_{0}},T_{x_{1}x_{2}})\right|\geq 1$.

Let  $p\in T_{x_{0}}$, $d\in T_{x_{1}x_{2}}$ and $p\sim d$.

\begin{fact}\label{case9.1-f1}
$G[E(T_{x_{1}x_{2}},T_{x_{0}})]\subset aK_{1}\cup bK_{2}$ ($b\geq1$ and $a+b\leq 2$).
\end{fact}
\begin{proof}
By Lemma \ref{lem-13}\ref{lem-13:item-1} and Lemma \ref{lem-13}\ref{lem-13:item-4}, we get $K_{3}\not\subset G[E(T_{x_{1}x_{2}},T_{x_{0}})]$.
If $P_{3}\subset G[E(T_{x_{1}x_{2}},T_{x_{0}})]$, then either $G[\mathcal{K}\cup\{p,\allowbreak p_{0,2},\allowbreak d\}]\cong R_{29}$ ($E(\{p,\allowbreak p_{0,2},\allowbreak d\})=\{pd,\allowbreak p_{0,2}d\}$) or $G$ contains a subdivision of $K_{3,3}$, a contradiction.
Thus, $G[E(T_{x_{1}x_{2}},T_{x_{0}})]\subset aK_{1}\cup bK_{2}$.
If $\left|T_{x_{1}x_{2}}-\{d\}\right|\geq 2$, then  $G[\mathcal{K}\cup\{p,\allowbreak d,\allowbreak q_{12,2},\allowbreak q_{12,3}\}]\cong R_{41}$, a contradiction.
If $a+b\geq 3$, then $G[\mathcal{K}\cup\{p,\allowbreak p_{0,2},\allowbreak q_{12,2}\}]\cong R_{24}$ by Lemma \ref{lem-13}\ref{lem-13:item-1}, a contradiction.
\end{proof}

\begin{fact}\label{case9.1-f3}
$\left|T_{x_{3}}\right|\leq 1$, $\left|E(T_{x_{3}},d)\right|=1$ (if $\left|T_{x_{3}}\right|= 1$) and $\left|E(T_{x_{1}x_{2}}-\{d\},T_{x_{3}})\right|=0$.
\end{fact}
\begin{proof}
By (F\ref{case2-f1}), $\left|T_{x_{3}}\right|\leq 1$.
If $\left|E(T_{x_{3}},d)\right|=0$, then $R_{2}\subset G$, a contradiction.
If $\left|E(T_{x_{1}x_{2}}-\{d\},T_{x_{3}})\right|\geq 1$, then  $G$ contains a subdivision of $K_{3,3}$, a contradiction.
\end{proof}

\begin{fact}\label{case9.1-f4}
At most two vertices in $T_{x_{1}}$ (or $T_{x_{2}}$) are not adjacent to $d$; at most two vertices in $T_{x_{1}}$ (or $T_{x_{2}}$) are adjacent to $d$; if  $S_{1}\subset T_{x_{1}}$ (or $S_{1}\subset T_{x_{2}}$), $S_{2}\subset T_{x_{2}}$ (or $S_{2}\subset T_{x_{1}}$), $\left|E(S_{1},\{p,d\})\right|=0$, and $\left|E(S_{2},d)\right|=\left|S_{2}\right|$, then $\left|S_{1}\cup S_{2}\right| \leq 2$.
\end{fact}
\begin{proof}
If at least three vertices in $T_{x_{1}}$ (or $T_{x_{2}}$) are not adjacent to $d$, then $G[\{x_{0},\allowbreak x_{1},\allowbreak x_{2},\allowbreak p,\allowbreak d,\allowbreak p_{1,1},\allowbreak p_{1,2},\allowbreak p_{1,3}\}]\cong R_{38}$, a contradiction.
If at least three vertices in $T_{x_{1}}$ (or $T_{x_{2}}$) are adjacent to $d$, then $G[\{x_{0},\allowbreak x_{2},\allowbreak x_{3},\allowbreak d,\allowbreak p_{1,1},\allowbreak p_{1,2},\allowbreak p_{1,3}\}]\cong R_{21}$, a contradiction.
If $\left|S_{1}\cup S_{2}\right| \geq 3$, then either $G[\{x_{0},\allowbreak x_{1},\allowbreak x_{3},\allowbreak d,\allowbreak p_{1,1},\allowbreak p_{2,1},\allowbreak p_{2,2}\}]\cong R_{21}$ ($E(\{d,\allowbreak p_{1,1},\allowbreak p_{2,1},\allowbreak p_{2,2}\})=\{dp_{2,1},\allowbreak dp_{2,2}\}$) or $G[\{x_{1},\allowbreak x_{3},\allowbreak p,\allowbreak d,\allowbreak p_{1,1},\allowbreak p_{1,2},\allowbreak p_{2,1}\}]\cong R_{21}$ ($E(\{p,\allowbreak d,\allowbreak p_{1,1},\allowbreak p_{1,2},\allowbreak p_{2,1}\})=\{pd,\allowbreak dp_{2,1}\}$), a contradiction.
\end{proof}

\begin{fact}\label{case9.1-f6}
$\left|T_{x_{0}x_{1}}\right|\leq 2$, $\left|T_{x_{0}x_{2}}\right|\leq 2$ and $\left|T_{x_{0}x_{3}}\right|\leq 1$.
\end{fact}
\begin{proof}
If $\left|T_{x_{0}x_{1}}\right|\geq 3$ or $\left|T_{x_{0}x_{2}}\right|\geq 3$, then $G[\{x_{0},\allowbreak x_{2},\allowbreak x_{3},\allowbreak d,\allowbreak q_{01,1},\allowbreak q_{01,2},\allowbreak q_{01,3}\}]\cong R_{21}$ by Lemma \ref{lem-13}\ref{lem-13:item-2}, a contradiction.
If $\left|T_{x_{0}x_{3}}\right|\geq 2$, then either $G[\{x_{0},\allowbreak x_{3},\allowbreak p,\allowbreak d,\allowbreak q_{03,1},\allowbreak q_{03,2}\}]\cong R_{2}$ ($q_{03,1}\nsim q_{03,2}$) or $G[\{x_{0},\allowbreak x_{1},\allowbreak x_{2},\allowbreak d,\allowbreak q_{03,1},\allowbreak q_{03,2}\}]\cong R_{13}$ ($q_{03,1}\sim q_{03,2}$) by Lemma \ref{lem-13}\ref{lem-13:item-2}, a contradiction.
\end{proof}

\begin{fact}\label{case9.1-f7}
$\left|T_{x_{1}x_{3}}\cup T_{x_{2}x_{3}}\cup T_{3}\right|=0$.
\end{fact}
\begin{proof}
If $\left|T_{x_{1}x_{3}}\cup T_{x_{2}x_{3}}\right|\geq 1$, then either $G[\mathcal{K}\cup\{p,\allowbreak q_{13,1}\}]\cong R_{14}$ or $G$ contains a subdivision of $K_{3,3}$, a contradiction.
If $\left|T_{3}\right|\geq 1$, then either $G[\{x_{0},\allowbreak x_{1},\allowbreak x_{3},\allowbreak p,\allowbreak d,\allowbreak w_{013,1}\}]\cong R_{12}$ or $G$ contains a subdivision of $K_{3,3}$, a contradiction.
\end{proof}

\begin{fact}\label{case9.1-f10}
$E(T_{0},T)=E(T_{0},d)$ and $\left|T_{0}\right|\leq 1$.
\end{fact}
\begin{proof}
If $\left|E(T_{0},T_{x_{1}x_{2}}\cup T_{x_{0}x_{3}}-\{d\})\right|\geq 1$, then $G[\{x_{0},\allowbreak x_{1},\allowbreak p,\allowbreak d,\allowbreak q_{12,2},\allowbreak t_{1}\}]\cong R_{2}$ by Lemma \ref{lem-13}\ref{lem-13:item-8}, a contradiction.
Thus, by (F\ref{case9.1-f7}) and Lemma \ref{lem-13}\ref{lem-13:item-8}, $E(T_{0},T)=E(T_{0},d)$ and $\left|T_{0}\right|\leq 1$.
\end{proof}

Suppose that $\left|E(T_{x_{0}},T_{x_{1}x_{2}})\right|\geq 2$.
By (F\ref{case9.1-f1}),  $\left|E(T_{x_{0}},T_{x_{1}x_{2}})\right|= \left|T_{x_{1}x_{2}}\right|=\left|T_{x_{0}}\right|=2$ and $G[E(T_{x_{1}x_{2}},T_{x_{0}})]= 2K_{2}$.
Then $\left|T_{x_{3}}\right|=0$ by (F\ref{case9.1-f3}) and   $G$ contains no subdivision of $K_{3,3}$.
Then $\left|T_{x_{1}}\cup T_{x_{2}}\right|=0$; otherwise, $R_{14}\subset G$ (if $\left|E(T_{x_{1}},T_{x_{1}x_{2}})\right|= 1$) or $R_{22}\subset G$ (if $\left|E(T_{x_{1}},T_{x_{1}x_{2}})\right|= 0$) or $G$ contains a subdivision of $K_{3,3}$ (if $\left|E(T_{x_{1}},T_{x_{1}x_{2}})\right|= 2$), a contradiction.
Then $\left|T_{x_{0}x_{3}}\right|=0$ for $G$ contains no subdivision of $K_{3,3}$.
Then $\left|T_{x_{0}x_{1}}\cup T_{x_{0}x_{2}}\right|=0$ by Lemma \ref{lem-13}\ref{lem-13:item-2}; otherwise, $G[\mathcal{K}\cup\{p,\allowbreak q_{12,2},\allowbreak q_{01,1}(q_{02,1})\}]\cong R_{28}$, a contradiction.
Then $\left|T_{0}\right|=0$ by (F\ref{case9.1-f10}); otherwise, $G[\{x_{0},\allowbreak x_{1},\allowbreak p_{0,2},\allowbreak d,\allowbreak q_{12,2},\allowbreak t_{1}\}]\cong R_{2}$, a contradiction.
Thus, $G\subset H_{138}$ by (F\ref{case9.1-f7}).
Let $\left|E(T_{x_{0}},T_{x_{1}x_{2}})\right|=1$.

Suppose that $\left|T_{0}\right|\geq 1$.
By (F\ref{case9.1-f10}), $E(T_{0},T)=E(T_{0},d)$ and $\left|T_{0}\right|= 1$.
Then $\left|T_{x_{1}x_{2}}-\{d\}\right|=0$; otherwise, $G[\{x_{1},\allowbreak x_{2},\allowbreak x_{3},\allowbreak p,\allowbreak d,\allowbreak q_{12,2},\allowbreak t_{1}\}]\cong R_{24}$, a contradiction.
Then $\left|T_{x_{1}}\cup T_{x_{2}}\right|=0$; otherwise, either $G[\{x_{1},\allowbreak x_{2},\allowbreak x_{3},\allowbreak p,\allowbreak d,\allowbreak p_{1,1},\allowbreak t_{1}\}]\cong R_{22}$ or $G[\mathcal{K}\cup\{d,\allowbreak p_{1,1},\allowbreak t_{1}\}]\cong R_{28}$, a contradiction.
Then $\left|T_{x_{0}x_{1}}\cup T_{x_{0}x_{2}}\right|=0$ by Lemma \ref{lem-13}\ref{lem-13:item-2}; otherwise, $G[\{x_{0},\allowbreak x_{1},\allowbreak x_{2},\allowbreak d,\allowbreak q_{01,1},\allowbreak t_{1}\}]\cong R_{14}$, a contradiction.
Then $\left|T_{x_{0}x_{3}}\right|=0$ by Lemma \ref{lem-13}\ref{lem-13:item-2}; otherwise, $G[\{x_{0},\allowbreak x_{1},\allowbreak x_{3},\allowbreak d,\allowbreak q_{03,1},\allowbreak t_{1}\}]\cong R_{2}$, a contradiction.
Then $\left|T_{x_{0}}-\{p\}\right|=0$ by (F\ref{case9.1-f1}) and  $d(G)\leq 3$.
Thus, $G\subset H_{137}$ by (F\ref{case9.1-f7}).
By symmetry and a similar argument, if $\left|T_{x_{0}}-\{p\}\right|\geq1$, then we can obtain the desired conclusion.
Let $\left|T_{0}\cup T_{x_{0}}-\{p\}\right|=0$.

Suppose that $\left|T_{x_{1}x_{2}}-\{d\}\right|\geq 1$.
By (F\ref{case9.1-f1}), $\left|T_{x_{1}x_{2}}-\{d\}\right|= 1$.
Then $\left|T_{x_{1}}\cup T_{x_{2}}\right|=0$; otherwise, $R_{14}\subset G$ (if $\left|E(T_{x_{1}},T_{x_{1}x_{2}}-\{d\})\right|= 1$) or $R_{40}\subset G$ (if $\left|E(T_{x_{1}},T_{x_{1}x_{2}})\right|= 0$) or $R_{28}\subset G$ (the remaining case) , a contradiction.
Then $\left|T_{x_{0}x_{1}}\cup T_{x_{0}x_{2}}\right|=0$; otherwise, $G[\mathcal{K}\cup\{p,\allowbreak q_{12,2},\allowbreak q_{01,1}\}]\cong R_{28}$, a contradiction.
Then $\left|T_{x_{3}}\cup T_{x_{0}x_{3}}\right|\leq 1$ by (F\ref{case9.1-f3}), (F\ref{case9.1-f6}) and   $G$ contains no subdivision of $K_{3,3}$.
Thus, $G\subset H_{138}$ by (F\ref{case9.1-f7}).
Let $\left|T_{x_{1}x_{2}}-\{d\}\right|=0$.

Suppose that $\left|T_{x_{0}x_{1}}\cup T_{x_{0}x_{2}}\right|\geq 1$.
Without loss of generality, let $\left|T_{x_{0}x_{1}}\right|\geq 1$ and $q_{01,1}\in T_{x_{0}x_{1}}$.
Then $\left|T_{x_{3}}\right|=0$ by (F\ref{case9.1-f3}); otherwise, either $R_{14}\subset G$ (if $\left|E(T_{x_{3}},T_{x_{0}x_{1}})\right|=0$) or $G$ contains a subdivision of $K_{3,3}$ (if $\left|E(T_{x_{3}},T_{x_{0}x_{1}})\right|\geq 1$), a contradiction.
Let $S_{3}=T_{x_{2}}-(T_{x_{2}}\cap N(\{p,d\}))$ and $S_{4}=T_{x_{1}}\cap N(d)$.
By (F\ref{case9.1-f4}), $\left|S_{3}\cup S_{4}\right| \leq 2$.
Then $\left|E(T_{x_{2}},d)\right|=0$, and $\left|E(T_{x_{1}},d)\right|=0$ (if  $\left|T_{x_{0}x_{2}}\right|\geq 1$); otherwise, either $G[\{x_{0},\allowbreak x_{1},\allowbreak x_{2},\allowbreak d,\allowbreak q_{01,1},\allowbreak p_{2,1}\}]\cong R_{18}$ or $G$ contains a subdivision of $K_{3,3}$, a contradiction.
Since  $R_{12}\not\subset G$, $\left|E(T_{x_{1}}-S_{4},T_{x_{0}x_{1}})\right|=0$.
Since $G$ contains no subdivision of $K_{3,3}$, $\left|E(S_{4},T_{x_{0}x_{1}})\right|\leq 1$.
If $p_{1,1}\in S_{4}$ and $\left|E(p_{1,1},T_{x_{0}x_{1}})\right|=1$, then its conclusion is covered by Case $2$ of Lemma \ref{lem-43} by taking $\{x_{0},x_{1},x_{2},d\}$ as the vertices of a new $K_{4}-e$ ($G\subset H_{87}$).
Let $\left|E(T_{x_{1}},T_{x_{0}x_{1}})\right|=0$.
Then $\left|S_{4}\cup T_{x_{0}x_{1}}\right| \leq 2$ by (F\ref{case9.1-f4}) and (F\ref{case9.1-f6}); otherwise, $G[\{x_{0},\allowbreak x_{1},\allowbreak d,\allowbreak q_{01,1},\allowbreak q_{01,2},\allowbreak p_{1,1}\}]\cong R_{16}$ ($E(\{d,\allowbreak q_{01,1},\allowbreak q_{01,2},\allowbreak p_{1,1}\})=\{dp_{1,1}\}$) or $G[\{x_{0},\allowbreak x_{1},\allowbreak d,\allowbreak q_{01,1},\allowbreak p_{1,1},\allowbreak p_{1,2}\}]\cong R_{16}$ ($E(\{d,\allowbreak q_{01,1},\allowbreak p_{1,1},\allowbreak p_{1,2}\})=\{dp_{1,1},\allowbreak dp_{1,2}\}$), a contradiction.
If  $\left|S_{4}\right| =1$, then  $\left|E(T_{x_{2}},T_{x_{0}x_{1}})\right|=0$; otherwise, $G[\{x_{0},\allowbreak x_{1},\allowbreak d,\allowbreak q_{01,1},\allowbreak p_{1,1},\allowbreak p_{2,1}\}]\cong R_{10}$, a contradiction.
Thus,  if  $\left|S_{4}\right| =1$, then $G\subset H_{142}$ by (F\ref{case9.1-f7}).
Let $\left|S_{4}\right| =0$.
Then $\left|E(T_{x_{2}},T_{x_{0}x_{1}})\right|+\left|E(T_{x_{1}},T_{x_{0}x_{2}})\right| \leq 1$; otherwise, either $G[\{x_{1},\allowbreak p,\allowbreak d,\allowbreak q_{01,1},\allowbreak p_{2,1},\allowbreak p_{2,2}\}]\cong R_{1}$ or $G$ contains a subdivision of $K_{3,3}$, a contradiction.
Then $\left|T_{x_{2}}\cup T_{x_{0}x_{1}}\right|\leq 3$ and $\left|T_{x_{1}}\cup T_{x_{0}x_{2}}\right|\leq 3$; otherwise, $G[\mathcal{K}\cup\{q_{01,2}(q_{01,1},q_{02,1},q_{02,2}),\allowbreak p_{2,1}(p_{1,1},p_{1,2}),\allowbreak p_{2,2}\}]\cong R_{25}$, a contradiction.
If  $\left|T_{x_{2}}\cup T_{x_{0}x_{1}}\right|=3$ (or $\left|T_{x_{1}}\cup T_{x_{0}x_{2}}\right|=3$), then $\left|E(T_{x_{2}},T_{x_{0}x_{1}})\right|= 1$ (or  $\left|E(T_{x_{1}},T_{x_{0}x_{2}})\right|\leq 1$); otherwise, $G[\{x_{0},\allowbreak x_{2},\allowbreak x_{3},\allowbreak p,\allowbreak q_{01,1},\allowbreak p_{2,1},\allowbreak p_{2,2}\}]\cong R_{21}$, a contradiction.
Thus,  $G\subset H_{73}$ (if $\left|T_{x_{1}}\right| =0$ and $\left|T_{x_{2}}\right| =2$) or $G\subset H_{74}$ (if $\left|T_{x_{1}}\cup T_{x_{2}}\right| =3$) or $G\subset H_{75}$ (if $\left|T_{x_{1}}\cup T_{x_{2}}\right| =4$) or $G\subset H_{76}$ (if $\left|T_{x_{1}}\cup T_{x_{2}}\right| \leq 1$) or $G\subset H_{77}$ (if $\left|T_{x_{1}}\right| =1$ and $\left|T_{x_{2}}\right| =1$) by (F\ref{case9.1-f4}) and (F\ref{case9.1-f7}).
Let  $\left|T_{x_{0}x_{1}}\cup T_{x_{0}x_{2}}\right|=0$.

Then $\left|T_{x_{3}}\cup T_{x_{0}x_{3}}\right|\leq 1$ by (F\ref{case9.1-f3}), (F\ref{case9.1-f6}) and   $G$ contains no subdivision of $K_{5}$.
If $\left|T_{x_{0}x_{3}}\right|=1$, then $\left|E(T_{x_{1}}\cup T_{x_{2}},d)\right|=0$; otherwise, either $R_{14}\subset G$ (if $\left|E(T_{x_{1}},T_{x_{0}x_{3}})\right|=0$) or $G$ contains a subdivision of $K_{3,3}$ (if $\left|E(T_{x_{1}},T_{x_{0}x_{3}})\right|\geq 1$), a contradiction.
Thus, if $\left|T_{x_{0}x_{3}}\right|\geq 1$, then $G\subset H_{140}$ by (F\ref{case9.1-f4}) and (F\ref{case9.1-f7}).
Without loss of generality, let $\left|T_{x_{1}}\right|\geq\left|T_{x_{2}}\right|$, and let $S_{1}=T_{x_{1}}-(T_{x_{1}}\cap N(\{p,d\}))$ and $S_{2}=T_{x_{2}}\cap N(d)$ as in (F\ref{case9.1-f4}).
Thus, if $\left|T_{x_{0}x_{3}}\right|=0$, then $G\subset H_{84}$ (if $\left|T_{x_{1}}\right|=4$ or $\left|T_{x_{2}}\right|=4$) or $G\subset H_{139}$ (if $\left|T_{x_{1}}\right|=3$ and $\left|S_{1}\right|=2$) or $G\subset H_{82}$ (if $\left|T_{x_{1}}\right|=3$ and $\left|S_{1}\right|=1$) or $G\subset H_{140}$ (if $\left|T_{x_{1}}\right|=2$ and $\left|S_{1}\right|=2$) or $G\subset H_{141}$ (if $\left|T_{x_{1}}\right|=2$ and $\left|S_{1}\right|=1$) or $G\subset H_{79}$ (if $\left|T_{x_{1}}\right|=2$ and $\left|S_{1}\right|=0$) or $G\subset H_{139}$ (if $\left|T_{x_{1}}\right|=1$ and $\left|S_{1}\right|=1$) or $G\subset H_{82}$ (if $\left|T_{x_{1}}\right|=1$ and $\left|S_{1}\right|=0$) by (F\ref{case9.1-f4}) and (F\ref{case9.1-f7}).

If  $\left|E(T_{x_{3}},T_{x_{1}x_{2}})\right|\geq 1$, then we can obtain the desired conclusion by symmetry and an analogous proof.
In the following cases and subcases, let $\left|E(T_{x_{0}}\cup T_{x_{3}},T_{x_{1}x_{2}})\right|=0$.

\noindent$\textbf{Subcase 1.2}$ Suppose that $\left|E(T_{x_{0}},T_{x_{1}x_{3}})\right|\geq 1$.

Let  $p\in T_{x_{0}}$, $e\in T_{x_{1}x_{3}}$ and $p\sim e$.

\begin{fact}\label{case9.2-f1}
$\left|T_{x_{0}}\right|=\left|T_{x_{1}x_{3}}\right|=1$.
\end{fact}
\begin{proof}
If $\left|T_{x_{0}}\right|\geq 2$, then either $G[\mathcal{K}\cup\{p_{0,2},\allowbreak e\}]\cong R_{14}$ ($E(\{p_{0,2},\allowbreak e\})=\varnothing$) or $G[\mathcal{K}\cup\{p,\allowbreak p_{0,2},\allowbreak e\}]\cong R_{83}$ ($E(\{p,\allowbreak p_{0,2},\allowbreak e\})=\{pe,\allowbreak p_{0,2}e\}$) by Lemma \ref{lem-13}\ref{lem-13:item-1}, a contradiction.
If $\left|T_{x_{1}x_{3}}\right|\geq 2$, then either $G[\mathcal{K}\cup\{p,\allowbreak q_{13,2}\}]\cong R_{14}$ ($E(\{p,\allowbreak q_{13,2}\})=\varnothing$) or $G$ contains a subdivision of $K_{3,3}$ by Lemma \ref{lem-13}\ref{lem-13:item-4}, a contradiction.
\end{proof}

\begin{fact}\label{case9.2-f3}
At most three vertices in $T_{x_{1}}$ are not adjacent to $e$; at most one vertex in $T_{x_{1}}$ is adjacent to $e$; at most one vertex in $T_{x_{2}}$ is not adjacent to $e$; at most three vertices in $T_{x_{2}}$ are adjacent to $e$.
\end{fact}
\begin{proof}
If at least four vertices in $T_{x_{1}}$ are not adjacent to $e$, then $G[\{x_{0},\allowbreak x_{1},\allowbreak p,\allowbreak e,\allowbreak p_{1,1},\ldots,p_{1,4}\}]\cong R_{37}$ by Lemma \ref{lem-13}\ref{lem-13:item-1}, a contradiction.
If at least two vertices in $T_{x_{1}}$ are adjacent to $e$, then $G[\{x_{0},\allowbreak x_{1},\allowbreak x_{2},\allowbreak e,\allowbreak p_{1,1},\allowbreak p_{1,2}\}]\cong R_{16}$ by Lemma \ref{lem-13}\ref{lem-13:item-1}, a contradiction.
If at least two vertices in $T_{x_{2}}$ are  not adjacent to $e$, then $G[\{x_{0},\allowbreak x_{2},\allowbreak p,\allowbreak e,\allowbreak p_{2,1},\allowbreak p_{2,2}\}]\cong R_{1}$ by Lemma \ref{lem-13}\ref{lem-13:item-1}, a contradiction.
If at least four vertices in $T_{x_{2}}$ are adjacent to $e$, then $G[\{x_{0},\allowbreak x_{1},\allowbreak p,\allowbreak e,\allowbreak p_{2,1},\ldots,p_{2,4}\}]\cong R_{37}$ by Lemma \ref{lem-13}\ref{lem-13:item-1}, a contradiction.
\end{proof}

\begin{fact}\label{case9.2-f7}
$\left|T_{x_{3}}\right|\leq 1$, $\left|T_{x_{2}x_{3}}\right|\leq 1$ and $\left|E(T_{x_{0}},T_{x_{2}x_{3}})\right|=1$ (if $\left|T_{x_{2}x_{3}}\right|= 1$).
\end{fact}
\begin{proof}
If $\left|T_{x_{3}}\right|\geq 2$, then $G[\{x_{0},\allowbreak x_{1},\allowbreak x_{3},\allowbreak p,\allowbreak p_{3,1},\allowbreak p_{3,2}\}]\cong R_{1}$ by Lemma \ref{lem-13}\ref{lem-13:item-1} and \ref{lem-13:item-2}, a contradiction.
If $\left|T_{x_{2}x_{3}}\right|\geq 1$ and $\left|E(T_{x_{0}},T_{x_{2}x_{3}})\right|=0$, then $G[\mathcal{K}\cup\{p,\allowbreak q_{23,1}\}]\cong R_{14}$ ($E(\{p,\allowbreak q_{23,1}\})=\varnothing$), a contradiction.
If $\left|T_{x_{2}x_{3}}\right|\geq 2$, then $G$ contains a subdivision of $K_{3,3}$, a contradiction.
\end{proof}

\begin{fact}\label{case9.2-f9}
$\left|T_{x_{0}x_{2}}\cup T_{x_{1}x_{2}}\cup T_{x_{0}x_{1}x_{3}}\right|=0$.
\end{fact}
\begin{proof}
By $E(T_{2})=E(T_{x_{0}x_{3}})$ and  Lemma \ref{lem-13}\ref{lem-13:item-6}, $\left|T_{x_{0}x_{2}}\right|=0$.
If $\left|T_{x_{1}x_{2}}\right|\geq 1$, then $G[\{x_{0},\allowbreak x_{1},\allowbreak x_{2},\allowbreak p,\allowbreak e,\allowbreak q_{12,1}\}]\cong R_{12}$, a contradiction.
If $\left|T_{x_{0}x_{1}x_{3}}\right|\geq 1$, then $G$ contains a subdivision of $K_{3,3}$, a contradiction.
\end{proof}

If $\left|T_{x_{0}x_{1}}\right|\geq 1$ or $\left|E(T_{x_{1}x_{3}},T_{x_{1}})\right|\geq 1$, then  its conclusion is covered by Case $2$ of Lemma \ref{lem-42} by taking $\{x_{0},x_{1},x_{2},s\}$ or $\{x_{1},x_{3},e,s\}$ ($s\in T_{x_{0}x_{1}}\cup T_{x_{1}}$) as the vertices of a new $K_{4}-e$ by Lemma \ref{lem-13}\ref{lem-13:item-2}.
Let $\left|T_{x_{0}x_{1}}\right|=\left|E(T_{x_{1}x_{3}},T_{x_{1}})\right|=0$.

Suppose that $\left|T_{x_{0}x_{2}x_{3}}\right|\geq 1$.
By Lemma \ref{lem-13}\ref{lem-13:item-7}, $\left|T_{x_{0}x_{2}x_{3}}\right|= 1$.
$\left|E(T_{x_{0}x_{2}x_{3}},\{p,e\})\right|\geq 1$; otherwise, $G[\{x_{0},\allowbreak x_{1},\allowbreak x_{2},\allowbreak p,\allowbreak e,\allowbreak w_{023,1}\}]\cong R_{12}$, a contradiction.
If $\left|E(T_{x_{0}x_{2}x_{3}},p)\right|=1$ and $\left|E(T_{x_{0}x_{2}x_{3}},e)\right|=0$,  then  its conclusion is covered by Case $1$ of Lemma \ref{lem-43} by taking $\{x_{1},x_{2},x_{3},e\}$ as the vertices of a new $K_{4}-e$ ($G\subset H_{88}$).
Let $\left|E(T_{x_{0}x_{2}x_{3}},e)\right|=1$.
By (F\ref{case9.2-f7}) and   $G$ contains no subdivision of $K_{3,3}$, $\left|T_{x_{2}x_{3}}\cup T_{x_{0}x_{3}}\right|=0$.
Since $G$ contains no subdivision of $K_{5}$ or $K_{3,3}$, $\left|E(T_{x_{1}},T_{x_{0}x_{2}x_{3}})\right|= \left|E(T_{x_{2}},e)\right|= 0$.
Thus, $\left|T_{x_{2}}\right|\leq 1$ by (F\ref{case9.2-f3}).
If $\left|E(T_{x_{0}x_{2}x_{3}},p)\right|=1$, then $\left|T_{x_{1}}\right|=0$ by Lemma \ref{lem-13}\ref{lem-13:item-2}; otherwise, $G[\{x_{0},\allowbreak x_{1},\allowbreak x_{2},\allowbreak p,\allowbreak w_{023,1},\allowbreak p_{1,1}\}]\cong R_{14}$, a contradiction.
If $\left|E(T_{x_{0}x_{2}x_{3}},p)\right|=1$, then $\left|T_{0}\right|=0$ by Lemma \ref{lem-13}\ref{lem-13:item-8} and (F\ref{case9.2-f9}); otherwise, $G[\mathcal{K}\cup\{p,\allowbreak w_{023,1},\allowbreak t_{1}\}]\cong R_{30}$, a contradiction.
If $\left|E(T_{x_{0}x_{2}x_{3}},p)\right|=1$, then $\left|E(T_{x_{3}},T_{x_{0}x_{2}x_{3}})\right|= 1$  by Lemma \ref{lem-13}\ref{lem-13:item-2} and (F\ref{case9.2-f7}); otherwise, $R_{7}\subset G$, a contradiction.
If $\left|E(T_{x_{0}x_{2}x_{3}},p)\right|=1$, then  $\left|E(T_{x_{2}},T_{x_{0}x_{2}x_{3}})\right|= 1$  for $R_{14}\not\subset G$.
Thus, if  $\left|E(T_{x_{0}x_{2}x_{3}},p)\right|=1$, then $G\subset H_{102}$  by (F\ref{case9.2-f1}) and (F\ref{case9.2-f9}).
If $\left|E(T_{x_{0}x_{2}x_{3}},p)\right|= 0$, then $\left|E(T_{x_{2}},T_{x_{0}x_{2}x_{3}})\right|= 0$; otherwise, $G[\{x_{0},\allowbreak x_{2},\allowbreak p,\allowbreak e,\allowbreak w_{023,1},\allowbreak p_{2,1}\}]\cong R_{12}$, a contradiction.
If $\left|E(T_{x_{0}x_{2}x_{3}},p)\right|= 0$, then  $\left|E(T_{x_{3}},T_{x_{0}x_{2}x_{3}})\right|= 0$ by Lemma \ref{lem-13}\ref{lem-13:item-2}; otherwise, $G[\{x_{0},\allowbreak x_{1},\allowbreak x_{3},\allowbreak p,\allowbreak w_{023,1},\allowbreak p_{3,1}\}]\cong R_{7}$, a contradiction.
If $\left|E(T_{x_{0}x_{2}x_{3}},p)\right|= 0$, then  $\left|T_{x_{1}}\right|\leq 1$; otherwise, $G[\mathcal{K}\cup\{w_{023,1},\allowbreak p_{1,1},\allowbreak p_{1,2}\}]\cong R_{27}$, a contradiction.
Thus, if $\left|E(T_{x_{0}x_{2}x_{3}},p)\right|= 0$, then  $G\subset A_{22}$ by  (F\ref{case9.2-f1}) and (F\ref{case9.2-f9}).
Let  $\left|T_{x_{0}x_{2}x_{3}}\right|=0$.

Suppose that $\left|T_{x_{2}x_{3}}\right|\geq  1$.
By (F\ref{case9.2-f7}), $\left|T_{x_{2}x_{3}}\right|= 1$ and $\left|E(T_{x_{0}},T_{x_{2}x_{3}})\right|=1$.
Since  $G$ contains no subdivision of $K_{3,3}$, $\left|T_{x_{0}x_{3}}\right|=0$.
Thus, $\left|T_{0}\right|=0$  by Lemma \ref{lem-13}\ref{lem-13:item-8} and (F\ref{case9.2-f9}).
Then $\left|E(T_{x_{1}}\cup T_{x_{2}},T_{x_{1}x_{3}}\cup T_{x_{2}x_{3}})\right|=0$; otherwise, either $G[\{x_{0},\allowbreak x_{1},\allowbreak x_{2},\allowbreak e,\allowbreak q_{23,1},\allowbreak p_{1,1}\}]\cong R_{10}$ or $G$ contains a subdivision of $K_{3,3}$, a contradiction.
Then $\left|T_{x_{1}}\right|\leq 1$ and $\left|T_{x_{2}}\right|\leq 1$; otherwise, $G[\{x_{0},\allowbreak x_{1},\allowbreak p,\allowbreak p_{1,1},\allowbreak p_{1,2},\allowbreak q_{23,1}\}]\cong R_{1}$, a contradiction.
Thus, $G\subset H_{143}$ by Lemma \ref{lem-13}\ref{lem-13:item-2}, (F\ref{case9.2-f1}), (F\ref{case9.2-f7}) and (F\ref{case9.2-f9}).
Let $\left|T_{x_{2}x_{3}}\right|=0$.

Suppose that $\left|T_{x_{0}x_{3}}\right|\geq  1$.
Then $\left|E(T_{x_{0}x_{3}})\right|=0$ by Lemma \ref{lem-13}\ref{lem-13:item-2}; otherwise, $G[\{x_{0},\allowbreak x_{1},\allowbreak x_{2},\allowbreak e,\allowbreak q_{03,1},\allowbreak q_{03,2}\}]\cong R_{10}$, a contradiction.
Then $\left|T_{x_{0}x_{3}}\right|\leq 3$ by Lemma \ref{lem-13}\ref{lem-13:item-2}; otherwise, $G[\{x_{0},\allowbreak x_{1},\allowbreak p,\allowbreak e,\allowbreak q_{03,1},\ldots,q_{03,4}\}]\cong R_{37}$, a contradiction.
Then $\left|T_{0}\right|=0$ by Lemma \ref{lem-13}\ref{lem-13:item-2}, (F\ref{case9.2-f9}); otherwise, $G[\{x_{1},\allowbreak x_{3},\allowbreak p,\allowbreak e,\allowbreak q_{03,1},\allowbreak t_{1}\}]\cong R_{3}$, a contradiction.
Since $G$ contains no subdivision of $K_{3,3}$, $\left|E(T_{x_{1}}\cup T_{x_{2}},T_{x_{1}x_{3}}\cup T_{x_{0}x_{3}})\right|= 0$; otherwise, $G[\{x_{1},\allowbreak x_{2},\allowbreak x_{3},\allowbreak e,\allowbreak p_{2,1},\allowbreak q_{03,1}\}]\cong R_{12}$, a contradiction.
Thus, $\left|T_{x_{2}}\right|\leq 1$ by (F\ref{case9.2-f3}).
Then $\left|T_{x_{1}}\cup T_{x_{0}x_{3}}\right|\leq 3$ by (F\ref{case9.2-f3}); otherwise, $G[\{x_{0},\allowbreak x_{1},\allowbreak e,\allowbreak p_{1,1},\allowbreak q_{03,1},\allowbreak q_{03,2},\allowbreak q_{03,3}\}]\cong R_{21}$  or $G[\{x_{0},\allowbreak x_{1},\allowbreak p,\allowbreak p_{1,1},\allowbreak p_{1,2},\allowbreak q_{03,1},\allowbreak q_{03,2}\}]\cong R_{21}$  or $G[\{x_{0},\allowbreak x_{1},\allowbreak p,\allowbreak p_{1,1},\allowbreak p_{1,2},\allowbreak p_{1,3},\allowbreak q_{03,1}\}]\cong R_{21}$, a contradiction.
Thus,   $G\subset H_{146}$ (if $\left|T_{x_{0}x_{3}}\right|=1$) or $G\subset H_{145}$ (if $\left|T_{x_{0}x_{3}}\right|=2$) or $G\subset H_{144}$ (if $\left|T_{x_{0}x_{3}}\right|=3$) by Lemma \ref{lem-13}\ref{lem-13:item-2}, (F\ref{case9.2-f1}), (F\ref{case9.2-f7}) and (F\ref{case9.2-f9}).
Let $\left|T_{x_{0}x_{3}}\right|=0$.

Let $S_{1}=T_{x_{1}}-(T_{x_{1}}\cap N(e))$ and $S_{2}=T_{x_{2}}\cap N(e)$, $\left|E(S_{1},e)\right|=0$, and all vertices of $S_{2}$ are adjacent to $e$.
Then $\left|S_{1}\cup S_{2}\right| \leq 3$ by (F\ref{case9.2-f3}); otherwise, $G[\{x_{0},\allowbreak x_{1},\allowbreak e,\allowbreak p_{1,1},\allowbreak p_{2,1},\allowbreak p_{2,2},\allowbreak p_{2,3}\}]\cong R_{21}$ ($E(\{e,\allowbreak p_{1,1},\allowbreak p_{2,1},\allowbreak p_{2,2},\allowbreak p_{2,3}\})=\{ep_{2,1},\allowbreak ep_{2,2},\allowbreak ep_{2,3}\}$) or $G[\{x_{0},\allowbreak x_{1},\allowbreak e,\allowbreak p_{1,1},\allowbreak p_{1,2},\allowbreak p_{2,1},\allowbreak p_{2,2}\}]\cong R_{21}$ ($E(\{e,\allowbreak p_{1,1},\allowbreak p_{1,2},\allowbreak p_{2,1},\allowbreak p_{2,2}\})=\{ep_{2,1},\allowbreak ep_{2,2}\}$) or $G[\{x_{1},\allowbreak p,\allowbreak e,\allowbreak p_{1,1},\allowbreak p_{1,2},\allowbreak p_{1,3},\allowbreak p_{2,1}\}]\cong R_{21}$ ($E(\{p,\allowbreak e,\allowbreak p_{1,1},\allowbreak p_{1,2},\allowbreak p_{1,3},\allowbreak p_{2,1}\})=\{pe,\allowbreak ep_{2,1}\}$), a contradiction.
Then $\left|(T_{x_{1}}\cup T_{x_{2}})-(S_{1}\cup S_{2})\right| \leq 1$ by (F\ref{case9.2-f3}), and $R_{10}\not\subset G$.
Thus,  $G\subset H_{144}$ (if $\left|S_{1}\right|=0$ and  $\left|T_{x_{1}}-S_{1}\right|=0$) or $G\subset H_{145}$ (if $\left|S_{1}\right|=1$ and  $\left|T_{x_{1}}-S_{1}\right|=0$) or $G\subset H_{146}$ (if $\left|S_{1}\right|=2$ and  $\left|T_{x_{1}}-S_{1}\right|=0$) or $G\subset H_{147}$ (if $\left|S_{1}\right|=3$ and  $\left|T_{x_{1}}-S_{1}\right|=0$) or $G\subset H_{92}$ (if $\left|S_{1}\right|=0$ and  $\left|T_{x_{1}}-S_{1}\right|=1$) or $G\subset H_{97}$ (if $\left|S_{1}\right|=1$ and  $\left|T_{x_{1}}-S_{1}\right|=1$) or $G\subset H_{96}$ (if $\left|S_{1}\right|=2$ and  $\left|T_{x_{1}}-S_{1}\right|=1$) or $G\subset H_{91}$ (if $\left|S_{1}\right|=3$ and  $\left|T_{x_{1}}-S_{1}\right|=1$) by Lemma \ref{lem-13}\ref{lem-13:item-2}, (F\ref{case9.2-f1}), (F\ref{case9.2-f7}) and (F\ref{case9.2-f9}).

If  $\left|E(T_{x_{0}},T_{x_{2}x_{3}})\right|\geq 1$ or $\left|E(T_{x_{3}},T_{x_{0}x_{1}}\cup T_{x_{0}x_{2}})\right|\geq 1$, then we can obtain the desired conclusion by symmetry and an analogous proof.
In the following cases and subcases, let $\left|E(T_{x_{0}}\cup T_{x_{3}},T_{2})\right|=0$ by Lemma \ref{lem-13}\ref{lem-13:item-2}.

\noindent$\textbf{Subcase 1.3}$ Suppose that $\left|E(T_{x_{1}},T_{x_{0}x_{1}})\right|\geq 1$.

Let  $q\in T_{x_{1}}$, $a\in T_{x_{0}x_{1}}$ and $q\sim a$.

\begin{fact}\label{case9.3-f2}
$\left|T_{x_{0}}\right|\leq 1$.
\end{fact}
\begin{proof}
If $\left|T_{x_{0}}\right|\geq 2$, then $G[\mathcal{K}\cup\{q,\allowbreak p_{0,1},\allowbreak p_{0,2}\}]\cong R_{22}$ by Lemma \ref{lem-13}\ref{lem-13:item-2}, a contradiction.
\end{proof}

\begin{fact}\label{case9.3-f3}
$\left|T_{x_{3}}\cup T_{x_{2}x_{3}}\right|=0$.
\end{fact}
\begin{proof}
If $\left|T_{x_{3}}\right|\geq 1$, then $G[\mathcal{K}\cup\{a,\allowbreak p_{3,1}\}]\cong R_{14}$, a contradiction.
If $\left|T_{x_{2}x_{3}}\right|\geq 1$, then $G[\mathcal{K}\cup\{a,\allowbreak q_{23,1}\}]\cong R_{18}$, a contradiction.
\end{proof}

\begin{fact}\label{case9.3-f4}
$\left|E(T_{x_{2}},a)\right|= \left|T_{x_{2}}\right|$ and if $\left|T_{x_{2}}\right|\geq 1$, then $\left|T_{x_{1}}\cup T_{x_{2}}\right|\leq 4$.
\end{fact}
\begin{proof}
If $\left|E(T_{x_{2}},a)\right|<\left|T_{x_{2}}\right|$, then $G[\{x_{0},\allowbreak x_{1},\allowbreak x_{2},\allowbreak q,\allowbreak a,\allowbreak p_{2,1}\}]\cong R_{14}$, a contradiction.
If $\left|T_{x_{2}}\right|\geq 4$, then $G[\mathcal{K}\cup\{a,\allowbreak p_{2,1},\ldots,p_{2,4}\}]\cong R_{79}$, a contradiction.
$\left|E(T_{x_{1}},a)\right|=1$; otherwise, $G[\{x_{1},\allowbreak x_{2},\allowbreak x_{3},\allowbreak q,\allowbreak p_{1,2},\allowbreak a\}]\cong R_{16}$, a contradiction.
If $\left|T_{x_{1}}\cup T_{x_{2}}\right|\geq 5$,  then either $R_{21}\subset G$ or $R_{38}\subset G$, a contradiction.
\end{proof}

\begin{fact}\label{case9.3-f6}
$\left|T_{x_{1}x_{3}}\right|\leq 1$, $\left|T_{x_{0}x_{3}}\right|\leq 1$ and $\left|E(T_{x_{1}x_{3}},q)\right|= \left|E(T_{x_{0}x_{3}},q)\right|=1$.
\end{fact}
\begin{proof}
If $\left|E(T_{x_{1}x_{3}},q)\right|= 0$, then $R_{16}\subset G$, a contradiction.
If $\left|T_{x_{1}x_{3}}\right|\geq 2$, then $G$ contains a subdivision of $K_{3,3}$, a contradiction.
If $\left|E(T_{x_{0}x_{3}},q)\right|= 0$, then $R_{1}\subset G$, a contradiction.
If $\left|T_{x_{0}x_{3}}\right|\geq 2$, then $G$ contains a subdivision of $K_{3,3}$, a contradiction.
\end{proof}

If $\left|E(T_{x_{0}x_{2}},q)\right|= 0$, then $R_{18}\subset G$, a contradiction.
Thus, if $\left|T_{x_{0}x_{2}}\right|\geq 1$, then  its conclusion is covered by Case $3$ of Lemma \ref{lem-43} by taking $\{x_{0},x_{1},x_{2},q_{02,1}\}$ ($q_{02,1}\in T_{x_{0}x_{2}}$) as the vertices of a new $K_{4}-e$.
If $\left|E(T_{x_{1}x_{2}},q)\right|= 0$, then $R_{1}\subset G$, a contradiction.
Thus, if $\left|T_{x_{1}x_{2}}\right|\geq 1$, then  its conclusion is covered by Case $2$ of Lemma \ref{lem-43} by taking $\{x_{0},x_{1},x_{2},a\}$  as the vertices of a new $K_{4}-e$.
If $\left|E(T_{x_{0}x_{1}x_{3}},T_{x_{0}x_{1}})\right|\geq 1$, then $\omega=4$, a contradiction.
If $\left|E(T_{x_{0}x_{1}x_{3}},q)\right|= 0$, then $R_{2}\subset G$, a contradiction.
Thus, if $\left|T_{x_{0}x_{1}x_{3}}\right|\geq 1$, then  its conclusion is covered by Case $4$ of Lemma \ref{lem-43} by taking $\{x_{1},x_{3},q,w_{013,1}\}$ ($w_{013,1}\in T_{x_{0}x_{1}x_{3}}$) as the vertices of a new $K_{4}-e$.
Let $\left|T_{x_{0}x_{2}}\cup T_{x_{1}x_{2}}\cup T_{x_{0}x_{1}x_{3}}\right|=0$.

Suppose that $\left|E(T_{x_{1}},T_{x_{0}x_{1}})\right|\geq 2$.
Then $G[E(T_{x_{1}},T_{x_{0}x_{1}})]\subset aK_{2}$; otherwise, $G[\{x_{1},\allowbreak x_{2},\allowbreak x_{3},\allowbreak q,\allowbreak p_{1,2},\allowbreak a\}]\cong R_{16}$ ($E(\{q,\allowbreak p_{1,2},\allowbreak a\})=\{qa,\allowbreak p_{1,2}a\}$) or $G[\{x_{1},\allowbreak x_{2},\allowbreak x_{3},\allowbreak q,\allowbreak a,\allowbreak q_{01,2}\}]\cong R_{16}$ ($E(\{q,\allowbreak a,\allowbreak q_{01,2}\})=\{qa,\allowbreak qq_{01,2}\}$), a contradiction.
By Lemma \ref{lem-13}\ref{lem-13:item-1} and \ref{lem-13:item-4}, $a\leq 8$.
Since $R_{88}\not\subset G$, $R_{89}\not\subset G$  and $R_{96}\not\subset G$, $G[T_{x_{1}}\cup T_{x_{0}x_{1}}]\subset 8K_{2}$.
By (F\ref{case9.3-f4}), symmetry  and   $G$ contains no subdivision of $K_{3,3}$, $\left|T_{x_{2}}\right|= 0$.
By (F\ref{case9.3-f6}), symmetry  and   $G$ contains no subdivision of $K_{3,3}$, $\left|T_{x_{1}x_{3}}\cup T_{x_{0}x_{3}}\right|= 0$.
If $\left|E(T_{x_{0}x_{2}x_{3}},\{q,a\})\right|<2$, then either $R_{12}\subset G$ or $R_{17}\subset G$, a contradiction.
Thus, by Lemma \ref{lem-13}\ref{lem-13:item-7} and   $G$ contains no subdivision of $K_{3,3}$, $\left|T_{3}\right|= 0$.
If $\left|E(T_{x_{1}},T_{x_{0}x_{1}})\right|\geq 5$, then $\left|T_{x_{0}}\right|=0$; otherwise, $G[\{x_{0},\allowbreak x_{1},\allowbreak x_{2},\allowbreak p_{1,2},\ldots,p_{1,5},\allowbreak a,\allowbreak p_{0,1}\}]\cong R_{47}$, a contradiction.
Thus, either  $G\subset H_{151}$ (if $\left|E(T_{x_{1}},T_{x_{0}x_{1}})\right|\geq 5$) or $G\subset H_{149}$ (if $\left|E(T_{x_{1}},T_{x_{0}x_{1}})\right|\leq 4$) by  Lemma \ref{lem-13}\ref{lem-13:item-2} and (F\ref{case9.3-f3}).
Let $\left|E(T_{x_{1}},T_{x_{0}x_{1}})\right|=1$.

Suppose that $\left|T_{x_{0}x_{2}x_{3}}\right|\geq 1$.
Then $\left|T_{x_{0}x_{2}x_{3}}\right|=1$ by Lemma \ref{lem-13}\ref{lem-13:item-7}.
If $\left|E(T_{x_{0}x_{2}x_{3}},\{q,a\})\right|<2$, then either $G[\{x_{0},\allowbreak x_{1},\allowbreak x_{3},\allowbreak q,\allowbreak a,\allowbreak w_{023,1}\}]\cong R_{12}$ ($E(\{q,\allowbreak a,\allowbreak w_{023,1}\})=\{qa\}$) or $G[\{x_{1},\allowbreak x_{2},\allowbreak x_{3},\allowbreak q,\allowbreak a,\allowbreak w_{023,1}\}]\cong R_{17}$ ($E(\{q,\allowbreak a,\allowbreak w_{023,1}\})=\{qa,\allowbreak qw_{023,1}\}$) or $G[\{x_{1},\allowbreak x_{2},\allowbreak x_{3},\allowbreak q,\allowbreak a,\allowbreak w_{023,1}\}]\cong R_{17}$ ($E(\{q,\allowbreak a,\allowbreak w_{023,1}\})=\{qa,\allowbreak aw_{023,1}\}$), a contradiction.
Thus, $\left|E(T_{x_{0}x_{2}x_{3}},\{q,a\})\right|= 2$.
Since $G$ contains no subdivision of $K_{3,3}$, $\left|E(T_{x_{0}x_{1}}-\{a\},T_{x_{0}x_{2}x_{3}})\right|= 0$.
Then $\left|T_{x_{0}x_{1}}-\{a\}\right|\leq 1$; otherwise, $G[\mathcal{K}\cup\{q_{01,2},\allowbreak q_{01,3},\allowbreak w_{023,1}\}]\cong R_{33}$, a contradiction.
By (F\ref{case9.3-f4}) and   $G$ contains no subdivision of $K_{3,3}$, $\left|T_{x_{2}}\right|= 0$.
By (F\ref{case9.3-f6})  and   $G$ contains no subdivision of $K_{3,3}$, $\left|T_{x_{0}x_{3}}\right|= 0$.
By (F\ref{case9.3-f6})  and $R_{17}\not\subset G$, $\left|E(T_{x_{1}x_{3}},T_{x_{0}x_{2}x_{3}})\right|= 1$.
Then at most one vertex in $T_{x_{1}}$ is not adjacent to $T_{x_{0}x_{2}x_{3}}$; otherwise, $G[\{x_{0},\allowbreak x_{1},\allowbreak x_{2},\allowbreak p_{1,2},\allowbreak p_{1,3},\allowbreak a,\allowbreak w_{023,1}\}]\cong R_{27}$, a contradiction.
If $\left|T_{x_{1}x_{3}}\right|= 1$, then  $\left|T_{x_{0}}\right|=0$; otherwise, $G[\mathcal{K}\cup\{q_{13,1},\allowbreak p_{0,1}\}]\cong R_{14}$, a contradiction.
If $\left|T_{x_{1}x_{3}}\right|= 1$, then  $\left|E(T_{x_{1}}-\{q\},T_{x_{0}x_{2}x_{3}})\right|= 0$ for $G$ contains no subdivision of $K_{3,3}$.
If $\left|T_{x_{1}x_{3}}\right|= 1$, then  $\left|T_{x_{0}x_{1}}-\{a\}\right|=0$; otherwise, $G[\{x_{0},\allowbreak x_{1},\allowbreak x_{3},\allowbreak a,\allowbreak q_{01,2},\allowbreak q_{13,1}\}]\cong R_{16}$, a contradiction.
Thus, if $\left|T_{x_{1}x_{3}}\right|= 1$, then $G\subset A_{18}$ by  Lemma \ref{lem-13}\ref{lem-13:item-8} and (F\ref{case9.3-f3}).
Let $\left|T_{x_{1}x_{3}}\right|= 0$.
If $\left|T_{x_{0}x_{1}}-\{a\}\right|= 1$, then no vertex in $T_{x_{1}}$ is not adjacent to $T_{x_{0}x_{2}x_{3}}$; otherwise, $G[\mathcal{K}\cup\{p_{1,2},\allowbreak q_{01,2},\allowbreak w_{023,1}\}]\cong R_{30}$, a contradiction.
If $\left|T_{x_{0}x_{1}}-\{a\}\right|= 1$, then at most two vertices in $T_{x_{1}}-\{q\}$ are adjacent to $T_{x_{0}x_{2}x_{3}}$; otherwise, $G[\{x_{0},\allowbreak x_{1},\allowbreak x_{2},\allowbreak q,\allowbreak p_{1,2},\allowbreak p_{1,3},\allowbreak p_{1,4},\allowbreak q_{01,2},\allowbreak w_{023,1}\}]\cong R_{79}$, a contradiction.
If $\left|T_{x_{0}x_{1}}-\{a\}\right|= 1$, then $\left|E(T_{x_{0}},T_{x_{0}x_{2}x_{3}})\right|= 1$ for $R_{30}\not\subset G$.
Thus, if $\left|T_{x_{0}x_{1}}-\{a\}\right|= 1$, then $G\subset H_{112}$ by   Lemma \ref{lem-13}\ref{lem-13:item-2}\ref{lem-13:item-8} and (F\ref{case9.3-f3}).
Let $\left|T_{x_{0}x_{1}}-\{a\}\right|= 0$.
If $\left|E(T_{x_{0}},T_{x_{0}x_{2}x_{3}})\right|= 1$, then no vertex in $T_{x_{1}}$ is not adjacent to $T_{x_{0}x_{2}x_{3}}$; otherwise, $G[\{x_{0},\allowbreak x_{1},\allowbreak x_{2},\allowbreak p_{1,2},\allowbreak w_{023,1},\allowbreak p_{0,1}\}]\cong R_{14}$, a contradiction.
If $\left|E(T_{x_{0}},T_{x_{0}x_{2}x_{3}})\right|= 1$, then at most two vertices in $T_{x_{1}}-\{q\}$ are adjacent to $T_{x_{0}x_{2}x_{3}}$; otherwise, $G[\{x_{0},\allowbreak x_{1},\allowbreak x_{2},\allowbreak q,\allowbreak p_{1,2},\allowbreak p_{1,3},\allowbreak p_{1,4},\allowbreak w_{023,1},\allowbreak p_{0,1}\}]\cong R_{79}$, a contradiction.
Thus, if $\left|E(T_{x_{0}},T_{x_{0}x_{2}x_{3}})\right|= 1$, then $G\subset H_{112}$  by   Lemma \ref{lem-13}\ref{lem-13:item-2}\ref{lem-13:item-8} and  (F\ref{case9.3-f3}).
Let $\left|E(T_{x_{0}},T_{x_{0}x_{2}x_{3}})\right|= 0$.
If $\left|T_{x_{0}}\right|= 1$, then at most two vertices in $T_{x_{1}}-\{q\}$ are adjacent to $T_{x_{0}x_{2}x_{3}}$ for $R_{55}\not\subset G$.
Since $R_{71}\not\subset G$, at most six vertices in $T_{x_{1}}-\{q\}$ are adjacent to $T_{x_{0}x_{2}x_{3}}$.
Thus, either  $G\subset A_{19}$ (if $\left|T_{x_{0}}\right|= 1$) or $G\subset A_{20}$ (if $\left|T_{x_{0}}\right|= 0$) by    Lemma \ref{lem-13}\ref{lem-13:item-2}\ref{lem-13:item-8} and (F\ref{case9.3-f3}).
Let  $\left|T_{x_{0}x_{2}x_{3}}\right|=0$.

Suppose that $\left|T_{x_{0}x_{3}}\right|\geq 1$.
By (F\ref{case9.3-f6}), $\left|T_{x_{0}x_{3}}\right|=1$ and $\left|E(T_{x_{0}x_{3}},q)\right|= 1$.
By (F\ref{case9.3-f3}), $\left|T_{x_{3}}\cup T_{x_{2}x_{3}}\right|=0$.
By (F\ref{case9.3-f4}) and   $G$ contains no subdivision of $K_{3,3}$, $\left|T_{x_{2}}\right|= 0$.
Then $\left|T_{0}\right|= 0$ by  Lemma \ref{lem-13}\ref{lem-13:item-8}; otherwise, $G[\{x_{0},\allowbreak x_{1},\allowbreak x_{2},\allowbreak q,\allowbreak q_{03,1},\allowbreak t_{1}\}]\cong R_{7}$, a contradiction.
Then $\left|E(T_{x_{0}},T_{x_{0}x_{3}})\right|=0$; otherwise, $G[\mathcal{K}\cup\{q_{03,1},\allowbreak p_{0,1}\}]\cong R_{17}$, a contradiction.
Let $S_{1}=T_{x_{1}}\cap N(T_{x_{0}x_{3}})-\{q\}$ and $S_{2}=T_{x_{1}}-(S_{1}\cup \{q\})$.
Since $R_{38}\not\subset G$, $\left|S_{2}\right|\leq 2$.
Since $R_{26}\not\subset G$, $\left|T_{x_{0}x_{1}}\cup S_{1}-\{a\}\right|\leq 1$.
Thus, either  $G\subset H_{101}$ (if $\left|S_{1}\right|=1$) or $G\subset H_{86}$ (if $\left|S_{1}\right|=0$) by  Lemma \ref{lem-13}\ref{lem-13:item-2} and (F\ref{case9.3-f3}).
Let $\left|T_{x_{0}x_{3}}\right|=0$.
Thus, $\left|T_{0}\right|=0$  Lemma \ref{lem-13}\ref{lem-13:item-8}.

Suppose that $\left|T_{x_{1}x_{3}}\right|\geq 1$.
By (F\ref{case9.3-f6}), $\left|T_{x_{1}x_{3}}\right|=1$ and $\left|E(T_{x_{1}x_{3}},q)\right|= 1$.
Then $\left|T_{x_{0}x_{1}}-\{a\}\right|=0$; otherwise, either $G[\{x_{0},\allowbreak x_{1},\allowbreak x_{3},\allowbreak a,\allowbreak q_{01,2},\allowbreak q_{13,1}\}]\cong R_{16}$ or $G$ contains a subdivision of $K_{3,3}$, a contradiction.
Then $\left|T_{x_{0}}\right|=0$ for $R_{14}\not\subset G$.
By (F\ref{case9.3-f4}), $R_{10}\not\subset G$ and   $G$ contains no subdivision of $K_{3,3}$, $\left|E(T_{x_{2}},T_{x_{1}x_{3}}\cup \{a\})\right|= 2$ and $\left|T_{x_{2}}\right|\leq 1$.
Then $\left|E(T_{x_{1}}-\{q\},T_{x_{1}x_{3}})\right|= 0$; otherwise, $G[\{x_{0},\allowbreak x_{1},\allowbreak x_{2},\allowbreak p_{1,2},\allowbreak a,\allowbreak q_{13,1}\}]\cong R_{16}$, a contradiction.
If $\left|T_{x_{2}}\right|= 1$, then $\left|T_{x_{1}}\right|\leq 3$  for $R_{38}\not\subset G$.
Thus, either $G\subset H_{150}$ (if $\left|T_{x_{2}}\right|=0$) or $G\subset H_{148}$ (if $\left|T_{x_{2}}\right|=1$) by  Lemma \ref{lem-13}\ref{lem-13:item-1} and (F\ref{case9.3-f3}).
Let $\left|T_{x_{1}x_{3}}\right|=0$.

Suppose that $\left|T_{x_{0}x_{1}}-\{a\}\right|\geq 1$.
Then $\left|E(T_{x_{2}},T_{x_{0}x_{1}}-\{a\})\right|= 0$ for $G$ contains no subdivision of $K_{3,3}$.
Then $\left|T_{x_{2}}\right|\leq 1$; otherwise, $G[\mathcal{K}\cup\{q_{01,2},\allowbreak p_{2,1},\allowbreak p_{2,2}\}]\cong R_{25}$, a contradiction.
If $\left|T_{x_{2}}\right|=1$, then $\left|T_{x_{0}x_{1}}-\{a\}\right|=1$; otherwise, $G[\mathcal{K}\cup\{q_{01,2},\allowbreak q_{01,3},\allowbreak p_{2,1}\}]\cong R_{28}$, a contradiction.
If $\left|T_{x_{2}}\right|=1$, then $\left|T_{x_{1}}\right|\leq 3$ for $R_{38}\not\subset G$.
Thus, if $\left|T_{x_{2}}\right|=1$, then $G\subset H_{142}$ by  Lemma \ref{lem-13}\ref{lem-13:item-2} and (F\ref{case9.3-f3}).
Let $\left|T_{x_{2}}\right|=0$.
If $\left|T_{x_{0}}\right|=1$, then $\left|T_{x_{0}x_{1}}\cup T_{x_{1}}-\{a\}\right|\leq 4$ for $R_{47}\not\subset G$ and $R_{87}\not\subset G$.
Thus, if $\left|T_{x_{0}}\right|=1$, then 
$G\subset H_{87}$ (if $\left|T_{x_{1}}\right|=4$) or 
$G\subset H_{80}$ (if $\left|T_{x_{1}}\right|=3$) or 
$G\subset H_{111}$ (if $\left|T_{x_{1}}\right|=2$) or 
$G\subset H_{114}$ (if $\left|T_{x_{1}}\right|=1$) 
by  Lemma \ref{lem-13}\ref{lem-13:item-2} and (F\ref{case9.3-f3}).
Note that $G[T_{x_{1}}\cup T_{x_{0}x_{1}}]\subset 8K_{2}$.
Thus, if $\left|T_{x_{0}}\right|=0$, then $G\subset H_{151}$ by (F\ref{case9.3-f3}).
Let $\left|T_{x_{0}x_{1}}-\{a\}\right|=0$.

Then $G\subset H_{151}$ (if $\left|T_{x_{2}}\right|=0$) or $G\subset H_{85}$ (if $\left|T_{x_{2}}\right|=1$) or $G\subset H_{95}$ (if $\left|T_{x_{2}}\right|=2$) or $G\subset H_{92}$ (if $\left|T_{x_{2}}\right|=3$) by  (F\ref{case9.3-f3}) and (F\ref{case9.3-f4}).

If  $\left|E(T_{x_{1}},T_{x_{1}x_{3}})\right|\geq 1$ or $\left|E(T_{x_{2}},T_{x_{0}x_{2}}\cup T_{x_{2}x_{3}})\right|\geq 1$, then we can obtain the desired conclusion by symmetry and an analogous proof.
In the following cases and subcases, let $\left|E(T_{x_{1}},T_{x_{0}x_{1}}\cup T_{x_{1}x_{3}})\right|=0$ and $\left|E(T_{x_{2}},T_{x_{0}x_{2}}\cup T_{x_{2}x_{3}})\right|=0$.

\noindent$\textbf{Subcase 1.4}$ Suppose that $\left|E(T_{x_{1}},T_{x_{0}x_{2}})\right|\geq 1$.

Let  $q\in T_{x_{1}}$, $b\in T_{x_{0}x_{2}}$ and $q\sim b$.
If $\left|T_{x_{0}x_{2}}-\{b\}\right|\geq 1$, then  its conclusion is covered by Case $1$ by taking $\{x_{0},x_{1},x_{2},q_{02,1}\}$ ($q_{02,1}\in T_{x_{0}x_{2}}-\{b\}$) as the vertices of a new $K_{4}-e$.
If $\left|E(T_{x_{1}x_{2}},q)\right|\geq 1$, then  its conclusion is covered by Case $3$ of Lemma \ref{lem-43} by taking $\{x_{1},x_{2},q_{12,1},q\}$ ($q_{12,1}\in T_{x_{1}x_{2}}$) as the vertices of a new $K_{4}-e$.
Let $\left|T_{x_{0}x_{2}}-\{b\}\right|=\left|E(T_{x_{1}x_{2}},q)\right|=0$.

\begin{fact}\label{case9.4-f1}
At most one vertex in $T_{x_{1}}-\{q\}$ is not adjacent to $b$; at most two vertices in $T_{x_{1}}-\{q\}$ are adjacent to $b$; if $p_{1,1}\in T_{x_{1}}$, then $\left|E(p_{1,1},T_{x_{0}x_{1}})\right| \leq 1$.
\end{fact}
\begin{proof}
If at least two vertices in $T_{x_{1}}-\{q\}$ are not adjacent to $b$, then $G[\mathcal{K}\cup\{p_{1,2},\allowbreak p_{1,3},\allowbreak b\}]\cong R_{25}$ by Lemma \ref{lem-13}\ref{lem-13:item-1}, a contradiction.
If at least three vertices in $T_{x_{1}}-\{q\}$ are adjacent to $b$, then $G[\mathcal{K}\cup\{q,\allowbreak p_{1,2},\allowbreak p_{1,3},\allowbreak p_{1,4},\allowbreak b\}]\cong R_{79}$ by Lemma \ref{lem-13}\ref{lem-13:item-1}, a contradiction.
If $\left|E(p_{1,1},T_{x_{0}x_{1}})\right|\geq 2$, then  $G$ contains a subdivision of $K_{3,3}$, a contradiction.
\end{proof}

\begin{fact}\label{case9.4-f4}
$\left|T_{x_{1}x_{2}}\cup T_{x_{0}}\right|\leq 1$, $\left|T_{x_{2}}\right|+\left|T_{x_{1}}\cap N(b)-\{q\}\right|\leq 2$, $\left|E(T_{x_{1}x_{2}},T_{x_{1}}\cup T_{x_{2}}-\{q\})\right|=0$ and $\left|T_{x_{0}x_{1}}\right|\leq 2$.
Furthermore, if $\left|E(T_{x_{2}},T_{x_{0}x_{1}})\right|=0$, then  $\left|T_{x_{2}}\cup T_{x_{0}x_{1}}\right|\leq 2$.
\end{fact}
\begin{proof}
If $\left|T_{x_{0}}\right|\geq 2$, then $G[\mathcal{K}\cup\{b,\allowbreak p_{0,1},\allowbreak p_{0,2}\}]\cong R_{25}$ by  Lemma \ref{lem-13}\ref{lem-13:item-2}, a contradiction.
If $\left|T_{x_{1}x_{2}}\right|\geq 2$, then $G[\mathcal{K}\cup\{q,\allowbreak b,\allowbreak q_{12,1},\allowbreak q_{12,2}\}]\cong R_{42}$, a contradiction.
If $\left|T_{x_{1}x_{2}}\cup T_{x_{0}}\right|\geq 2$, then $G[\mathcal{K}\cup\{b,\allowbreak q_{12,1},\allowbreak p_{0,1}\}]\cong R_{28}$, a contradiction.
If $\left|T_{x_{2}}\right|+\left|T_{x_{1}}\cap N(b)-\{q\}\right|\geq 3$, then either $R_{38}\subset G$ (if $\left|T_{x_{2}}\right|\geq 3$) or $R_{21}\subset G$ (the remaining case) by (F\ref{case9.4-f1}), a contradiction.
If $\left|E(T_{x_{1}x_{2}},T_{x_{1}}\cup T_{x_{2}}-\{q\})\right|\geq 1$, then either $G[\{x_{1},\allowbreak x_{2},\allowbreak q,\allowbreak b,\allowbreak q_{12,1},\allowbreak p_{2,1}\}]\cong R_{12}$ or $G$ contains a subdivision of $K_{3,3}$, a contradiction.
If $\left|T_{x_{0}x_{1}}\right|\geq 3$, then $G[\{x_{0},\allowbreak x_{1},\allowbreak x_{2},\allowbreak q,\allowbreak b,\allowbreak q_{01,1},\allowbreak q_{01,2},\allowbreak q_{01,3}\}]\cong R_{42}$, a contradiction.
If $\left|E(T_{x_{2}},T_{x_{0}x_{1}})\right|=0$ and $\left|T_{x_{2}}\cup T_{x_{0}x_{1}}\right|\geq 3$, then either $G[\{x_{0},\allowbreak x_{1},\allowbreak x_{2},\allowbreak b,\allowbreak p_{2,1},\allowbreak p_{2,2},\allowbreak q_{01,1}\}]\cong R_{25}$ ($E(\{b,\allowbreak p_{2,1},\allowbreak p_{2,2},\allowbreak q_{01,1}\})=\varnothing$) or $G[\{x_{0},\allowbreak x_{1},\allowbreak x_{2},\allowbreak b,\allowbreak p_{2,1},\allowbreak q_{01,1},\allowbreak q_{01,2}\}]\cong R_{28}$ ($E(\{b,\allowbreak p_{2,1},\allowbreak q_{01,1},\allowbreak q_{01,2}\})=\varnothing$), a contradiction.
\end{proof}

\begin{fact}\label{case9.4-f6}
$\left|T_{x_{3}}\cup T_{x_{1}x_{3}}\cup T_{3}\cup T_{0}\right|=0$.
\end{fact}
\begin{proof}
If $\left|T_{x_{3}}\right|\geq 1$, then $G[\mathcal{K}\cup\{b,\allowbreak p_{3,1}\}]\cong R_{14}$, a contradiction.
If $\left|T_{x_{1}x_{3}}\right|\geq 1$, then $G[\mathcal{K}\cup\{b,\allowbreak q_{13,1}\}]\cong R_{18}$, a contradiction.
If $\left|T_{x_{0}x_{1}x_{3}}\right|\geq 1$, then either $G[\{x_{1},\allowbreak x_{2},\allowbreak x_{3},\allowbreak q,\allowbreak b,\allowbreak w_{013,1}\}]\cong R_{12}$ ($E(\{q,\allowbreak b,\allowbreak w_{013,1}\})=\{qb\}$) or $G$ contains a subdivision of $K_{3,3}$, a contradiction.
If $\left|T_{x_{0}x_{2}x_{3}}\right|\geq 1$, then either $G[\{x_{0},\allowbreak x_{2},\allowbreak x_{3},\allowbreak q,\allowbreak b,\allowbreak w_{023,1}\}]\cong R_{14}$ ($E(\{q,\allowbreak b,\allowbreak w_{023,1}\})=\{qb\}$) or $G$ contains a subdivision of $K_{3,3}$, a contradiction.
Thus, $\left|T_{3}\right|=0$ by Lemma \ref{lem-13}\ref{lem-13:item-7}.
If $\left|E(T_{x_{1}x_{2}},T_{0})\right|\geq 1$, then $G[\{x_{0},\allowbreak x_{1},\allowbreak x_{2},\allowbreak b,\allowbreak q_{12,1},\allowbreak t_{1}\}]\cong R_{14}$ by Lemma \ref{lem-13}\ref{lem-13:item-8}, a contradiction.
If $\left|E(T_{x_{0}x_{3}},T_{0})\right|\geq 1$, then $G[\{x_{0},\allowbreak x_{1},\allowbreak x_{2},\allowbreak q,\allowbreak q_{03,1},\allowbreak t_{1}\}]\cong R_{3}$ by Lemma \ref{lem-13}\ref{lem-13:item-8} and (F\ref{case9.4-f10}), a contradiction.
Thus, $\left|T_{0}\right|=0$ by Lemma \ref{lem-13}\ref{lem-13:item-8}.
\end{proof}

\begin{fact}\label{case9.4-f10}
$\left|E(T_{x_{2}x_{3}},q)\right|=\left|T_{x_{2}x_{3}}\right|$, $\left|E(T_{x_{0}x_{3}},q)\right|=0$ and $\left|T_{x_{2}x_{3}}\cup T_{x_{0}x_{3}}\right| \leq1$.
\end{fact}
\begin{proof}
If $\left|E(T_{x_{2}x_{3}},q)\right|=0$, then $R_{2}\subset G$, a contradiction.
If $\left|E(T_{x_{0}x_{3}},q)\right|\geq 1$, then  $G$ contains a subdivision of $K_{3,3}$, a contradiction.
If $\left|T_{x_{2}x_{3}}\right|\geq 2$, then $G$ contains a subdivision of $K_{3,3}$, a contradiction.
If $\left|T_{x_{0}x_{3}}\right|\geq 2$, then either $G[\{x_{0},\allowbreak x_{3},\allowbreak q,\allowbreak b,\allowbreak q_{03,1},\allowbreak q_{03,2}\}]\cong R_{2}$ ($E(\{q,\allowbreak b,\allowbreak q_{03,1},\allowbreak q_{03,2}\})=\{qb\}$) or $G[\{x_{0},\allowbreak x_{1},\allowbreak x_{2},\allowbreak b,\allowbreak q_{03,1},\allowbreak q_{03,2}\}]\cong R_{16}$ ($E(\{b,\allowbreak q_{03,1},\allowbreak q_{03,2}\})=\{q_{03,1}q_{03,2}\}$), a contradiction.
If $\left|T_{x_{2}x_{3}}\cup T_{x_{0}x_{3}}\right| \geq 2$, then  $G$ contains a subdivision of $K_{3,3}$, a contradiction.
\end{proof}

Suppose that $\left|T_{x_{2}x_{3}}\right|\geq 1$.
By (F\ref{case9.4-f10}), $\left|T_{x_{2}x_{3}}\right| = 1$,  $\left|T_{x_{0}x_{3}}\right| = 0$ and $\left|E(T_{x_{2}x_{3}},q)\right|=1$.
By (F\ref{case9.4-f6}) and symmetry, $\left|T_{x_{0}}\cup T_{x_{3}}\cup T_{x_{0}x_{1}}\cup T_{x_{1}x_{3}}\cup T_{3}\cup T_{0}\right|=0$.
Then $\left|E(T_{x_{1}}-\{q\},T_{x_{0}x_{2}}\cup T_{x_{2}x_{3}})\right|=0$; otherwise, either $G[\{x_{0},\allowbreak x_{1},\allowbreak x_{3},\allowbreak p_{1,2},\allowbreak b,\allowbreak q_{23,1}\}]\cong R_{2}$ ($E(\{p_{1,2},\allowbreak b,\allowbreak q_{23,1}\})=\{p_{1,2}b\}$ or $\{p_{1,2}q_{23,1}\}$) or  $G$ contains a subdivision of $K_{3,3}$, a contradiction.
Thus, $G\subset H_{86}$ by   (F\ref{case9.4-f1}) and (F\ref{case9.4-f4}).
Let $\left|T_{x_{2}x_{3}}\right|=0$.

Suppose that $\left|T_{x_{0}x_{3}}\right|\geq 1$.
By (F\ref{case9.4-f10}), $\left|T_{x_{0}x_{3}}\right|=1$ and $\left|E(T_{x_{0}x_{3}},q)\right|=0$.
Then $\left|E(T_{x_{0}x_{3}},T_{x_{1}}\cup T_{x_{2}}-\{q\})\right|=0$; otherwise, $G[\{x_{1},\allowbreak x_{3},\allowbreak q,\allowbreak p_{1,2},\allowbreak b,\allowbreak q_{03,1}\}]\cong R_{2}$ ($E(\{q,\allowbreak p_{1,2},\allowbreak b,\allowbreak q_{03,1}\})=\{qb,\allowbreak p_{1,2}q_{03,1}\}$) or $G[\{x_{2},\allowbreak x_{3},\allowbreak q,\allowbreak b,\allowbreak q_{03,1},\allowbreak p_{2,1}\}]\cong R_{2}$ ($E(\{q,\allowbreak b,\allowbreak q_{03,1},\allowbreak p_{2,1}\})=\{qb,\allowbreak q_{03,1}p_{2,1}\}$) or  $G$ contains a subdivision of $K_{3,3}$, a contradiction.
Since $G$ contains no subdivision of $K_{3,3}$,  $\left|E(T_{x_{2}},T_{x_{0}x_{1}})\right|=0$.
Thus,  $G\subset H_{152}$ (if  $\left|T_{x_{0}}\right|=1$, $\left|T_{x_{1}}-\{q\}\right|=1$ and $\left|T_{x_{2}}\right|=2$)  or  $G\subset H_{153}$ (if  $\left|T_{x_{0}}\right|=1$, $\left|T_{x_{1}}-\{q\}\right|=1$ and $\left|T_{x_{2}}\right|=1$) or  $G\subset H_{81}$ (if  $\left|T_{x_{0}}\right|=1$, $\left|T_{x_{1}}-\{q\}\right|=1$ and $\left|T_{x_{2}}\right|=0$) or $G\subset H_{75}$ (if  $\left|T_{x_{0}}\right|=1$, $\left|T_{x_{1}}-\{q\}\right|=0$ and $\left|T_{x_{2}}\right|=2$) or $G\subset H_{74}$ (if  $\left|T_{x_{0}}\right|=1$, $\left|T_{x_{1}}-\{q\}\right|=0$ and $\left|T_{x_{2}}\right|=1$) or $G\subset H_{73}$ (if  $\left|T_{x_{0}}\right|=1$, $\left|T_{x_{1}}-\{q\}\right|=0$ and $\left|T_{x_{2}}\right|=0$) or $G\subset H_{75}$ (if  $\left|T_{x_{0}}\right|=0$, $\left|T_{x_{1}}-\{q\}\right|=1$ and $\left|T_{x_{2}}\right|=2$) or $G\subset H_{74}$ (if  $\left|T_{x_{0}}\right|=0$, $\left|T_{x_{1}}-\{q\}\right|=1$ and $\left|T_{x_{2}}\right|=1$) or $G\subset H_{73}$ (if  $\left|T_{x_{0}}\right|=0$, $\left|T_{x_{1}}-\{q\}\right|=1$ and $\left|T_{x_{2}}\right|=0$) or $G\subset H_{78}$ (if  $\left|T_{x_{0}}\right|=0$, $\left|T_{x_{1}}-\{q\}\right|=0$ and $\left|T_{x_{2}}\right|=2$) or $G\subset H_{77}$ (if  $\left|T_{x_{0}}\right|=0$, $\left|T_{x_{1}}-\{q\}\right|=0$ and $\left|T_{x_{2}}\right|=1$) or $G\subset H_{76}$ (if  $\left|T_{x_{0}}\right|=0$, $\left|T_{x_{1}}-\{q\}\right|=0$ and $\left|T_{x_{2}}\right|=0$) by (F\ref{case9.4-f1}), (F\ref{case9.4-f4}) and (F\ref{case9.4-f6}).
Let $\left|T_{x_{0}x_{3}}\right|=0$.

Suppose that $\left|T_{x_{0}x_{1}}\right|\geq 1$.
Then $\left|E(b,T_{x_{1}}-\{q\})\right| \leq 1$; otherwise, $G[\mathcal{K}\cup\{q,\allowbreak p_{1,2},\allowbreak p_{1,3},\allowbreak b,\allowbreak q_{01,1}\}]\cong R_{78}$, a contradiction.
If $\left|E(b,T_{x_{1}}-\{q\})\right|=1$, then  $\left|T_{x_{0}x_{1}}\right|=1$; otherwise, $G[\{x_{0},\allowbreak x_{1},\allowbreak q,\allowbreak p_{1,2},\allowbreak b,\allowbreak q_{01,1},\allowbreak q_{01,2}\}]\cong R_{26}$, a contradiction.
If $\left|E(b,T_{x_{1}}-\{q\})\right|=1$, then   $\left|E(T_{x_{2}},T_{x_{0}x_{1}})\right|=0$; otherwise, $G[\{x_{2},\allowbreak q,\allowbreak p_{1,2},\allowbreak b,\allowbreak q_{01,1},\allowbreak p_{2,1}\}]\cong R_{1}$, a contradiction.
Thus, if $\left|E(b,T_{x_{1}}-\{q\})\right|=1$, then either $G\subset H_{4}$ (if  $\left|T_{x_{0}}\right|=1$) or $G\subset H_{80}$ (if  $\left|T_{x_{0}}\right|=0$) by (F\ref{case9.4-f1}), (F\ref{case9.4-f4}) and (F\ref{case9.4-f6}).
Let $\left|E(b,T_{x_{1}}-\{q\})\right|=0$.
Then $\left|E(T_{x_{2}},T_{x_{0}x_{1}})\right|\leq 1$; otherwise, $G[\{x_{2},\allowbreak q,\allowbreak b,\allowbreak q_{01,1},\allowbreak p_{2,1},\allowbreak p_{2,2}\}]\cong R_{2}$ ($E(\{q,\allowbreak b,\allowbreak q_{01,1},\allowbreak p_{2,1},\allowbreak p_{2,2}\})=\{qb,\allowbreak q_{01,1}p_{2,1},\allowbreak q_{01,1}p_{2,2}\}$) or $G[\{x_{0},\allowbreak q,\allowbreak b,\allowbreak q_{01,1},\allowbreak q_{01,2},\allowbreak p_{2,1}\}]\cong R_{2}$ ($E(\{q,\allowbreak b,\allowbreak q_{01,1},\allowbreak q_{01,2},\allowbreak p_{2,1}\})=\{qb,\allowbreak q_{01,1}p_{2,1},\allowbreak q_{01,2}p_{2,1}\}$) or  $G$ contains a subdivision of $K_{3,3}$, a contradiction.
Thus, $G\subset H_{73}$ (if  $\left|T_{x_{0}}\right|=0$ and $\left|T_{x_{2}}\right|\leq 1$) or $G\subset H_{74}$ (if  $\left|T_{x_{0}}\right|=1$ and $\left|T_{x_{2}}\right|\leq 1$) or $G\subset H_{81}$ (if  $\left|T_{x_{0}}\right|=0$ and $\left|T_{x_{2}}\right|=2$) or $G\subset H_{153}$ (if  $\left|T_{x_{0}}\right|=1$ and $\left|T_{x_{2}}\right|=2$) by (F\ref{case9.4-f1}), (F\ref{case9.4-f4}) and (F\ref{case9.4-f6}).
Let $\left|T_{x_{0}x_{1}}\right|=0$.

If  $\left|T_{x_{1}}\cap N(b)-\{q\}\right|\geq 1$, then $\left|T_{x_{1}x_{2}}\right|=0$ by (F\ref{case9.4-f4}); otherwise, $G[\{x_{1},\allowbreak x_{2},\allowbreak x_{3},\allowbreak q,\allowbreak p_{1,2},\allowbreak b,\allowbreak q_{12,1}\}]\cong R_{26}$, a contradiction.
Thus,  
$G\subset A_{24}$ (if  $\left|T_{x_{1}}\cap N(b)-\{q\}\right|\geq 1$) or 
$G\subset A_{30}$ (if  $\left|T_{x_{1}}\cap N(b)-\{q\}\right|=0$ and $\left|T_{x_{0}}\right|=1$) or $G\subset H_{11}$ (if  $\left|T_{x_{1}}\cap N(b)-\{q\}\right|=0$ and $\left|T_{x_{0}}\right|=0$) by (F\ref{case9.4-f1}), (F\ref{case9.4-f4}) and (F\ref{case9.4-f6}).

If  $\left|E(T_{x_{1}},T_{x_{2}x_{3}})\right|\geq 1$ or $\left|E(T_{x_{2}},T_{x_{0}x_{1}}\cup T_{x_{1}x_{3}})\right|\geq 1$, then we can obtain the desired conclusion by symmetry and an analogous proof.
In the following cases and subcases, let $\left|E(T_{x_{1}},T_{x_{0}x_{2}}\cup T_{x_{2}x_{3}})\right|=0$ and $\left|E(T_{x_{2}},T_{x_{0}x_{1}}\cup T_{x_{1}x_{3}})\right|=0$.

\noindent$\textbf{Subcase 1.5}$ Suppose that $\left|E(T_{x_{1}},T_{x_{0}x_{3}})\right|\geq 1$.

Let  $q\in T_{x_{1}}$, $c\in T_{x_{0}x_{3}}$ and $q\sim c$.
If $\left|T_{x_{0}x_{1}}\cup T_{x_{1}x_{3}}\right|\geq 1$, then  its conclusion is covered by Subcase $1.4$ by taking $\{x_{0}(x_{3}),x_{1},x_{2},s\}$ ($s\in T_{x_{0}x_{1}}(T_{x_{1}x_{3}})$) as the vertices of a new $K_{4}-e$.
Since $G$ contains no subdivision of $K_{5}$, $\left|E(T_{x_{0}x_{2}x_{3}},T_{x_{1}})\right| =0$.
Since $R_{12}\not\subset G$, $\left|E(T_{x_{0}x_{2}x_{3}},c)\right| =1$.
Thus, if $\left|T_{x_{0}x_{2}x_{3}}\right|\geq 1$, then  its conclusion is covered by Subcase $1.2$ by taking $\{x_{0},x_{1},x_{2},w_{023,1}\}$ ($w_{023,1}\in T_{x_{0}x_{2}x_{3}}$) as the vertices of a new $K_{4}-e$.
Let $\left|T_{x_{0}x_{1}}\cup T_{x_{1}x_{3}}\cup T_{x_{0}x_{2}x_{3}}\right|=0$.

\begin{fact}\label{case9.5-f1}
$\left|E(T_{x_{1}},T_{x_{0}x_{3}}-\{c\})\right| =0$; at most one vertex in $T_{x_{1}}-\{q\}$ is adjacent to $c$; at most two vertices in $T_{x_{1}}-\{q\}$ are not adjacent to $c$.
\end{fact}
\begin{proof}
If $\left|E(T_{x_{1}},T_{x_{0}x_{3}}-\{c\})\right|  \geq1$, then $G$ contains a subdivision of $K_{3,3}$, a contradiction.
If at least two vertices in $T_{x_{1}}-\{q\}$ are adjacent to $c$, then $G[\mathcal{K}\cup\{q,\allowbreak p_{1,2},\allowbreak p_{1,3},\allowbreak c\}]\cong R_{81}$, a contradiction.
If at least three vertices in $T_{x_{1}}-\{q\}$ are not adjacent to $c$, then $G[\mathcal{K}\cup\{p_{1,2},\allowbreak p_{1,3},\allowbreak p_{1,4},\allowbreak c\}]\cong R_{38}$, a contradiction.
\end{proof}

\begin{fact}\label{case9.5-f4}
$\left|E(T_{x_{0}x_{3}}-\{c\})\right| =0$ and at most two vertices in $T_{x_{0}x_{3}}-\{c\}$ are not adjacent to $c$.
\end{fact}
\begin{proof}
If $\left|E(T_{x_{0}x_{3}}-\{c\})\right|  \geq1$, then $G[\{x_{0},\allowbreak x_{1},\allowbreak x_{2},\allowbreak q,\allowbreak q_{03,2},\allowbreak q_{03,3}\}]\cong R_{10}$ by Lemma \ref{lem-13}\ref{lem-13:item-4}, a contradiction.
If at least three vertices in $T_{x_{0}x_{3}}-\{c\}$ are not adjacent to $c$, then $G[\{x_{0},\allowbreak x_{1},\allowbreak x_{3},\allowbreak q,\allowbreak c,\allowbreak q_{03,2},\allowbreak q_{03,3},\allowbreak q_{03,4}\}]\cong R_{39}$, a contradiction.
\end{proof}

\begin{fact}\label{case9.5-f7}
At most one vertex in $T_{x_{2}}$ is not adjacent to $c$ and at most two vertices in $T_{x_{2}}$ are adjacent to $c$.
\end{fact}
\begin{proof}
If at least two vertices in $T_{x_{2}}$ are not adjacent to $c$, then $G[\{x_{0},\allowbreak x_{2},\allowbreak q,\allowbreak c,\allowbreak p_{2,1},\allowbreak p_{2,2}\}]\cong R_{1}$, a contradiction.
If at least three vertices in $T_{x_{2}}$ are adjacent to $c$, then $G[\mathcal{K}\cup\{c,\allowbreak p_{2,1},\allowbreak p_{2,2},\allowbreak p_{2,3}\}]\cong R_{81}$, a contradiction.
\end{proof}

\begin{fact}\label{case9.5-f9}
$\left|T_{x_{0}}\right| \leq 1$ and $\left|T_{x_{3}}\right| \leq 1$.
\end{fact}
\begin{proof}
If $\left|T_{x_{0}}\right|  \geq2$, then $G[\mathcal{K}\cup\{q,\allowbreak p_{0,1},\allowbreak p_{0,2}\}]\cong R_{22}$, a contradiction.
\end{proof}

\begin{fact}\label{case9.5-f11}
$\left|T_{x_{0}x_{2}}\cup T_{x_{2}x_{3}}\cup T_{x_{1}x_{2}}\cup T_{x_{0}x_{1}x_{3}}\cup T_{0}\right|=0$.
\end{fact}

\begin{proof}
If $\left|T_{x_{0}x_{2}}\cup T_{x_{2}x_{3}}\right|  \geq1$, then $G[\{x_{0},\allowbreak x_{1},\allowbreak x_{2},\allowbreak q,\allowbreak c,\allowbreak q_{02,1}(q_{23,1})\}]\cong R_{12}$, a contradiction.
If $\left|T_{x_{1}x_{2}}\right|  \geq1$, then either $G[\{x_{0},\allowbreak x_{1},\allowbreak x_{2},\allowbreak q,\allowbreak c,\allowbreak q_{12,1}\}]\cong R_{12}$ ($E(\{q,\allowbreak c,\allowbreak q_{12,1}\})=\{qc\}$) or $G$ contains a subdivision of $K_{3,3}$, a contradiction.
If $\left|T_{x_{0}x_{1}x_{3}}\right|  \geq1$, then  $G$ contains a subdivision of $K_{3,3}$, a contradiction.
If $\left|E(T_{0},T_{x_{0}x_{3}})\right|  \geq1$, then either $G[\{x_{0},\allowbreak x_{1},\allowbreak x_{2},\allowbreak q,\allowbreak q_{03,2},\allowbreak t_{1}\}]\cong R_{3}$ ($E(\{q,\allowbreak q_{03,2},\allowbreak t_{1}\})=\{q_{03,2}t_{1}\}$) or $G[\{x_{0},\allowbreak x_{1},\allowbreak x_{2},\allowbreak q,\allowbreak c,\allowbreak t_{1}\}]\cong R_{7}$ ($E(\{q,\allowbreak c,\allowbreak t_{1}\})=\{qc,\allowbreak ct_{1}\}$) by Lemma \ref{lem-13}\ref{lem-13:item-8}, a contradiction.
Thus, by Lemma \ref{lem-13}\ref{lem-13:item-8}, $\left|T_{0}\right|=0$.
\end{proof}

Suppose that $\left|E(T_{x_{0}x_{3}})\right|\geq 1$.
By  Lemma \ref{lem-13}\ref{lem-13:item-4} and (F\ref{case9.5-f4}), $\left|E(T_{x_{0}x_{3}})\right|= 1$   and $\left|E(T_{x_{0}x_{3}}-\{c\})\right| =0$.
Let $q_{03,1}\in T_{x_{0}x_{3}}-\{c\}$ and $c\sim q_{03,1}$.
Then $\left|E(T_{x_{2}},q_{03,1})\right|=\left|T_{x_{2}}\right|$ and $\left|E(T_{x_{2}},c)\right|=0$; otherwise, either $G[\{x_{0},\allowbreak x_{1},\allowbreak x_{2},\allowbreak c,\allowbreak q_{03,1},\allowbreak p_{2,1}\}]\cong R_{10}$ ($E(\{c,\allowbreak q_{03,1},\allowbreak p_{2,1}\})=\{cq_{03,1}\}$) or $G$ contains a subdivision of $K_{3,3}$, a contradiction.
Then $\left|T_{x_{2}}\right| \leq 1$; otherwise, $G[\{x_{0},\allowbreak x_{2},\allowbreak q,\allowbreak c,\allowbreak p_{2,1},\allowbreak p_{2,2}\}]\cong R_{1}$, a contradiction.
Then $\left|E(T_{x_{1}},c)\right| =\left|T_{x_{1}}\right|$ and $\left|T_{x_{1}}-\{q\}\right| \leq 1$ by (F\ref{case9.5-f1}); otherwise, $G[\{x_{0},\allowbreak x_{1},\allowbreak x_{2},\allowbreak p_{1,2},\allowbreak c,\allowbreak q_{03,1}\}]\cong R_{10}$, a contradiction.
Thus, $\left|T_{x_{1}}\cup T_{x_{2}}-\{q\}\right|\leq 1$; otherwise, $G[\{x_{1},\allowbreak x_{2},\allowbreak q,\allowbreak p_{1,2},\allowbreak q_{03,1},\allowbreak p_{2,1}\}]\cong R_{1}$, a contradiction.
If $\left|T_{x_{1}}\cup T_{x_{2}}-\{q\}\right|=1$, then $\left|T_{x_{0}x_{3}}-\{c\}\right| =1$; otherwise, either $G[\{x_{0},\allowbreak x_{1},\allowbreak x_{3},\allowbreak q,\allowbreak p_{1,2},\allowbreak q_{03,1},\allowbreak q_{03,3}\}]\cong R_{77}$ ($E(\{q,\allowbreak p_{1,2},\allowbreak q_{03,1},\allowbreak q_{03,3}\})=\varnothing$) or   $G$ contains a subdivision of $K_{3,3}$, a contradiction.
Thus,  $G\subset A_{44}$ (if $\left|T_{x_{2}}\right| =1$) or $G\subset A_{45}$ (if $\left|T_{x_{1}}-\{q\}\right| =1$) or $G\subset A_{46}$ (if $\left|T_{x_{1}}\cup T_{x_{2}}-\{q\}\right|=0$) by (F\ref{case9.5-f4}), (F\ref{case9.5-f9}) and (F\ref{case9.5-f11}).
Let $\left|E(T_{x_{0}x_{3}})\right|=0$.

Then $\left|E(T_{x_{2}},c)\right|+\left|T_{x_{1}}-T_{x_{1}}\cap N(c)\right|\leq 2$ and $\left|E(T_{x_{1}},c)\right|+\left|T_{x_{2}}-T_{x_{2}}\cap N(c)\right|\leq 2$ by (F\ref{case9.5-f1}) and (F\ref{case9.5-f7}); otherwise, $G[\{x_{1},\allowbreak x_{2},\allowbreak p_{1,2},\allowbreak p_{1,3},\allowbreak c,\allowbreak p_{2,1}\}]\cong R_{1}$ ($E(\{p_{1,2},\allowbreak p_{1,3},\allowbreak c,\allowbreak p_{2,1}\})=\{cp_{2,1}\}$) or $G[\{x_{0},\allowbreak x_{1},\allowbreak p_{1,2},\allowbreak c,\allowbreak p_{2,1},\allowbreak p_{2,2}\}]\cong R_{1}$ ($E(\{p_{1,2},\allowbreak c,\allowbreak p_{2,1},\allowbreak p_{2,2}\})=\{cp_{2,1},\allowbreak cp_{2,2}\}$) or $G[\{x_{0},\allowbreak x_{2},\allowbreak q,\allowbreak p_{1,2},\allowbreak c,\allowbreak p_{2,1}\}]\cong R_{1}$ ($E(\{q,\allowbreak p_{1,2},\allowbreak c,\allowbreak p_{2,1}\})=\{qc,\allowbreak p_{1,2}c\}$) or $G[\{x_{0},\allowbreak x_{2},\allowbreak q,\allowbreak c,\allowbreak p_{2,1},\allowbreak p_{2,2}\}]\cong R_{1}$ ($E(\{q,\allowbreak c,\allowbreak p_{2,1},\allowbreak p_{2,2}\})=\{qc\}$), a contradiction.
If  $\left|E(T_{x_{2}},c)\right|\geq 1$, then $\left|T_{x_{0}x_{3}}-\{c\}\right| =0$ for $G$ contains no subdivision of $K_{5}$.
Thus, if  $\left|E(T_{x_{2}},c)\right|\geq 1$, then  $G\subset A_{38}$ (if $\left|E(T_{x_{2}},c)\right|=2$ and $\left|E(T_{x_{1}},c)\right|=1$) or $G\subset A_{37}$ (if $\left|E(T_{x_{2}},c)\right|=2$ and $\left|E(T_{x_{1}},c)\right|=2$) or $G\subset A_{39}$ (if $\left|E(T_{x_{2}},c)\right|=1$ and $\left|E(T_{x_{1}},c)\right|=1$) or $G\subset A_{38}$ (if $\left|E(T_{x_{2}},c)\right|=1$ and $\left|E(T_{x_{1}},c)\right|=2$) by  (F\ref{case9.5-f9}) and (F\ref{case9.5-f11}).
Let $\left|E(T_{x_{2}},c)\right|=0$.
If $\left|E(T_{x_{1}},c)\right|=2$, then $\left|T_{x_{0}x_{3}}-\{c\}\right| =0$; otherwise, $G[\{x_{0},\allowbreak x_{2},\allowbreak x_{3},\allowbreak q,\allowbreak p_{1,2},\allowbreak c,\allowbreak q_{03,2}\}]\cong R_{77}$, a contradiction.
Thus, if  $\left|E(T_{x_{1}},c)\right|=2$, then $G\subset A_{41}$ by  (F\ref{case9.5-f9}) and (F\ref{case9.5-f11}).
Let $\left|E(T_{x_{1}},c)\right|=1$.
Then $\left|T_{x_{0}x_{3}}-\{c\}\right|+\left|T_{x_{1}}-T_{x_{1}}\cap N(c)\right|\leq 2$; otherwise, $G[\{x_{0},\allowbreak x_{1},\allowbreak x_{3},\allowbreak p_{1,2},\allowbreak p_{1,3},\allowbreak c,\allowbreak q_{03,2}\}]\cong R_{77}$, a contradiction.
Thus,  $G\subset A_{43}$ (if $\left|T_{x_{0}x_{3}}-\{c\}\right|=2$) or $G\subset A_{42}$ (if $\left|T_{x_{0}x_{3}}-\{c\}\right|=1$) or $G\subset A_{40}$ (if $\left|T_{x_{0}x_{3}}-\{c\}\right|=0$) by (F\ref{case9.5-f4}), (F\ref{case9.5-f9}) and (F\ref{case9.5-f11}).

If  $\left|E(T_{x_{2}},T_{x_{0}x_{3}})\right|\geq 1$, then we can obtain the desired conclusion by symmetry and an analogous proof.
In the following cases and subcases, let $\left|E(T_{x_{1}}\cup T_{x_{2}},T_{x_{0}x_{3}})\right|=0$.

\noindent$\textbf{Subcase 1.6}$ Suppose that $\left|E(T_{x_{1}},T_{x_{1}x_{2}})\right|\geq 1$.

Let  $q\in T_{x_{1}}$, $d\in T_{x_{1}x_{2}}$ and $q\sim d$.
If $\left|T_{x_{0}x_{1}}\right|\geq 1$, then  its conclusion is covered by Subcase $1.3$ by taking $\{x_{0},x_{1},x_{2},q_{01,1}\}$ ($q_{01,1}\in T_{x_{0}x_{1}}$) as the vertices of a new $K_{4}-e$.
If $\left|E(T_{x_{1}x_{2}}-\{d\},q)\right|\geq 1$, then  its conclusion is covered by Case $4$ of Lemma \ref{lem-43} by taking $\{x_{1},x_{2},x_{3},d\}$  as the vertices of a new $K_{4}-e$.
If $p_{1,1}\in T_{x_{1}}-\{q\}$, $q_{12,1}\in T_{x_{1}x_{2}}-\{d\}$, $p_{1,1}\sim q_{12,1}$ and $p_{1,1}\nsim d$, then  its conclusion is covered by Subcase $1.3$ by taking $\{x_{1},x_{2},q,d\}$  as the vertices of a new $K_{4}-e$.
If $p_{1,1}\in T_{x_{1}}-\{q\}$, $q_{12,1}\in T_{x_{1}x_{2}}-\{d\}$, $p_{1,1}\sim q_{12,1}$ and $p_{1,1}\sim d$, then  its conclusion is covered by Case $4$ of Lemma \ref{lem-43} by taking $\{x_{1},x_{2},q,d\}$  as the vertices of a new $K_{4}-e$.
By  symmetry, let $\left|T_{x_{0}x_{1}}\cup T_{x_{1}x_{3}}\right|=\left|E(T_{x_{1}x_{2}}-\{d\},T_{x_{1}})\right|=0$.

\begin{fact}\label{case9.6-f1}
At most six  vertices ($S_{1}$) in $T_{x_{1}}$ are not adjacent to $d$; at most one vertex  ($S_{2}$) in $T_{x_{1}}-\{q\}$ is adjacent to $d$; $\left|T_{x_{1}x_{2}}\cup S_{1}-\{d\}\right|\leq 6$; if $\left|S_{2}\right|=1$, then $\left|T_{x_{1}x_{2}}-\{d\}\right|=0$.
\end{fact}
\begin{proof}
If at least seven vertices in $T_{x_{1}}-\{q\}$  are not adjacent to $d$, then $G[\mathcal{K}\cup\{p_{1,2},\ldots,p_{1,8},\allowbreak d\}]\cong R_{86}$, a contradiction.
If at least two vertices in $T_{x_{1}}-\{q\}$ are adjacent to $d$, then $G[\{x_{0},\allowbreak x_{2},\allowbreak x_{3},\allowbreak q,\allowbreak p_{1,2},\allowbreak p_{1,3},\allowbreak d\}]\cong R_{21}$, a contradiction.
If $\left|T_{x_{1}x_{2}}-\{d\}\right|\geq 7$, then $G[\mathcal{K}\cup\{q,\allowbreak d,\allowbreak q_{12,2},\ldots,q_{12,8}\}]\cong R_{70}$, a contradiction.
If $\left|T_{x_{1}x_{2}}\cup S_{1}-\{d\}\right|\geq 7$, then  $R_{90}\subset G$ (if $1\leq \left|T_{x_{1}x_{2}}-\{d\}\right|\leq 2$) or $R_{96}\subset G$ (if $\left|T_{x_{1}x_{2}}-\{d\}\right|=3$) or $R_{89}\subset G$ (if $\left|T_{x_{1}x_{2}}-\{d\}\right|=4$) or $R_{88}\subset G$ (if $\left|T_{x_{1}x_{2}}-\{d\}\right|=5$) or $R_{92}\subset G$ (if $\left|T_{x_{1}x_{2}}-\{d\}\right|=6$), a contradiction.
If $\left|T_{x_{1}x_{2}}-\{d\}\right|\geq 1$, then  $R_{21}\subset G$, a contradiction.
\end{proof}

\begin{fact}\label{case9.6-f5}
$\left|T_{x_{0}}\cup T_{x_{3}}\cup T_{x_{0}x_{3}}\cup T_{x_{0}x_{2}}\cup T_{x_{2}x_{3}}\cup T_{3}\cup T_{0}\right|= 0$.
\end{fact}
\begin{proof}
If $\left|T_{x_{0}}\cup T_{x_{3}}\cup T_{x_{0}x_{3}}\right|  \geq1$, then $G[\{x_{0},\allowbreak x_{1},\allowbreak x_{2},\allowbreak q,\allowbreak d,\allowbreak p_{0,1}(p_{3,1},q_{03,1})\}]\cong R_{14}$, a contradiction.
If $\left|T_{x_{0}x_{2}}\cup T_{x_{2}x_{3}}\right|  \geq1$, then $G[\{x_{0},\allowbreak x_{1},\allowbreak x_{2},\allowbreak q,\allowbreak d,\allowbreak q_{02,1}(q_{23,1})\}]\cong R_{18}$, a contradiction.
If $\left|T_{x_{0}x_{1}x_{3}}\right|  \geq1$, then either $G[\{x_{0},\allowbreak x_{1},\allowbreak x_{3},\allowbreak q,\allowbreak d,\allowbreak w_{013,1}\}]\cong R_{16}$ ($E(\{q,\allowbreak d,\allowbreak w_{013,1}\})=\{qd\}$) or $G$ contains a subdivision of $K_{3,3}$, a contradiction.
If $\left|T_{x_{0}x_{2}x_{3}}\right|  \geq1$, then either $G[\{x_{0},\allowbreak x_{1},\allowbreak x_{2},\allowbreak q,\allowbreak d,\allowbreak w_{023,1}\}]\cong R_{18}$ ($E(\{q,\allowbreak d,\allowbreak w_{023,1}\})=\{qd\}$) or $G$ contains a subdivision of $K_{3,3}$, a contradiction.
If $\left|E(T_{0},d)\right|  \geq1$, then $G[\mathcal{K}\cup\{q,\allowbreak d,\allowbreak t_{1}\}]\cong R_{28}$ by Lemma \ref{lem-13}\ref{lem-13:item-8}, a contradiction.
If $\left|E(T_{0},T_{x_{1}x_{2}}-\{d\})\right|  \geq1$, then $G[\{x_{1},\allowbreak x_{2},\allowbreak q,\allowbreak d,\allowbreak q_{12,2},\allowbreak t_{1}\}]\cong R_{14}$ by Lemma \ref{lem-13}\ref{lem-13:item-8}, a contradiction.
Thus, by Lemma \ref{lem-13}\ref{lem-13:item-8}, $\left|T_{0}\right|= 0$.
\end{proof}

\begin{fact}\label{case9.6-f6}
At most one vertex in $T_{x_{2}}$ is not adjacent to $d$; at most two vertices in $T_{x_{2}}$ are adjacent to $d$.
\end{fact}
\begin{proof}
If at least two vertices in $T_{x_{2}}$ are not adjacent to $d$, then $G[\{x_{0},\allowbreak x_{1},\allowbreak x_{2},\allowbreak q,\allowbreak d,\allowbreak p_{2,1},\allowbreak p_{2,2}\}]\cong R_{25}$, a contradiction.
If at least three vertices in $T_{x_{2}}$  are adjacent to $d$, then $G[\{x_{0},\allowbreak x_{1},\allowbreak x_{3},\allowbreak d,\allowbreak p_{2,1},\allowbreak p_{2,2},\allowbreak p_{2,3}\}]\cong R_{21}$, a contradiction.
\end{proof}

If $\left|T_{x_{2}}\right|=0$, then $G\subset H_{120}$ (if $\left|T_{x_{1}x_{2}}-\{d\}\right|=0$) or $G\subset H_{131}$ (if $\left|T_{x_{1}x_{2}}-\{d\}\right|=1$) or $G\subset H_{151}$ (if $\left|T_{x_{1}x_{2}}-\{d\}\right|\geq 2$) by (F\ref{case9.6-f1}) and (F\ref{case9.6-f5}).
Let $\left|T_{x_{2}}\right|\geq 1$ and $p_{2,1}\in T_{x_{2}}$.
Let $S_{1}=T_{x_{1}}-(N(d)\cap T_{x_{1}})$ and $S_{2}=N(d)\cap T_{x_{1}}-\{q\}$.
If $\left|S_{2}\right| =1$, then $\left|E(T_{x_{2}},T_{x_{1}x_{2}})\right| =\left|T_{x_{2}}\right|$ for $R_{28}\not\subset G$.
If $\left|S_{2}\right| =1$, then $\left|S_{1}\cup T_{x_{2}}\right| \leq 2$ by (F\ref{case9.6-f6}); otherwise either $G[\{x_{1},x_{2},x_{3}\}\cup S_{1}\cup T_{x_{2}}\cup T_{x_{1}x_{2}}]\cong R_{14}$ (if $\left|S_{1}\right|=2$) or $G[\{x_{1},x_{2},x_{3}\}\cup S_{1}\cup T_{x_{2}}\cup T_{x_{1}x_{2}}]\cong R_{28}$ (if $\left|S_{1}\right|=1$), a contradiction.
Thus, if $\left|S_{2}\right| =1$, then either $G\subset A_{2}$ (if $\left|S_{1}\right|=0$) or $G\subset A_{3}$ (if $\left|S_{1}\right|=1$) by (F\ref{case9.6-f1}), (F\ref{case9.6-f5}) and (F\ref{case9.6-f6}).
Let $\left|S_{2}\right|  =0$.
If $\left|E(p_{2,1},T_{x_{1}x_{2}}-\{d\})\right|\geq 1$, then $p_{2,1}\sim d$; otherwise, $G[\{x_{1},\allowbreak x_{2},\allowbreak q,\allowbreak d,\allowbreak q_{12,2},\allowbreak p_{2,1}\}]\cong R_{18}$, a contradiction.
Thus, $\left|E(T_{x_{2}},T_{x_{1}x_{2}}-\{d\})\right|\leq 1$ for $G$ contains no subdivision of $K_{3,3}$.
If $\left|E(T_{x_{2}},T_{x_{1}x_{2}}-\{d\})\right|=1$, then $\left|T_{x_{1}}-\{q\}\right| =0$; otherwise, $G[\{x_{0},\allowbreak x_{1},\allowbreak x_{2},\allowbreak q,\allowbreak p_{1,2},\allowbreak q_{12,2},\allowbreak p_{2,1}\}]\cong R_{25}$, a contradiction.
If $\left|E(T_{x_{2}},T_{x_{1}x_{2}}-\{d\})\right|=1$, then  $\left|E(T_{x_{2}},d)\right|+\left|T_{x_{1}x_{2}}-\{d\}\right| \leq 4$ by (F\ref{case9.6-f6}); otherwise, either $R_{21}\subset G$ (if  $\left|E(T_{x_{2}},d)\right|=2$) or $R_{37}\subset G$ (if  $\left|E(T_{x_{2}},d)\right|=1$), a contradiction.
Thus, if $\left|E(T_{x_{2}},T_{x_{1}x_{2}}-\{d\})\right|=1$, then either $G\subset H_{115}$ (if  $\left|E(T_{x_{2}},d)\right|=2$) or $G\subset H_{10}$ (if  $\left|E(T_{x_{2}},d)\right|=1$) by (F\ref{case9.6-f1}), (F\ref{case9.6-f5}) and (F\ref{case9.6-f6}).
Let $\left|E(T_{x_{2}},T_{x_{1}x_{2}}-\{d\})\right|=0$.
Then $\left|E(T_{x_{2}},d)\right|+\left|T_{x_{1}}-\{q\}\right|\leq 2$  by (F\ref{case9.6-f6}); otherwise, either $R_{25}\subset G$ (if  $\left|E(T_{x_{2}},d)\right|\geq 1$) or $R_{47}\subset G$ (if  $\left|E(T_{x_{2}},d)\right|=0$), a contradiction.
Then $\left|E(T_{x_{2}},d)\right|+\left|T_{x_{1}x_{2}}-\{d\}\right|\leq 2$  by (F\ref{case9.6-f6}); otherwise, $R_{91}\subset G$ (if  $\left|E(T_{x_{2}},d)\right|=0$) or $R_{61}\subset G$ (if  $\left|E(T_{x_{2}},d)\right|=1$) or $R_{21}\subset G$ (if  $\left|E(T_{x_{2}},d)\right|=2$), a contradiction.
Then $\left|T_{x_{1}}-\{q\}\right|+\left|T_{x_{1}x_{2}}-\{d\}\right|\leq 2$  by (F\ref{case9.6-f6}); otherwise, $R_{47}\subset G$ (if  $\left|T_{x_{1}x_{2}}-\{d\}\right|=0$) or $R_{93}\subset G$ (if $\left|T_{x_{1}x_{2}}-\{d\}\right|=1$) or $R_{87}\subset G$ (if  $\left|T_{x_{1}x_{2}}-\{d\}\right|=2$) or $R_{91}\subset G$ (if  $\left|T_{x_{1}x_{2}}-\{d\}\right|\geq 3$), a contradiction.
Thus,  $G\subset A_{3}$ (if $\left|E(T_{x_{2}},d)\right| =2$) or $G\subset H_{11}$ (if $\left|E(T_{x_{2}},d)\right| =1$) or $G\subset H_{10}$ (if $\left|E(T_{x_{2}},d)\right| =0$ and $\left|T_{x_{1}x_{2}}-\{d\}\right|\geq 1$) or $G\subset A_{6}$ (if $\left|E(T_{x_{2}},d)\right| =0$ and $\left|T_{x_{1}x_{2}}-\{d\}\right|=0$) or by  (F\ref{case9.6-f5}) and (F\ref{case9.6-f6}).

If  $\left|E(T_{x_{2}},T_{x_{1}x_{2}})\right|\geq 1$, then we can obtain the desired conclusion by symmetry and an analogous proof.
In the following cases and subcases, let $\left|E(T_{1},T_{2})\right|=0$.

\noindent$\textbf{Case 2.}$ Suppose that $\left|E(T_{1},T_{3})\right|\geq 1$.

If $\left|E(T_{x_{0}},T_{x_{0}x_{1}x_{3}})\right|\geq 1$, then  its conclusion is covered by Case $1$ of Lemma \ref{lem-43} by taking $\{x_{0},x_{1},p_{0,1},w_{013,1}\}$ ($p_{0,1}\in T_{x_{0}}$, $w_{013,1}\in T_{x_{0}x_{1}x_{3}}$ and $p_{0,1}\sim w_{013,1}$) as the vertices of a new $K_{4}-e$.
If $\left|E(T_{x_{1}},T_{x_{0}x_{1}x_{3}})\right|\geq 1$, then its conclusion is covered by Case $4$ of Lemma \ref{lem-43} by taking $\{x_{1},x_{3},p_{1,1},w_{013,1}\}$ ($p_{1,1}\in T_{x_{1}}$, $w_{013,1}\in T_{x_{0}x_{1}x_{3}}$ and $p_{1,1}\sim w_{013,1}$) as the vertices of a new $K_{4}-e$.
By  Lemma \ref{lem-13}\ref{lem-13:item-7} and symmetry, let $\left|E(T_{x_{0}}\cup T_{x_{3}},T_{3})\right|=\left|E(T_{x_{1}},T_{x_{0}x_{1}x_{3}})\right|=\left|E(T_{x_{2}},T_{x_{0}x_{2}x_{3}})\right|=0$.
Thus, without loss of generality, let $\left|E(T_{x_{2}},T_{x_{0}x_{1}x_{3}})\right|\geq 1$  by  Lemma \ref{lem-13}\ref{lem-13:item-7}.
Let  $y\in T_{x_{2}}$, $u\in T_{x_{0}x_{1}x_{3}}$ and $y\sim u$.

If $\left|E(T_{x_{0}x_{2}},T_{x_{0}x_{1}x_{3}})\right|<\left|T_{x_{0}x_{2}}\right|$, then  its conclusion is covered by Case $1$ of Lemma \ref{lem-43} by taking $\{x_{0},x_{1},x_{2},q_{02,1}\}$ ($q_{02,1}\in T_{x_{0}x_{2}}$, $w_{013,1}\in T_{x_{0}x_{1}x_{3}}$ and $q_{02,1}\nsim w_{013,1}$) as the vertices of a new $K_{4}-e$.
Since $K_{4}\not\subset G$, $\left|E(T_{x_{0}x_{1}},T_{x_{0}x_{1}x_{3}})\right|=0$.
Thus, if $\left|T_{x_{0}x_{1}}\right|\geq 1$, then  its conclusion is covered by Subcase $1.6$ by taking $\{x_{0},x_{1},x_{2},q_{01,1}\}$ ($q_{01,1}\in T_{x_{0}x_{1}}$) as the vertices of a new $K_{4}-e$.
Since $G$ contains no subdivision of $K_{5}$, $\left|E(T_{x_{1}x_{2}},T_{x_{0}x_{1}x_{3}})\right|=0$.
Thus, if $\left|T_{x_{1}x_{2}}\right|\geq 1$, then  its conclusion is covered by Case $4$ of Lemma \ref{lem-43} by taking $\{x_{0},x_{1},x_{2},q_{12,1}\}$ ($q_{12,1}\in T_{x_{1}x_{2}}$) as the vertices of a new $K_{4}-e$.
By symmetry, let $\left|T_{x_{0}x_{1}}\cup T_{x_{1}x_{3}}\cup T_{x_{1}x_{2}}\right|=0$, $\left|E(T_{x_{0}x_{2}},T_{x_{0}x_{1}x_{3}})\right|=\left|T_{x_{0}x_{2}}\right|$ and $\left|E(T_{x_{2}x_{3}},T_{x_{0}x_{1}x_{3}})\right|=\left|T_{x_{2}x_{3}}\right|$.

\begin{fact}\label{case10.2-f1}
$\left|T_{x_{0}x_{3}}\cup T_{x_{0}x_{2}x_{3}}\right|=0$.
\end{fact}
\begin{proof}
If $\left|T_{x_{0}x_{3}}\cup T_{x_{0}x_{2}x_{3}}\right| \geq1$, then $G$ contains a subdivision of $K_{5}$, a contradiction.
\end{proof}

\begin{fact}\label{case10.2-f2}
At most one vertex in $T_{x_{2}}$ is not adjacent to $u$; at most seven vertices in $T_{x_{2}}$ are adjacent to $u$.
\end{fact}
\begin{proof}
If at least two vertices in $T_{x_{2}}$ are not adjacent to $u$, then $G[\mathcal{K}\cup\{p_{2,2},\allowbreak p_{2,3},\allowbreak u\}]\cong R_{27}$, a contradiction.
If at least eight vertices in $T_{x_{2}}$  are adjacent to $u$, then $G[\mathcal{K}\cup\{y,\allowbreak p_{2,2},\ldots,p_{2,8},\allowbreak u\}]\cong R_{71}$, a contradiction.
\end{proof}

\begin{fact}\label{case10.2-f4}
$\left|T_{x_{0}}\right|\leq 1$, $\left|T_{x_{3}}\right|\leq 1$, $\left|T_{x_{1}}\right|\leq 2$, $\left|T_{x_{0}x_{2}}\cup T_{x_{3}} \right|\leq 1$   and  $\left|T_{x_{2}x_{3}}\cup T_{x_{0}}\right|\leq 1$.
\end{fact}
\begin{proof}
If $\left|T_{x_{0}}\right| \geq2$ or $\left|T_{x_{3}}\right| \geq2$, then $G[\mathcal{K}\cup\{u,\allowbreak p_{0,1},\allowbreak p_{0,2}\}]\cong R_{27}$, a contradiction.
If $\left|T_{x_{1}}\right| \geq3$, then $G[\{x_{0},\allowbreak x_{1},\allowbreak x_{2},\allowbreak y,\allowbreak u,\allowbreak p_{1,1},\allowbreak p_{1,2},\allowbreak p_{1,3}\}]\cong R_{38}$, a contradiction.
If $\left|T_{x_{0}x_{2}}\right| \geq2$ or $\left|T_{x_{2}x_{3}}\right| \geq2$, then $G$ contains a subdivision of $K_{3,3}$, a contradiction.
If $\left|T_{x_{0}x_{2}}\cup T_{x_{3}} \right| \geq2$, then $G[\mathcal{K}\cup\{q_{02,1},\allowbreak p_{3,1}\}]\cong R_{14}$, a contradiction.
\end{proof}

Suppose that $\left|T_{x_{0}x_{2}}\right|\geq1$. 
By (F\ref{case10.2-f4}),  $\left|T_{x_{0}x_{2}}\right|=1$ and $\left|T_{x_{3}}\right|=0$.
Then at most six vertices in $T_{x_{2}}$ are adjacent to $u$; otherwise, $G[\mathcal{K}\cup\{y,\allowbreak p_{2,2},\ldots,p_{2,7},\allowbreak u,\allowbreak q_{02,1}\}]\cong R_{73}$, a contradiction.
Then $\left|T_{x_{1}}\right|\leq 1$; otherwise, $G[\mathcal{K}\cup\{q_{02,1},\allowbreak p_{1,1},\allowbreak p_{1,2}\}]\cong R_{25}$, a contradiction.
If  $\left|E(T_{x_{2}},u)\right| \geq 3$, then  $\left|T_{x_{0}}\cup T_{x_{1}}\right|=0$; otherwise, $G[\{x_{0},\allowbreak x_{2},\allowbreak x_{3},\allowbreak y,\allowbreak p_{2,2},\allowbreak p_{2,3},\allowbreak u,\allowbreak p_{0,1},\allowbreak q_{02,1}\}]\cong R_{55}$  or $G[\{x_{1},\allowbreak x_{2},\allowbreak x_{3},\allowbreak y,\allowbreak p_{2,2},\allowbreak p_{2,3},\allowbreak u,\allowbreak p_{1,1},\allowbreak q_{02,1}\}]\cong R_{55}$, a contradiction.
Thus, either  $G\subset A_{20}$ (if $\left|E(T_{x_{2}},u)\right| \geq 3$) or $G\subset A_{23}$ (if $\left|E(T_{x_{2}},u)\right| \leq 2$ and $\left|T_{x_{2}x_{3}}\right|=0$) or $G\subset A_{19}$ (if $\left|E(T_{x_{2}},u)\right| \leq 2$ and $\left|T_{x_{2}x_{3}}\right|=1$) by  Lemma \ref{lem-13}\ref{lem-13:item-8}, (F\ref{case10.2-f1}), (F\ref{case10.2-f2}) and (F\ref{case10.2-f4}).
By symmetry, let $\left|T_{x_{0}x_{2}}\cup T_{x_{2}x_{3}}\right|=0$.

If $\left|E(T_{x_{2}},u)\right| \geq 4$, then  $\left|T_{x_{0}}\cup T_{x_{3}}\right|=0$; otherwise, $G[\{x_{0},\allowbreak x_{1},\allowbreak x_{2},\allowbreak y,\allowbreak p_{2,2},\allowbreak p_{2,3},\allowbreak p_{2,4},\allowbreak u,\allowbreak p_{0,1}\}]\cong R_{55}$ ($E(\{y,\allowbreak p_{2,2},\allowbreak p_{2,3},\allowbreak p_{2,4},\allowbreak u,\allowbreak p_{0,1}\})=\{yu,\allowbreak p_{2,2}u,\allowbreak p_{2,3}u,\allowbreak p_{2,4}u\}$) or $G[\{x_{1},\allowbreak x_{2},\allowbreak x_{3},\allowbreak y,\allowbreak p_{2,2},\allowbreak p_{2,3},\allowbreak p_{2,4},\allowbreak u,\allowbreak p_{3,1}\}]\cong R_{55}$ ($E(\{y,\allowbreak p_{2,2},\allowbreak p_{2,3},\allowbreak p_{2,4},\allowbreak u,\allowbreak p_{3,1}\})=\{yu,\allowbreak p_{2,2}u,\allowbreak p_{2,3}u,\allowbreak p_{2,4}u\}$), a contradiction.
If $\left|E(T_{x_{2}},u)\right| \geq 4$, then  $\left|T_{x_{1}}\right|=0$; otherwise, $G[\mathcal{K}\cup\{y,\allowbreak p_{2,2},\allowbreak p_{2,3},\allowbreak p_{2,4},\allowbreak p_{1,1}\}]\cong R_{47}$, a contradiction.
Thus, if $\left|E(T_{x_{2}},u)\right| \geq 4$, then $G\subset A_{20}$ by  Lemma \ref{lem-13}\ref{lem-13:item-8}, (F\ref{case10.2-f1}) and (F\ref{case10.2-f2}).
If $\left|E(T_{x_{2}},u)\right| \geq 2$, then $\left|T_{x_{1}}\right|\leq 1$; otherwise, $G[\{x_{1},\allowbreak x_{2},\allowbreak y,\allowbreak p_{2,2},\allowbreak u,\allowbreak p_{1,1},\allowbreak p_{1,2}\}]\cong R_{77}$, a contradiction.
If $\left|E(T_{x_{2}},u)\right| = 3$ and $\left|T_{x_{1}}\right|=1$, then $\left|T_{x_{0}}\cup T_{x_{3}}\cup T_{0}\right|=0$  and $\left|T_{x_{2}}\right|=3$; otherwise, $G[\mathcal{K}\cup\{y,\allowbreak p_{2,2},\allowbreak p_{2,3},\allowbreak p_{2,4},\allowbreak p_{1,1}\}]\cong R_{47}$ ($E(\{y,\allowbreak p_{2,2},\allowbreak p_{2,3},\allowbreak p_{2,4},\allowbreak p_{1,1}\})=\varnothing$) or $G[\{x_{0},\allowbreak x_{1},\allowbreak x_{3},\allowbreak y,\allowbreak p_{2,2},\allowbreak p_{2,3},\allowbreak u,\allowbreak p_{1,1},\allowbreak t_{1}\}]\cong R_{47}$ ($E(\{y,\allowbreak p_{2,2},\allowbreak p_{2,3},\allowbreak u,\allowbreak p_{1,1},\allowbreak t_{1}\})=\{yu,\allowbreak p_{2,2}u,\allowbreak p_{2,3}u,\allowbreak ut_{1}\}$) or $G[\{x_{0},\allowbreak x_{1},\allowbreak x_{2},\allowbreak y,\allowbreak p_{2,2},\allowbreak p_{2,3},\allowbreak u,\allowbreak p_{0,1},\allowbreak p_{1,1}\}]\cong R_{94}$ ($E(\{y,\allowbreak p_{2,2},\allowbreak p_{2,3},\allowbreak u,\allowbreak p_{0,1},\allowbreak p_{1,1}\})=\{yu,\allowbreak p_{2,2}u,\allowbreak p_{2,3}u\}$) or $G[\{x_{1},\allowbreak x_{2},\allowbreak x_{3},\allowbreak y,\allowbreak p_{2,2},\allowbreak p_{2,3},\allowbreak u,\allowbreak p_{3,1},\allowbreak p_{1,1}\}]\cong R_{94}$ ($E(\{y,\allowbreak p_{2,2},\allowbreak p_{2,3},\allowbreak u,\allowbreak p_{3,1},\allowbreak p_{1,1}\})=\{yu,\allowbreak p_{2,2}u,\allowbreak p_{2,3}u\}$), a contradiction.
Thus, if $\left|E(T_{x_{2}},u)\right| = 3$,  then either $G\subset A_{26}$ (if $\left|T_{x_{1}}\right|=0$) or $G\subset H_{109}$ (if $\left|T_{x_{1}}\right|=1$) by  Lemma \ref{lem-13}\ref{lem-13:item-8}, (F\ref{case10.2-f1}) and (F\ref{case10.2-f2}).
If $\left|E(T_{x_{2}},u)\right| \leq 2$, then either $G\subset A_{29}$ (if $\left|E(T_{x_{2}},u)\right|= 2$) or $G\subset A_{32}$ (if $\left|E(T_{x_{2}},u)\right|\leq 1$) by  Lemma \ref{lem-13}\ref{lem-13:item-8}, (F\ref{case10.2-f1}), (F\ref{case10.2-f2}) and (F\ref{case10.2-f4}).

In the following cases and subcases, let $\left|E(T_{1},T_{3})\right|=0$.

\noindent$\textbf{Case 3.}$ Suppose that $\left|E(T_{2},T_{3})\right|\geq 1$.

Since $\omega\leq 3$, $\left|E(T_{x_{0}x_{1}}\cup T_{x_{1}x_{3}},T_{x_{0}x_{1}x_{3}})\right|=0$.
Since  $G$ contains no subdivision of $K_{3,3}$, $\left|E(T_{x_{1}x_{2}},T_{x_{0}x_{1}x_{3}})\right|=0$.
By Lemma \ref{lem-13}\ref{lem-13:item-7} and symmetry, we distinguish the following two subcases.

Suppose that $\left|E(T_{x_{0}x_{2}},T_{x_{0}x_{1}x_{3}})\right|\geq 1$.
Let  $b\in T_{x_{0}x_{2}}$, $u\in T_{x_{0}x_{1}x_{3}}$ and $b\sim u$.
If $\left|T_{x_{0}x_{1}}\cup T_{x_{1}x_{3}}\cup T_{x_{1}x_{2}}\right|\geq 1$, then  its conclusion is covered by Case $4$ of Lemma \ref{lem-43} by taking $\{x_{0},x_{1},x_{2},s\}$ ($s\in T_{x_{0}x_{1}}\cup T_{x_{1}x_{2}}$) as the vertices of a new $K_{4}-e$.
Then $\left|E(T_{x_{2}x_{3}},T_{x_{0}x_{1}x_{3}})\right|=\left|T_{x_{2}x_{3}}\right|$; otherwise, $G[\{x_{0},\allowbreak x_{1},\allowbreak x_{3},\allowbreak b,\allowbreak u,\allowbreak q_{23,1}\}]\cong R_{14}$, a contradiction.
Thus, if $\left|T_{x_{2}x_{3}}\right|\geq 1$, then  its conclusion is covered by Case $2$ by taking $\{x_{2},x_{3},q_{23,1},u\}$ ($q_{23,1}\in T_{x_{2}x_{3}}$) as the vertices of a new $K_{4}-e$.
Since $G$ contains no subdivision of $K_{3,3}$, $\left|E(T_{x_{0}x_{2}},u)\right|= 1$.
Thus, if $\left|T_{x_{0}x_{2}}\right|\geq 2$, then  its conclusion is covered by Case $1$ of Lemma \ref{lem-43} by taking $\{x_{0},x_{1},x_{2},q_{02,1}\}$ ($q_{02,1}\in T_{x_{0}x_{2}}$ and $q_{02,1}\nsim u$) as the vertices of a new $K_{4}-e$.
Since $G$ contains no subdivision of $K_{3,3}$, $\left|T_{x_{0}x_{3}}\right|=0$.
Let  $\left|T_{2}-\{b\}\right|=0$.
Then $\left|T_{x_{3}}\right|=0$; otherwise, $G[\mathcal{K}\cup\{b,\allowbreak p_{3,1}\}]\cong R_{14}$, a contradiction.
Then $\left|T_{x_{0}}\right|\leq 1$; otherwise, $G[\mathcal{K}\cup\{u,\allowbreak p_{0,1},\allowbreak p_{0,2}\}]\cong R_{27}$, a contradiction.
Then $\left|T_{x_{1}}\right| \leq 1$; otherwise, $G[\mathcal{K}\cup\{b,\allowbreak p_{1,1},\allowbreak p_{1,2}\}]\cong R_{25}$, a contradiction.
Then $\left|T_{x_{2}}\right| \leq 1$; otherwise, $G[\{x_{0},\allowbreak x_{1},\allowbreak x_{2},\allowbreak b,\allowbreak u,\allowbreak p_{2,1},\allowbreak p_{2,2}\}]\cong R_{27}$, a contradiction.
Thus,  $G\subset A_{21}$ by Lemma \ref{lem-13}\ref{lem-13:item-8}.
Let $\left|E(T_{x_{0}x_{2}}\cup T_{x_{2}x_{3}},T_{3})\right|=0$.

Suppose that $\left|E(T_{x_{0}x_{3}},T_{x_{0}x_{1}x_{3}})\right|\geq 1$.
Let  $c\in T_{x_{0}x_{3}}$, $u\in T_{x_{0}x_{1}x_{3}}$ and $c\sim u$.
If $\left|T_{x_{0}x_{1}}\right|\geq 1$, then  its conclusion is covered by Case $4$ of Lemma \ref{lem-43} by taking $\{x_{0},x_{1},x_{2},q_{01,1}\}$ ($q_{01,1}\in T_{x_{0}x_{1}}$) as the vertices of a new $K_{4}-e$.
Let  $\left|T_{x_{0}x_{1}}\cup T_{x_{1}x_{3}}\right|=0$.
Then $\left|T_{x_{2}}\cup T_{x_{0}x_{2}}\cup T_{x_{2}x_{3}}\right|=0$; otherwise, $G[\{x_{0},\allowbreak x_{1},\allowbreak x_{2},\allowbreak c,\allowbreak u,\allowbreak p_{2,1}(q_{02,1},q_{23,1})\}]\cong R_{14}$, a contradiction.
Then $\left|T_{x_{1}x_{2}}\right|=0$; otherwise, $G[\{x_{0},\allowbreak x_{1},\allowbreak x_{2},\allowbreak c,\allowbreak u,\allowbreak q_{12,1}\}]\cong R_{18}$, a contradiction.
Then $\left|T_{x_{0}}\right|\leq 1$ and $\left|T_{x_{3}}\right|\leq 1$; otherwise, $G[\mathcal{K}\cup\{u,\allowbreak p_{0,1},\allowbreak p_{0,2}\}]\cong R_{27}$, a contradiction.
Then $\left|T_{x_{1}}\right|\leq 1$; otherwise, $G[\{x_{0},\allowbreak x_{1},\allowbreak x_{2},\allowbreak c,\allowbreak u,\allowbreak p_{1,1},\allowbreak p_{1,2}\}]\cong R_{25}$, a contradiction.
If $\left|E(T_{0},c)\right|\geq 1$, then either $R_{14}\subset G$ (if $\left|E(T_{0},u)\right|= 0$) or $R_{18}\subset G$ (the remaining case) by Lemma \ref{lem-13}\ref{lem-13:item-8}, a contradiction.
If $\left|E(T_{0},T_{x_{0}x_{3}}-\{c\})\right|\geq 1$, then either $R_{9}\subset G$ (if $\left|E(T_{0},u)\right|= 0$) or $R_{12}\subset G$ (the remaining case) by Lemma \ref{lem-13}\ref{lem-13:item-8}, a contradiction.
Thus, $\left|E(T_{0},T_{x_{0}x_{3}})\right|=0$.
Since $G$ contains no subdivision of $K_{3,3}$, $\left|E(T_{x_{0}x_{3}}-\{c\},u)\right| =0$.
Then $\left|E(T_{x_{0}x_{3}}-\{c\})\right| =0$ and  $\left|E(T_{x_{0}x_{3}}-\{c\}, c)\right| \leq 1$ by  Lemma \ref{lem-13}\ref{lem-13:item-4}; otherwise, $G[\{x_{0},\allowbreak x_{1},\allowbreak c,\allowbreak q_{03,2},\allowbreak q_{03,3},\allowbreak u\}]\cong R_{16}$ ($E(\{c,\allowbreak q_{03,2},\allowbreak q_{03,3},\allowbreak u\})=\{cu,\allowbreak q_{03,2}q_{03,3}\}$) or $G[\{x_{0},\allowbreak x_{1},\allowbreak x_{2},\allowbreak c,\allowbreak q_{03,2},\allowbreak q_{03,3}\}]\cong R_{16}$ ($E(\{c,\allowbreak q_{03,2},\allowbreak q_{03,3}\})=\{cq_{03,2},\allowbreak cq_{03,3}\}$), a contradiction.
Since $R_{73}\not\subset G$, at most six vertices ($S_{1}$) in $T_{x_{0}x_{3}}-\{c\}$ are not adjacent to $c$.
If $\left|T_{0}\cup T_{x_{1}}\right|\geq 1$, then $\left|S_{1}\right| \leq 2$ by Lemma \ref{lem-13}\ref{lem-13:item-8} and $R_{55}\not\subset G$.
If $\left|T_{x_{1}}\right|\geq 1$, then $\left|E(T_{x_{0}x_{3}}-\{c\}, c)\right| =0$; otherwise, $G[\{x_{0},\allowbreak x_{1},\allowbreak x_{2},\allowbreak c,\allowbreak q_{03,2},\allowbreak p_{1,1}\}]\cong R_{10}$, a contradiction.
Thus, either $G\subset A_{19}$ (if $\left|T_{0}\cup T_{x_{1}}\right| \geq 1$) or $G\subset A_{20}$ (if $\left|T_{0}\cup T_{x_{1}}\right| =0$) by Lemma \ref{lem-13}\ref{lem-13:item-8}.
\end{proof}

\begin{thm}\label{thm-41}
If $G \in \mathcal{G}_3$ and contains $K_{4}-e$ as an induced subgraph, then $G$ is an induced subgraph of a graph in $\{H_{4},\allowbreak H_{10},\allowbreak H_{11},\allowbreak H_{71},\allowbreak \ldots,\allowbreak H_{153},\allowbreak A_{3},\allowbreak A_{6},\allowbreak A_{15},\allowbreak \ldots,\allowbreak A_{32},\allowbreak A_{34}, \allowbreak A_{35},\allowbreak A_{37},\allowbreak \ldots,\allowbreak A_{47},\allowbreak  D_1,\allowbreak D_2\}$.
\end{thm}
\begin{proof}
Let $K_{4}-e$ be an induced subgraph of $G$ with $\mathcal{K}=\{x_{0},x_{1},x_{2},x_{3}\}$ and $x_{0}\nsim x_{3}$.
If $E(T-T_{0})-E(T_{x_{0}x_{3}})\neq\varnothing$, then the conclusion follows from Lemmas \ref{lem-42}--\ref{lem-45}.
Let $E(T-T_{0})-E(T_{x_{0}x_{3}})=\varnothing$.
We distinguish the following four cases.

\noindent$\textbf{Case 1.}$ $\left|T_{3}\right|\geq 1$.

Without loss of generality, let $\left|T_{x_{0}x_{1}x_{3}}\right|\geq 1$ and $u\in T_{x_{0}x_{1}x_{3}}$.
By  Lemma \ref{lem-13}\ref{lem-13:item-7}, $\left|T_{x_{0}x_{1}x_{3}}\right|= 1$, $\left|T_{x_{0}x_{1}x_{2}}\cup T_{x_{1}x_{2}x_{3}}\right|= 0$ and $\left|T_{x_{0}x_{2}x_{3}}\right|\leq 1$.
If $\left|T_{x_{0}x_{2}x_{3}}\right|=1$, then  its conclusion is covered by Case $3$ of Lemma \ref{lem-45} by taking $\{x_{0},x_{1},x_{3},u\}$ ($u\in T_{x_{0}x_{1}x_{3}}$) as the vertices of a new $K_{4}-e$.
If $\left|T_{x_{0}x_{1}}\cup T_{x_{1}x_{2}}\right|\geq 1$, then  its conclusion is covered by Case $4$ of Lemma \ref{lem-43} by taking $\{x_{0},x_{1},x_{2},s\}$ ($s\in T_{x_{0}x_{1}}\cup T_{x_{1}x_{2}}$) as the vertices of a new $K_{4}-e$.
If $\left|T_{x_{0}x_{2}}\right|\geq 1$, then  its conclusion is covered by Case $1$ of Lemma \ref{lem-43} by taking $\{x_{0},x_{1},x_{2},q_{02,1}\}$ ($q_{02,1}\in T_{x_{0}x_{2}}$) as the vertices of a new $K_{4}-e$.
Let  $\left|(T_{2}- T_{x_{0}x_{3}})\cup T_{x_{0}x_{2}x_{3}}\right|=0$.
Then $\left|T_{x_{0}}\right|\leq 1$, $\left|T_{x_{3}}\right|\leq 1$ and $\left|T_{x_{2}}\right|\leq 1$; otherwise, $G[\mathcal{K}\cup\{u,\allowbreak p_{0,1}(p_{3,1},p_{2,1}),\allowbreak p_{0,2}(p_{3,2},p_{2,2})\}]\cong R_{27}$, a contradiction.
Then $\left|T_{x_{1}}\right|\leq 5$; otherwise, $G[\mathcal{K}\cup\{u,\allowbreak p_{1,1},\ldots,p_{1,6}\}]\cong R_{84}$, a contradiction.
Then $\left|E(T_{x_{0}x_{3}})\right|=0$; otherwise, $G[\{x_{0},\allowbreak x_{1},\allowbreak x_{2},\allowbreak u,\allowbreak q_{03,1},\allowbreak q_{03,2}\}]\cong R_{16}$, a contradiction.
Then $\left|T_{x_{0}x_{3}}\right|\leq 7$; otherwise, $G[\mathcal{K}\cup\{u,\allowbreak q_{03,1},\ldots,q_{03,8}\}]\cong R_{71}$, a contradiction.
Then $\left|E(T_{x_{0}x_{3}},T_{0})\right|=0$; otherwise, either $G[\{x_{0},\allowbreak x_{1},\allowbreak x_{2},\allowbreak u,\allowbreak q_{03,1},\allowbreak t_{1}\}]\cong R_{9}$ ($E(\{u,\allowbreak q_{03,1},\allowbreak t_{1}\})=\{q_{03,1}t_{1}\}$) or $G[\{x_{0},\allowbreak x_{1},\allowbreak x_{2},\allowbreak u,\allowbreak q_{03,1},\allowbreak t_{1}\}]\cong R_{12}$ ($E(\{u,\allowbreak q_{03,1},\allowbreak t_{1}\})=\{ut_{1},\allowbreak q_{03,1}t_{1}\}$), a contradiction.
If $\left|T_{x_{1}}\right|\geq 4$, then $\left|T_{x_{0}}\cup T_{x_{2}}\cup T_{x_{3}}\cup T_{0}\right|=0$ by Lemma \ref{lem-13}\ref{lem-13:item-8}; otherwise, $G[\{x_{0},\allowbreak x_{1},\allowbreak x_{2},\allowbreak u,\allowbreak p_{1,1},\ldots,p_{1,4},\allowbreak p_{0,1}\}]\cong R_{47}$ or $G[\mathcal{K}\cup\{p_{1,1},\ldots,p_{1,4},\allowbreak p_{2,1}\}]\cong R_{47}$ or $G[\{x_{1},\allowbreak x_{2},\allowbreak x_{3},\allowbreak u,\allowbreak p_{1,1},\ldots,p_{1,4},\allowbreak p_{3,1}\}]\cong R_{47}$ or $G[\{x_{0},\allowbreak x_{1},\allowbreak x_{3},\allowbreak u,\allowbreak p_{1,1},\ldots,p_{1,4},\allowbreak t_{1}\}]\cong R_{47}$, a contradiction.
If $\left|T_{x_{0}x_{3}}\right|\geq 4$, then $\left|T_{x_{2}}\cup T_{0}\cup T_{x_{1}}\right|=0$  by Lemma \ref{lem-13}\ref{lem-13:item-8}; otherwise, $G[\mathcal{K}\cup\{q_{03,1},\ldots,q_{03,4},\allowbreak p_{2,1}(p_{1,1})\}]\cong R_{55}$ or $G[\{x_{0},\allowbreak x_{1},\allowbreak x_{3},\allowbreak u,\allowbreak q_{03,1},\ldots,q_{03,4},\allowbreak t_{1}\}]\cong R_{55}$ ($u\sim t_{1}$), a contradiction.
If $\left|T_{x_{0}x_{3}}\right|\geq 2$, then $\left|T_{x_{1}}\right|\leq 1$; otherwise, $G[\{x_{0},\allowbreak x_{1},\allowbreak x_{3},\allowbreak q_{03,1},\allowbreak q_{03,2},\allowbreak p_{1,1},\allowbreak p_{1,2}\}]\cong R_{77}$, a contradiction.
If $\left|T_{x_{0}x_{3}}\right|=3$ and $\left|T_{x_{1}}\right|=1$, then $\left|T_{x_{2}}\cup T_{0}\cup T_{x_{3}}\cup T_{x_{0}}\right|=0$; otherwise, either $R_{47}\subset G$ (if $\left|T_{x_{0}}\cup T_{x_{3}}\right|\geq 1$) or  $R_{94}\subset G$ (if $\left|T_{x_{2}}\cup T_{0}\right|\geq 1$), a contradiction.
If $\left|T_{x_{0}x_{3}}\right|\geq 1$, then $\left|T_{x_{1}}\right|\leq 2$; otherwise, $G[\mathcal{K}\cup\{q_{03,1},\allowbreak p_{1,1},\allowbreak p_{1,2},\allowbreak p_{1,3}\}]\cong R_{38}$, a contradiction.
Thus, $G\subset H_{133}$ (if $\left|T_{x_{0}x_{3}}\right|= 0$ and $\left|T_{x_{1}}\right|\geq 4$) or $G\subset A_{35}$ (if $\left|T_{x_{0}x_{3}}\right|= 0$ and $\left|T_{x_{1}}\right|\leq 3$) or $G\subset A_{32}$ (if $\left|T_{x_{0}x_{3}}\right|= 1$) or $G\subset A_{29}$ (if $\left|T_{x_{0}x_{3}}\right|=2$) or $G\subset A_{26}$ (if $\left|T_{x_{0}x_{3}}\right|=3$ and $\left|T_{x_{1}}\right|=0$) or $G\subset H_{109}$ (if $\left|T_{x_{0}x_{3}}\right|= 3$ and $\left|T_{x_{1}}\right|\geq 1$) or $G\subset A_{20}$ (if $\left|T_{x_{0}x_{3}}\right|\geq 4$) by Lemma \ref{lem-13}\ref{lem-13:item-8}.

In the following cases, let $\left|T_{3}\right|=0$.

\noindent$\textbf{Case 2.}$ $\left|T_{x_{0}x_{1}}\right| \neq 0$.

If $\left|T_{x_{0}x_{3}}\right|\geq 1$, then  its conclusion is covered by Subcase $1.4$ of Lemma \ref{lem-45} by taking $\{x_{0},x_{1},x_{2},q_{01,1}\}$ ($q_{01,1}\in T_{x_{0}x_{1}}$) as the vertices of a new $K_{4}-e$.
If $\left|T_{x_{1}x_{2}}\right|\geq 1$, then  its conclusion is covered by Subcase $1.6$ of Lemma \ref{lem-45} by taking $\{x_{0},x_{1},x_{2},q_{12,1}\}$ ($q_{12,1}\in T_{x_{1}x_{2}}$) as the vertices of a new $K_{4}-e$.
Let  $\left|T_{x_{0}x_{3}}\cup T_{x_{1}x_{2}}\right|=0$.
Then $\left|T_{x_{2}x_{3}}\right|=0$; otherwise, $G[\mathcal{K}\cup\{q_{01,1},\allowbreak q_{23,1}\}]\cong R_{18}$, a contradiction.
Then $\left|T_{x_{3}}\right|=0$; otherwise, $G[\mathcal{K}\cup\{q_{01,1},\allowbreak p_{3,1}\}]\cong R_{14}$, a contradiction.
Then $\left|T_{x_{0}}\right| \leq 1$ and $\left|T_{x_{2}}\right| \leq 1$; otherwise, $G[\mathcal{K}\cup\{q_{01,1},\allowbreak p_{0,1}(p_{2,1}),\allowbreak p_{0,2}(p_{2,2})\}]\cong R_{25}$, a contradiction.
Then $\left|T_{x_{1}}\cup T_{x_{0}x_{1}}\right|\leq 8$ by Lemma \ref{lem-13}\ref{lem-13:item-1} and \ref{lem-13:item-4}; otherwise, $R_{95}\subset G$ (if $\left|T_{x_{0}x_{1}}\right|=1$) or $R_{86}\subset G$ (if $\left|T_{x_{0}x_{1}}\right|=2$ or $3$) or  $R_{90}\subset G$ (if $\left|T_{x_{0}x_{1}}\right|=4$)  or $R_{96}\subset G$ (if $\left|T_{x_{0}x_{1}}\right|=5$)  or $R_{89}\subset G$ (if $\left|T_{x_{0}x_{1}}\right|=6$)  or $R_{88}\subset G$ (if $\left|T_{x_{0}x_{1}}\right|=7$) or $R_{92}\subset G$ (if $\left|T_{x_{0}x_{1}}\right|=8$), a contradiction.
If $\left|T_{x_{0}x_{1}}\right| \geq5$, then $\left|T_{x_{0}}\right| =0$; otherwise, $G[\mathcal{K}\cup\{q_{01,1},\ldots,q_{01,5},\allowbreak p_{0,1}\}]\cong R_{91}$, a contradiction.
If $\left|T_{x_{0}x_{1}}\right| \geq2$, then $\left|T_{x_{2}}\right| =0$; otherwise, $R_{28}\subset G$, a contradiction.
If $\left|T_{x_{0}x_{2}}\right| \geq1$, then $\left|T_{x_{0}x_{1}}\cup T_{x_{0}x_{2}}\right| \leq 4$; otherwise, $G[\{x_{1},\allowbreak x_{2},\allowbreak q_{01,1},\allowbreak q_{01,2},\allowbreak q_{02,1},\allowbreak q_{02,2},\allowbreak q_{02,3}\}]\cong R_{21}$ or $G[\{x_{1},\allowbreak x_{2},\allowbreak q_{01,1},\allowbreak q_{01,2},\allowbreak q_{01,3},\allowbreak q_{02,1},\allowbreak q_{02,2}\}]\cong R_{21}$ or $G[\mathcal{K}\cup\{q_{01,1},\allowbreak q_{02,1},\ldots,q_{02,4}\}]\cong R_{61}$ or $G[\mathcal{K}\cup\{q_{01,1},\ldots,q_{01,4},\allowbreak q_{02,1}\}]\cong R_{61}$, a contradiction.
If $\left|T_{x_{1}x_{3}}\right| \geq1$, then $\left|T_{x_{0}x_{1}}\cup T_{x_{1}x_{3}}\right| =2$; otherwise, $G[\{x_{0},\allowbreak x_{1},\allowbreak x_{3},\allowbreak q_{01,1},\allowbreak q_{01,2},\allowbreak q_{13,1}\}]\cong R_{16}$ or $G[\{x_{0},\allowbreak x_{1},\allowbreak x_{3},\allowbreak q_{01,1},\allowbreak q_{13,1},\allowbreak q_{13,2}\}]\cong R_{16}$, a contradiction.
Thus, if $\left|T_{x_{0}x_{1}}\right| \geq 5$, then $G\subset H_{151}$.
Let $\left|T_{x_{0}x_{1}}\right| \leq 4$.
If $\left|T_{x_{0}x_{1}}\right| =4$ and $\left|T_{x_{0}}\right| =1$, then $\left|T_{x_{1}}\right| =0$; otherwise, $G[\mathcal{K}\cup\{q_{01,1},\ldots,q_{01,4},\allowbreak p_{0,1},\allowbreak p_{1,1}\}]\cong R_{87}$, a contradiction.
Thus, if $\left|T_{x_{0}x_{1}}\right| =4$, then either $G\subset H_{151}$ (if $\left|T_{x_{0}}\right| =0$) or $G\subset H_{114}$ (if $\left|T_{x_{0}}\right| =1$).
Let $\left|T_{x_{0}x_{1}}\right| \leq 3$.
If $\left|T_{x_{0}x_{1}}\right| =3$ and $\left|T_{x_{0}}\right| =1$, then $\left|T_{x_{1}}\right| \leq 1$; otherwise, $G[\mathcal{K}\cup\{q_{01,1},\allowbreak q_{01,2},\allowbreak q_{01,3},\allowbreak p_{0,1},\allowbreak p_{1,1},\allowbreak p_{1,2}\}]\cong R_{93}$, a contradiction.
If $\left|T_{x_{0}x_{2}}\right| \geq 1$, then $\left|T_{x_{1}}\right| \leq 1$; otherwise, $G[\{x_{0},\allowbreak x_{1},\allowbreak x_{2},\allowbreak q_{01,1},\allowbreak q_{02,1},\allowbreak p_{1,1},\allowbreak p_{1,2}\}]\cong R_{25}$, a contradiction.
Thus, if $\left|T_{x_{0}x_{1}}\right| =3$, then either $G\subset H_{133}$ (if $\left|T_{x_{0}}\cup T_{x_{0}x_{2}}\right| =0$) or $G\subset H_{114}$ (if $\left|T_{x_{0}}\cup T_{x_{0}x_{2}}\right| \geq 1$).
Let $\left|T_{x_{0}x_{1}}\right| \leq 2$.
If $\left|T_{x_{0}x_{2}}\right| \geq 2$, then $\left|T_{x_{1}}\right| =0$; otherwise, $G[\{x_{0},\allowbreak x_{1},\allowbreak x_{2},\allowbreak q_{01,1},\allowbreak q_{02,1},\allowbreak q_{02,2},\allowbreak p_{1,1}\}]\cong R_{28}$, a contradiction.
If $\left|T_{x_{0}x_{1}}\right| =2$ and $\left|T_{x_{0}}\right| =1$, then $\left|T_{x_{1}}\right| \leq 2$; otherwise, $G[\{x_{0},\allowbreak x_{1},\allowbreak x_{3},\allowbreak q_{01,1},\allowbreak q_{01,2},\allowbreak p_{0,1},\allowbreak p_{1,1},\allowbreak p_{1,2},\allowbreak p_{1,3}\}]\cong R_{47}$, a contradiction.
Thus, if $\left|T_{x_{0}x_{1}}\right| =2$, then  $G\subset H_{79}$ (if $\left|T_{x_{0}x_{2}}\right| =2$) or $G\subset H_{80}$ (if $\left|T_{x_{0}x_{2}}\right| =1$) or $G\subset H_{120}$ (if $\left|T_{x_{0}x_{2}}\right| =0$ and $\left|T_{x_{0}}\right| =0$) or $G\subset H_{78}$ (if $\left|T_{x_{0}x_{2}}\right| =0$ and $\left|T_{x_{0}}\right| =1$).
Let $\left|T_{x_{0}x_{1}}\right| =1$.
By symmetry, if $\left|T_{x_{0}x_{2}}\right| \geq 2$, we can obtain a similar conclusion.
Let $\left|T_{x_{0}x_{2}}\right| \leq 1$.
Then $\left|T_{x_{0}x_{2}}\cup T_{x_{1}x_{3}}\right|\leq 1$; otherwise, $G[\mathcal{K}\cup\{q_{02,1},\allowbreak q_{13,1}\}]\cong R_{18}$, a contradiction.
Then $\left|T_{x_{0}}\cup T_{x_{1}x_{3}}\right|\leq 1$; otherwise, $G[\mathcal{K}\cup\{p_{0,1},\allowbreak q_{13,1}\}]\cong R_{14}$, a contradiction.
Thus, if $\left|T_{x_{0}x_{2}}\right| =1$, then  $G\subset H_{77}$.
Let $\left|T_{x_{0}x_{2}}\right| =0$.
If $\left|T_{x_{0}}\cup T_{x_{2}}\right| \geq 1$, then $\left|T_{x_{1}}\right| \leq 3$; otherwise, $G[\{x_{0},\allowbreak x_{1},\allowbreak x_{2},\allowbreak q_{01,1},\allowbreak p_{0,1},\allowbreak p_{1,1},\ldots,p_{1,4}\}]\cong R_{47}$ or $G[\mathcal{K}\cup\{p_{2,1},\allowbreak p_{1,1},\ldots,p_{1,4}\}]\cong R_{47}$, a contradiction.
Thus, $G\subset H_{85}$ (if $\left|T_{x_{0}}\right| =1$) or $G\subset H_{87}$ (if $\left|T_{x_{0}}\right| =0$  and $\left|T_{x_{2}}\right| =1$) or $G\subset H_{120}$ (if $\left|T_{x_{0}}\right| =0$  and $\left|T_{x_{2}}\right| =0$).

In the following cases, by symmetry,  let  $\left|T_{2}-(T_{x_{1}x_{2}}\cup T_{x_{0}x_{3}})\right| = 0$.

\noindent$\textbf{Case 3.}$ $\left|T_{x_{1}x_{2}}\right| \neq 0$.

If $\left|T_{x_{0}x_{3}}\right|\geq 1$, then  its conclusion is covered by Subcase $1.1$ of Lemma \ref{lem-45} by taking $\{x_{0},x_{1},x_{2},q_{12,1}\}$ ($q_{12,1}\in T_{x_{1}x_{2}}$) as the vertices of a new $K_{4}-e$.
Let  $\left|T_{x_{0}x_{3}}\right|=0$.
If $t\in T_{0}$, $q_{12,1},q_{12,2}\in T_{x_{1}x_{2}}$ and $t\sim q_{12,1},q_{12,2}$, then $\left|T_{1}\right| =0$ and $\left|T_{x_{1}x_{2}}\right|=2$; otherwise, $R_{2}\subset G$ (if $\left|T_{x_{0}}\cup T_{x_{3}}\right|\geq 1$) or $R_{40}\subset G$ (if $\left|T_{x_{1}}\cup T_{x_{2}}\right|\geq 1$) or $R_{41}\subset G$ (if $\left|T_{x_{1}x_{2}}\right|\geq 3$), a contradiction.
Thus, if $t\sim q_{12,1},q_{12,2}$, then $G\subset H_{138}$ by Lemma \ref{lem-13}\ref{lem-13:item-8}.
Let $\left|N(t)\cap T_{x_{1}x_{2}}\right|\leq1$ for any vertex $t\in T_{0}$.
If $\left|T_{x_{1}}\right| \geq 1$ and $\left|T_{x_{2}}\right| \geq 1$, then $\left|T_{x_{1}}\cup T_{x_{2}}\right| \leq 4$; otherwise, either $G[\{x_{1},\allowbreak x_{2},\allowbreak p_{1,1}(p_{2,1}),\allowbreak p_{1,2}(p_{2,2}),\allowbreak p_{2,1}(p_{1,1}),\allowbreak p_{2,2}(p_{1,2}),\allowbreak p_{2,3}(p_{1,3})\}]\cong R_{21}$ or $G[\mathcal{K}\cup\{p_{1,1}(p_{2,1}),\allowbreak p_{2,1},\ldots,p_{2,4}(p_{1,1},\ldots,p_{1,4})\}]\cong R_{47}$, a contradiction.
If $\left|T_{x_{1}x_{2}}\right| \geq 4$, then $\left|T_{x_{1}}\cup T_{x_{2}}\right| \leq 4$; otherwise, $G[\mathcal{K}\cup\{q_{12,1},\ldots,q_{12,4},\allowbreak p_{2,1},\ldots,p_{2,5}(p_{1,1},\ldots,p_{1,5})\}]\cong R_{67}$, a contradiction.
If $\left|T_{x_{1}x_{2}}\right| \geq 2$, then $\left|T_{x_{1}}\right| \leq 5$ and $\left|T_{x_{2}}\right| \leq 5$; otherwise, $G[\mathcal{K}\cup\{q_{12,1},\allowbreak q_{12,2},\allowbreak p_{1,1},\ldots,p_{1,6}\}]\cong R_{90}$, a contradiction.
If $\left|T_{x_{1}x_{2}}\right| =1$, then $\left|T_{x_{1}}\right| \leq 6$ and $\left|T_{x_{2}}\right| \leq 6$; otherwise, $G[\mathcal{K}\cup\{q_{12,1},\allowbreak p_{1,1},\ldots,p_{1,7}\}]\cong R_{86}$, a contradiction.
Then $\left|T_{x_{0}}\cup T_{x_{3}}\cup T_{0}\right| \leq 1$ by $d(G) <4$; otherwise, $G[\mathcal{K}\cup\{q_{12,1},\allowbreak t_{1},\allowbreak t_{2}\}]\cong R_{24}$, a contradiction.
If $\left|T_{x_{0}}\cup T_{x_{3}}\cup T_{0}\right| =1$, then $\left|T_{x_{1}}\right| \leq 4$ and $\left|T_{x_{2}}\right| \leq 4$; otherwise, $G[\{x_{0},\allowbreak x_{1},\allowbreak x_{2},\allowbreak p_{0,1}(p_{3,1}),\allowbreak p_{1,1},\ldots,p_{1,5}(p_{2,1},\ldots,p_{2,5})\}]\cong R_{46}$  or $G[\{x_{1},\allowbreak x_{2},\allowbreak q_{12,1},\allowbreak t_{1},\allowbreak p_{1,1},\ldots,p_{1,5}(p_{2,1},\ldots,p_{2,5})\}]\cong R_{46}$ ($E(\{q_{12,1},\allowbreak t_{1},\allowbreak p_{1,1},\ldots,p_{1,5}\})=\{q_{12,1}t_{1}\}$), a contradiction.
Thus, if $\left|T_{x_{1}x_{2}}\right| \geq 4$ or $\left|T_{x_{0}}\cup T_{x_{3}}\cup T_{0}\right| =1$, then $G\subset D_{2}$.
Let $\left|T_{x_{1}x_{2}}\right| \leq 3$ and $\left|T_{x_{0}}\cup T_{x_{3}}\cup T_{0}\right| =0$.
If $\left|T_{x_{1}}\right| =6$, then $G\subset H_{120}$.
If $\left|T_{x_{1}}\right| =5$, then $G\subset H_{134}$.
If $\left|T_{x_{1}}\right| \leq 4$ and $\left|T_{x_{2}}\right| \leq 4$, then $G\subset D_{2}$.

In the following cases, let  $\left|T_{x_{1}x_{2}}\right| = 0$.

\noindent$\textbf{Case 4.}$ The remaining case.

If $\left|T_{x_{1}}\right| \geq 1$ and $\left|T_{x_{2}}\right| \geq 1$, then $\left|T_{x_{1}}\cup T_{x_{2}}\right| \leq 4$; otherwise, either $G[\{x_{1},\allowbreak x_{2},\allowbreak p_{1,1}(p_{2,1}),\allowbreak p_{1,2}(p_{2,2}),\allowbreak p_{2,1}(p_{1,1}),\allowbreak p_{2,2}(p_{1,2}),\allowbreak p_{2,3}(p_{1,3})\}]\cong R_{21}$ or  $G[\mathcal{K}\cup\{p_{1,1},\allowbreak p_{2,1},\ldots,p_{2,4}\}(\{p_{1,1},\ldots,p_{1,4},\allowbreak p_{2,1}\})]\cong R_{47}$, a contradiction.
Then $\left|T_{x_{0}}\cup T_{x_{3}}\right| \leq 2$; otherwise, either $G[\{x_{0},\allowbreak x_{1},\allowbreak x_{3},\allowbreak p_{0,1}(p_{3,1}),\allowbreak p_{0,2}(p_{3,2}),\allowbreak p_{3,1}(p_{0,1})\}]\cong R_{1}$  or $G[\mathcal{K}\cup\{p_{0,1},\allowbreak p_{0,2},\allowbreak p_{0,3}(\{p_{3,1},\allowbreak p_{3,2},\allowbreak p_{3,3}\})\}]\cong R_{23}$, a contradiction.
If $\left|T_{x_{0}}\cup T_{x_{3}}\right| \geq 1$, then $\left|T_{x_{1}}\right| \leq 4$ and $\left|T_{x_{2}}\right| \leq 4$; otherwise, $G[\{x_{0},\allowbreak x_{1},\allowbreak x_{2},\allowbreak p_{0,1},\allowbreak p_{1,1},\ldots,p_{1,5}(p_{2,1},\ldots,p_{2,5})\}]\cong R_{46}$, a contradiction.
If $\left|T_{x_{1}}\cup T_{x_{2}}\right| \geq 1$, then $\left|E(T_{x_{0}x_{3}})\right| =0$ and $\left|T_{x_{0}x_{3}}\right| \leq 3$; otherwise, either $G[\{x_{0},\allowbreak x_{1},\allowbreak x_{2},\allowbreak p_{1,1}(p_{2,1}),\allowbreak q_{03,1},\allowbreak q_{03,2}\}]\cong R_{10}$ ($q_{03,1}\sim q_{03,2}$) or $G[\mathcal{K}\cup\{p_{1,1}(p_{2,1}),\allowbreak q_{03,1},\ldots,q_{03,4}\}]\cong R_{55}$, a contradiction.
If $\left|T_{x_{1}}\cup T_{x_{2}}\right| \geq 1$, then $\left|T_{x_{0}}\right| \leq 1$ and $\left|T_{x_{3}}\right| \leq 1$; otherwise, $G[\mathcal{K}\cup\{p_{1,1}(p_{2,1}),\allowbreak p_{0,1}(p_{3,1}),\allowbreak p_{0,2}(p_{3,2})\}]\cong R_{22}$, a contradiction.
If $\left|T_{x_{1}}\right| \geq 2$ or $\left|T_{x_{2}}\right| \geq 2$, then $\left|T_{x_{0}x_{3}}\right| \leq 1$; otherwise, $G[\{x_{0},\allowbreak x_{1},\allowbreak x_{3},\allowbreak p_{1,1}(p_{2,1}),\allowbreak p_{1,2}(p_{2,2}),\allowbreak q_{03,1},\allowbreak q_{03,2}\}]\cong R_{77}$, a contradiction.
If $\left|T_{x_{1}}\right| \geq 3$ or $\left|T_{x_{2}}\right| \geq 3$, then $\left|T_{x_{0}x_{3}}\right| =0$; otherwise, $G[\mathcal{K}\cup\{p_{1,1}(p_{2,1}),\allowbreak p_{1,2}(p_{2,2}),\allowbreak p_{1,3}(p_{2,3}),\allowbreak q_{03,1}\}]\cong R_{38}$, a contradiction.
If $\left|T_{x_{1}}\cup T_{x_{2}}\right|=0$, then $G\subset D_{1}$ by Lemma \ref{lem-13}\ref{lem-13:item-8}.
Let $\left|T_{x_{1}}\cup T_{x_{2}}\right|\geq 1$.
If $\left|T_{x_{1}}\right| \geq 5$ or $\left|T_{x_{2}}\right| \geq 5$, then $\left|T_{x_{0}}\cup T_{x_{3}}\right| =0$; otherwise, $G[\{x_{0},\allowbreak x_{1},\allowbreak x_{2},\allowbreak p_{0,1}(p_{3,1}),\allowbreak p_{1,1},\ldots,p_{1,5}(p_{2,1},\ldots,p_{2,5})\}]\cong R_{46}$, a contradiction.
Thus, if  $\left|T_{x_{1}}\right| \geq 5$ or $\left|T_{x_{2}}\right| \geq 5$, then $G\subset H_{120}$ by Lemma \ref{lem-13}\ref{lem-13:item-1}.
Let $\left|T_{x_{1}}\right|  \leq4$ and $\left|T_{x_{2}}\right| \leq4$.
If $\left|T_{x_{1}}\right| \geq 3$ or $\left|T_{x_{2}}\right| \geq 3$, then either $G\subset H_{82}$ (if $\left|T_{x_{1}}\right| = 3$) or $G\subset H_{84}$ (if $\left|T_{x_{1}}\right| = 4$).
Let $\left|T_{x_{1}}\right|  \leq2$ and $\left|T_{x_{2}}\right| \leq2$.
Thus, if  $\left|T_{x_{1}}\right| =2$ or $\left|T_{x_{2}}\right| =2$, then $G\subset A_{47}$.
Let $\left|T_{x_{1}}\right|  \leq1$ and $\left|T_{x_{2}}\right| \leq1$.
If $\left|T_{x_{1}}\right|=1$ and $\left|T_{x_{2}}\right|=1$, then $\left|T_{x_{0}x_{3}}\right| \leq 2$; otherwise, $G[\mathcal{K}\cup\{p_{1,1},\allowbreak p_{2,1},\allowbreak q_{03,1},\allowbreak q_{03,2},\allowbreak q_{03,3}\}]\cong R_{94}$, a contradiction.
Thus, either $G\subset A_{23}$ (if $\left|T_{x_{1}}\cup T_{x_{2}}\right|=2$) or $G\subset A_{19}$ (if $\left|T_{x_{1}}\cup T_{x_{2}}\right|\leq 1$).

This completes the proof.
\end{proof}

\section{Planar graphs with diameter at least 4}\label{s-diam4}

We prove the following theorem in this section.

\begin{thm}\label{thm-51}
If $G \in \mathcal{G}$ and $d(G) \geq 4$, then $G$ is a connected induced subgraph of a graph in $\left\{A_{1},\ldots,A_{51},D_1,D_2\right\}$.
\end{thm}

Let $G\in\mathcal{G}$ and suppose that $d(G)\geq4$.
Since a diametral path is induced and $\lambda_2(P_6)\approx1.24698>1$, Corollary \ref{cor-subgraph} implies that $d(G)\leq4$.
Hence it remains to consider $d(G)=4$. 
Since $P_5$ is an induced subgraph of $G$ and $\lambda_2(P_5)=1$, eigenvalue interlacing and the definition of $\mathcal{G}$ give $\lambda_2(G)=1$.

\subsection{Proof strategy}

Let $P=v_1v_2v_3v_4v_5$ be a diametral path. We characterize $G$ by adding vertices to $P$ and analyzing their adjacencies both to $P$ and to vertices outside $P$.
We first show that no additional vertex can be adjacent to $v_1$ or $v_5$.
Let $X=N(v_2)\cap N(v_4)$ and $Y=V(G)\setminus(X\cup V(P))$.
Then $|X|\geq1$ and $N_{V(G)\setminus Y}(Y)\subseteq X$.
We determine the possible structures of $G[X]$ and $G[Y]$ in Lemmas \ref{diam4:Xstructure} and \ref{diam4:Ystructure}, and show that $|X|\leq2$ whenever $G[Y]$ contains an edge.
The case $|X|\leq2$ is treated in Lemma \ref{diam4:Xsmall}; the case $|X|\geq3$ is treated in Lemmas \ref{diam4:Xcycle} and \ref{diam4:Xforest}.

\subsection{Structural characterization}

For vertex sets $S,T\subseteq V(G)$, write $N_S(T):=N(T)\cap S$.
Throughout this section, let $G\in\mathcal{G}$ with $d(G)=4$, and let $P=v_1v_2v_3v_4v_5$ be a diametral path of $G$.

We claim that $N(v_1)=\{v_2\}$, $N(v_5)=\{v_4\}$, and $N(v_2)\setminus \{v_1\}=N(v_4)\setminus \{v_5\}$.
First, suppose that there exists a vertex $v\neq v_2$ in $N(v_1)$.
Since $P_6\not\subset G$, we have $N(v)\cap V(P)\setminus\{v_1\}\neq \varnothing$ and $N(v)\cap V(P)\setminus\{v_1\}\subseteq\{v_2,v_3\}$.
If $N(v)\cap V(P)=\{v_1,v_2\}$, then $R_{98}\subset G$, a contradiction.
If $N(v)\cap V(P)=\{v_1,v_3\}$, then $R_{2}\subset G$, a contradiction.
If $N(v)\cap V(P)=\{v_1,v_2,v_3\}$, then $R_{6}\subset G$, a contradiction.
Hence $N(v_1)=\{v_2\}$.
The proof for $N(v_5)=\{v_4\}$ is similar and thus omitted.
Next, suppose that there exists a vertex $v\in N(v_2)\setminus\{v_1\}$ such that $v\notin N(v_4)$.
Since $d(G)=4$, we have $N(v)\cap V(P)\setminus \{v_2\}\subseteq\{v_1,v_3\}$.
By the previous proof, we only need to consider the following two cases.
If $N(v)\cap V(P)=\{v_2,v_3\}$, then $R_{3}\subset G$, a contradiction.
If $N(v)\cap V(P)=\{v_2\}$, then $R_{1}\subset G$, a contradiction.
Therefore, if $v\in N(v_2)\backslash \{v_1\}$, then $v\in N(v_4)$.
Similarly, we obtain that if $v\in N(v_4)\backslash \{v_5\}$, then $v\in N(v_2)$.
Thus, $N(v_2)\setminus \{v_1\}=N(v_4)\setminus \{v_5\}$.

Denote by $X$ the vertex set $N(v_2)\cap N(v_4)$ and by $Y$ the set of vertices in $V(G)\setminus (X\cup V(P))$.
Then $|X|\geq 1$ and $N_{V(G)\setminus Y}(Y)\subseteq X$.

To characterize $G$, we first determine the possible structures of $G[X]$ and $G[Y]$ (see Lemmas \ref{diam4:Xstructure} and \ref{diam4:Ystructure}) and show that $|X|\leq2$ if $G[Y]$ contains an edge.
We then consider the cases $|X|\leq2$ and $|X|\geq3$. The former is treated in Lemma \ref{diam4:Xsmall}, and the latter in Lemmas \ref{diam4:Xcycle} and \ref{diam4:Xforest}.

We interpret a union indexed from $1$ to $k$ as empty when $k=0$.

\begin{lem}\label{diam4:Xstructure}
  The induced subgraph $G[X]$ is isomorphic to $\bigcup_{i=1}^s C_{a_i}\cup \bigcup_{j=1}^t P_{b_j}$, where $0\leq s\leq 1$, $t\geq 0$, $3\leq a_i\leq 6$ and $1\leq b_j\leq 5$. Moreover, if $s=1$, then $G[X]$ is isomorphic to $C_{a_1}$.
\end{lem}

\begin{proof}
Suppose that there exists a vertex $u\in X$ such that $d_{G[X]}(u)\geq3$.
  Since every vertex in $X$ is adjacent to both $v_2$ and $v_4$, the graph $G$ then contains $K_{3,3}$ as a subgraph, a contradiction.
  Therefore, $d_{G[X]}(u)\leq2$ for every $u\in X$. Hence $G[X]=\bigcup_{i=1}^s C_{a_i}\cup \bigcup_{j=1}^t P_{b_j}$, where $a_i\geq3$, $b_j\geq1$, and $s,t\geq0$.
  If $G[X]$ contained two cycles, then $G$ would contain a subdivision of $K_5$, a contradiction.
  Therefore, $0\leq s\leq 1$.
  Similarly, if $s=1$ and $G[X]$ contained a path component $P_b$ with $b\geq1$, then $G$ would contain a subdivision of $K_5$, a contradiction.
  Therefore, $s=1$ implies $t=0$ and $G[X]=C_{a_1}$.
  Finally, $P_6$ is forbidden for $\lambda_2(P_6)>1$. Thus $3\leq a_i\leq6$ and $1\leq b_j\leq5$.
\end{proof}

\begin{lem}\label{diam4:Ystructure}
The induced subgraph $G[Y]$ is isomorphic to $sK_1\cup tK_2$ for some integers $s,t\geq 0$.
\end{lem}

\begin{proof}
Suppose that there exist three vertices in $Y$, say $u_1,u_2,u_3$, such that $G[\{u_1,u_2,u_3\}]\cong P_3$ or $G[\{u_1,u_2,u_3\}]\cong K_3$.
Let $V_1=\{u_1,u_2,u_3\}$ and $V_2=N_X(V_1)$.
If $V_2=\varnothing$, then $G[V(P)\cup V_1]\cong P_5\cup P_3$ or $P_5\cup K_3$.
Since $\lambda_2(P_5\cup P_3)=\sqrt{2}>1$ and $\lambda_2(P_5\cup K_3)=\sqrt{3}>1$, we have $\lambda_2(G)>1$ by Lemma \ref{lem-2}, a contradiction.
Thus, $V_2\neq\varnothing$.
Let $w$ be a vertex in $V_2$. Then $1\leq |N_{V_1}(w)|\leq 3$.
If $|N_{V_1}(w)|=1$ and $G[V_1]=P_3$, then $R_{1}\subset G$ or $P_6\subset G$, a contradiction.
If $G[V_1]=K_3$, then $R_{98}\subset G$, a contradiction.
If $|N_{V_1}(w)|=2$ and $G[V_1]=P_3$, then $R_{3}\subset G$ or $R_{2}\subset G$, a contradiction. If $G[V_1]=K_3$, then $R_{6}\subset G$, a contradiction.
If $|N_{V_1}(w)|=3$, then $R_{9}\subset G$ or $R_{15}\subset G$, a contradiction.
Therefore, for any three vertices in $Y$, say $u_1,u_2,u_3$, we have $G[\{u_1,u_2,u_3\}]=3K_1$ or $G[\{u_1,u_2,u_3\}]=K_1\cup K_2$.
This completes the proof.
\end{proof}

\begin{lem}\label{diam4:NXK2}
Let $H=K_2\subset G[Y]$. Then $1\leq |N_X(V(H))|\leq 2$.
Moreover, if $|N_X(V(H))|= 2$, then $G[N_X(V(H))\cup V(H)]$ is isomorphic to $Q_i$, $i=1,2,5,6,7,8$ (see Fig.~\ref{diam4:Qfig}) and $|X|=2$.
If $|N_X(V(H))|= 1$, then $G[N_X(V(H))\cup V(H)]$ is isomorphic to $Q_i$, $i=9,10$ (see Fig.~\ref{diam4:Qfig}) and $|X|=1$.
\end{lem}

\par\addvspace{2mm}
\noindent\begin{minipage}{\linewidth}
\captionsetup{skip=3pt}
\centering
\begin{tikzpicture}[x=0.5mm, y=0.5mm, inner xsep=0pt, inner ysep=0pt, outer xsep=0pt, outer ysep=0pt]
\path[line width=0mm] (-3.31,147.28) rectangle +(237.32,97.28);
\definecolor{L}{rgb}{0,0,0}
\definecolor{F}{rgb}{0,0,0}
\path[line width=0.15mm, draw=L, fill=F] (0.00,240.00) circle (0.4mm);
\path[line width=0.15mm, draw=L, fill=F] (30.13,239.95) circle (0.4mm);
\path[line width=0.15mm, draw=L, fill=F] (-0.31,209.98) circle (0.4mm);
\path[line width=0.15mm, draw=L, fill=F] (30.13,209.98) circle (0.4mm);
\path[line width=0.15mm, draw=L] (-0.08,239.95) -- (29.90,210.21);
\path[line width=0.15mm, draw=L] (-0.31,240.18) -- (-0.08,209.98);
\path[line width=0.45mm, draw=L] (30.14,239.77) -- (29.93,209.96);
\path[line width=0.15mm, draw=L] (-0.10,209.75) -- (30.14,209.96);
\path[line width=0.15mm, draw=L] (-0.31,209.96) -- (30.14,239.77);
\path[line width=0.15mm, draw=L, fill=F] (0.00,240.00) circle (0.4mm);
\path[line width=0.15mm, draw=L, fill=F] (30.13,239.95) circle (0.4mm);
\path[line width=0.15mm, draw=L, fill=F] (-0.31,209.98) circle (0.4mm);
\path[line width=0.15mm, draw=L, fill=F] (30.13,209.98) circle (0.4mm);
\path[line width=0.15mm, draw=L, fill=F] (0.00,240.00) circle (0.4mm);
\path[line width=0.15mm, draw=L, fill=F] (30.13,239.95) circle (0.4mm);
\path[line width=0.15mm, draw=L, fill=F] (-0.31,209.98) circle (0.4mm);
\path[line width=0.15mm, draw=L, fill=F] (30.13,209.98) circle (0.4mm);
\path[line width=0.15mm, draw=L, fill=F] (0.00,240.00) circle (0.4mm);
\path[line width=0.15mm, draw=L, fill=F] (30.13,239.95) circle (0.4mm);
\path[line width=0.15mm, draw=L, fill=F] (-0.31,209.98) circle (0.4mm);
\path[line width=0.15mm, draw=L, fill=F] (30.13,209.98) circle (0.4mm);
\path[line width=0.15mm, draw=L, fill=F] (0.00,240.00) circle (0.4mm);
\path[line width=0.15mm, draw=L, fill=F] (30.13,239.95) circle (0.4mm);
\path[line width=0.15mm, draw=L, fill=F] (-0.31,209.98) circle (0.4mm);
\path[line width=0.15mm, draw=L, fill=F] (30.13,209.98) circle (0.4mm);
\path[line width=0.15mm, draw=L, fill=F] (0.00,240.00) circle (0.4mm);
\path[line width=0.45mm, draw=L, fill=F] (30.13,239.95) circle (0.8mm);
\path[line width=0.15mm, draw=L, fill=F] (-0.31,209.98) circle (0.4mm);
\path[line width=0.45mm, draw=L, fill=F] (30.13,209.98) circle (0.8mm);
\path[line width=0.15mm, draw=L] (0.16,240.05) -- (30.05,240.05);
\path[line width=0.15mm, draw=L, fill=F] (50.14,240.00) circle (0.4mm);
\path[line width=0.15mm, draw=L, fill=F] (80.27,239.95) circle (0.4mm);
\path[line width=0.15mm, draw=L, fill=F] (49.83,209.98) circle (0.4mm);
\path[line width=0.15mm, draw=L, fill=F] (80.27,209.98) circle (0.4mm);
\path[line width=0.15mm, draw=L] (50.06,239.95) -- (80.04,210.21);
\path[line width=0.45mm, draw=L] (80.28,239.77) -- (80.07,209.96);
\path[line width=0.15mm, draw=L] (50.04,209.75) -- (80.28,209.96);
\path[line width=0.15mm, draw=L] (49.83,209.96) -- (80.28,239.77);
\path[line width=0.15mm, draw=L, fill=F] (50.14,240.00) circle (0.4mm);
\path[line width=0.15mm, draw=L, fill=F] (80.27,239.95) circle (0.4mm);
\path[line width=0.15mm, draw=L, fill=F] (49.83,209.98) circle (0.4mm);
\path[line width=0.15mm, draw=L, fill=F] (80.27,209.98) circle (0.4mm);
\path[line width=0.15mm, draw=L, fill=F] (50.14,240.00) circle (0.4mm);
\path[line width=0.15mm, draw=L, fill=F] (80.27,239.95) circle (0.4mm);
\path[line width=0.15mm, draw=L, fill=F] (49.83,209.98) circle (0.4mm);
\path[line width=0.15mm, draw=L, fill=F] (80.27,209.98) circle (0.4mm);
\path[line width=0.15mm, draw=L, fill=F] (50.14,240.00) circle (0.4mm);
\path[line width=0.15mm, draw=L, fill=F] (80.27,239.95) circle (0.4mm);
\path[line width=0.15mm, draw=L, fill=F] (49.83,209.98) circle (0.4mm);
\path[line width=0.15mm, draw=L, fill=F] (80.27,209.98) circle (0.4mm);
\path[line width=0.15mm, draw=L, fill=F] (50.14,240.00) circle (0.4mm);
\path[line width=0.15mm, draw=L, fill=F] (80.27,239.95) circle (0.4mm);
\path[line width=0.15mm, draw=L, fill=F] (49.83,209.98) circle (0.4mm);
\path[line width=0.15mm, draw=L, fill=F] (80.27,209.98) circle (0.4mm);
\path[line width=0.15mm, draw=L, fill=F] (50.14,240.00) circle (0.4mm);
\path[line width=0.45mm, draw=L, fill=F] (80.27,239.95) circle (0.8mm);
\path[line width=0.15mm, draw=L, fill=F] (49.83,209.98) circle (0.4mm);
\path[line width=0.45mm, draw=L, fill=F] (80.27,209.98) circle (0.8mm);
\path[line width=0.15mm, draw=L] (50.30,240.05) -- (80.19,240.05);
\path[line width=0.15mm, draw=L, fill=F] (100.29,240.19) circle (0.4mm);
\path[line width=0.15mm, draw=L, fill=F] (130.41,240.15) circle (0.4mm);
\path[line width=0.15mm, draw=L, fill=F] (99.98,210.17) circle (0.4mm);
\path[line width=0.15mm, draw=L, fill=F] (130.41,210.17) circle (0.4mm);
\path[line width=0.15mm, draw=L] (99.98,240.38) -- (100.21,210.17);
\path[line width=0.45mm, draw=L] (130.42,239.97) -- (130.21,210.15);
\path[line width=0.15mm, draw=L] (100.19,209.94) -- (130.42,210.15);
\path[line width=0.15mm, draw=L, fill=F] (100.29,240.19) circle (0.4mm);
\path[line width=0.15mm, draw=L, fill=F] (130.41,240.15) circle (0.4mm);
\path[line width=0.15mm, draw=L, fill=F] (99.98,210.17) circle (0.4mm);
\path[line width=0.15mm, draw=L, fill=F] (130.41,210.17) circle (0.4mm);
\path[line width=0.15mm, draw=L, fill=F] (100.29,240.19) circle (0.4mm);
\path[line width=0.15mm, draw=L, fill=F] (130.41,240.15) circle (0.4mm);
\path[line width=0.15mm, draw=L, fill=F] (99.98,210.17) circle (0.4mm);
\path[line width=0.15mm, draw=L, fill=F] (130.41,210.17) circle (0.4mm);
\path[line width=0.15mm, draw=L, fill=F] (100.29,240.19) circle (0.4mm);
\path[line width=0.15mm, draw=L, fill=F] (130.41,240.15) circle (0.4mm);
\path[line width=0.15mm, draw=L, fill=F] (99.98,210.17) circle (0.4mm);
\path[line width=0.15mm, draw=L, fill=F] (130.41,210.17) circle (0.4mm);
\path[line width=0.15mm, draw=L, fill=F] (100.29,240.19) circle (0.4mm);
\path[line width=0.15mm, draw=L, fill=F] (130.41,240.15) circle (0.4mm);
\path[line width=0.15mm, draw=L, fill=F] (99.98,210.17) circle (0.4mm);
\path[line width=0.15mm, draw=L, fill=F] (130.41,210.17) circle (0.4mm);
\path[line width=0.15mm, draw=L, fill=F] (100.29,240.19) circle (0.4mm);
\path[line width=0.45mm, draw=L, fill=F] (130.41,240.15) circle (0.8mm);
\path[line width=0.15mm, draw=L, fill=F] (99.98,210.17) circle (0.4mm);
\path[line width=0.45mm, draw=L, fill=F] (130.41,210.17) circle (0.8mm);
\path[line width=0.15mm, draw=L] (100.44,240.24) -- (130.34,240.24);
\path[line width=0.15mm, draw=L, fill=F] (150.04,240.19) circle (0.4mm);
\path[line width=0.15mm, draw=L, fill=F] (180.17,240.15) circle (0.4mm);
\path[line width=0.15mm, draw=L, fill=F] (149.74,210.17) circle (0.4mm);
\path[line width=0.15mm, draw=L, fill=F] (180.17,210.17) circle (0.4mm);
\path[line width=0.45mm, draw=L] (180.18,239.97) -- (179.97,210.15);
\path[line width=0.15mm, draw=L] (149.95,209.94) -- (180.18,210.15);
\path[line width=0.15mm, draw=L, fill=F] (150.04,240.19) circle (0.4mm);
\path[line width=0.15mm, draw=L, fill=F] (180.17,240.15) circle (0.4mm);
\path[line width=0.15mm, draw=L, fill=F] (149.74,210.17) circle (0.4mm);
\path[line width=0.15mm, draw=L, fill=F] (180.17,210.17) circle (0.4mm);
\path[line width=0.15mm, draw=L, fill=F] (150.04,240.19) circle (0.4mm);
\path[line width=0.15mm, draw=L, fill=F] (180.17,240.15) circle (0.4mm);
\path[line width=0.15mm, draw=L, fill=F] (149.74,210.17) circle (0.4mm);
\path[line width=0.15mm, draw=L, fill=F] (180.17,210.17) circle (0.4mm);
\path[line width=0.15mm, draw=L, fill=F] (150.04,240.19) circle (0.4mm);
\path[line width=0.15mm, draw=L, fill=F] (180.17,240.15) circle (0.4mm);
\path[line width=0.15mm, draw=L, fill=F] (149.74,210.17) circle (0.4mm);
\path[line width=0.15mm, draw=L, fill=F] (180.17,210.17) circle (0.4mm);
\path[line width=0.15mm, draw=L, fill=F] (150.04,240.19) circle (0.4mm);
\path[line width=0.15mm, draw=L, fill=F] (180.17,240.15) circle (0.4mm);
\path[line width=0.15mm, draw=L, fill=F] (149.74,210.17) circle (0.4mm);
\path[line width=0.15mm, draw=L, fill=F] (180.17,210.17) circle (0.4mm);
\path[line width=0.15mm, draw=L, fill=F] (150.04,240.19) circle (0.4mm);
\path[line width=0.45mm, draw=L, fill=F] (180.17,240.15) circle (0.8mm);
\path[line width=0.15mm, draw=L, fill=F] (149.74,210.17) circle (0.4mm);
\path[line width=0.45mm, draw=L, fill=F] (180.17,210.17) circle (0.8mm);
\path[line width=0.15mm, draw=L] (150.20,240.24) -- (180.09,240.24);
\node[anchor=base,font=\small] at (165,201) {$\smash{Q_{4}}$};
\node[anchor=base,font=\small] at (65,201) {$\smash{Q_{2}}$};
\node[anchor=base,font=\small] at (115,201) {$\smash{Q_{3}}$};
\node[anchor=base,font=\small] at (15,201) {$\smash{Q_{1}}$};
\path[line width=0.15mm, draw=L, fill=F] (199.89,240.61) circle (0.4mm);
\path[line width=0.15mm, draw=L, fill=F] (230.01,240.57) circle (0.4mm);
\path[line width=0.15mm, draw=L, fill=F] (199.58,210.59) circle (0.4mm);
\path[line width=0.15mm, draw=L, fill=F] (230.01,210.59) circle (0.4mm);
\path[line width=0.15mm, draw=L] (199.58,240.80) -- (199.81,210.59);
\path[line width=0.45mm, draw=L] (230.02,240.39) -- (229.81,210.57);
\path[line width=0.15mm, draw=L, fill=F] (199.89,240.61) circle (0.4mm);
\path[line width=0.15mm, draw=L, fill=F] (230.01,240.57) circle (0.4mm);
\path[line width=0.15mm, draw=L, fill=F] (199.58,210.59) circle (0.4mm);
\path[line width=0.15mm, draw=L, fill=F] (230.01,210.59) circle (0.4mm);
\path[line width=0.15mm, draw=L, fill=F] (199.89,240.61) circle (0.4mm);
\path[line width=0.15mm, draw=L, fill=F] (230.01,240.57) circle (0.4mm);
\path[line width=0.15mm, draw=L, fill=F] (199.58,210.59) circle (0.4mm);
\path[line width=0.15mm, draw=L, fill=F] (230.01,210.59) circle (0.4mm);
\path[line width=0.15mm, draw=L, fill=F] (199.89,240.61) circle (0.4mm);
\path[line width=0.15mm, draw=L, fill=F] (230.01,240.57) circle (0.4mm);
\path[line width=0.15mm, draw=L, fill=F] (199.58,210.59) circle (0.4mm);
\path[line width=0.15mm, draw=L, fill=F] (230.01,210.59) circle (0.4mm);
\path[line width=0.15mm, draw=L, fill=F] (199.89,240.61) circle (0.4mm);
\path[line width=0.15mm, draw=L, fill=F] (230.01,240.57) circle (0.4mm);
\path[line width=0.15mm, draw=L, fill=F] (199.58,210.59) circle (0.4mm);
\path[line width=0.15mm, draw=L, fill=F] (230.01,210.59) circle (0.4mm);
\path[line width=0.15mm, draw=L, fill=F] (199.89,240.61) circle (0.4mm);
\path[line width=0.45mm, draw=L, fill=F] (230.01,240.57) circle (0.8mm);
\path[line width=0.15mm, draw=L, fill=F] (199.58,210.59) circle (0.4mm);
\path[line width=0.45mm, draw=L, fill=F] (230.01,210.59) circle (0.8mm);
\path[line width=0.15mm, draw=L] (200.04,240.66) -- (229.94,240.66);
\node[anchor=base,font=\small] at (215,201) {$\smash{Q_{5}}$};
\path[line width=0.15mm, draw=L, fill=F] (50.60,189.77) circle (0.4mm);
\path[line width=0.15mm, draw=L, fill=F] (80.73,189.73) circle (0.4mm);
\path[line width=0.15mm, draw=L, fill=F] (50.29,159.75) circle (0.4mm);
\path[line width=0.15mm, draw=L, fill=F] (80.73,159.75) circle (0.4mm);
\path[line width=0.45mm, draw=L] (80.74,189.55) -- (80.53,159.74);
\path[line width=0.15mm, draw=L, fill=F] (50.60,189.77) circle (0.4mm);
\path[line width=0.15mm, draw=L, fill=F] (80.73,189.73) circle (0.4mm);
\path[line width=0.15mm, draw=L, fill=F] (50.29,159.75) circle (0.4mm);
\path[line width=0.15mm, draw=L, fill=F] (80.73,159.75) circle (0.4mm);
\path[line width=0.15mm, draw=L, fill=F] (50.60,189.77) circle (0.4mm);
\path[line width=0.15mm, draw=L, fill=F] (80.73,189.73) circle (0.4mm);
\path[line width=0.15mm, draw=L, fill=F] (50.29,159.75) circle (0.4mm);
\path[line width=0.15mm, draw=L, fill=F] (80.73,159.75) circle (0.4mm);
\path[line width=0.15mm, draw=L, fill=F] (50.60,189.77) circle (0.4mm);
\path[line width=0.15mm, draw=L, fill=F] (80.73,189.73) circle (0.4mm);
\path[line width=0.15mm, draw=L, fill=F] (50.29,159.75) circle (0.4mm);
\path[line width=0.15mm, draw=L, fill=F] (80.73,159.75) circle (0.4mm);
\path[line width=0.15mm, draw=L, fill=F] (50.60,189.77) circle (0.4mm);
\path[line width=0.15mm, draw=L, fill=F] (80.73,189.73) circle (0.4mm);
\path[line width=0.15mm, draw=L, fill=F] (50.29,159.75) circle (0.4mm);
\path[line width=0.15mm, draw=L, fill=F] (80.73,159.75) circle (0.4mm);
\path[line width=0.15mm, draw=L, fill=F] (50.60,189.77) circle (0.4mm);
\path[line width=0.45mm, draw=L, fill=F] (80.73,189.73) circle (0.8mm);
\path[line width=0.15mm, draw=L, fill=F] (50.29,159.75) circle (0.4mm);
\path[line width=0.45mm, draw=L, fill=F] (80.73,159.75) circle (0.8mm);
\path[line width=0.15mm, draw=L] (50.76,190.00) -- (80.65,190.00);
\path[line width=0.15mm, draw=L, fill=F] (100.75,189.97) circle (0.4mm);
\path[line width=0.15mm, draw=L, fill=F] (130.87,189.92) circle (0.4mm);
\path[line width=0.15mm, draw=L, fill=F] (100.44,159.94) circle (0.4mm);
\path[line width=0.15mm, draw=L, fill=F] (130.87,159.94) circle (0.4mm);
\path[line width=0.45mm, draw=L] (130.88,189.74) -- (130.67,159.93);
\path[line width=0.15mm, draw=L, fill=F] (100.75,189.97) circle (0.4mm);
\path[line width=0.15mm, draw=L, fill=F] (130.87,189.92) circle (0.4mm);
\path[line width=0.15mm, draw=L, fill=F] (100.44,159.94) circle (0.4mm);
\path[line width=0.15mm, draw=L, fill=F] (130.87,159.94) circle (0.4mm);
\path[line width=0.15mm, draw=L, fill=F] (100.75,189.97) circle (0.4mm);
\path[line width=0.15mm, draw=L, fill=F] (130.87,189.92) circle (0.4mm);
\path[line width=0.15mm, draw=L, fill=F] (100.44,159.94) circle (0.4mm);
\path[line width=0.15mm, draw=L, fill=F] (130.87,159.94) circle (0.4mm);
\path[line width=0.15mm, draw=L, fill=F] (100.75,189.97) circle (0.4mm);
\path[line width=0.15mm, draw=L, fill=F] (130.87,189.92) circle (0.4mm);
\path[line width=0.15mm, draw=L, fill=F] (100.44,159.94) circle (0.4mm);
\path[line width=0.15mm, draw=L, fill=F] (130.87,159.94) circle (0.4mm);
\path[line width=0.15mm, draw=L, fill=F] (100.75,189.97) circle (0.4mm);
\path[line width=0.15mm, draw=L, fill=F] (130.87,189.92) circle (0.4mm);
\path[line width=0.15mm, draw=L, fill=F] (100.44,159.94) circle (0.4mm);
\path[line width=0.15mm, draw=L, fill=F] (130.87,159.94) circle (0.4mm);
\path[line width=0.15mm, draw=L, fill=F] (100.75,189.97) circle (0.4mm);
\path[line width=0.45mm, draw=L, fill=F] (130.87,189.92) circle (0.8mm);
\path[line width=0.15mm, draw=L, fill=F] (100.44,159.94) circle (0.4mm);
\path[line width=0.45mm, draw=L, fill=F] (130.87,159.94) circle (0.8mm);
\path[line width=0.15mm, draw=L] (100.40,189.76) -- (130.29,189.76);
\node[anchor=base,font=\small] at (165,151) {$\smash{Q_{9}}$};
\node[anchor=base,font=\small] at (65,151) {$\smash{Q_{7}}$};
\node[anchor=base,font=\small] at (115,151) {$\smash{Q_{8}}$};
\node[anchor=base,font=\small] at (15,151) {$\smash{Q_{6}}$};
\node[anchor=base,font=\small] at (215,151) {$\smash{Q_{10}}$};
\path[line width=0.15mm, draw=L] (50.34,189.94) -- (50.34,159.93);
\path[line width=0.15mm, draw=L, fill=F] (149.97,174.90) circle (0.4mm);
\path[line width=0.15mm, draw=L] (100.02,240.28) -- (130.22,210.21);
\path[line width=0.15mm, draw=L] (150.05,240.25) -- (179.95,210.14);
\path[line width=0.15mm, draw=L] (199.67,210.57) -- (229.99,210.46);
\path[line width=0.15mm, draw=L, fill=F] (0.19,189.97) circle (0.4mm);
\path[line width=0.15mm, draw=L, fill=F] (30.32,189.92) circle (0.4mm);
\path[line width=0.15mm, draw=L, fill=F] (-0.12,159.94) circle (0.4mm);
\path[line width=0.15mm, draw=L, fill=F] (30.32,159.94) circle (0.4mm);
\path[line width=0.45mm, draw=L] (30.33,189.74) -- (30.12,159.93);
\path[line width=0.15mm, draw=L] (0.10,159.72) -- (30.33,159.93);
\path[line width=0.15mm, draw=L, fill=F] (0.19,189.97) circle (0.4mm);
\path[line width=0.15mm, draw=L, fill=F] (30.32,189.92) circle (0.4mm);
\path[line width=0.15mm, draw=L, fill=F] (-0.12,159.94) circle (0.4mm);
\path[line width=0.15mm, draw=L, fill=F] (30.32,159.94) circle (0.4mm);
\path[line width=0.15mm, draw=L, fill=F] (0.19,189.97) circle (0.4mm);
\path[line width=0.15mm, draw=L, fill=F] (30.32,189.92) circle (0.4mm);
\path[line width=0.15mm, draw=L, fill=F] (-0.12,159.94) circle (0.4mm);
\path[line width=0.15mm, draw=L, fill=F] (30.32,159.94) circle (0.4mm);
\path[line width=0.15mm, draw=L, fill=F] (0.19,189.97) circle (0.4mm);
\path[line width=0.15mm, draw=L, fill=F] (30.32,189.92) circle (0.4mm);
\path[line width=0.15mm, draw=L, fill=F] (-0.12,159.94) circle (0.4mm);
\path[line width=0.15mm, draw=L, fill=F] (30.32,159.94) circle (0.4mm);
\path[line width=0.15mm, draw=L, fill=F] (0.19,189.97) circle (0.4mm);
\path[line width=0.15mm, draw=L, fill=F] (30.32,189.92) circle (0.4mm);
\path[line width=0.15mm, draw=L, fill=F] (-0.12,159.94) circle (0.4mm);
\path[line width=0.15mm, draw=L, fill=F] (30.32,159.94) circle (0.4mm);
\path[line width=0.15mm, draw=L, fill=F] (0.19,189.97) circle (0.4mm);
\path[line width=0.45mm, draw=L, fill=F] (30.32,189.92) circle (0.8mm);
\path[line width=0.15mm, draw=L, fill=F] (-0.12,159.94) circle (0.4mm);
\path[line width=0.45mm, draw=L, fill=F] (30.32,159.94) circle (0.8mm);
\path[line width=0.15mm, draw=L] (-0.15,189.76) -- (29.74,189.76);
\path[line width=0.15mm, draw=L] (80.62,190.05) -- (50.16,159.97);
\path[line width=0.15mm, draw=L] (130.83,189.67) -- (100.25,159.97);
\path[line width=0.45mm, draw=L, fill=F] (180.01,190.18) circle (0.8mm);
\path[line width=0.45mm, draw=L, fill=F] (179.88,160.10) circle (0.8mm);
\path[line width=0.45mm, draw=L] (179.88,160.10) -- (180.01,190.05);
\path[line width=0.15mm, draw=L] (180.01,190.05) -- (149.93,175.01);
\path[line width=0.15mm, draw=L] (149.93,175.01) -- (179.88,160.10);
\path[line width=0.15mm, draw=L, fill=F] (199.93,175.16) circle (0.4mm);
\path[line width=0.45mm, draw=L, fill=F] (229.96,190.43) circle (0.8mm);
\path[line width=0.45mm, draw=L, fill=F] (229.84,160.35) circle (0.8mm);
\path[line width=0.45mm, draw=L] (229.84,160.35) -- (229.96,190.30);
\path[line width=0.15mm, draw=L] (229.96,190.30) -- (199.89,175.26);
\end{tikzpicture}%

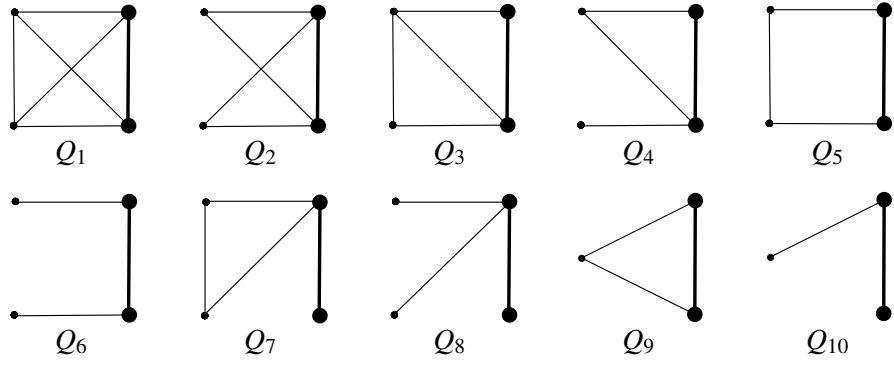
\captionof{figure}{Some possible substructures $Q_1,Q_2,\ldots,Q_{10}$.}
\label{diam4:Qfig}
\end{minipage}\par
\addvspace{2mm}

\begin{proof}
Since $G$ is connected, $E(Y)$ is a matching, and vertices of $Y$ have no neighbors in $V(P)\setminus X$, the edge $H$ has a neighbor in $X$. Hence $|N_X(V(H))|\geq1$.
If $|N_X(V(H))|\geq3$, then $G$ contains a subdivision of $K_{3,3}$, a contradiction.
In Fig.~\ref{diam4:Qfig}, the large vertices and thick edge represent $H$ in $Y$.
If $|N_X(V(H))|= 2$, then $G[N_X(V(H))\cup V(H)]=Q_i$, $1\leq i\leq 8$.
Since $R_{14}\not\subset G$ $(R_{7}\not\subset G)$, we have $G[N_X(V(H))\cup V(H)]\neq Q_3$ $(Q_4)$, respectively.
Moreover, we claim that $|X|=2$.
Otherwise, suppose that $|X|\geq3$.
If $G[N_X(V(H))\cup V(H)]=Q_i$, $i=1,2$, then $R_{4}\subset G$, a contradiction.
If $G[N_X(V(H))\cup V(H)]=Q_i$, $i=5,6,7,8$, then $R_{1}\subset G$, a contradiction.
If $|N_X(V(H))|= 1$, then $G[N_X(V(H))\cup V(H)]=Q_i$, $i=9,10$.
Similarly, since $R_{1},R_{4}\not\subset G$, we have $|X|=1$.
\end{proof}

Let $T_{w_{i_1}\cdots w_{i_k}}$ be the set of isolated vertices in $G[Y]$
which are adjacent precisely to the vertices $w_{i_1},\ldots, w_{i_k}$ of $X$.
Since $G$ does not contain $K_{3,3}$ as a subgraph, we have $1\leq k\leq 2$.

\par\addvspace{2mm}
\noindent\begin{minipage}{\linewidth}
\captionsetup{skip=3pt}
\centering
\begin{tikzpicture}[x=0.35mm, y=0.35mm, inner xsep=0pt, inner ysep=0pt, outer xsep=0pt, outer ysep=0pt]
\definecolor{L}{rgb}{0,0,0}
\definecolor{F}{rgb}{0,0,0}
\node[anchor=base,font=\small] at (47.5,118) {$\smash{W_{1}}$};
\path[line width=0.15mm, draw=L, fill=F] (208.90,156.83) circle (0.4mm);
\path[line width=0.15mm, draw=L, fill=F] (218.98,157.08) circle (0.4mm);
\path[line width=0.15mm, draw=L] (208.65,156.83) -- (218.98,157.08);
\path[line width=0.15mm, draw=L, fill=F] (248.78,177.18) circle (0.4mm);
\path[line width=0.15mm, draw=L] (248.78,177.18) -- (219.31,156.78);
\path[line width=0.15mm, draw=L] (219.31,156.78) -- (248.78,137.13);
\path[line width=0.15mm, draw=L, fill=F] (248.78,137.13) circle (0.4mm);
\path[line width=0.15mm, draw=L] (248.78,137.13) -- (278.76,157.03);
\path[line width=0.15mm, draw=L, fill=F] (278.76,157.03) circle (0.4mm);
\path[line width=0.15mm, draw=L] (278.76,157.03) -- (248.78,176.93);
\path[line width=0.15mm, draw=L, fill=F] (248.78,157.03) circle (0.4mm);
\path[line width=0.15mm, draw=L] (218.81,156.78) -- (278.25,157.28);
\path[line width=0.15mm, draw=L, fill=F] (288.84,157.28) circle (0.4mm);
\path[line width=0.15mm, draw=L] (288.84,157.28) -- (278.51,157.03);
\node[anchor=base,font=\small] at (148.04,118) {$\smash{W_{2}}$};
\path[line width=0.15mm, draw=L, fill=F] (306.97,163.47) circle (0.4mm);
\path[line width=0.15mm, draw=L, fill=F] (317.05,163.72) circle (0.4mm);
\path[line width=0.15mm, draw=L] (306.72,163.47) -- (317.05,163.72);
\path[line width=0.15mm, draw=L, fill=F] (346.85,183.82) circle (0.4mm);
\path[line width=0.15mm, draw=L] (346.85,183.82) -- (317.38,163.42);
\path[line width=0.15mm, draw=L] (317.38,163.42) -- (346.85,143.78);
\path[line width=0.15mm, draw=L, fill=F] (346.85,143.78) circle (0.4mm);
\path[line width=0.15mm, draw=L] (346.85,143.78) -- (376.82,163.68);
\path[line width=0.15mm, draw=L, fill=F] (376.82,163.68) circle (0.4mm);
\path[line width=0.15mm, draw=L] (376.82,163.68) -- (346.85,183.57);
\path[line width=0.15mm, draw=L, fill=F] (346.85,163.68) circle (0.4mm);
\path[line width=0.15mm, draw=L] (316.88,163.42) -- (376.32,163.93);
\path[line width=0.15mm, draw=L, fill=F] (386.91,163.93) circle (0.4mm);
\path[line width=0.15mm, draw=L] (386.91,163.93) -- (376.58,163.68);
\path[line width=0.15mm, draw=L, fill=F] (347.18,173.66) circle (0.4mm);
\path[line width=0.15mm, draw=L] (347.18,173.66) -- (316.92,163.96);
\path[line width=0.15mm, draw=L] (346.95,173.66) -- (376.74,163.73);
\path[line width=0.15mm, draw=L, fill=F] (346.86,194.75) circle (0.4mm);
\path[line width=0.15mm, draw=L] (346.86,194.75) -- (317.01,163.68);
\path[line width=0.15mm, draw=L] (346.86,194.75) -- (377.31,163.99);
\path[line width=0.15mm, draw=L, fill=F] (35.50,67.50) circle (0.4mm);
\path[line width=0.15mm, draw=L, fill=F] (45.58,67.75) circle (0.4mm);
\path[line width=0.15mm, draw=L] (35.25,67.50) -- (45.58,67.75);
\path[line width=0.15mm, draw=L] (45.91,67.45) -- (75.38,47.81);
\path[line width=0.15mm, draw=L, fill=F] (75.38,47.81) circle (0.4mm);
\path[line width=0.15mm, draw=L] (75.38,47.81) -- (105.35,67.70);
\path[line width=0.15mm, draw=L, fill=F] (105.35,67.70) circle (0.4mm);
\path[line width=0.15mm, draw=L, fill=F] (115.44,67.96) circle (0.4mm);
\path[line width=0.15mm, draw=L] (115.44,67.96) -- (105.11,67.70);
\path[line width=0.15mm, draw=L, fill=F] (75.38,98.78) circle (0.4mm);
\path[line width=0.15mm, draw=L] (75.38,98.78) -- (45.54,67.71);
\path[line width=0.15mm, draw=L] (75.38,98.78) -- (105.84,68.02);
\path[line width=0.15mm, draw=L, fill=F] (192.45,76.56) circle (0.4mm);
\path[line width=0.15mm, draw=L, dash pattern=on 2.00mm off 1.00mm] (179.54,67.34);
\path[line width=0.15mm, draw=L, fill=F] (216.75,71.47) circle (0.15mm);
\path[line width=0.15mm, draw=L, fill=F] (216.75,76.13) circle (0.15mm);
\path[line width=0.15mm, draw=L, fill=F] (216.75,81.11) circle (0.15mm);
\path[line width=0.15mm, draw=L, fill=F] (178.54,66.19) circle (0.15mm);
\path[line width=0.15mm, draw=L, fill=F] (174.87,72.91) circle (0.15mm);
\path[line width=0.15mm, draw=L, fill=F] (176.67,69.33) circle (0.15mm);
\node[anchor=base,font=\small] at (248.78,118) {$\smash{W_{3}}$};
\node[anchor=base,font=\small] at (346.85,118) {$\smash{W_{4}}$};
\node[anchor=base,font=\small] at (75.4,35.42) {$\smash{W_{5}}$};
\path[line width=0.15mm, draw=L, fill=F] (248.95,196.20) circle (0.4mm);
\path[line width=0.15mm, draw=L] (248.95,196.20) -- (219.06,157.01);
\path[line width=0.15mm, draw=L] (248.95,196.20) -- (278.84,157.01);
\path[line width=0.15mm, draw=L] (192.52,76.89) -- (207.79,85.86);
\path[line width=0.15mm, draw=L, fill=F] (207.79,85.86) circle (0.4mm);
\path[line width=0.15mm, draw=L] (207.79,85.86) -- (226.06,86.19);
\path[line width=0.15mm, draw=L, fill=F] (226.06,86.19) circle (0.4mm);
\path[line width=0.15mm, draw=L] (192.85,76.56) -- (207.79,66.60);
\path[line width=0.15mm, draw=L, fill=F] (207.79,66.60) circle (0.4mm);
\path[line width=0.15mm, draw=L] (207.79,66.60) -- (226.06,66.60);
\path[line width=0.15mm, draw=L, fill=F] (226.06,66.60) circle (0.4mm);
\path[line width=0.15mm, draw=L] (192.52,76.56) -- (189.20,51.98);
\path[line width=0.15mm, draw=L, fill=F] (189.20,51.98) circle (0.4mm);
\path[line width=0.15mm, draw=L] (189.20,51.98) -- (178.90,55.97);
\path[line width=0.15mm, draw=L, fill=F] (178.90,55.97) circle (0.4mm);
\path[line width=0.15mm, draw=L] (178.90,55.97) -- (192.52,76.89);
\path[line width=0.15mm, draw=L] (192.52,76.23) -- (168.94,88.85);
\path[line width=0.15mm, draw=L, fill=F] (168.94,88.85) circle (0.4mm);
\path[line width=0.15mm, draw=L] (168.94,88.85) -- (164.62,77.56);
\path[line width=0.15mm, draw=L, fill=F] (164.62,77.56) circle (0.4mm);
\path[line width=0.15mm, draw=L] (164.62,77.56) -- (192.52,76.56);
\path[line width=0.15mm, draw=L, fill=F] (285.59,81.71) circle (0.4mm);
\path[line width=0.15mm, draw=L] (275.43,81.48) -- (285.36,81.48);
\path[line width=0.15mm, draw=L] (285.36,81.48) -- (315.39,101.34);
\path[line width=0.15mm, draw=L] (315.15,101.34) -- (346.10,81.48);
\path[line width=0.15mm, draw=L, fill=F] (346.10,81.48) circle (0.4mm);
\path[line width=0.15mm, draw=L] (346.10,81.48) -- (356.04,81.71);
\path[line width=0.15mm, draw=L, fill=F] (356.04,81.71) circle (0.4mm);
\node[anchor=base,font=\small] at (315,35.42) {$\smash{W_{7}}$};
\path[line width=0.15mm, draw=L, fill=F] (275.30,81.38) circle (0.4mm);
\path[line width=0.15mm, draw=L] (315.95,82.30) -- (285.92,81.61);
\path[line width=0.15mm, draw=L] (315.95,82.07) -- (345.97,81.38);
\path[line width=0.15mm, draw=L, fill=F] (315.11,101.21) circle (0.4mm);
\path[line width=0.15mm, draw=L] (315.11,101.21) -- (314.90,92.16);
\path[line width=0.15mm, draw=L, fill=F] (314.90,92.16) circle (0.4mm);
\path[line width=0.15mm, draw=L] (314.90,92.16) -- (285.82,81.39);
\path[line width=0.15mm, draw=L] (314.90,92.16) -- (346.12,81.39);
\path[line width=0.15mm, draw=L, fill=F] (314.68,72.13) circle (0.4mm);
\path[line width=0.15mm, draw=L] (314.68,72.13) -- (346.12,81.61);
\path[line width=0.15mm, draw=L] (314.90,72.13) -- (286.04,81.61);
\path[line width=0.15mm, draw=L, fill=F] (314.90,82.04) circle (0.4mm);
\path[line width=0.15mm, draw=L] (314.90,82.04) -- (314.90,72.13);
\path[line width=0.15mm, draw=L, fill=F] (314.72,87.11) circle (0.15mm);
\path[line width=0.15mm, draw=L, fill=F] (314.72,89.74) circle (0.15mm);
\path[line width=0.15mm, draw=L, fill=F] (314.72,84.57) circle (0.15mm);
\path[line width=0.15mm, draw=L, fill=F] (314.73,62.12) circle (0.4mm);
\path[line width=0.15mm, draw=L] (314.73,62.12) -- (346.20,81.50);
\path[line width=0.15mm, draw=L] (314.73,62.00) -- (285.61,81.50);
\path[line width=0.15mm, draw=L, fill=F] (314.72,56.58) circle (0.15mm);
\path[line width=0.15mm, draw=L, fill=F] (314.72,59.20) circle (0.15mm);
\path[line width=0.15mm, draw=L, fill=F] (314.72,54.04) circle (0.15mm);
\path[line width=0.15mm, draw=L, fill=F] (314.73,51.44) circle (0.4mm);
\path[line width=0.15mm, draw=L] (314.73,51.44) -- (346.09,81.50);
\path[line width=0.15mm, draw=L] (314.73,51.55) -- (285.49,81.61);
\path[line width=0.15mm, draw=L] (314.73,62.12) -- (321.07,62.12);
\path[line width=0.15mm, draw=L, fill=F] (321.07,62.12) circle (0.4mm);
\node[anchor=base,font=\small] at (192.45,35.42) {$\smash{W_{6}}$};
\path[line width=0.15mm, draw=L] (314.89,51.44) -- (321.45,51.44);
\path[line width=0.15mm, draw=L, fill=F] (320.96,51.44) circle (0.4mm);
\path[line width=0.15mm, draw=L] (248.79,196.13) -- (248.79,177.25);
\path[line width=0.15mm, draw=L, fill=F] (259.97,196.13) circle (0.4mm);
\path[line width=0.15mm, draw=L] (259.97,196.13) -- (248.79,196.13);
\path[line width=0.15mm, draw=L, fill=F] (258.43,177.25) circle (0.4mm);
\path[line width=0.15mm, draw=L] (258.43,177.25) -- (248.79,177.25);
\path[line width=0.15mm, draw=L] (347.03,194.71) -- (347.03,183.49);
\path[line width=0.15mm, draw=L, fill=F] (357.79,194.94) circle (0.4mm);
\path[line width=0.15mm, draw=L] (347.03,194.94) -- (357.79,194.94);
\path[line width=0.15mm, draw=L, fill=F] (75.44,79.99) circle (0.4mm);
\path[line width=0.15mm, draw=L] (75.44,79.99) -- (45.43,67.62);
\path[line width=0.15mm, draw=L] (75.44,79.99) -- (105.68,67.85);
\path[line width=0.15mm, draw=L, fill=F] (80.48,86.17) circle (0.4mm);
\path[line width=0.15mm, draw=L] (80.48,86.17) -- (75.44,98.77);
\path[line width=0.15mm, draw=L] (80.48,86.17) -- (75.44,79.76);
\path[line width=0.15mm, draw=L, fill=F] (75.38,64.31) circle (0.4mm);
\path[line width=0.15mm, draw=L] (75.38,64.31) -- (45.72,67.58);
\path[line width=0.15mm, draw=L] (75.57,64.12) -- (105.43,67.78);
\coordinate (qw1v0) at (17.50000,177.50000);
\coordinate (qw1v1) at (77.50000,177.50000);
\coordinate (qw1v2) at (7.50000,177.50000);
\coordinate (qw1v3) at (87.50000,177.50000);
\coordinate (qw1v4) at (47.50000,197.50000);
\coordinate (qw1v5) at (47.50000,184.50000);
\coordinate (qw1v6) at (47.50000,178.50000);
\coordinate (qw1v7) at (47.50000,172.50000);
\coordinate (qw1v8) at (47.50000,165.00000);
\coordinate (qw1v9) at (47.50000,159.00000);
\coordinate (qw1v10) at (47.50000,153.00000);
\coordinate (qw1v11) at (47.50000,147.00000);
\coordinate (qw1v12) at (47.50000,141.00000);
\coordinate (qw1v13) at (47.50000,135.00000);
\draw[line width=0.15mm,draw=L] (qw1v0) -- (qw1v2);
\draw[line width=0.15mm,draw=L] (qw1v1) -- (qw1v3);
\draw[line width=0.15mm,draw=L] (qw1v0) -- (qw1v4);
\draw[line width=0.15mm,draw=L] (qw1v1) -- (qw1v4);
\draw[line width=0.15mm,draw=L] (qw1v0) -- (qw1v5);
\draw[line width=0.15mm,draw=L] (qw1v1) -- (qw1v5);
\draw[line width=0.15mm,draw=L] (qw1v0) -- (qw1v6);
\draw[line width=0.15mm,draw=L] (qw1v1) -- (qw1v6);
\draw[line width=0.15mm,draw=L] (qw1v0) -- (qw1v7);
\draw[line width=0.15mm,draw=L] (qw1v1) -- (qw1v7);
\draw[line width=0.15mm,draw=L] (qw1v0) -- (qw1v8);
\draw[line width=0.15mm,draw=L] (qw1v1) -- (qw1v8);
\draw[line width=0.15mm,draw=L] (qw1v0) -- (qw1v9);
\draw[line width=0.15mm,draw=L] (qw1v1) -- (qw1v9);
\draw[line width=0.15mm,draw=L] (qw1v0) -- (qw1v10);
\draw[line width=0.15mm,draw=L] (qw1v1) -- (qw1v10);
\draw[line width=0.15mm,draw=L] (qw1v0) -- (qw1v11);
\draw[line width=0.15mm,draw=L] (qw1v1) -- (qw1v11);
\draw[line width=0.15mm,draw=L] (qw1v0) -- (qw1v12);
\draw[line width=0.15mm,draw=L] (qw1v1) -- (qw1v12);
\draw[line width=0.15mm,draw=L] (qw1v0) -- (qw1v13);
\draw[line width=0.15mm,draw=L] (qw1v1) -- (qw1v13);
\draw[line width=0.15mm,draw=L] (qw1v4) -- (qw1v5);
\draw[line width=0.15mm,draw=L] (qw1v5) -- (qw1v6);
\draw[line width=0.15mm,draw=L] (qw1v6) -- (qw1v7);
\path[line width=0.15mm,draw=L,fill=F] (qw1v0) circle (0.4mm);
\path[line width=0.15mm,draw=L,fill=F] (qw1v1) circle (0.4mm);
\path[line width=0.15mm,draw=L,fill=F] (qw1v2) circle (0.4mm);
\path[line width=0.15mm,draw=L,fill=F] (qw1v3) circle (0.4mm);
\path[line width=0.15mm,draw=L,fill=F] (qw1v4) circle (0.4mm);
\path[line width=0.15mm,draw=L,fill=F] (qw1v5) circle (0.4mm);
\path[line width=0.15mm,draw=L,fill=F] (qw1v6) circle (0.4mm);
\path[line width=0.15mm,draw=L,fill=F] (qw1v7) circle (0.4mm);
\path[line width=0.15mm,draw=L,fill=F] (qw1v8) circle (0.4mm);
\path[line width=0.15mm,draw=L,fill=F] (qw1v9) circle (0.4mm);
\path[line width=0.15mm,draw=L,fill=F] (qw1v10) circle (0.4mm);
\path[line width=0.15mm,draw=L,fill=F] (qw1v11) circle (0.4mm);
\path[line width=0.15mm,draw=L,fill=F] (qw1v12) circle (0.4mm);
\path[line width=0.15mm,draw=L,fill=F] (qw1v13) circle (0.4mm);
\coordinate (qw2v0) at (118.04000,177.50000);
\coordinate (qw2v1) at (178.04000,177.50000);
\coordinate (qw2v2) at (108.04000,177.50000);
\coordinate (qw2v3) at (188.04000,177.50000);
\coordinate (qw2v4) at (148.04000,197.50000);
\coordinate (qw2v5) at (148.04000,184.50000);
\coordinate (qw2v6) at (148.04000,172.50000);
\coordinate (qw2v7) at (148.04000,165.00000);
\coordinate (qw2v8) at (148.04000,160.00000);
\coordinate (qw2v9) at (148.04000,155.00000);
\coordinate (qw2v10) at (148.04000,150.00000);
\coordinate (qw2v11) at (148.04000,145.00000);
\coordinate (qw2v12) at (148.04000,140.00000);
\coordinate (qw2v13) at (148.04000,135.00000);
\draw[line width=0.15mm,draw=L] (qw2v0) -- (qw2v2);
\draw[line width=0.15mm,draw=L] (qw2v1) -- (qw2v3);
\draw[line width=0.15mm,draw=L] (qw2v0) -- (qw2v4);
\draw[line width=0.15mm,draw=L] (qw2v1) -- (qw2v4);
\draw[line width=0.15mm,draw=L] (qw2v0) -- (qw2v5);
\draw[line width=0.15mm,draw=L] (qw2v1) -- (qw2v5);
\draw[line width=0.15mm,draw=L] (qw2v0) -- (qw2v6);
\draw[line width=0.15mm,draw=L] (qw2v1) -- (qw2v6);
\draw[line width=0.15mm,draw=L] (qw2v0) -- (qw2v7);
\draw[line width=0.15mm,draw=L] (qw2v1) -- (qw2v7);
\draw[line width=0.15mm,draw=L] (qw2v0) -- (qw2v8);
\draw[line width=0.15mm,draw=L] (qw2v1) -- (qw2v8);
\draw[line width=0.15mm,draw=L] (qw2v0) -- (qw2v9);
\draw[line width=0.15mm,draw=L] (qw2v1) -- (qw2v9);
\draw[line width=0.15mm,draw=L] (qw2v0) -- (qw2v10);
\draw[line width=0.15mm,draw=L] (qw2v1) -- (qw2v10);
\draw[line width=0.15mm,draw=L] (qw2v0) -- (qw2v11);
\draw[line width=0.15mm,draw=L] (qw2v1) -- (qw2v11);
\draw[line width=0.15mm,draw=L] (qw2v0) -- (qw2v12);
\draw[line width=0.15mm,draw=L] (qw2v1) -- (qw2v12);
\draw[line width=0.15mm,draw=L] (qw2v0) -- (qw2v13);
\draw[line width=0.15mm,draw=L] (qw2v1) -- (qw2v13);
\draw[line width=0.15mm,draw=L] (qw2v4) -- (qw2v5);
\draw[line width=0.15mm,draw=L] (qw2v5) -- (qw2v6);
\path[line width=0.15mm,draw=L,fill=F] (qw2v0) circle (0.4mm);
\path[line width=0.15mm,draw=L,fill=F] (qw2v1) circle (0.4mm);
\path[line width=0.15mm,draw=L,fill=F] (qw2v2) circle (0.4mm);
\path[line width=0.15mm,draw=L,fill=F] (qw2v3) circle (0.4mm);
\path[line width=0.15mm,draw=L,fill=F] (qw2v4) circle (0.4mm);
\path[line width=0.15mm,draw=L,fill=F] (qw2v5) circle (0.4mm);
\path[line width=0.15mm,draw=L,fill=F] (qw2v6) circle (0.4mm);
\path[line width=0.15mm,draw=L,fill=F] (qw2v7) circle (0.4mm);
\path[line width=0.15mm,draw=L,fill=F] (qw2v8) circle (0.4mm);
\path[line width=0.15mm,draw=L,fill=F] (qw2v9) circle (0.4mm);
\path[line width=0.15mm,draw=L,fill=F] (qw2v10) circle (0.4mm);
\path[line width=0.15mm,draw=L,fill=F] (qw2v11) circle (0.4mm);
\path[line width=0.15mm,draw=L,fill=F] (qw2v12) circle (0.4mm);
\path[line width=0.15mm,draw=L,fill=F] (qw2v13) circle (0.4mm);
\end{tikzpicture}%

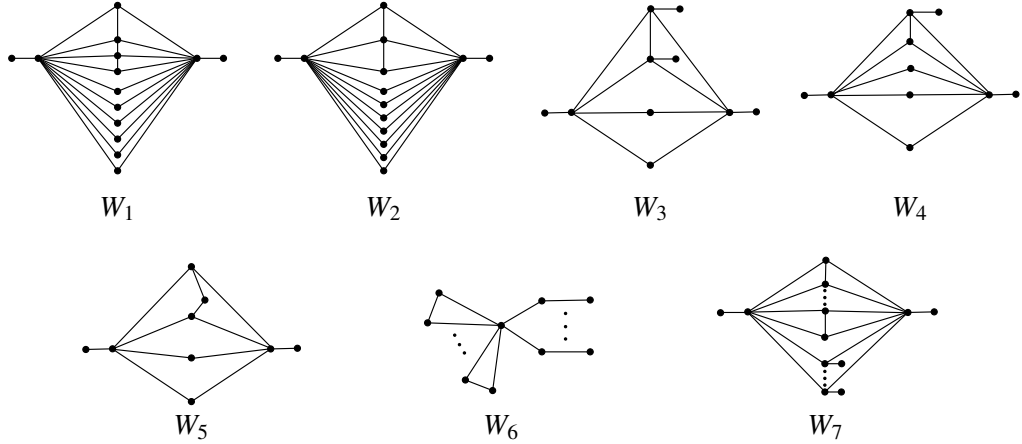
\captionof{figure}{The graphs $W_1,W_2,\ldots,W_{7}$.}
\label{diam4:Wfig}
\end{minipage}\par
\addvspace{2mm}

\begin{lem}\label{diam4:Xsmall}
If $|X|\leq2$, then $G$ is an induced subgraph of $D_1$ or $D_2$.
\end{lem}

\begin{proof}
Suppose first that $|X|=1$. If there exists $H=K_2\subset G[Y]$, then Lemma \ref{diam4:NXK2} gives $G[X\cup V(H)]=Q_i$ for some $i\in\{9,10\}$.
If $G[Y]$ contains no edges, then $G[X\cup Y]$ is a star.
Therefore, we have $G\subset W_6\subset D_1$ (see Fig.~\ref{diam4:Wfig}).
Next, suppose that $|X|=2$.
Let $X=\{w_1,w_2\}$.
When $G[X]=2K_1$, Lemma \ref{diam4:NXK2} implies that every edge $H$ in $G[Y]$ satisfies $G[X\cup V(H)]=Q_i$ for some $i\in\{2,6,8\}$.
Let $E_1$ be the set of edges $H$ in $G[Y]$ such that $G[X\cup V(H)]=Q_6$.
Since $R_{1},R_{21}\not\subset G$, we have $|T_{w_1}|+|T_{w_2}|+|E_1|\leq 2$.
Then $G\subset D_1$.
When $G[X]=K_2$, Lemma \ref{diam4:NXK2} implies that every edge $H$ in $G[Y]$ satisfies $G[X\cup V(H)]=Q_i$ for some $i\in\{1,5,7\}$.
Let $E_2$ be the set of edges $H$ in $G[Y]$ such that $G[X\cup V(H)]=Q_5$.
Since $R_{21},R_{45},R_{46}\not\subset G$, we have $|T_{w_1}|+|T_{w_2}|+|E_2|\leq 4$.
Then $G\subset D_2$.
\end{proof}

%

\begin{lem}\label{diam4:Xcycle}
If $G[X]\cong C_a$ for some $3\leq a\leq6$, then $G$ is an induced subgraph of one of $A_1,\ldots,A_{18}$.
\end{lem}

\begin{proof}
Let $C_a=w_1w_2\cdots w_aw_1$.
Since $|X|\geq3$, Lemma \ref{diam4:NXK2} gives $G[Y]=sK_1$ for some $s\geq0$.
We distinguish the following two cases.

\noindent\textbf{Case 1:} $G[X]=C_3$.

First, by Lemma \ref{lem-1} and Corollary \ref{cor-subgraph}, we have the following facts.

\begin{fact}\label{diam4:f1}
$|T_{w_iw_j}|\leq3$ and $|T_{w_i}|\leq3$ for $1\leq i<j\leq3$. ($R_{62},R_{54}\not\subset G$)\footnote{This notation records the forbidden induced subgraphs from which the stated bounds follow; the same convention is used below.}.
\end{fact}

\begin{fact}\label{diam4:f2}
If $|T_{w_iw_j}|\geq 1$ for some $1\leq i<j\leq 3$, then $|T_{w_i}|\leq 2$ and $|T_{w_j}|\leq 2$. ($R_{51}\not\subset G$).
\end{fact}

\begin{fact}\label{diam4:f3}
If $|T_{w_iw_j}|\geq1$ and $|T_{w_iw_k}|\geq1$, then $|T_{w_j}|\leq1$ and $|T_{w_k}|\leq1$. ($R_{25}\not\subset G$).
\end{fact}

\begin{fact}\label{diam4:f4}
If $|T_{w_iw_j}|\geq2$ and $|T_{w_jw_k}|\geq1$, then $|T_{w_k}|=0$. ($R_{28}\not\subset G$).
\end{fact}

Without loss of generality, we may assume that $|T_{w_1w_2}|\geq |T_{w_2w_3}|\geq|T_{w_1w_3}|$.
By (F\ref{diam4:f1}), we distinguish the following four cases.

\noindent\textbf{Subcase 1.1:} $|T_{w_1w_2}|=3$.

Then $|T_{w_1w_3}|=|T_{w_2w_3}|=0$ for $R_{59}\not\subset G$.
Since $R_{56},R_{58}\not\subset G$, we have $|T_{w_i}|=0$ for $i=1,2,3$.
Hence $G=A_1$.

\noindent\textbf{Subcase 1.2:} $|T_{w_1w_2}|=2$.

If $|T_{w_1w_2}|=|T_{w_1w_3}|=|T_{w_2w_3}|=2$, then by (F\ref{diam4:f4}), we have $|T_{w_i}|=0$ for $i=1,2,3$. Hence $G=A_2$.
If $|T_{w_1w_2}|=|T_{w_2w_3}|=2$ and $|T_{w_1w_3}|=1$, then by (F\ref{diam4:f3}) and (F\ref{diam4:f4}), we have $|T_{w_2}|\leq 1$ and $|T_{w_1}|=|T_{w_3}|=0$. Hence $G\subset A_3$.
If $|T_{w_1w_2}|=|T_{w_2w_3}|=2$ and $|T_{w_1w_3}|=0$, then by (F\ref{diam4:f2}) and (F\ref{diam4:f4}), we have $|T_{w_2}|\leq 2$ and $|T_{w_1}|=|T_{w_3}|=0$. Hence $G\subset A_4$.
If $|T_{w_1w_2}|=2$ and $|T_{w_1w_3}|=|T_{w_2w_3}|=1$, then by (F\ref{diam4:f3}) and (F\ref{diam4:f4}), we have $|T_{w_1}|\leq 1$, $|T_{w_2}|\leq 1$ and $|T_{w_3}|=0$. Hence $G\subset A_5$.
If $|T_{w_1w_2}|=2$, $|T_{w_2w_3}|=1$ and $|T_{w_1w_3}|=0$, then by (F\ref{diam4:f2}), (F\ref{diam4:f3}), and (F\ref{diam4:f4}), we have $|T_{w_3}|=0$, $|T_{w_1}|\leq 1$ and $|T_{w_2}|\leq 2$. Hence $G\subset A_6$.
If $|T_{w_1w_2}|=2$, $|T_{w_2w_3}|=|T_{w_1w_3}|=0$, since $R_{24}\not\subset G$, we have $|T_{w_3}|\leq 1$.
When $|T_{w_3}|= 1$, since $R_{52}\not\subset G$, we have $|T_{w_1}|=|T_{w_2}|=0$. Hence $G=A_7$.
When $|T_{w_3}|= 0$, by (F\ref{diam4:f2}), we have $|T_{w_1}|\leq 2$ and $|T_{w_2}|\leq 2$. Hence $G\subset A_8$.

\noindent\textbf{Subcase 1.3:} $|T_{w_1w_2}|=1$.

If $|T_{w_1w_2}|=|T_{w_1w_3}|=|T_{w_2w_3}|=1$, then by (F\ref{diam4:f3}), we have $|T_{w_i}|\leq 1\,(1\leq i\leq 3)$. Hence $G\subset A_9$.
If $|T_{w_1w_2}|=|T_{w_2w_3}|=1$ and $|T_{w_1w_3}|=0$, then by (F\ref{diam4:f2}) and (F\ref{diam4:f3}), we have $|T_{w_1}|\leq 1$, $|T_{w_2}|\leq 2$ and $|T_{w_3}|\leq 1$. Hence $G\subset A_{10}$.
If $|T_{w_1w_2}|=1$ and $|T_{w_1w_3}|=|T_{w_2w_3}|=0$, since $R_{23}\not\subset G$, we have $|T_{w_3}|\leq 2$.
When $|T_{w_3}|= 2$, since $R_{22}\not\subset G$, we have $|T_{w_1}|=|T_{w_2}|=0$. Hence $G=A_{11}$.
When $|T_{w_3}|\leq 1$, by (F\ref{diam4:f2}), we have $|T_{w_1}|\leq 2$, $|T_{w_2}|\leq 2$. Hence $G\subset A_{12}$.

\noindent\textbf{Subcase 1.4:} $|T_{w_1w_2}|=0$.

By (F\ref{diam4:f1}), if $|T_{w_i}|=3$ for some $w_i$, say $w_1$, then $|T_{w_2}|=|T_{w_3}|=0$, for $R_{48}\not\subset G$. Hence $G=A_{13}$.
Otherwise, $G\subset A_{14}$.

\noindent\textbf{Case 2:} $G[X]=C_a\,(4\leq a\leq 6)$.

For any $w_i\in X$, there exists a $P_3\subset G[X]$ such that $w_i$ is a pendant vertex of $P_3$.
Therefore, since $R_{22}\not\subset G$, we have $|T_{w_i}|\leq 1\,(1\leq i\leq a)$.
If $|T_{w_iw_j}|\neq 0$, then $j=i+1 \bmod a$.
Otherwise, $G$ contains a subdivision of $K_5$, a contradiction.
When $a=4$, since $R_{22}\not\subset G$, we have $|T_{w_iw_{i+1}}|\leq 1$.
Moreover, since $R_{14}\not\subset G$, if $|T_{w_iw_{i+1}}|\neq 0$, then $|T_{w_{i+2}w_{i+3}}|= 0$.
If $|T_{w_i}|=1$ for some $w_i$, then $|T_{w_sw_t}|=0\,(1\leq s<t\leq a)$.
Otherwise, $R_{22}$ or $R_{7}\subset G$, a contradiction.
Therefore, if $|T_{w_i}|=1$ for some $w_i$, then $G\subset A_{15}$.
Otherwise, $G\subset A_{16}$.
When $a=5$, since $R_{14}\not\subset G$, we have $|T_{w_i}|=0\,(1\leq i\leq 5)$.
Since $R_{82}\not\subset G$, we have $|T_{w_iw_j}|=0\,(1\leq i<j\leq 5)$.
Hence $G=A_{17}$.
When $a=6$, since $P_{6}\not\subset G$, we have $|T_{w_i}|=0\,(1\leq i\leq 6)$ and $|T_{w_iw_j}|=0\,(1\leq i<j\leq 6)$.
Hence $G=A_{18}$.
This completes the proof.
\end{proof}

\begin{lem}\label{diam4:Xforest}
If $|X|\geq3$ and $G[X]$ contains no cycle, then $G$ is an induced subgraph of $A_i$ $(i=19,\ldots,51)$ or $D_1$.
\end{lem}

\begin{proof}
By Lemma \ref{diam4:Xstructure}, if $G[X]$ contains a $P_b$, then $1\leq b\leq 5$.
Since $|X|\geq3$, by Lemma \ref{diam4:NXK2}, we have $G[Y]=tK_1\,(t\geq0)$.
For any $u\in Y $, we have $d(u)\leq 2$. Otherwise, $G$ contains $K_{3,3}$ as a subgraph, a contradiction.
If there exists a $P_b\,(b\geq 3)\subset G[X]$, since $R_{16}\not\subset G$, we have $G[X]=P_b\cup sP_1\,(b\geq 3,s\geq 0)$.
We now distinguish the following five cases.

\noindent\textbf{Case 1:} $G[X]=P_5\cup sP_1\,(s\geq0)$.

Let $P_5={w_1w_2w_3w_4w_5}$ and $V(sP_1)=\{w_6,\ldots,w_{5+s}\}$.
We claim that $d(u)=1$ for any $u\in Y$.
Suppose that there exists a vertex $u$ such that $d(u)=2$.
Let $N_X(u)=\{w_i,w_j\}, i<j$.
If $G[\{u,w_i,w_j\}]=K_3$, then $1\leq i<j\leq 5$.
Then $R_{18}\subset G$, a contradiction.
Next, suppose that $G[\{u,w_i,w_j\}]=P_3$.
If $1\leq i<j\leq 5$, then $i=1$, $j=5$.
Otherwise, $G$ has a subdivision of $K_5$, a contradiction.
Then $R_{12}\subset G$, a contradiction.
If $1\leq i\leq 5$ and $6\leq j\leq 5+s$, then $R_{12}\subset G$, a contradiction.
If $6\leq i<j\leq 5+s$, then $R_{8}\subset G$, a contradiction.
Thus $d(u)=1$ for any $u\in Y$.
Next we consider the order of $T_{w_i}\,(i\geq 1)$.
Since $R_{1},P_{6}\not\subset G$, we have $|T_{w_i}|=0$, $1\leq i\leq 5$ and $i\neq 3$.
Since $R_{9}\not\subset G$, we have $|T_{w_i}|=0$ for $6\leq i\leq s+5$.
Since $R_{21}\not\subset G$, we have $|T_{w_3}|\leq 1$.
If $|T_{w_3}|=1$, then $s\leq 2$.
Otherwise, $R_{46}\subset G$, a contradiction.
Then $G\subset A_{19}$.
If $|T_{w_3}|=0$, since $R_{95}\not\subset G$, we have $s\leq 6$. Then $G\subset A_{20}$.

\noindent\textbf{Case 2:} $G[X]=P_4\cup sP_1\,(s\geq0)$.

Let $P_4={w_1w_2w_3w_4}$ and $V(sP_1)=\{w_5,\ldots,w_{4+s}\}$.
Since $R_{9}\not\subset G$, we have $|T_{w_i}|=0$ for $5\leq i\leq s+4$.
Since $R_{14}\not\subset G$, we have $|T_{w_i}|=0$, $i=1,4$.
Moreover, since $R_{22}\not\subset G$, we have $|T_{w_2}|\leq 1$ and $|T_{w_3}|\leq 1$.

First, suppose that there exists a vertex $u\in Y$ such that $d(u)=2$.
Let $N_X(u)= \{w_i,w_j\}$, where $i<j$.
First we claim that if $G[\{u,w_i,w_j\}]=K_3$, then $i=2$, $j=3$ and $|T_{w_iw_j}|=1$.
When $N_X(u)=\{w_1,w_2\}$ or $N_X(u)=\{w_3,w_4\}$,
then $R_{18}\subset G$, a contradiction.
When $N_X(u)=\{w_2,w_3\}$, since $R_{28}\not\subset G$, we have $|T_{w_iw_j}|=1$.
Next we claim that if $G[\{u,w_i,w_j\}]=P_3$, then $i=1$, $j=4$ and $|T_{w_iw_j}|=1$.
If $1\leq i<j\leq4$, then $i=1$ and $j=4$.
Otherwise, $G$ has a subdivision of $K_5$, a contradiction.
Moreover, since $R_{22}\not\subset G$, we have $|T_{w_1w_4}|=1$.
If $5\leq i<j\leq s+4$, then $R_{8}\subset G$, a contradiction.
If $1\leq i\leq4$ and $5\leq j\leq4+s$, then $R_{12}\subset G$, a contradiction.
Since $R_{82}\not\subset G$, we have $|T_{w_2w_3}|+|T_{w_1w_4}|\leq 1$.
If $|T_{w_2w_3}|=1$, then $s=0$. Otherwise, $R_{25}\subset G$, a contradiction.
Then $G\subset A_{21}$.
If $|T_{w_1w_4}|=1$, then $s=0$. Otherwise, $G$ has a subdivision of $K_5$, a contradiction.
Then $G\subset A_{22}$.

Next, suppose that $d(u)=1$ for any $u\in Y$.
If $|T_{w_2}|=1$ or $|T_{w_3}|=1$, then $s\leq 2$. Otherwise, $R_{46}\subset G$, a contradiction.
Therefore, we have $G\subset A_{23}$.
If $|T_{w_2}|=|T_{w_3}|=0$, since $R_{95}\not\subset G$, we have $s\leq 6$.
Then $G\subset W_1\subset A_{20}$.

\noindent\textbf{Case 3:} $G[X]=P_3\cup sP_1\,(s\geq0)$.

Let $P_3={w_1w_2w_3}$ and $V(sP_1)=\{w_4,\ldots,w_{3+s}\}$.
First we claim that if $d(u)=2$ for $u\in Y$, then $N_X(u)\subseteq \{w_1,w_2,w_3\}$.
On the contrary, if $N_X(u)\subseteq V(sP_1)$, then $R_{8}\subset G$, a contradiction.
If $|N_X(u)\cap V(sP_1)|=|N_X(u)\cap V(P_3)|=1$,  then we consider the following two cases.
If $N_X(u)\cap V(P_3)=\{w_1\}$ or $\{w_3\}$, then $R_{12}\subset G$, a contradiction.
If $N_X(u)\cap V(P_3)=\{w_2\}$, then $G$ contains a subdivision of $K_{3,3}$, a contradiction.
Since $R_{9}\not\subset G$, we have $|T_{w_i}|=0$, $4\leq i\leq s+3$.
Since $R_{22}\not\subset G$, we have $|T_{w_1}|\leq 1$ and $|T_{w_3}|\leq 1$.
Since $R_{28}\not\subset G$, we have $|T_{w_1w_2}|\leq1$ and $|T_{w_2w_3}|\leq1$.

First we have the following facts.

\begin{fact}\label{diam4:f5}
If $|T_{w_1w_2}|\neq 0$, then $|T_{w_3}|=0$. Similarly, if $|T_{w_2w_3}|\neq 0$, then $|T_{w_1}|=0$ ($R_{14}\not\subset G$).
\end{fact}

\begin{fact}\label{diam4:f6}
$|T_{w_1w_3}|+|T_{w_2}|\leq 3$. ($R_j\not\subset G$ for $j\in\{38,47,55,77,94\}$).
\end{fact}

\noindent\textbf{Subcase 3.1:} There exists a vertex $u\in Y$ such that $d(u)=2$.

First we claim that $s=0$.
Assume on the contrary that $s\geq1$.
If $|T_{w_1w_2}|\neq 0$ or $|T_{w_2w_3}|\neq 0$, we have $R_{25}\subset G$, a contradiction. If $|T_{w_1w_3}|\neq 0$, then $G$ contains a subdivision of $K_5$ as a subgraph, a contradiction.
Thus $s=0$.
By (F\ref{diam4:f6}), we consider the order of $T_{w_1w_3}$.
When $|T_{w_1w_3}|=3$, by (F\ref{diam4:f6}), we have $|T_{w_2}|=0$.
If $|T_{w_1w_2}|=1$ and $|T_{w_2w_3}|=0$ (or $|T_{w_1w_2}|=0$ and $|T_{w_2w_3}|=1$), by (F\ref{diam4:f5}), then $G\subset A_{24}$.
If $|T_{w_1w_2}|=|T_{w_2w_3}|=1$, by (F\ref{diam4:f5}), then $G=A_{25}$.
If $|T_{w_1w_2}|=|T_{w_2w_3}|=0$, then $G\subset A_{26}$.
When $|T_{w_1w_3}|=2$, by (F\ref{diam4:f6}), we have $|T_{w_2}|\leq 1$.
Similarly, by (F\ref{diam4:f5}), $G\subset A_{27},A_{28}$ or $A_{29}$.
When $|T_{w_1w_3}|=1$, by (F\ref{diam4:f6}), we have $|T_{w_2}|\leq 2$.
Similarly, by (F\ref{diam4:f5}), $G\subset A_{30},A_{31}$ or $A_{32}$.
When $|T_{w_1w_3}|=0$, by (F\ref{diam4:f6}), we have $|T_{w_2}|\leq 3$.
Similarly, by (F\ref{diam4:f5}), $G\subset A_{33}$ or $A_{34}$.

\noindent\textbf{Subcase 3.2:} For any $u\in Y$, $d(u)=1$.

If $|T_{w_2}|\neq 0$, since $R_{21}, R_{47}\not\subset G$, we have $|T_{w_2}|+s\leq 3$.
Then $G\subset A_{29},A_{32}$ or $A_{35}$.
Next we consider the case $|T_{w_2}|=0$.
Since $R_{46}\not\subset G$, if $|T_{w_1}|= 1$ or $|T_{w_3}|= 1$, then $s\leq 3$.
Then $G\subset A_{26}$.
If $|T_{w_1}|= |T_{w_3}|= 0$, since $R_{97}\not\subset G$, we have $s\leq 7$.
Then $G\subset W_2\subset A_{20}$.

\noindent\textbf{Case 4:} $G[X]=aP_2\cup sP_1\,(a\geq1,s\geq0)$.

If $Y=\varnothing$, then $G\subset D_{1}$.
Next we let $Y\neq \varnothing$ and we distinguish the following two cases.

\noindent\textbf{Subcase 4.1:} There exists $u\in Y$ such that $d(u)=2$.

Since $R_{8}\not \subset G$, we have $N_X(u)\not \subseteq V(sP_1)$.
Suppose that there exists $u\in Y$ such that $G[N_X(u)]=P_2$.
Then $|X|=3$.
Otherwise, we have $R_{23}\subset G$ or $R_{13}\subset G$, a contradiction.
Let $X=\{w_1,w_2,w_3\}$ and $N_X(u)=\{w_1,w_2\}$.
Since $R_{6},R_{22}\not\subset G$, we have $|T_{w_1}|=|T_{w_2}|=|T_{w_3}|=0$.
Since $R_{24}\not\subset G$, we have $|T_{w_1w_2}|=1$.
Since $R_{12}\not\subset G$, we have $|T_{w_1w_3}|=|T_{w_2w_3}|=0$.
Then $G=A_{36}$.

Next assume that there does not exist $u\in Y$ such that $G[N_X(u)]=P_2$.
Suppose that there exists $u\in Y$ (let $N_X(u)=\{w_2,w_3\}$) such that $|N_X(w_2)|=1$ and $|N_X(w_3)|=0$.
We assert that $a=1$ and $s\leq 3$.
On the contrary, if $a\geq 2$, then $R_{8}\subset G$, a contradiction.
If $s\geq 4$, then $R_{37}\subset G$, a contradiction.
Let $X=\{w_1,w_2,\ldots,w_{2+s}\}$, where $N_X(w_2)=\{w_1\}$.
First we have the following facts.

\begin{fact}\label{diam4:f7}
$|T_{w_2w_3}|+|T_{w_1}|\leq 2$. ($R_{1}$, $R_{21}\not\subset G$).
\end{fact}

\begin{fact}\label{diam4:f8}
$s+|T_{w_2}|\leq 3$. ($R_{21}\not \subset G$).
\end{fact}

\begin{fact}\label{diam4:f9}
$|T_{w_i}|=0$, $3\leq i\leq 2+s$. ($R_{7}\not\subset G$ if $|T_{w_3}|\geq 1$, $R_{2}\not\subset G$ if $|T_{w_i}|\geq 1\,(4\leq i\leq 2+s)$).
\end{fact}

 Next we claim that if $|T_{w_iw_j}|\neq0$, then $(i,j)=(1,3)$ or $(2,3)$.
Since $R_{2}\not \subset G$, we have $|T_{w_1w_j}|=0$, $j\geq 4$.
Since $|T_{w_2w_3}|\neq 0$, we have $|T_{w_2w_j}|=0, j\geq 4$.
Otherwise, $G$ contains a subdivision of $K_{3,3}$ as a subgraph.
Since $R_{8}\not\subset G$, we have $|T_{w_iw_j}|=0, j>i\geq 3$.
Moreover, if $|T_{w_1w_3}|\neq 0$, then $s=1$.
Otherwise, $G$ contains a subdivision of $K_5$ as a subgraph.
Since $R_{1}$, $R_{21}\not \subset G$, we have $|T_{w_1w_3}|+|T_{w_2}|\leq 2$.
If $|T_{w_1w_3}|=2$, then $s=1$ and $|T_{w_2}|=0$; hence, by (F\ref{diam4:f7}) and (F\ref{diam4:f9}), $G$ is an induced subgraph of $A_{37}$ or $A_{38}$.
If $|T_{w_1w_3}|=1$, then $s=1$ and $|T_{w_2}|\leq1$; hence, by (F\ref{diam4:f7}) and (F\ref{diam4:f9}), $G$ is an induced subgraph of $A_{38}$ or $A_{39}$.
Finally, suppose that $|T_{w_1w_3}|=0$.
When $s=1$, by (F\ref{diam4:f8}), we have $|T_{w_2}|\leq 2$.
By (F\ref{diam4:f7}) and (F\ref{diam4:f9}), we have either $G\subset A_{40}$ or $G\subset A_{41}$.
When $s\geq 2$, since $R_{77}\not\subset G$, we have $|T_{w_2w_3}|=1$.
By (F\ref{diam4:f7}), we have $|T_{w_1}|\leq 1$.
Then by (F\ref{diam4:f8}) and (F\ref{diam4:f9}), we have $G\subset A_{42}$ if $s=2$ and $G\subset A_{43}$ if $s=3$.

Finally, suppose that every vertex $u\in Y$ of degree $2$ has the following property: after writing $N_X(u)=\{w_2,w_3\}$, if $w_2\nsim w_3$, then $|N_X(w_2)|=|N_X(w_3)|=1$.
We assert that $a=2$ and $s\leq 2$.
Otherwise, if $a\geq 3$, then $R_{8}\subset G$, a contradiction.
If $s\geq 3$, then $R_{37}\subset G$, a contradiction.
Let $X=\{w_1,w_2,\ldots,w_{4+s}\}$, where $N_X(w_2)=\{w_1\}$ and $N_X(w_3)=\{w_4\}$.
Since $R_{1},P_6\not \subset G$, we have $|T_{w_i}|=0$, $i=1,2,3,4$.
Since $G$ contains no subdivision of $K_5$, we have $|T_{w_1w_3}|=|T_{w_2w_4}|=0$.
Suppose that $|T_{w_1w_4}|\neq0$.
If $|T_{w_1w_4}|\geq 2$ or $|T_{w_2w_3}|\geq 2$, then $R_{1}\subset G$, a contradiction.
If $s\geq 1$, then $G$ contains the subdivision of $K_5$ as a subgraph, a contradiction.
Therefore, we have $|T_{w_1w_4}|=1$, $|T_{w_2w_3}|=1$ and $s=0$.
Then $G=A_{44}$.
Next, suppose that $|T_{w_1w_4}|=0$.
Then $|T_{w_2w_3}|\leq 2$;
otherwise, we have $R_{21}\subset G$, a contradiction.
If $|T_{w_2w_3}|=2$, then $s=0$ since $R_{21}\not\subset G$, and hence $G= A_{45}$.
If $|T_{w_2w_3}|=1$, since $R_{2}\not\subset G$, we have $|T_{w_i}|=0\,(i\geq 5)$.
Then $G\subset A_{46}$.

\noindent\textbf{Subcase 4.2:} For any $u\in Y$, $d(u)=1$.

If there exists a vertex $u\in Y$ such that $N_X(u)\subseteq V(aP_2)$, then $a=1$ since $R_{10}\not\subset G$.
Let $X=\{w_1,w_2,\ldots,w_{2+s}\}$ and $G[\{w_1,w_2\}]=P_2$.
Since $|X|\geq3$ and $R_{21}\not\subset G$, we have $|T_{w_i}|\leq2$ for $i=1,2$.
Since $R_{3}\not\subset G$, we have $|T_{w_i}|=0$, $i=3,\ldots, 2+s$.
Without loss of generality, suppose that $|T_{w_1}|\geq |T_{w_2}|$.
If $|T_{w_1}|=2$ or $|T_{w_2}|=2$, since $R_{21}\not\subset G$, we have $s= 1$ and then $G\subset A_{47}$.
If $|T_{w_1}|=|T_{w_2}|=1$, since $R_{45}\not\subset G$, we have $s\leq 2$ and then $G\subset W_3 \subset A_{43}$.
If $|T_{w_1}|=1$ and $|T_{w_2}|=0$, since $R_{46}\not\subset G$, we have $s\leq 3$ and then $G\subset W_4\subset A_{43}$.

Next, suppose that $N_X(u)\subseteq V(sP_1)$ for any $u\in Y$.
Since $R_{4}\not \subset G$, we have $|T_{w}|\leq 1$ for $w\in V(sP_1)$.
Then $G\subset W_7\subset D_{1}$.

\noindent\textbf{Case 5:} $G[X]=sP_1\,(s\geq3)$.

Let $X=\{w_1,w_2,\ldots,w_{s}\}$.

\noindent\textbf{Subcase 5.1:} There exists a vertex $u\in Y$ such that $d(u)=2$.

Since $R_{77}\not\subset G$ and $s\geq3$, we have $|T_{w_iw_j}|\leq1$ for $1\leq i<j\leq s$.
Moreover, if $|T_{w_iw_j}|=1$ and $|T_{w_kw_\ell}|=1$, then $\{i,j\}\cap\{k,\ell\}\neq\varnothing$.
Otherwise, we have $R_{2}\subset G$, a contradiction.
If there exist three vertices of degree 2 in $Y$, say $u_1,u_2,u_3$, such that $|N_X(u_1)\cap N_X(u_2)\cap N_X(u_3)|=1$, then $G$ contains a subdivision of $K_{3,3}$ as a subgraph, a contradiction.
Hence, there are at most three vertices of degree $2$ in $Y$.
Since $R_{2}\not\subset G$, we have $|T_w|=0$ for every $w\in X$ that is not adjacent to a vertex of degree $2$ in $Y$.

If there are three vertices of degree $2$ in $Y$, say $u_1,u_2,u_3$, then, without loss of generality, we may assume that $N_X(u_1)=\{w_1,w_2\}$, $N_X(u_2)=\{w_2,w_3\}$, and $N_X(u_3)=\{w_1,w_3\}$.
Then $|X|=3$; otherwise, $G$ contains a subdivision of $K_5$ as a subgraph, a contradiction.
Since $R_{1}\not \subset G$, we have $|T_{w_i}|= 0$, $i=1,2,3$.
Then $G=A_{48}$.

If there are two vertices of degree $2$ in $Y$, say $u_1,u_2$, then, without loss of generality, we may assume that $N_X(u_1)=\{w_1,w_2\}$ and $N_X(u_2)=\{w_2,w_3\}$.
Since $R_{1}\not\subset G$, we have $|T_{w_1}|=|T_{w_3}|=0$.
Since $R_{21}\not \subset G$, we have $|T_{w_2}|+s\leq 4$.
If $s=3$, then $G\subset A_{49}$.
If $s=4$, then $G=A_{50}$.

If there is exactly one vertex of degree $2$ in $Y$, say $u_1$, then, without loss of generality, we may assume that $N_X(u_1)=\{w_1,w_2\}$ and $|T_{w_1}|\geq|T_{w_2}|$.
Since $R_{21},R_{37}\not \subset G$, we have $|T_{w_1}|+s\leq 4$.
If $s=3$, then $G\subset A_{51}$.
If $s=4$, then $G=W_5\subset A_{46}$.

\noindent\textbf{Subcase 5.2:} For any $u\in Y$, $d(u)=1$.

Since $s\geq3$ and $R_{21}\not\subset G$, we have $|T_{w}|\leq 1$ for $w\in X$.
Then $G\subset W_7\subset D_1$.
\end{proof}

\begin{proof}[Proof of Theorem \ref{thm-51}]
The result follows from Lemmas \ref{diam4:Xstructure}, \ref{diam4:Xsmall}, \ref{diam4:Xcycle}, and \ref{diam4:Xforest}.
\end{proof}

\section{Proof of Theorem \ref{thm-6}}
We  now combine the results of the preceding sections to prove the main theorem.

\begin{proof}[Proof]

Suppose first that $G\in\mathcal{G}$. If $d(G)\geq4$, then Theorem \ref{thm-51} implies that $G$ is a connected induced subgraph of one of $A_1,\ldots,A_{51},D_1,D_2$. Now suppose that $d(G)<4$. Since $G$ is planar, $\omega(G)\leq4$. If $\omega(G)=4$, the conclusion follows from Theorem \ref{thm-31}. If $\omega(G)=3$ and $G$ contains $K_4-e$ as an induced subgraph, the conclusion follows from Theorem \ref{thm-41}. In all remaining cases, $G$ is $K_4-e$-free, so the conclusion follows from Lemma \ref{lem-5}. Together, these four results show that $G$ is a connected induced subgraph of a graph in
$\{A_1,\ldots,A_{91},H_1,\ldots,H_{153},D_1,D_2\}$.

Conversely, suppose that $G$ is a connected induced subgraph of a graph $F$ in the displayed family. Every such $F$ is planar and satisfies $\lambda_2(F)\leq1$ by Lemma \ref{lem-3}. Since planarity is inherited by subgraphs and eigenvalue interlacing (Lemma \ref{lem-2}) gives $\lambda_2(G)\leq\lambda_2(F)$, the graph $G$ is connected and planar and satisfies $\lambda_2(G)\leq1$. Hence $G\in\mathcal{G}$.
\end{proof}

\renewcommand{\refname}{References}

\end{document}